\documentclass[10pt]{amsart}

\usepackage{amsmath, amssymb, amsthm, latexsym, amscd}
\usepackage[all]{xy}
\usepackage{graphicx}

 \newlength{\baseunit}               
 \newcount{\numlines}                
\newtheorem*{thm}{Theorem}
\newtheorem{theorem}{Theorem}
\newtheorem{prop}{Proposition}
\newtheorem{lemma}{Lemma}
\newtheorem{corollary}{Corollary}

\theoremstyle{definition}
\newtheorem{define}{Definition}

\newtheorem{remark}{Remark}

\DeclareMathOperator{\Hom}{Hom}

\DeclareMathOperator{\End}{End}

\DeclareMathOperator{\Id}{Id}

\DeclareMathOperator{\Ext}{Ext}

\begin{document}
\pagestyle{plain} 
\title{Category $ \mathcal{O} $ and $ \mathfrak{sl}_k  $ Link Invariants}
\author{Joshua Sussan}
\date{\today}
\maketitle

\setcounter{tocdepth}{1}
\tableofcontents

\section{Introduction}
The program of categorification via category $ \mathcal{O} $ was introduced by J. Bernstein, I. Frenkel, and M. Khovanov in [BFK] where tensor powers of the standard two dimensional representation of 
$ \mathfrak{sl}_2 $ were recognized as  Grothendieck groups of certain subcategories of $ \mathcal{O} $ for various $ \mathfrak{gl}_n. $  They had two different constructions.  One was based on studying certain blocks with singular generalized central characters.  The other was based on examining the trivial regular block but by considering various parabolic subcategories.  In the first case the action of $ \mathfrak{sl}_2 $ was categorified by projective functors acting on these singular categories.  The intertwiners were categorified by derived Zuckerman functors acting on the derived category $ \mathcal{O}. $  In the latter case, the action of the Lie algebra was lifted to Zuckerman functors and the intertwiners became projective functors.  These two constructions are related by Koszul duality where the projective functors get exchanged for Zuckerman functors and visa-versa [Rh], [MOS].

Several conjectures posed in [BFK] were solved in [Str2].  The most important result of that work was that the Jones polynomial was lifted to a functorial invariant of tangles.
Natural transformations between these functors became invariants of 2-tangles acting as cobordisms between two different tangles [Str4]. One of the principle goals of categorification is to obtain invariants for $ n+1 $ dimensional topological objects by homological realizations of classical $ n $ dimensional invariants.

There is a parallel approach to categorification of link invariants due to Khovanov.  In [Khov], a bi-graded homology theory was constructed whose graded Euler characteristic is the Jones polynomial.  
M. Khovanov and L. Rozansky extended this work to a homology theory categorifying the HOMFLYPT polynomial [KR].  The theory of matrix factorization is their primary tool.
H. Murakami, T. Ohtsuki, and S. Yamada have an intepretation of the HOMFLYPT polynomial through a polynomial assigned to certain planar graphs [MOY].
Khovanov and Rozansky categorified these polynomials which is a crucial step in the categorification of the HOMFLYPT polynomial.

The goal of this paper is to categorify the $ \mathfrak{sl}_k $ tangle invariant via category $ \mathcal{O}. $  It is important to understand all of the planar graphs representation theoretically.  The planar graphs give a graphical interpretation of intertwiners between various tensor products $ \Lambda^{i_1} V_{k-1} \otimes \cdots \otimes \Lambda^{i_r} V_{k-1} $ of fundamental $ \mathfrak{sl}_k $ representations.  The first step is to categorify tensor powers of the standard $ k- $ dimensional representation of $ \mathfrak{sl}_k. $ In order to accomplish this, we look at more general singular blocks of category $ \oplus_{\bf d} \mathcal{O}_{\bf d}(\mathfrak{gl}_n). $  Here, $ {\bf d} $ denotes a $ k- $ tuple $ (d_{k-1}, \ldots, d_0). $  Then $ \mathcal{O}_{\bf d}(\mathfrak{gl}_n) $ is a block corresponding to an integral dominant weight whose stabalizer under the dot action of the Weyl group is $ S_{d_{k-1}} \times \cdots \times S_{d_0}. $ Parabolic subcategories of these blocks,  $ \mathcal{O}_{\bf d}^{\mathfrak{p}}(\mathfrak{gl}_n) $ will provide a categorification of tensor products of fundamental representations.  The action of the Lie algebra $ \mathfrak{sl}_k $ will be lifted to projective functors.  These are functors given by tensoring with a finite dimensional representation and then projecting onto a certain block.  Morphisms between various $ \Lambda^{i_1} V_{k-1} \otimes \cdots \otimes \Lambda^{i_r} V_{k-1} $ will become inclusion and derived Zuckerman functors between various parabolic subcategories of these blocks.  Relations between the intertwiners become functorial isomorphisms between various compositions of inclusion and derived Zuckerman functors.  It should be possible to do this construction in the Koszul dual situation as well.  In this case projective functors categorify the intertwiners and derived Zuckerman functors categorify the action of the Lie algebra.

This construction only gives rise to a categorification of $ \mathfrak{sl}_k- $ modules.  In order to get a categorification of the quantum group, we consider the above categories to be a category of graded modules.  Then in the Grothendieck group the shift functor will descend to the action of multiplication by the quantum parameter $ q. $  The theory of graded category $ \mathcal{O} $ is due to W. Soergel and has been studied extensively by C. Stroppel [Str1].   Projective functors have graded lifts.  This allows us to give a categorification of the quantum Serre relations.  Zuckerman functors  and inclusion functors also have graded lifts. 
Graded functors may now be assigned to all flat tangles.  In order to extend this construction to tangle with crossings, we consider adjunction morphisms between shifted compositions of inclusion and Zuckerman functors for a minimal parabolic and the identity functor.  Cones of these morphisms are assigned to the crossings.
In this graded setup it is easy to see that these cones satisfy the $ \mathfrak{sl}_k $ skein relation.  A cobordism between two tangles should give rise to a natural transformation between the functors assigned to the tangles.  This natural transformation should become an invariant of 2-tangles as in [Str4], [Khov3], [Jac].  

In section 2 we review the necessary finite dimensional representation theory of $ \mathfrak{sl}_k. $  A categorification of the fundamental representations will be given in section 3.  An equivalence of categories between certain parabolic subcategories which generalizes corollary 5 in [BFK] is also given in this section.  This equivalence is used in the definition of the functors assigned to the flat tangles. In section 4 we provide functorial analogues of the graphic relations presented in section 2.  Each of the graphical relations will give rise to a functorial isomorphism.  Natural transformations are constructed and shown that when restricted to generalized Verma modules, they are isomorphisms.  Most of the section is devoted to calculations of the relevant functors on generalized Verma modules.  In section 5 graded category $ \mathcal{O} $ will be discussed and all of the results in the previous sections will be lifted to the graded case.
We define in section 6 our functor valued tangle invariant $ \mathbb{F} $ and prove that it satisfies the Reidemeister relations.  

Suppose that tangles $ T $ and $ T' $ are morphisms from $ r $ points to $ r' $ points labeled from the set $ \lbrace 1, k-1 \rbrace. $  Denote by $ n $ and $ n' $ the sum of the labels for the $ r $ and $ r' $ points respectively.  These compositions of $ n $ and $ n' $ naturally give rise to parabolic subalgebras $ \mathfrak{p} $ and $ \mathfrak{p}' $ of $ \mathfrak{gl}_n $ and $ \mathfrak{gl}_{n'}. $ Theorem ~\ref{main} is the main result of this work.

\begin{thm}
If tangles $ T $ and $ T' $ are ambient isotopic tangles, then $ \mathbb{F}(T) $ and $ \mathbb{F}(T') $ are isomorphic as derived functors from $ D^b(\oplus_{\bf d} \mathcal{O}_{\bf d}^{\mathfrak{p}}(\mathfrak{gl}_n)) $ to $ D^b(\oplus_{\bf d'} \mathcal{O}_{\bf d'}^{\mathfrak{p'}}(\mathfrak{gl}_{n'})). $ 
\end{thm}
The functorial isomorphisms of section 4 play a crucial role in proving this theorem.  When restricted to $ (0,0)- $ tangles, the invariant is a complex of graded vector spaces and thus homology groups may be assigned to links.

It follows immediately from the definition of $ \mathbb{F} $ that on the Grothendieck group, the following equality holds:

\begin{thm}
Let $ \Omega_i $ and $ \Pi_i $ be the functors assigned to the crossings.  Then
$ q^{k}[\Omega_i]-q^{-k}[\Pi_i]=(-1)^{k-1}(q^{1}-q^{-1})[\Id] $
\end{thm}

An explicit relationship between the functorial tangle invariant introduced by Khovanov in [Khov5] and the category $ \mathcal{O} $ invariant of [Str2] was given in [Str5].  Stroppel  considered the subcategory of the trivial block of $ \mathcal{O}(\mathfrak{sl}_{2n}) $ of modules locally finite with respect to a parabolic subalgebra whose reductive part is $ \mathfrak{gl}_n \oplus \mathfrak{gl}_n $ and whose projective presentation only consists of projective-injective modules.  This subcategory categorifies the space of invariants in tensor products of the standard two dimensional of $ \mathfrak{sl}_2. $ This subcategory is equivalent to a category of modules over an associative algebra which is isomorphic to the algebra constructed in [Khov5].
A natural question is to find a connection between the functorial $ \mathfrak{sl}_k $ invariant in this work and the invariant constructed via matrix factorization in [KR] or Soergel bimodules in [Khov4].  One should probably look for a subcategory of $ \mathcal{O} $ that categorifies the space of invariants in tensor products of the standard $ k- $ dimensional representation of $ \mathfrak{sl}_k $ and then study the associative algebra governing that subcategory.   

{\it Acknowledgements: } The author is grateful to his advisor, Igor Frenkel for his encouragement and support throughout the development of this project.  In addition, the author is very thankful to Mikhail Khovanov, Raphael Rouquier, Catharina Stroppel, and Gregg Zuckerman for helpful conversations and comments on preliminary versions of this work.

\section{Representation Theory of $ \mathfrak{sl}_k $}

\subsection{Basic Definitions}
We begin by reviewing the finite dimensional representation theory of $ \mathfrak{sl}_k.  $  
Recall that $ \mathfrak{sl}_k $ is the Lie algebra of $ k \times k $ matrices with entries in $ \mathbb{C}. $
Denote the matrix with only 1 in the $ (i,j) $ entry by $ e_{i,j}. $
There is the triangular decomposition $ \mathfrak{sl}_k = {\mathfrak{n}}^{-}+{\mathfrak{h}}+{\mathfrak{n}}^{+} $
where $ {\mathfrak{n}}^{-} $ is the subalgebra of lower triangular matrices, $ {\mathfrak{n}}^{+} $ is
the subalgebra of upper triangular matrices and $ {\mathfrak{h}} $ is the Cartan subalgebra of
diagonal matrices.

The dual $ {\mathfrak{h}}^{*} $ of the Cartan subalgebra has a basis
$ \alpha_1=e_{1,1}^{*}-e_{2,2}^{*}, \ldots, \alpha_i=e_{i,i}^{*}-e_{i+1,i+1}^{*}, \ldots,
\alpha_{k-1}=e_{k-1,k-1}^{*} - e_{k,k}^{*} $ where $ e_{i,j}^{*}(e_{r,s}) = \delta_{i,r} 
\delta_{j,s}. $
These $ \alpha_i $ are called the positive simple roots for $ \mathfrak{sl}_k. $

A linear map $ \lambda \colon {\mathfrak{h}} \rightarrow \mathbb{C} $
may be written in coordinates $ \lambda = \lambda_{1} e_{1,1}^* + \cdots +
\lambda_{k} e_{k,k}^*. $

The finite dimensional irreducible representations of $ \mathfrak{sl}_k $ are indexed by these weights 
$ \lambda. $ There exists a unique irreducible representation of highest weight 
$ \lambda = \lambda_{1} e_{1,1}^* + \cdots + \lambda_{k} e_{k,k}^* $ when 
$ \lambda_1 \geq \lambda_2 \geq \cdots \geq \lambda_k. $
For our purposes, one of the most important representations is the standard $ k- $ dimensional vector space
$ \mathbb{C}^{k} $ where $ \mathfrak{sl}_k $ acts on it naturally as matrices.  Call this space $ V_{k-1}. $
The basis for this vector space is given by $ e_1, \ldots, e_k $ and $ e_{ij}(e_m) = 
\delta_{j,m} e_i. $  Notice that the highest weight of this representation is $ e_{1,1}^*. $
We denote the $ i\text{th} $ exterior power of a module $ V $ by $ \Lambda^{i} V. $

\begin{prop}
Let $ V_{k-1} $ be the module defined above.  Then $ \Lambda^{i} V_{k-1} $ is an irreducible module
of highest weight $ e_{1,1}^* + \cdots + e_{i,i}^*. $
\end{prop}

\begin{proof}
See [FH] page 221.
\end{proof}

Now we define the algebra $ \mathcal{U}_{q}(\mathfrak{sl}_k) $ and its fundamental representations.

\begin{define}
The quantum group $ \mathcal{U}_{q}(\mathfrak{sl}_k) $ is the associative algebra over $ \mathbb{C}(q) $ with generators $ E_i, F_i, K_i, K_{i}^{-1} $ for $ i=1, \ldots, k-1 $ satisfying the following conditions:
\begin{enumerate}
\item $ K_i K_i^{-1} = K_i^{-1} K_i = 1$
\item $ K_i K_j = K_j K_i $
\item $ K_i E_j = q^{c_{i,j}} E_j K_i $
\item $ K_i F_j = q^{-c_{i,j}} F_j K_i $
\item $ E_i F_j - F_j E_i = \delta_{i,j} \frac{K_i - K_i^{-1}}{q-q^{-1}} $
\item $ E_i E_j = E_j E_i $ if $ |i-j|>1 $
\item $ F_i F_j = F_j F_i $ if $ |i-j|>1 $
\item $ E_i^2 E_{i \pm 1} - (q+q^{-1}) E_i E_{i \pm 1} E_i + E_{i \pm 1} E_i^2 = 0 $
\item $ F_i^2 F_{i \pm 1} - (q+q^{-1}) F_i F_{i \pm 1} F_i + F_{i \pm 1} F_i^2 = 0 $
\end{enumerate}
where 
\begin{eqnarray*} 
c_{i,j} = 2 & \textrm{ if } &j=i \\
-1 & \textrm{ if } &j=i \pm1\\
0 & \textrm{if} &|i-j|>1.
\end{eqnarray*}
\end{define}

The most basic representation of $ \mathcal{U}(\mathfrak{sl}_k) $ is $ V_{k-1}. $  It is the k-dimensional vector space with basis $ v_0, \ldots, v_{k-1}. $  
The algebra acts on this space as follows:
\begin{eqnarray*}
E_i v_{j} = 0 &\textrm{if}	&j \neq i-1\\
E_i v_{j}= v_i &\textrm{if}	&j=i-1\\
F_i v_j = 0 &\textrm{if} & j \neq i\\ 
F_i v_j = v_{i-1} &\textrm{if} &j=i\\ 
K_i^{\pm 1} v_j = q^{\pm 1} v_i &\textrm{if}	&j=i\\
K_i^{\pm 1} v_{j} = q^{\mp 1} v_{i-1} &\textrm{if} &j=i-1\\ 
K_i^{\pm 1} v_j = v_j &\textrm{if} &j \neq i-1, i.
\end{eqnarray*}
There are several intertwiners between various representations that will be important for later.
There is the map $ \Lambda^k V_{k-1} \rightarrow V_{k-1} \otimes \Lambda^{k-1} V_{k-1} $ given by
$$ v_{k-1} \wedge \cdots \wedge v_0 \rightarrow \sum_{j=0}^{k-1} (-1)^{j} q^{k-j-1} v_j \otimes (v_{k-1} \wedge \cdots \wedge \hat{v_j} \wedge \cdots \wedge v_0), $$ 
where $ \hat{v_j} $ means that the term is omitted from the expression.

There is the map $ \Lambda^k V_{k-1} \rightarrow \Lambda^{k-1} V_{k-1} \otimes V_{k-1} $ given by
$$ v_{k-1} \wedge \cdots \wedge v_0 \rightarrow \sum_{j=0}^{k-1} (-1)^{k-1-j} q^{j} (v_{k-1} \wedge \cdots \wedge \hat{v_j} \wedge \cdots \wedge v_0 \otimes v_j). $$

There is a map in the other direction: $ V_{k-1} \otimes \Lambda^{k-1} V_{k-1} \rightarrow \Lambda^k V_{k-1} $ given by 
$$ v_j \otimes (v_{k-1} \wedge \cdots \hat{v_j} \wedge \cdots \wedge v_0) \rightarrow (-1)^{k-j-1} q^{-j} v_{k-1} \wedge \cdots \wedge v_0. $$

There is also the map $ \Lambda^{k-1} V_{k-1} \otimes V_{k-1} \rightarrow \Lambda^k V_{k-1} $ given by 
$$ (v_{k-1} \wedge \cdots \hat{v_j} \wedge \cdots \wedge v_0) \otimes v_j \rightarrow (-1)^{j} q^{j-k+1} v_{k-1} \wedge \cdots \wedge v_0. $$

The inclusion map $ \Lambda^2 V_{k-1} \rightarrow V_{k-1} \otimes V_{k-1} $ is determined by
$$ v_i \wedge v_j \rightarrow v_i \otimes v_j - q^{1} v_j \otimes v_i $$ where $ i>j. $

The projection map $ V_{k-1} \otimes V_{k-1} \rightarrow \Lambda^2 V_{k-1} $ is given by
\begin{eqnarray*} 
v_i \otimes v_j \rightarrow q^{-1} v_i \wedge v_j& \textrm{ if } &i>j \\
-v_i \wedge v_j & \textrm{ if } &i<j.
\end{eqnarray*}

More generally,
let 
$$ \pi_{r_1} \otimes \cdots \otimes \pi_{r_t} \colon V_{k-1}^{\otimes (r_1+\cdots + r_t)} \rightarrow \Lambda^{r_1}  V_{k-1} \otimes \cdots \otimes \Lambda^{r_t}  V_{k-1} $$ 
be the canonical projection map and
$$ i_{r_1} \otimes \cdots \otimes i_{r_t} \colon \Lambda^{r_1}  V_{k-1} \otimes \cdots \otimes \Lambda^{r_t}  V_{k-1} \rightarrow V_{k-1}^{\otimes (r_1+\cdots +r_t)} $$ 
be the canonical inclusion map.
Since we will not need these more general intertwiners, we will omit the precise formulas.

\subsection{Graphical Calculus} There is a graphical description of the intertwiners between tensor products of fundamental representations via colored trivalent graphs.  
A line segment labeled by an integer $ i $ where $ 1 \leq i \leq k, $ will depict the representation $ \Lambda^{i}V_{k-1}. $   Note that an edge labeled by $ k $ depicts the trivial representation and thus such an edge may be added or removed at will.
The figure below depicts the canonical projection map $ \Lambda^{i} V_{k-1} \otimes \Lambda^{j} V_{k-1} \rightarrow \Lambda^{i+j} V_{k-1}. $ 

\begin{figure}[htb]
  \centering
  \includegraphics{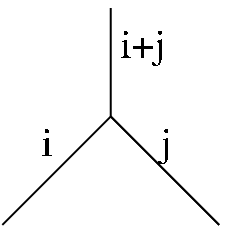}
  \caption{}
  \label{fig1}
\end{figure}

The inclusion $ \Lambda^{i+j} V_{k-1} \rightarrow \Lambda^{i} V_{k-1} \otimes \Lambda^{j} V_{k-1} $ is given by the diagram below.

\begin{figure}[htb]
  \centering
  \includegraphics{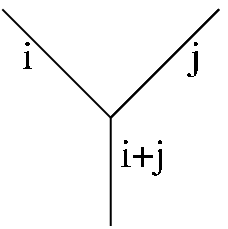}
  \caption{}
  \label{fig2}
\end{figure}

There are relations between these maps.  Graphically, the relations are depicted by the following five diagrams.

\begin{figure}[htb]
  \centering
  \includegraphics{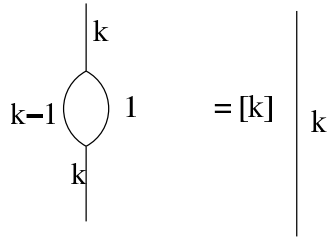}
  \caption{Diagram Relation 1}
  \label{relation1}
\end{figure}

\begin{figure}[htb]
  \centering
  \includegraphics{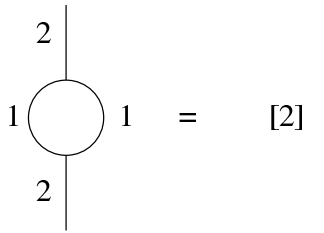}
  \caption{Diagram Relation 2}
  \label{relation2}
\end{figure}

\begin{figure}[htb]
  \centering
  \includegraphics{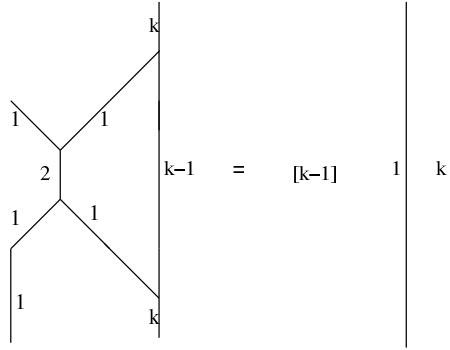}
  \caption{Diagram Relation 3}
  \label{relation3}
\end{figure}

\begin{figure}[htb]
  \centering
  \includegraphics{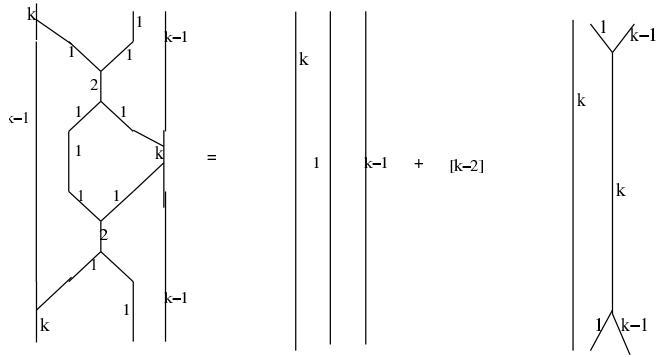}
  \caption{Diagram Relation 4}
  \label{relation4}
\end{figure}

\begin{figure}[htb]
  \centering
  \includegraphics{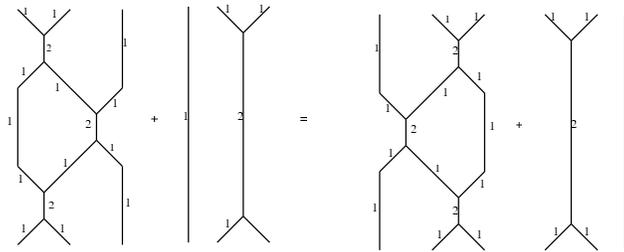}
  \caption{Diagram Relation 5}
  \label{relation5}
\end{figure}

If the trivalent graph has external edges labeled only by $ k, $ then it represents a composition of intertwiners from the the trivial representation to itself so a Laurent polynomial may be assigned to it.  We will say that such a graph is closed.

\begin{theorem}
There is a Laurent polynomial $ \langle D \rangle_k $ in $ \mathbb{Z}[q, q^{-1}] $ which may be assigned to closed, colored trivalent graphs $ D $ which satisfy the five relations above.
It is invariant under ambient isotopy of $ \mathbb{R}^{2}. $
\end{theorem}

\begin{proof}
See sections 1 and 2 of [MOY].
\end{proof}

The Laurent polynomial in the theorem above coincides with the representation-theoretic polynomial.

A version of the HOMFLYPT for links may then be defined.  Each crossing may be resolved in the following two ways:

\begin{figure}[htb]
  \centering
  \includegraphics{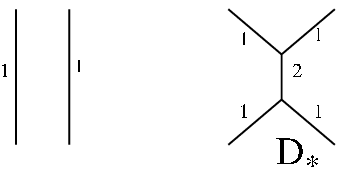}
  \caption{}
  \label{resolve}
\end{figure}

Let $ D_{+}, D_{-}, D_{0} $ be identical link diagrams except near a crossing as given by the figure below.

\begin{figure}[htb]
  \centering
  \includegraphics{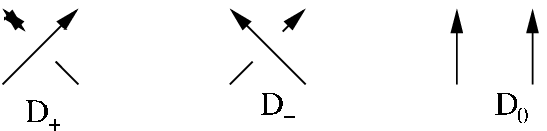}
  \caption{}
  \label{skein}
\end{figure}

\begin{define}
Let $ D $ be a planar projection of a link. 
Let
\begin{enumerate}
\item $ \langle D_+ \rangle_k = q^{1} \langle D_0 \rangle_k - \langle D_{*} \rangle_k, $
\item $ \langle D_{-} \rangle_k = q^{-1} \langle D_0 \rangle_k - \langle D_{*} \rangle_k, $
\item $ P_k(D) = q^{(k)(-w(D))} \langle D \rangle_{k} $ where $ w(D) $ is the difference between the number of positive and negative crossings.
\end{enumerate}
\end{define}

\begin{theorem}
\begin{enumerate}
\item $ P_k(D) $ satisfies the Reidemeister moves and thus is an invariant of a link $ L $ which has planar projection $ D. $ 
\item The invariant above satisfies the skein relation for the one variable $ \mathfrak{sl}_k $ specialization of the HOMFLYPT polynomial:
$$ q^{k} P_k(D_{+}) - q^{-k}P_k(D_{-}) = (q^{1}-q^{-1})P_k(D_{0}). $$
\end{enumerate}
\end{theorem}

\begin{proof}
This is theorem 3.2 of [MOY].
\end{proof}
 
\section{Categorification of $ \mathfrak{sl}_k $ Modules}

The most important representations to categorify are the tensor powers $ V_{k-1}^{\otimes n}. $
When $ k=2, $ the categorifications were constructed in [BFK]. The authors worked on the most singular blocks of category $ \mathcal{O} $ which correspond to maximal subgroups of the symmetric group.  They also suggested
what to do for the $ \mathfrak{sl}_{k} $ case: consider other blocks of category $ \mathcal{O}. $

\subsection{Categorification of $ \Lambda^{i_1} V_{k-1} \otimes \cdots \otimes \Lambda^{i_r} V_{k-1} $ }

\begin{define}
\begin{enumerate}
\item Denote the Grothendieck group of an abelian or triangulated category $ \mathcal{C} $ by $ [\mathcal{C}]. $  It is the free abelian group generated by the symbols $ [M] $ where $ M $ is an object of $ \mathcal{C}. $  The only relations in this group are of the form $ [N] = [M]+[P] $ when there is a short exact sequence or distinguished triangle of the form 
$$ 0 \rightarrow M \rightarrow N \rightarrow P \rightarrow 0. $$
\item The image of an object $ M $ or an exact functor $ F $ in the Grothendieck group will be indicated by $ [M] $ and $ [F] $ respectively.
\item Given an abelian category $ \mathcal{C}, $ let $ D^b(\mathcal{C}) $ denote the corresponding bounded derived category.
\item Denote by $ H^j \colon D^b(\mathcal{C}) \rightarrow \mathcal{C}, $ the $ j \text{th} $ cohomology functor.  If $ F $ is a left exact functor, the $ j \text{th} $ derived functor 
$ R^j F $ is defined to be $ H^j \circ RF. $  If $ F $ is right exact, the $ j \text{th} $ derived functor $ L_j F $ is defined to be $ H^{-j} \circ LF. $
\end{enumerate}
\end{define}

The truncation functors $ \tau^{\leq n} $ and $ \tau^{\geq n} $ are well defined on the derived category. 

\begin{prop}
Let $ X $ be an object in $ D^b(\mathcal{C}). $ The following triangles exist in the derived category:
\begin{eqnarray*}
&\tau^{\leq n} X \rightarrow X \rightarrow \tau^{\geq n+1} X\\
&\tau^{\leq n-1} X \rightarrow \tau^{\leq n} X \rightarrow H^n(X)[-n]
\end{eqnarray*}
\end{prop}

\begin{proof}
See proposition 4.1.8 of [Sch].
\end{proof}

\begin{define}
Let $ \mathcal{O}(\mathfrak{gl}_n) $ be the category of $ \mathfrak{gl}_n $ modules which satisfy the following properties:
\begin{enumerate}
\item Finitely generated as $ \mathcal{U}(\mathfrak{gl}_n)-  $ modules.
\item Diagonalizable under the action of the Cartan subalgebra $ \mathfrak{h}. $
\item Locally finite under the action of the Borel subalgebra $ {\mathfrak{b}} = {\mathfrak{h}} + {\mathfrak{n}}^{+}. $
\end{enumerate}
\end{define}

This category decomposes into a direct sum of subcategories corresponding to the generalized central characters.

\begin{define}
\begin{enumerate}
\item Let 
$ \mathcal{O}_{(d_{k-1}, d_{k-2}, \ldots, d_{0})} = \mathcal{O}_{\bf d} $
be the block of $ \mathcal{O}(\mathfrak{gl}_n) $ for the central character corresponding to the weight 
$ \sum_{i=0}^{k-1} \sum_{j=1}^{d_i} i e_{d_0 + \cdots + d_{i-1} + j}^{*} - \rho, $
where 
$$ \rho = \frac{n-1}{2}e_1 + \frac{n-3}{2}e_2 + \cdots + \frac{1-n}{2}e_n $$
is half the sum of the positive roots.
\item Let $ M(a_1, \ldots, a_n) $ be the Verma module with highest weight $ a_1 e_1+ \cdots + a_n e_n - \rho. $
\item Let $ L(a_1, \ldots, a_n) $ be the Verma module with highest weight $ a_1 e_1+ \cdots + a_n e_n - \rho. $
\end{enumerate}
\end{define}

There are $ d_i $ terms in the weight from the definition above with coefficient $ i. $  Note that $ \sum_{j=0}^{k-1} d_j = n. $

\begin{prop}
Assume that the following direct sum is over all $ \bf{d} $ such that the entries are non-negative integers and the sum of the entries is n.  Then
$ \mathbb{C} \otimes_{\mathbb{Z}} [\oplus_{\bf{d}} \mathcal{O}_{{\bf d}}] \cong V_{k-1}^{\otimes n}. $
\end{prop}

\begin{proof}
The image of the Verma module $ [M(a_, \ldots, a_n)] $ gets mapped to $ v_{a_1} \otimes \cdots \otimes v_{a_n}. $
\end{proof}

This proposition is the first step towards categorification of $ \mathfrak{sl}_k- $ modules.
Next we would like to categorify the action of the Lie algebra.  The desired functors come directly from [BFK].  It is essentially the projective functor of
tensoring with the $ n-$ dimensional representation $ V_{n-1}. $  One only has to be careful about projecting onto the various blocks.
This is done next.

\begin{define}
\begin{enumerate}
	\item Let 
		$ \mathcal{E}_{i} \colon \mathcal{O}_{(d_{k-1},d_{k-2}, \ldots, d_{0})} \rightarrow
		\mathcal{O}_{(d_{k-1},\ldots, d_i +1, d_{i-1}-1, \ldots, d_{0})} $ be the functor defined by
		$$ \mathcal{E}_{i} M = \text{proj}_{(d_{k-1}, \ldots, d_i+1, d_{i-1}-1, \ldots, d_{0})}
		(V_{n-1} \otimes M). $$
	\item Let 
		$ \mathcal{F}_{i} \colon \mathcal{O}_{(d_{k-1}, d_{k-2}, \ldots, d_{0})} \rightarrow
		\mathcal{O}_{(d_{k-1}, \ldots, d_{i}-1, d_{i-1}+1, \ldots, d_{0})} $ by
		$$ \mathcal{F}_{i} M = \text{proj}_{(d_{k-1}, \ldots, d_{i}-1, d_{i-1}+1, \ldots, d_{0})}
	 ({V_{n-1}^{*}} \otimes M). $$
	\item Let 
		$ \mathcal{H}_{i} \colon \mathcal{O}_{(d_{k-1}, d_{k-2}, \ldots, d_{0})} \rightarrow
		\mathcal{O}_{(d_{k-1}, d_{k-2}, \ldots, d_{0})} $ be 
		$ \Id^{\oplus (d_i - d_{i-1})}. $
\end{enumerate}
\end{define}

\begin{theorem}
Functorial Serre relations are satisfied on $ \oplus_{\bf d} \mathcal{O}(\mathfrak{gl}_n). $
\begin{enumerate}
	\item $ \mathcal{H}_i \mathcal{H}_j \cong \mathcal{H}_j \mathcal{H}_i. $ 
	\item If $ d_i > d_{i-1} $ then $ \mathcal{E}_i \mathcal{F}_i \cong \mathcal{F}_i \mathcal{E}_i \oplus \mathcal{H}_i. $
	\item If $ d_i = d_{i-1} $ then $ \mathcal{E}_i \mathcal{F}_i \cong \mathcal{F}_i \mathcal{E}_i. $
	\item If $ d_i < d_{i-1} $ then $ \mathcal{F}_i \mathcal{E}_i \cong \mathcal{E}_i \mathcal{F}_i \oplus \mathcal{H}_i. $
	\item If $ i \neq j, $ then $ \mathcal{E}_i \mathcal{F}_j \cong \mathcal{F}_j \mathcal{E}_i. $
	\item If $ d_i - d_{i-1} \geq 0, $ then $ \mathcal{H}_i \mathcal{E}_i \cong \mathcal{E}_i \mathcal{H}_i \oplus \mathcal{E}_i^{\oplus 2}. $
	\item If $ d_i - d_{i-1} \leq -2, $ then $ \mathcal{E}_i \mathcal{H}_i \cong \mathcal{H}_i \mathcal{E}_i \oplus \mathcal{E}_i^{\oplus 2}. $
	\item If $ d_i - d_{i-1} = -1, $ then $ \mathcal{H}_i \mathcal{E}_i \cong \mathcal{E}_i \mathcal{H}_i. $
	\item If $ d_i - d_{i-1} \geq 2, $ then $ \mathcal{F}_i \mathcal{H}_i \cong \mathcal{H}_i \mathcal{F}_i \oplus \mathcal{F}_i^{\oplus 2}. $
	\item If $ d_i - d_{i-1} \leq 0, $ then $ \mathcal{H}_i \mathcal{F}_i \cong \mathcal{F}_i \mathcal{H}_i \oplus \mathcal{F}_i^{\oplus 2}. $
	\item If $ d_i - d_{i-1} = 1, $ then $ \mathcal{H}_i \mathcal{F}_i \cong \mathcal{F}_i \mathcal{H}_i. $
	\item If $ |i-j|>1, $ then $ \mathcal{H}_i \mathcal{E}_j \cong \mathcal{E}_j \mathcal{H}_i. $  
	\item If $ |i-j|>1, $ then $ \mathcal{H}_i \mathcal{F}_j \cong \mathcal{F}_j \mathcal{H}_i. $  
	\item If $ |i-j|>1, $ then $ \mathcal{E}_i \mathcal{E}_j \cong \mathcal{E}_j \mathcal{E}_i. $  
	\item If $ |i-j|>1, $ then $ \mathcal{F}_i \mathcal{F}_j \cong \mathcal{F}_j \mathcal{F}_i. $  
	\item If $ j=i+1, $ then $ \mathcal{E}_i \mathcal{E}_i \mathcal{E}_j \oplus \mathcal{E}_j \mathcal{E}_i \mathcal{E}_i \cong
		\mathcal{E}_i \mathcal{E}_j \mathcal{E}_i \oplus \mathcal{E}_i \mathcal{E}_j \mathcal{E}_i. $
	\item If $ j=i+1, $ then $ \mathcal{F}_i \mathcal{F}_i \mathcal{F}_j \oplus \mathcal{F}_j \mathcal{F}_i \mathcal{F}_i \cong
		\mathcal{F}_i \mathcal{F}_j \mathcal{F}_i \oplus \mathcal{F}_i \mathcal{F}_j \mathcal{F}_i. $
	\item If $ d_i - d_{i-1}-1 \geq 0, $ then 
		$ \mathcal{H}_i \mathcal{E}_{i+1} \oplus \mathcal{E}_{i+1} \cong \mathcal{E}_{i+1} \mathcal{H}_i. $
	\item If $ d_i - d_{i-1} \leq 0, $ then 
		$ \mathcal{E}_{i+1} \mathcal{H}_{i} \oplus \mathcal{E}_{i+1} \cong \mathcal{H}_{i} \mathcal{E}_{i+1}. $
	\item If $ d_i - d_{i-1} -1  \geq 0, $ then 
		$ \mathcal{H}_{i} \mathcal{E}_{i-1} \oplus \mathcal{E}_{i-1} \cong \mathcal{E}_{i-1} \mathcal{H}_{i}. $
	\item If $ d_i - d_{i-1} \leq 0, $ then 
		$ \mathcal{E}_{i-1} \mathcal{H}_{i} \oplus \mathcal{E}_{i-1} \cong \mathcal{H}_{i} \mathcal{E}_{i-1}. $
	\item If $ d_i - d_{i-1} \geq 0, $ then 
		$ \mathcal{F}_{i+1} \mathcal{H}_{i} \oplus \mathcal{F}_{i+1} \cong \mathcal{H}_{i} \mathcal{F}_{i+1}. $
	\item If $ d_i - d_{i-1}+1 \geq 0, $ then 
		$ \mathcal{H}_i \mathcal{F}_{i+1} \oplus \mathcal{F}_{i+1} \cong \mathcal{F}_{i+1} \mathcal{H}_i. $
	\item If $ d_i - d_{i-1} \geq 0, $ then 
		$ \mathcal{F}_{i-1} \mathcal{H}_{i} \oplus \mathcal{F}_{i-1} \cong \mathcal{H}_{i} \mathcal{F}_{i-1}. $
	\item If $ d_i - d_{i-1} +1  \geq 0, $ then 
		$ \mathcal{H}_{i} \mathcal{F}_{i-1} \oplus \mathcal{F}_{i-1} \cong \mathcal{F}_{i-1} \mathcal{H}_{i}. $
\end{enumerate}
\end{theorem}

One only has to compute in the Grothendieck group to check these relations [BG].  Details will be provided later when a functorial version of the quantum Serre relations is proved.

Our next goal is to categorify tensor products of exterior powers of $ V_{k-1}. $  
The case $ k=2 $ was addressed in [BFK].  When $ k=2, $ $ \Lambda^{2}V_{k-1} $
is just a one dimensional vector space.  The main ingredient in the categorification of
$ V_{1}^{\otimes n} \rightarrow V_{1}^{\otimes (n-2)} $ is the Zuckerman functor and parabolic subcategories.  Thus for general
$ k, $ we expect a locally finite subcategory to categorify 
a tensor product of fundamental representations.
We now give a definition of various parabolic subalgebras of $ \mathfrak{sl}_k $ and the corresponding locally finite subcategories first introduced by A. Rocha-Caridi in [Ro].

\begin{define}
\begin{enumerate}
\item The subalgebra $ \mathfrak{p}_{(r_1, \ldots, r_t)} $ is the parabolic subalgebra whose reductive subalgebra is
$ \mathfrak{gl}_{r_1} \oplus \cdots \oplus \mathfrak{gl}_{r_t}, $  where $ r_1 + \cdots + r_t = n. $
\item Denote by $ \mathcal{O}_{\bf d}^{\mathfrak{p}} $ the full subcategory of $ \mathcal{O}_{\bf d} $ of modules locally finite with respect to the subalgebra
$ \mathfrak{p}. $ By abuse of notation, let
$ \mathcal{O}_{\bf d}^{(r_1, \ldots, r_t)} $ be the category $ \mathcal{O}_{\bf d}^{\mathfrak{p}_{(r_1, \ldots, r_t)}}. $
\end{enumerate}
\end{define}

It is important to know some objects which are in these locally finite categories.  
\begin{define}
\begin{enumerate}
\item Let S denote the subset of simple roots defining the parabolic subalgebra $ \mathfrak{p}. $
\item Let $ P_{\mathfrak{p}}^{+} = \lbrace \lambda \in \mathfrak{h}^{*} | \langle \lambda, \alpha \rangle \in \mathbb{N}, \forall \alpha \in S \rbrace. $
\end{enumerate}
\end{define}

Given such a $ \lambda \in P_{\mathfrak{p}}^{+}, $ we may define the generalized Verma module
$ M^{\mathfrak{p}}(\lambda) = \mathcal{U}(\mathfrak{g}) \otimes_{\mathcal{U}(\mathfrak{p})} E(\lambda), $ where $ E(\lambda) $
is the simple $ \mathfrak{p}- $ module with highest weight $ \lambda. $

Now it is easy to give conditions for which modules in these singular blocks are locally finite.
Let $ \lambda = a_1 e_{1}^{*} + \cdots + a_n e_{n}^{*}. $  We need to give a condition for
$ \lambda $ to be in this set.  Since $ \langle \rho, h_i \rangle =1, $ the condition that
$ \langle \lambda - \rho, h_i \rangle \geq 0 $ simply becomes $ a_i > a_{i+1}. $
We now could state which simple modules and generalized Verma modules are in the locally finite
subcategories.

\begin{lemma}
\begin{enumerate}
	\item $ L(a_1, \ldots, a_n) $ is a simple module in $ \mathcal{O}_{\bf d}^{\mathfrak{p}} $ if $ a_i > a_{i+1} $ whenever $ \alpha_i $ is a simple root of $ \mathfrak{p}. $
	\item $ M^{\mathfrak{p}}(a_1, \ldots, a_n) $ is a generalized Verma module in $ \mathcal{O}_{\bf d}^{\mathfrak{p}} $ if $ a_i > a_{i+1} $ whenever $ \alpha_i $ is a simple root of $ \mathfrak{p}. $
\end{enumerate}
\end{lemma}

\begin{proof}
This follows directly from above by evaluating at all simple roots $ \alpha_i $ which define the parabolic subalgebra $ \mathfrak{p}. $
\end{proof}

We will often group coefficients in the highest weight of a generalized module such as
$$ M^{\mathfrak{p}}(a_1, \ldots, a_{i-1}, \underbrace{a_i, \ldots, a_r}, \ldots, a_n) $$ 
to stress that $ a_i > \cdots > a_r $ so that it is locally finite with respect to a certain subalgebra $ \mathfrak{p}. $

\begin{prop}
$ \mathbb{C} \otimes_{\mathbb{Z}} [\oplus_{\bf d} \mathcal{O}_{\bf d}^{(r_1, \ldots, r_t)}] \cong \Lambda^{r_1} V_{k-1} \otimes \cdots \otimes \Lambda^{r_t} V_{k-1}. $
\end{prop}

\begin{proof}
It suffices to show that this is an isomorphism on the basis of generalized Verma modules.
Let $ \mathfrak{p} $ be the subalgebra given above.
The isomorphism sends $ [M^{\mathfrak{p}}(a_1, \ldots, a_n)] $ to
$$ (v_{a_1} \wedge \cdots \wedge v_{r_1}) \otimes \cdots \otimes
(v_{a_{r_1 + \cdots + r_{t-1} + 1}} \wedge \cdots \wedge v_{a_{r_1+ \cdots + r_t}}). $$
This is clearly a bijection.
\end{proof}

\begin{remark}
If any of the $ r_i $ above is larger than $ k, $ then the category contains no non-trivial objects.
\end{remark}

The Zuckerman functor plays an obvious role in this setup.  It categorifies projection maps
of these modules.

\begin{define}
\begin{enumerate}
	\item Let the Zuckerman functor $ \Gamma^{(r_1, \ldots, r_t)} \colon \mathcal{O}_{\bf d} 
		\rightarrow \mathcal{O}_{\bf d}^{(r_1, \ldots, r_t)} $ be the functor of taking the maximal locally $ \mathfrak{p}_{(r_1, \ldots, r_t)} $ finite submodule.
	\item Let the dual Zuckerman $ Z^{(r_1, \ldots, r_t)} \colon \mathcal{O}_{\bf d} 
		\rightarrow \mathcal{O}_{\bf d}^{(r_1, \ldots, r_t)} $ be the functor of taking the maximal locally $ \mathfrak{p}_{(r_1, \ldots, r_t)} $ finite quotient.
	\item Let $ \epsilon_{(r_1, \ldots, r_t)} \colon \mathcal{O}_{\bf d}^{(r_1, \ldots, r_t)} \rightarrow \mathcal{O}_{\bf d} $ be the natural inclusion functor.
\end{enumerate}
\end{define}

\begin{remark}
The dual Zuckerman functor is also known as the Bernstein functor.  See [KV].
\end{remark}

If there is an inclusion of parabolic subalgebras $ \mathfrak{q} \subset \mathfrak{p}, $ there is an obvious generalization of these definitions for categories locally finite with respect to these algebras.
For example, one may take a module locally finite with respect to $ \mathfrak{p} $ and apply the Zuckerman functor $ \Gamma_{\mathfrak{q}}^{\mathfrak{p}}. $

The Zuckerman functor is left exact.  One usually studies its right derived functor and its cohomology functors.  Taking the
derived functor is important in categorification so that it becomes exact as a functor on the derived
category.  We denote its right derived functor shifted by $ j $  by $ R\Gamma[j]. $ 
The dual Zuckerman functor is right exact and similarly, one should consider its left derived functor.
On the Grothendieck group, $ LZ $ and $ R\Gamma $ descend to the canonical map from
$ V_{k-1}^{\otimes n} $ to $ \Lambda^{r_1} V_{k-1} \otimes \cdots \otimes \Lambda^{r_t} V_{k-1}. $

\begin{prop}
\begin{enumerate}
	\item $ [LZ^{(r_1, \ldots, r_t)}] = \pi_{r_1} \otimes \cdots \otimes \pi_{r_t}. $
	\item $ [\epsilon_{(r_1, \ldots, r_t)}[-\Sigma_{r=1}^t \frac{i_r(i_r-1)}{2}]] = i_{r_1} \otimes \cdots \otimes i_{r_t}. $
\end{enumerate}
\end{prop}

\begin{proof}
By proposition 5.5 of [ES], one may easily prove the first part by computing on the basis of Verma modules.
The second part follows from the generalized BGG resolution.
\end{proof}

Let $ d $ be the codimension of $ \mathfrak{q} $ in $ \mathfrak{p}. $
\begin{lemma}
The derived Zuckerman functor $ R\Gamma_{\mathfrak{q}}^{\mathfrak{p}} $ and the inclusion functor $ \epsilon_{\mathfrak{q}}^{\mathfrak{p}} $ satisfy the following adjointness properties in the derived category:
\begin{enumerate}
	\item $ \Hom(\epsilon_{\mathfrak{p}}^{\mathfrak{q}} X, Y) \cong \Hom(X, R\Gamma_{\mathfrak{q}}^{\mathfrak{p}} Y) $
	\item $ \Hom(X, \epsilon_{\mathfrak{p}}^{\mathfrak{q}}[2d] Y) \cong \Hom(R\Gamma_{\mathfrak{q}}^{\mathfrak{p}} X, Y). $
\end{enumerate}
\end{lemma}

\begin{proof}
See the proof of theorem 5 in [BFK].
\end{proof}

\begin{lemma}
\label{adjoint}
The derived dual Zuckerman functor $ LZ_q^p $ and the inclusion functor $ \epsilon_q^p $ satisfy the following adjointness properties in the derived category:
\begin{enumerate}
	\item $ \Hom(\epsilon_{\mathfrak{p}}^{\mathfrak{q}}[-2d] X, Y) \cong \Hom(X, LZ_{\mathfrak{q}}^{\mathfrak{p}} Y) $
	\item $ \Hom(X, \epsilon_{\mathfrak{p}}^{\mathfrak{q}} Y) \cong \Hom(LZ_{\mathfrak{q}}^{\mathfrak{p}} X, Y). $
\end{enumerate}
\end{lemma}

\begin{corollary}
The derived Zuckerman and dual Zuckerman functors are related by 
$ R\Gamma_{\mathfrak{q}}^{\mathfrak{p}}[2d] \cong LZ_{\mathfrak{q}}^{\mathfrak{p}}. $
\end{corollary}

\begin{proof}
By the previous two lemmas these functors are adjoints of the same functor so they are the same up to isomorphism.
\end{proof}

\subsection{Equivalences of Categories}
\label{sec3.2}

We look to generalize section 3.2.2 of [BFK] to the subcategories with generalized central character considered
here.  
In order to make the notation more compact, we will have the following notation for several important subalgebras.

\begin{define}
Let
\begin{enumerate}
\item $ \mathfrak{p}_i $ to be the parabolic subalgebra corresponding to 
$ (\underbrace{1, \ldots, 1,}_{i-1} k, 1, \ldots, 1). $
\item $ \mathfrak{q}_i $ to be the parabolic subalgebra corresponding to
$ (\underbrace{1, \ldots, 1,}_{i} k-1, 1, \ldots, 1). $
\item $ \mathfrak{r}_i $ be the parabolic subalgebra corresponding to $ (\underbrace{1, \ldots, 1}_{i-1}, k-2, 1, \ldots 1). $
\item $ \mathfrak{s}_i $ be the parabolic subalgebra corresponding to $ (\underbrace{1, \ldots, 1}_{i-1}, 2, 1, \ldots 1). $
\item $ \mathfrak{t}_i $ be the parabolic subalgebra corresponding to $ (\underbrace{1, \ldots, 1}_{i-1}, 3, 1, \ldots 1). $
\end{enumerate}
\end{define}

Clearly $ \mathfrak{p}_j \supset \mathfrak{q}_j $ and $ \mathfrak{p}_j \supset \mathfrak{q}_{j-1}. $  Recall from earlier the definitions of the functors $ \epsilon_{\mathfrak{p}_j}^{\mathfrak{q}_j}, $ $ \epsilon_{\mathfrak{p}_j}^{\mathfrak{q}_{j-1}}, $ $ LZ_{\mathfrak{q}_j}^{\mathfrak{p}_j}, $ and
$ LZ_{\mathfrak{q}_{j-1}}^{\mathfrak{p}_j}. $
Note that the derived Zuckerman functors in this setup could have non-zero cohomology functors only in degrees
0 through 2(k-1).  As in [BFK], the middle cohomology functor plays a significant role.  Given a generalized
Verma module $ M^{\mathfrak{p}_j}(\alpha), $ we would like to compute its image under the functor 
$ L_{(k-1)}Z_{\mathfrak{q}_j}^{\mathfrak{p}_{j+1}} \circ \epsilon_{\mathfrak{p}_j}^{\mathfrak{q}_j}. $

The main result of this subsection is there are equivalences of categories:
$$ \mathcal{O}_{(d_{k-1}, \ldots, d_{0})}^{\mathfrak{p}_j} \cong \mathcal{O}_{(d_{k-1}, \ldots, d_{0})}^{\mathfrak{p}_{j+1}}. $$
The equivalences are compositions of inclusions into larger categories and the middle cohomology of the derived
Zuckerman functor.  
The plan of the proof is exactly that from [BFK].  First the action of these functors on generalized Verma 
modules is computed.  Then after some exactness and adjointness statements are proved, the equivalence will
follow from lemma 2 of [BFK].  

First we recall the generalized BGG resolution. Let $ S $ be a subset of simple roots defining a parabolic subalgebra $ \mathfrak{p} $ and $ W_{\mathfrak{p}} $ the corresponding Weyl group.  Denote by $ W^{\mathfrak{p}} $ the set of shortest coset representatives in $ W/W_{\mathfrak{p}} $ and $ \rho_S $ half the sum of the positive roots of the subalgebra $ \mathfrak{p}. $ 

\begin{theorem}
Let $ \mu $ be dominant integral.  For all $ j = 0, \ldots, \dim u, $ define
$$ C_{j}^{S} = \oplus_{w \in W^{\mathfrak{p}}, l(w)=j} M^{\mathfrak{p}}(w(\mu + \rho) - \rho_S). $$
Then there is an exact sequence
$$ 0 \rightarrow C_{dim u}^{S} \rightarrow \cdots \rightarrow C_{0}^{S} \rightarrow L(\mu + \rho) \rightarrow 0  $$
with all maps nontrivial.
\end{theorem}

Let $ \alpha = (a_1, \ldots, a_n).  $  In order for the generalized Verma module 
$ M^{\mathfrak{p}_j}(\alpha) $ to be locally $ \mathcal{U}(\mathfrak{p}_j)- $ finite, 
$$ a_j = k-1, a_{j+1} = k-2, \ldots, a_{j+k-1} = 0. $$
From the previous theorem, it follows that there is a resolution of this generalized Verma module in terms of generalized Verma modules
in larger parabolic categories.

\begin{corollary}
There is an exact sequence in $ \mathcal{O}_{\bf{d}}^{q_j} $ of the form:
$$ 0 \rightarrow M^{\mathfrak{q}_j}(\sigma_1 \cdots \sigma_{k-1}. \alpha) \rightarrow M^{\mathfrak{q}_j}(\sigma_1 \cdots \sigma_{k-2}. \alpha) 
\rightarrow \cdots \rightarrow M^{\mathfrak{q}_j}(\alpha) \rightarrow M^{\mathfrak{p}_j}(\alpha) \rightarrow 0,  $$
where $ \sigma_i $ permutes the elements $ a_{j+i-1} $ and $ a_{j+i} $ of $ (a_j, \ldots, a_{j+k-1}). $
\end{corollary}

\begin{proof}
This follows from the theorem and parabolic induction.
\end{proof}

The cohomology functors on the modules $ M^{\mathfrak{q}_j}(\sigma_1 \cdots \sigma_{l}. \alpha) $ play a critical role.

\begin{lemma}
Suppose $ a_{k+j} = k-l-1. $ Then,
$ L_{i}Z_{\mathfrak{q}_j}^{\mathfrak{p}_{j+1}} M^{\mathfrak{q}_j}(\sigma_1 \cdots \sigma_{l}. \alpha) \cong $
\begin{eqnarray*} 
M^{\mathfrak{p}_{j+1}}(a_1, \ldots, a_{j-1}, k-l-1, k-1, \ldots, 0, a_{k+j-1}, \ldots, a_n) & \textrm{ if } & i=k-l-1 \\
0 & \textrm{ if } & i \neq k-l-1.
\end{eqnarray*}
\end{lemma}

\begin{proof}
It is easily seen that 
$$ M^{\mathfrak{q}_j}(\sigma_1 \cdots \sigma_{l}. \alpha) = 
M^{\mathfrak{q}_j}(a_1, \ldots, a_{j-1}, k-l-1, k-1, k-2, \ldots, k-l, k-l-2, \ldots, 0, a_{k+j}, \ldots, a_{n}). $$
Thus $ L_{i}Z_{\mathfrak{q}_j}^{\mathfrak{p}_{j+1}} M^{\mathfrak{q}_j}(\sigma_1 \cdots \sigma_{l}. \alpha) = 0 $ if $ a_{k+j} \neq k-l-1. $

Now suppose $ S_k $ permutes the elements $ (k-1, k-2, \ldots, k-l, k-l-2, \ldots, 0, a_{k+j}), $ by
\begin{eqnarray*}
&\sigma_{l+1} \cdots \sigma_{k-1}. (k-1, k-2, \ldots, k-l, k-l-2, \ldots, 0, a_{k+j}) =\\
&(k-1, k-2, \ldots, k-l, a_{k+j}, k-l-2, \ldots, 0). 
\end{eqnarray*} 
The length of this element in the symmetric group is $ k-l-1. $
Thus from [ES] proposition 5.5, 
$$ L_{i}Z_{\mathfrak{q}_j}^{\mathfrak{p}_{j+1}} M^{\mathfrak{q}_j}(\sigma_1 \cdots \sigma_{l}. \alpha) \cong
M^{\mathfrak{p}_{j+1}}(a_1, \ldots, a_{j-1}, k-l-1, k-1, \ldots, 0, a_{k+j-1}, \ldots, a_n),  $$
if $ i=(k-l-1). $
\end{proof}

Now we are in position to compute $ L_iZ_{\mathfrak{q}_j}^{\mathfrak{p}_{j+1}} \circ \epsilon_{\mathfrak{p}_j}^{\mathfrak{q}_j} M^{\mathfrak{p}_j}(\alpha). $
Let 
$$ \beta = (a_1, \ldots, a_{j-1}, a_{j+k}, a_j, \ldots, a_{j+k-1}, a_{j+k+1}, \ldots, a_n). $$

\begin{lemma}
There are isomorphisms: 
$ L_{(k-1)}Z_{\mathfrak{q}_j}^{\mathfrak{p}_{j+1}} \circ \epsilon_{\mathfrak{p}_j}^{\mathfrak{q}_j} M^{\mathfrak{p}_j}(\alpha) \cong $
\begin{eqnarray*} 
M^{\mathfrak{p}_{j+1}}(\beta) & \textrm{ if } & i=k-1 \\
0 & \textrm{ if } & i \neq k-1
\end{eqnarray*}
\end{lemma}

\begin{proof}
Consider the generalized BGG resolution.  This gives rise to short exact sequences:
\begin{eqnarray*} 
&0 \rightarrow K_0 \rightarrow M^{\mathfrak{q}_j}(\alpha) \rightarrow M^{\mathfrak{p}_j}(\alpha) \rightarrow 0\\
&0 \rightarrow K_1 \rightarrow M^{\mathfrak{q}_j}(\sigma_1. \alpha) \rightarrow K_0 \rightarrow 0\\
&0 \rightarrow K_2 \rightarrow M^{\mathfrak{q}_j}(\sigma_1 \sigma_2. \alpha) \rightarrow K_1 \rightarrow 0\\
&\cdots\\
&0 \rightarrow K_{k-2} \rightarrow M^{\mathfrak{q}_j}(\sigma_1 \cdots \sigma_{k-2}. \alpha) \rightarrow K_{k-3} \rightarrow 0\\
&0 \rightarrow K_{k-1} \rightarrow M^{\mathfrak{q}_j}(\sigma_1 \cdots \sigma_{k-2} \sigma_{k-1}. \alpha) \rightarrow K_{k-2} \rightarrow 0. 
\end{eqnarray*} 

Assume $ a_{k+j} = k-l-1. $  
By the previous lemma, for all $ s $ 
$$ L_s Z_{\mathfrak{q}_j}^{\mathfrak{p}_{j+1}} M^{\mathfrak{q}_j}(e) \cong L_s Z_{\mathfrak{q}_j}^{\mathfrak{p}_{j+1}} M^{\mathfrak{q}_j}(\sigma_1. \alpha) \cong \cdots \cong L_s Z_{\mathfrak{q}_j}^{\mathfrak{p}_{j+1}}M^{\mathfrak{q}_j}(\sigma_1 \ldots \sigma_{l-1}. \alpha)= 0. $$
Therefore
$$ L_s Z_{\mathfrak{q}_j}^{\mathfrak{p}_{j+1}} M^{\mathfrak{p}_j}(\alpha) \cong L_{s-1} Z_{\mathfrak{q}_j}^{\mathfrak{p}_{j+1}} K_0 \cong \cdots \cong L_{s-l} Z_{\mathfrak{q}_j}^{\mathfrak{p}_{j+1}} K_{l-1}. $$

Also by the previous lemma, for all $ s $,

$$ L_s Z_{\mathfrak{q}_j}^{\mathfrak{p}_{j+1}} M^{\mathfrak{q}_j}(\sigma_{l+1}. \alpha) \cong \cdots \cong L_s Z_{\mathfrak{q}_j}^{\mathfrak{p}_{j+1}} M^{\mathfrak{q}_j}(\sigma_1 \ldots \sigma_{k-1}. \alpha) \cong 0. $$
Therefore for all $ s, $ 
$$ L_{s} Z_{\mathfrak{q}_j}^{\mathfrak{p}_{j+1}} K_{k-2} \cong \cdots \cong L_{s+l} Z_{\mathfrak{q}_j}^{\mathfrak{p}_{j+1}} K_{l} \cong 0. $$

Now consider the short exact sequence
$$ 0 \rightarrow K_l \rightarrow M^{\mathfrak{q}_j}(\sigma_1 \ldots \sigma_l. \alpha) \rightarrow K_{l-1} \rightarrow 0. $$
Due to the above
$$ L_{s} Z_{\mathfrak{q}_j}^{\mathfrak{p}_{j+1}} M^{\mathfrak{q}_j}(\sigma_1 \ldots \sigma_l. \alpha) \cong L_{s} Z_{\mathfrak{q}_j}^{\mathfrak{p}_{j+1}} K_{l-1} \cong L_{s+l} Z_{\mathfrak{q}_j}^{\mathfrak{p}_{j+1}} \epsilon_{\mathfrak{p}_j}^{\mathfrak{q}_j} M^{\mathfrak{p}_j}(\alpha). $$

By the previous lemma 
$$ L_{s} Z_{\mathfrak{q}_j}^{\mathfrak{p}_{j+1}} M^{\mathfrak{q}_j}(\sigma_1 \ldots \sigma_l. \alpha) \cong M^{\mathfrak{p}_{j+1}}(\beta) $$ 
for $ s=k-l-1, $
and it is zero otherwise.

Thus $ L_{s} Z_{\mathfrak{q}_j}^{\mathfrak{p}_{j+1}} M^{\mathfrak{p}_j}(\alpha) \cong M^{\mathfrak{p}_{j+1}}(\beta) $ if $ s = k-1 $ and is zero otherwise.
\end{proof}

\begin{lemma}
Suppose $ i \neq k-1. $ Then for any module $ M $ in $ \mathcal{O}_{(d_{k-1},d_{k-2}, \ldots, d_{0})}^{\mathfrak{p}_j}, $ there is an isomorphism
$$ L_{i}Z_{\mathfrak{q}_j}^{\mathfrak{p}_{j+1}} \circ \epsilon_{\mathfrak{p}_j}^{\mathfrak{q}_j} M =0. $$ 
\end{lemma}

\begin{proof}
We know this is true by the previous lemma for generalized Verma modules.  It is then true for projective modules by induction on the length of their Verma
flag.
Now consider a simple module $ S. $  It has a projective resolution which gives rise to several short exact sequences:
$$ 0 \rightarrow K_1 \rightarrow P_1 \rightarrow S \rightarrow 0 $$
$$ 0 \rightarrow K_2 \rightarrow P_2 \rightarrow K_1 \rightarrow 0 $$
$$ \cdots $$
$$ 0 \rightarrow K_r \rightarrow P_r \rightarrow K_{r-1} \rightarrow 0. $$
The first short exact sequence implies 
$ L_{i}Z_{\mathfrak{q}_j}^{\mathfrak{p}_{j+1}} \circ \epsilon_{\mathfrak{p}_j}^{\mathfrak{q}_j} S \cong L_{(i-1)}Z_{\mathfrak{q}_j}^{\mathfrak{p}_{j+1}} \circ \epsilon_{\mathfrak{p}_j}^{\mathfrak{q}_j} K_1 $ for
$ i<(k-1). $
Continuing in this way we get 
$ L_{i}Z_{\mathfrak{q}_j}^{\mathfrak{p}_{j+1}} \circ \epsilon_{\mathfrak{p}_j}^{\mathfrak{q}_j} S \cong L_{i-r}Z_{\mathfrak{q}_j}^{\mathfrak{p}_{j+1}} \circ \epsilon_{\mathfrak{p}_j}^{q_j} P_r $ for $ i<(k-1). $
Since $ K_r \cong 0, $ 
$ L_{i}Z_{\mathfrak{q}_j}^{\mathfrak{p}_{j+1}} \circ \epsilon_{\mathfrak{p}_j}^{\mathfrak{q}_j} S \cong 0, $ for $ i<(k-1). $
By the duality theorem for Zuckerman functors, this isomorphism is true for all $ i \neq (k-1). $
Now the lemma follows for any $ M $ by induction on the length of its composition series. 
\end{proof}

\begin{prop}
There is an equivalence of categories 
$$ L_{k-1} Z_{\mathfrak{q}_j}^{\mathfrak{p}_{j+1}} \circ \epsilon_{\mathfrak{p}_j}^{\mathfrak{q}_j} \colon \mathcal{O}_{\bf d}^{\mathfrak{p}_j} \overset \cong \to \mathcal{O}_{\bf d}^{\mathfrak{p}_{j+1}}. $$
The inverse equivalence is given by
$ L_{k-1} Z_{\mathfrak{q}_j}^{\mathfrak{p}_j} \circ \epsilon_{\mathfrak{p}_{j+1}}^{\mathfrak{q}_j}. $
\end{prop}

\begin{proof}
These functors send generalized Verma modules to generalized Verma modules.  By examining the long exact sequence for the functor and using that the cohomology functors vanish in all but one degree, it is clear that this functor is exact exact. It suffices by lemma 2 of [BFK] to show that they are adjoint.
In the derived category $ LZ_{\mathfrak{q}_j}^{\mathfrak{p}_{j+1}} $ is left adjoint to $ \epsilon_{\mathfrak{p}_{j+1}}^{\mathfrak{q}_j} $ by lemma~\ref{adjoint}. Furthermore,
$ LZ_{\mathfrak{q}_j}^{\mathfrak{p}_{j+1}}[-2(k-1)]  $ is right adjoint to $ \epsilon_{\mathfrak{p}_{j+1}}^{\mathfrak{q}_j}. $
Thus we have 
\begin{eqnarray*}
\Hom_{\mathcal{O}} (L_{k-1}Z_{\mathfrak{q}_j}^{\mathfrak{p}_{j+1}} \circ \epsilon_{\mathfrak{p}_j}^{\mathfrak{q}_j}M,N) &\cong
& \Hom_{D^b(\mathcal{O})}(LZ_{\mathfrak{q}_j}^{\mathfrak{p}_{j+1}} \circ \epsilon_{\mathfrak{p}_j}^{\mathfrak{q}_j}[-(k-1)] M, N)\\
&\cong &\Hom_{D^b(\mathcal{O})}(\epsilon_{\mathfrak{p}_j}^{\mathfrak{q}_j} M, \epsilon_{\mathfrak{p}_{j+1}}^{\mathfrak{q}_j}[k-1] N)\\ 
&\cong &\Hom_{D^b(\mathcal{O})}(M, LZ_{\mathfrak{q}_j}^{\mathfrak{p}_j} \epsilon_{\mathfrak{p}_{j+1}}^{\mathfrak{q}_j}[-(k-1)] N)\\
& \cong &\Hom_{\mathcal{O}}(M, L_{k-1}Z_{\mathfrak{q}_j}^{\mathfrak{p}_j} \epsilon_{\mathfrak{p}_{j+1}}^{\mathfrak{q}_j} N). 
\end{eqnarray*} 
\end{proof}

\begin{remark}
The functor $ LZ_{\mathfrak{q}_j}^{\mathfrak{p}_{j+1}} \circ \epsilon_{\mathfrak{p}_j}^{\mathfrak{q}_j} $ gives an equivalence of the corresponding derived categories.
\end{remark}

Next assume that all of the $ a_i $ are 1 or k-1.

\begin{prop}
There is an equivalence of categories
$$ \mathcal{O}_{\bf d}^{(a_1, \ldots a_i, k, a_{i+1}, \ldots a_n)}(\mathfrak{gl}_n) \cong \mathcal{O}_{\bf d}^{(a_1, \ldots a_{i+1}, k, a_{i+2}, \ldots a_n)}(\mathfrak{gl}_n) $$ 
given by the functor
$$ L_{(k-1)} Z_{(a_1, \ldots a_i, a_{i+1}, k-a_{i+1}, a_{i+1} \ldots a_n)}^{(a_1, \ldots a_{i+1}, k, \ldots a_n)} \epsilon_{(a_1, \ldots a_i, k, a_{i+1}, \ldots a_n)}^{(a_1, \ldots a_i, a_{i+1}, k-a_{i+1}, a_{i+1} \ldots a_n)}. $$
\end{prop}

\begin{proof}
The proof is exactly the same as the previous proposition.  
\end{proof}

Composing the equivalences given by the proposition above, we get the following equivalence.
$$ \Delta_{(a_1, \ldots a_i, k, \ldots a_r)}^{(a_1, \ldots a_j, k, \ldots a_r)} \colon \mathcal{O}_{\bf d}^{(a_1, a_i, k, \ldots a_r)} \rightarrow  \mathcal{O}_{\bf d}^{(a_1, \ldots a_j, k, \ldots a_r)}. $$

Next we define functors between parabolic categories for Lie algebras of different ranks.

\begin{define}
\begin{enumerate}
\item Let $ Y $ be the one dimensional $ \mathfrak{gl}_k \oplus \mathfrak{gl}_{n}- $ module with weight
$$ -\frac{n}{2}(e_1 + \cdots + e_k)+\frac{k}{2}(e_{k+1}+\cdots + e_{k+n}). $$
\item Let $ Y' $ be the one dimensional $ \mathfrak{gl}_n- $ module with weight $ -\frac{k}{2}(e_1 + \cdots + e_n). $
\item Let $ W $ be the finite dimensional irreducible $ \mathfrak{gl}_k- $ module with highest weight $ (k-1)e_1 + \cdots + (0)e_k-\rho_k. $
\item The functor $ \zeta \colon \mathcal{O}_{\bf d}^{(a_1, \ldots, a_r)} (\mathfrak{gl}_n) \rightarrow \mathcal{O}_{\bf d'}^{(k, a_1, \ldots a_r)}(\mathfrak{gl}_{k+n}) $ is given by 
$$ \zeta(M) = \text{Ind}_{\mathcal{U}(\mathfrak{gl}_k + \mathfrak{gl}_n)}^{\mathcal{U}(\mathfrak{gl}_{k+n})}(Y \otimes W \otimes M) $$
where $ {\bf d'} = (d_{k-1}+1, \ldots, d_0+1). $
\item The functor $ \nu \colon \mathcal{O}_{\bf d'}^{(k, a_1, \ldots a_r)}(\mathfrak{gl}_{k+n}) \rightarrow \mathcal{O}_{\bf d}^{(a_1, \ldots a_r)}(\mathfrak{gl}_n) $ is given by taking the sum of the weight spaces of a module M for the
weights of the form 
$$ (k-1) e_1 + (k-2)e_2 + \cdots + (0)e_k + x_{k+1} e_{k+1} + \cdots x_{k+n} e_{k+n} - \rho, $$ 
tensor with $ Y', $ where the $ x_i \in \mathbb{Z}. $
\end{enumerate}
\end{define}

The next lemma is a direct generalization of proposition 17 of [BFK].

\begin{lemma}
The functors $ \zeta $ and $ \nu $ are inverse equivalences of categories.
\end{lemma}

Next we assign functors to various oriented tangles.  If an edge is labeled by $ k-1, $ then we give it a down arrow.  If an edge is labeled by $ 1, $ then it is given an up arrow.

\begin{figure}[htb]
  \centering
  \includegraphics{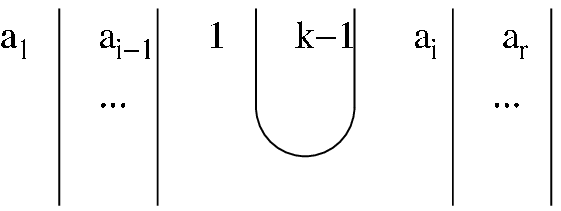}
  \caption{}
  \label{cup+}
\end{figure}

To the tangle in Figure ~\ref{cup+} we assign the functor 
$ \cup_{i,+, r} \colon D^b(\mathcal{O}_{\bf d}^{(a_1, \ldots, a_r)}) \rightarrow D^b(\mathcal{O}_{\bf d'}^{(a_1, \ldots, a_{i-1}, 1, k-1, a_i, \ldots, a_r)}). $  It is given by
$$ \cup_{i, +, r}(M) = \epsilon_{(a_1, \ldots, a_{i-1}, k, \ldots, a_r}^{(a_1, \ldots, a_{i-1}, 1, k-1, \ldots, a_r)}[-(k-1)] 
\Delta_{(k, a_1, \ldots, a_r)}^{(a_1, \ldots, a_{i-1}, k, a_r)} \zeta(M). $$

\begin{figure}[htb]
  \centering
  \includegraphics{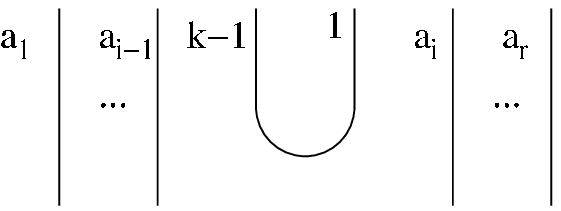}
  \caption{}
  \label{cup-}
\end{figure}

To the tangle in Figure ~\ref{cup-} we assign the functor 
$ \cup_{i,-, r} \colon D^b(\mathcal{O}_{\bf d}^{(a_1, \ldots, a_r)}) \rightarrow D^b(\mathcal{O}_{\bf d'}^{(a_1, \ldots, a_{i-1}, k-1, 1, a_i, \ldots, a_r)}). $  It is given by
$$ \cup_{i, -, r}(M) = \epsilon_{(a_1, \ldots, a_{i-1}, k, \ldots, a_r}^{(a_1, \ldots, a_{i-1}, k-1, 1, \ldots, a_r)}[-(k-1)] 
\Delta_{(k, a_1, \ldots, a_r)}^{(a_1, \ldots, a_{i-1}, k, a_r)} \zeta(M). $$

\begin{figure}[htb]
  \centering
  \includegraphics{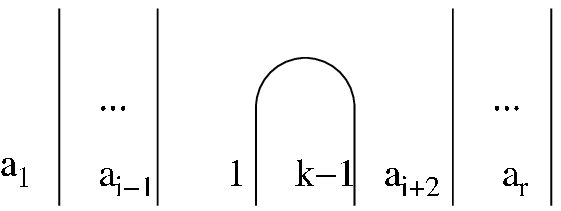}
  \caption{}
  \label{cap+}
\end{figure}

To the tangle in Figure ~\ref{cap+} we assign the functor
$ \cap_{i, +, r} \colon D^b(\mathcal{O}_{\bf d'}^{(a_1, \ldots, a_r)}) \rightarrow D^b(\mathcal{O}_{\bf d}^{(a_1, \ldots, a_{i-1}, a_{i+2}, \ldots, a_r)}). $  It is given by
$$ \cap_{i, +, r}(M) = \nu \Delta_{(a_1, \ldots, a_{i-1}, k, a_{i+2}, \ldots, a_r)}^{(k, a_1, \ldots, a_{i-1}, a_{i+2}, \ldots, a_r)} LZ_{(a_1, \ldots, a_r)}^{(a_1, \ldots, a_{i-1}, k, a_{i+2}, \ldots, a_r)}. $$

\begin{figure}[htb]
  \centering
  \includegraphics{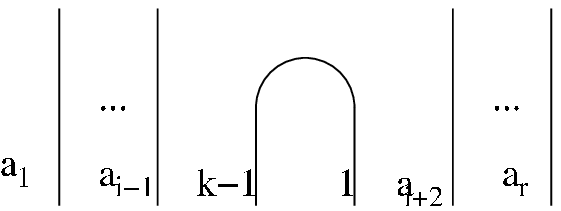}
  \caption{}
  \label{cap-}
\end{figure}

To the tangle in Figure ~\ref{cap-} we assign the functor
$ \cap_{i, -, r} \colon D^b(\mathcal{O}_{\bf d'}^{(a_1, \ldots, a_r)}) \rightarrow D^b(\mathcal{O}_{\bf d}^{(a_1, \ldots, a_{i-1}, a_{i+2}, \ldots, a_r)}). $  It is given by
$$ \cap_{i, +, r}(M) = \nu \Delta_{(a_1, \ldots, a_{i-1}, k, a_{i+2}, \ldots, a_r)}^{(k, a_1, \ldots, a_{i-1}, a_{i+2}, \ldots, a_r)} LZ_{(a_1, \ldots, a_r)}^{(a_1, \ldots, a_{i-1}, k, a_{i+2}, \ldots, a_r)}. $$

\section{Diagram Relations}
In this section we verify that the relations between various intertwiners which have a graphical interpretation given in figures ~\ref{relation1}, ~\ref{relation2}, ~\ref{relation3}, ~\ref{relation4}, and ~\ref{relation5} gives rise to relations between various inclusion and Zuckerman functors.   

\subsection{Diagram 1}
Our goal is to prove that  $ LZ_{\mathfrak{q}_j}^{\mathfrak{p}_j} \epsilon_{\mathfrak{p}_j}^{\mathfrak{q}_j} $ is a direct sum of shifted identity functors.
First we compute the cohomology functors on the generalized Verma module $ M^{\mathfrak{p}_j}(\alpha) $ where $ \alpha = (a_1, \ldots a_{j-1}, k-1, \ldots 0, a_{j+k}, \ldots a_n). $

\begin{lemma}
\label{lemmadiagram1}
$ L_{i}Z_{\mathfrak{q}_j}^{\mathfrak{p}_j} \epsilon_{\mathfrak{p}_j}^{\mathfrak{q}_j} M^{\mathfrak{p}_j}(\alpha) \cong M^{\mathfrak{p}_j}(\alpha) $ if $ i $ is even and $ 0 \leq i \leq 2(k-1). $
Otherwise it is 0.
\end{lemma}

\begin{proof}
First we note that $ L_s Z M^{\mathfrak{q}_j}(\sigma_1 \cdots \sigma_t.\alpha) \cong $
\begin{eqnarray*}
M^{\mathfrak{p}_{j}}(\alpha)	&\textrm{ if }		&s=t\\
0						&\textrm{ if } 		&s \neq t.
\end{eqnarray*}

Suppose $ i=2r. $  Then the above fact and the short exact sequences from the previous section imply
\begin{eqnarray*}
L_{i}Z_{\mathfrak{q}_j}^{\mathfrak{p}_j} \epsilon_{\mathfrak{p}_j}^{\mathfrak{q}_j} M^{\mathfrak{p}_j}(\alpha) &\cong &L_{i-1}Z_{\mathfrak{q}_j}^{\mathfrak{p}_j}K_0\\
L_{i-1}Z_{\mathfrak{q}_j}^{\mathfrak{p}_j}K_0 &\cong &L_{i-2}Z_{\mathfrak{q}_j}^{\mathfrak{p}_j}K_1\\
&\cdots &\\
L_{i-r+1}Z_{\mathfrak{q}_j}^{\mathfrak{p}_j}K_{r-2} &\cong &L_{i-r}Z_{\mathfrak{q}_j}^{p_j}K_{r-1}. 
\end{eqnarray*} 

Now consider $ 0 \rightarrow K_r \rightarrow M^{\mathfrak{q}_j}(\sigma_1 \cdots \sigma_r.\alpha) \rightarrow K_{r-1} \rightarrow 0. $
This gives a long exact sequence
\begin{align*}
0 \rightarrow &L_{i-r+1}Z_{q_j}^{p_j}K_{r-1} \rightarrow L_{i-r}Z_{\mathfrak{q}_j}^{\mathfrak{p}_j}K_r \rightarrow L_{i-r}Z_{\mathfrak{q}_j}^{\mathfrak{p}_j}M^{\mathfrak{q}_j}(\sigma_1 \ldots \sigma_r.\alpha) \cong M^{\mathfrak{p}_j}(\alpha) \rightarrow\\
&L_{i-r-1}Z_{\mathfrak{q}_j}^{\mathfrak{p}_j}K_{r-1} \rightarrow L_{i-r-1}Z_{\mathfrak{q}_j}^{\mathfrak{p}_j}K_r \rightarrow L_{i-r-1}Z_{\mathfrak{q}_j}^{\mathfrak{p}_j}M^{\mathfrak{q}_j}(\sigma_1 \cdots \sigma_r.\alpha) \cong 0. 
\end{align*}
Thus it suffices to show $ L_{i-r}Z_{\mathfrak{q}_j}^{\mathfrak{p}_j}K_r \cong L_{i-r-1}Z_{\mathfrak{q}_j}^{\mathfrak{p}_j}K_r \cong 0. $
Using the first fact of the proof we easily see
\begin{eqnarray*}
L_{i-r}Z_{\mathfrak{q}_j}^{\mathfrak{p}_j}K_r &\cong & L_{i-r-1}Z_{\mathfrak{q}_j}^{\mathfrak{p}_j}K_{r+1}\\
L_{i-r-1}Z_{\mathfrak{q}_j}^{\mathfrak{p}_j}K_{r+1} &\cong & L_{i-r-2}Z_{\mathfrak{q}_j}^{\mathfrak{p}_j}K_{r+2}\\
&\cdots &\\
L_{i-k+2}Z_{\mathfrak{q}_j}^{\mathfrak{p}_j}K_{k-2} &\cong &L_{i-k+1}Z_{\mathfrak{q}_j}^{\mathfrak{p}_j}K_{k-1} \cong 0. 
\end{eqnarray*} 
Thus $ L_{i-r}Z_{\mathfrak{q}_j}^{\mathfrak{p}_j}K_{r} \cong 0. $
Again using the first fact we easily see that
\begin{align*}
L_{i-r-1}Z_{\mathfrak{q}_j}^{\mathfrak{p}_j}K_r &\cong L_{i-r-2}Z_{\mathfrak{q}_j}^{\mathfrak{p}_j}K_{r+1}\\
L_{i-r-2}Z_{\mathfrak{q}_j}^{\mathfrak{p}_j}K_{r+1} &\cong L_{i-r-3}Z_{\mathfrak{q}_j}^{\mathfrak{p}_j}K_{r+2}\\
&\cdots\\ 
L_{i-k+1}Z_{\mathfrak{q}_j}^{\mathfrak{p}_j}K_{k-2} &\cong L_{i-k}Z_{\mathfrak{q}_j}^{\mathfrak{p}_j}K_{k-1} \cong 0. 
\end{align*} 
Thus $ L_{i-r-1}Z_{\mathfrak{q}_j}^{\mathfrak{p}_j}K_{r} \cong 0. $

Now we proceed for odd i.  Let $ i=2r+1. $  Then we get
\begin{align*}
L_{i}Z_{\mathfrak{q}_j}^{\mathfrak{p}_j} \epsilon_{\mathfrak{p}_j}^{\mathfrak{q}_j} M^{\mathfrak{p}_j}(\alpha) &\cong L_{i-1}Z_{\mathfrak{q}_j}^{\mathfrak{p}_j}K_0\\
L_{i-1}Z_{\mathfrak{q}_j}^{\mathfrak{p}_j}K_0 &\cong L_{i-2}Z_{\mathfrak{q}_j}^{\mathfrak{p}_j}K_1\\
&\cdots\\
L_{i-r+1}Z_{\mathfrak{q}_j}^{\mathfrak{p}_j}K_{r-2} &\cong L_{i-r}Z_{\mathfrak{q}_j}^{\mathfrak{p}_j}K_{r-1}. 
\end{align*}

Now consider $ 0 \rightarrow K_r \rightarrow M^{q_j}(\sigma_1 \cdots \sigma_r.\alpha) \rightarrow K_{r-1} \rightarrow 0. $
This gives a long exact sequence
$$ \rightarrow L_{i-r}Z_{\mathfrak{q}_j}^{\mathfrak{p}_j}M^{\mathfrak{q}_j}(\sigma_1 \cdots \sigma_r.\alpha) \cong 0 \rightarrow
L_{i-r}Z_{\mathfrak{q}_j}^{\mathfrak{p}_j}K_{r-1} \rightarrow L_{i-r-1}Z_{\mathfrak{q}_j}^{\mathfrak{p}_j}K_r \rightarrow L_{i-r-1}Z_{\mathfrak{q}_j}^{\mathfrak{p}_j}M^{\mathfrak{q}_j}(\sigma_1 \cdots \sigma_r.\alpha) \rightarrow . $$
Thus it suffices to prove $ L_{i-r-1}Z_{\mathfrak{q}_j}^{\mathfrak{p}_j}K_{r} \cong 0. $
Again using the first fact, it is easy to see that 
\begin{align*}
L_{i-r-1}Z_{\mathfrak{q}_j}^{\mathfrak{p}_j}K_r &\cong L_{i-r-2}Z_{\mathfrak{q}_j}^{\mathfrak{p}_j}K_{r+1}\\
L_{i-r-2}Z_{\mathfrak{q}_j}^{\mathfrak{p}_j}K_{r+1} &\cong L_{i-r-3}Z_{\mathfrak{q}_j}^{\mathfrak{p}_j}K_{r+2}\\
&\cdots\\ 
L_{i-k+1}Z_{\mathfrak{q}_j}^{\mathfrak{p}_j}K_{k-2} &\cong L_{i-k}Z_{\mathfrak{q}_j}^{\mathfrak{p}_j}K_{k-1} \cong 0. 
\end{align*} 
Thus $ L_{i-r-1}Z_{\mathfrak{q}_j}^{\mathfrak{p}_j}K_{r} \cong 0. $
\end{proof}

\begin{lemma}
Let $ X \overset f\to Y \overset g\to Z \overset h\to T(X) $ be a distinguished triangle in a triangulated category.
If $ h=0, $ then this triangle is isomorphic to
$ X \rightarrow X \oplus Z \rightarrow Z \rightarrow T(X). $
\end{lemma}

\begin{proof}
Since $ h $ is zero, using the axioms of a triangulated category one could construct a morphism from the triangle 
$ X \overset f\to Y \overset g\to Z \overset h\to T(X) $ to the triangle 
$ X \rightarrow X \oplus Z \rightarrow Z \rightarrow T(X). $ The map from $ Y $ to $ X \oplus Z $ is an isomorphism because the other two maps in the morphism of triangles are isomorphisms.
\end{proof}

\begin{prop}
\label{prop8}
$ LZ_{\mathfrak{q}_j}^{\mathfrak{p}_j} \epsilon_{\mathfrak{p}_j}^{\mathfrak{q}_j} \cong \oplus^{k-1}_{r=0} \Id[2r]. $
\end{prop}

\begin{proof}
We claim that for $ 0 \leq n \leq k-1, $
$$ \tau^{\leq -2n} LZ_{\mathfrak{q}_j}^{\mathfrak{p}_j} \epsilon_{\mathfrak{p}_j}^{\mathfrak{q}_j} M^{\mathfrak{p}_j}(\alpha) \cong \oplus_{r=0}^{(k-1)-n} M^{\mathfrak{p}_j}(\alpha)[2(k-1)-2r]. $$
We proceed by induction.  The base case is $ n=k-1. $  
Assume by induction 
$$ \tau^{\leq -2n} LZ_{\mathfrak{q}_j}^{\mathfrak{p}_j} \epsilon_{\mathfrak{p}_j}^{\mathfrak{q}_j} M^{\mathfrak{p}_j}(\alpha) \cong \oplus_{r=0}^{(k-1)-n} M^{\mathfrak{p}_j}[2(k-1)-2r]. $$
Consider the distinguished triangle
$$ \tau^{\leq(-2n+1)} LZ_{\mathfrak{q}_j}^{\mathfrak{p}_j} \epsilon_{\mathfrak{p}_j}^{\mathfrak{q}_j} M^{\mathfrak{p}_j}(\alpha) \rightarrow \tau^{\leq(-2n+2)} LZ_{\mathfrak{q}_j}^{\mathfrak{p}_j} \epsilon_{\mathfrak{p}_j}^{\mathfrak{q}_j} M^{\mathfrak{p}_j}(\alpha) \rightarrow
H^{-2n+2} LZ_{\mathfrak{q}_j}^{\mathfrak{p}_j} \epsilon_{\mathfrak{p}_j}^{\mathfrak{q}_j}M^{\mathfrak{p}_j}(\alpha) [2n-2]. $$
This triangle is isomorphic to
$$ \tau^{\leq(-2n)} LZ_{\mathfrak{q}_j}^{\mathfrak{p}_j} \epsilon_{\mathfrak{p}_j}^{\mathfrak{q}_j} M^{\mathfrak{p}_j}(\alpha) \rightarrow \tau^{\leq(-2n+2)} LZ_{\mathfrak{q}_j}^{\mathfrak{p}_j} \epsilon_{\mathfrak{p}_j}^{\mathfrak{q}_j} M^{\mathfrak{p}_j}(\alpha) \rightarrow
M^{\mathfrak{p}_j}(\alpha)[2n-2]. $$
We would like to compute
$$ \Hom(M^{\mathfrak{p}_j}(\alpha)[2n-2], \tau^{\leq -2n} LZ_{\mathfrak{q}_j}^{\mathfrak{p}_j} \epsilon_{\mathfrak{p}_j}^{\mathfrak{q}_j}M^{\mathfrak{p}_j}[1]). $$
By the induction hypothesis this is equal to
$$ \Hom(M^{\mathfrak{p}_j}(\alpha)[2n-2], \oplus_{r=0}^{2(k-1)-n} M^{\mathfrak{p}_j}(\alpha)[2(k-1)-2r+1]). $$
Since there are no non-trivial morphism between a generalized Verma module and itself in different degrees, the above must be zero.
Thus by the previous lemma
$$ \tau^{\leq -2n+2} LZ_{\mathfrak{q}_j}^{\mathfrak{p}_j} \epsilon_{\mathfrak{p}_j}^{\mathfrak{q}_j} M^{\mathfrak{p}_j}(\alpha) \cong \tau^{\leq -2n} LZ_{\mathfrak{q}_j}^{\mathfrak{p}_j} \epsilon_{\mathfrak{p}_j}^{\mathfrak{q}_j} M^{\mathfrak{p}_j}(\alpha) \oplus M^{\mathfrak{p}_j}(\alpha)[2n-2] \cong \oplus_{r=0}^{2(k-1)-(n-1)} M^{\mathfrak{p}_j}(\alpha)[2(k-1)-2r] $$
which completes the induction.

In the Koszul dual situation, the inclusion functor and derived Zuckerman functor become graded translation off an intersection of walls and graded translation back respectively.  For this fact see [MOS] or [Rh].  Then by theorem 4.12 of [J], the composition of these functors is a direct sum of identity functors with various shifts in the grading.  By Koszul duality these become identity functors with shifts in the grading and homological shifts.  The computation above for the generalized module indicate exactly what the homological shifts are.
\end{proof} 

\begin{corollary}
$ LZ_{\mathfrak{q}_j}^{\mathfrak{p}_j} \epsilon_{\mathfrak{p}_j}^{\mathfrak{q}_j}[-(k-1)] \cong \oplus^{k-1}_{r=0} \Id [2r-(k-1)]. $
\end{corollary}

\begin{remark}
Lemma 3.5.4 of [BGS] is a geometric version of this result.
\end{remark}

\subsection{Diagram 2}
This subsection is a restatement of proposition 16 of [BFK].  It could also be proven using the techniques for diagram 1, but here we just refer to [BFK].  
The following is a functorial version of the second relation in the graphical calculus of [MOY].

\begin{prop}
There is an isomorphism of functors
$ LZ^{\mathfrak{s}_i} \epsilon_{\mathfrak{s}_i} \cong \Id \oplus \Id[2] $.
\end{prop}

\begin{proof}
See [BFK] proposition 16.
\end{proof}

\begin{corollary}
$ \epsilon_{\mathfrak{s}_i}[-1]LZ^{\mathfrak{s}_i}\epsilon_{\mathfrak{s}_i}[-1]LZ^{\mathfrak{s}_i} \cong \epsilon_{\mathfrak{s}_i}[-1]LZ^{\mathfrak{s}_i}[1] \oplus \epsilon_{\mathfrak{s}_i}[-1]LZ^{\mathfrak{s}_i}[-1] $.
\end{corollary}

\begin{proof}
This is a direct consequence of the previous proposition.
\end{proof}

\subsection{Diagram 3}
\label{sec4.3}

We would like to have a functorial isomorphism corresponding to the third graphical relation.

Recall the definitions of the algebras, $ \mathfrak{p}_i, \mathfrak{q}_i, $ and $ \mathfrak{s}_i. $  There are obvious inclusions of subalgebras: $ \mathfrak{q}_i \subset \mathfrak{p}_i, $ $ \mathfrak{q}_{i+1} \subset \mathfrak{s}_i+\mathfrak{q}_{i+1}, $ $ \mathfrak{q}_{i-1} \subset \mathfrak{p}_i. $

We begin by studying the generalized Verma module with highest weight $ \alpha, $ $ M^{\mathfrak{p}_{i+1}}(\alpha)= $
$$ M^{\mathfrak{p}_{i+1}}(a_1 \ldots a_i, k-1 \ldots, 0, a_{i+k+1} \ldots, a_n). $$

For notational purposes, we will omit the components $ a_1, \ldots, a_{i-1}, a_{i+k+1}, \ldots, a_n $ from the weights in the generalized Verma modules in this subsection.
Denote by $ M^{\mathfrak{q}_{i+1}}(\sigma_1 \cdots \sigma_t.\alpha) $ the module
$$ M^{\mathfrak{q}_{i+1}}(a_i, k-1-t, k-1, \ldots, \widehat{k-1-t}, \ldots, 0). $$ 
(The term $ k-1-t $ is omitted.)

The generalized BGG resolution gives rise to the following long exact sequence
$$ 
0 \rightarrow M^{\mathfrak{q}_{i+1}}(\sigma_1 \cdots \sigma_{k-1}.\alpha) \rightarrow M^{\mathfrak{q}_{i+1}}(\sigma_1 \cdots \sigma_{k-2}.\alpha) \rightarrow \cdots \rightarrow
M^{\mathfrak{q}_{i+1}}(e.\alpha) \rightarrow M^{\mathfrak{p}_{i+1}}(a_i, k-1, \ldots, 0) \rightarrow 0. $$

Then the above sequence gives rise to the following set of short exact sequences:

\begin{align*}
&0 \rightarrow K_0 \rightarrow M^{\mathfrak{q}_{i+1}}(e.\alpha) \rightarrow M^{\mathfrak{p}_{i+1}}(a_i, k-1, \ldots, 0) \rightarrow 0\\
&0 \rightarrow K_1 \rightarrow M^{\mathfrak{q}_{i+1}}(\sigma_1.\alpha) \rightarrow K_0 \rightarrow 0\\
&\cdots\\
&0 \rightarrow K_{k-2} \rightarrow M^{\mathfrak{q}_{i+1}}(\sigma_1 \cdots \sigma_{k-2}.\alpha) \rightarrow K_{k-3} \rightarrow 0\\
&0 \rightarrow M^{\mathfrak{q}_{i+1}}(\sigma_1 \cdots \sigma_{k-1}.\alpha) \rightarrow K_{k-2} \rightarrow 0. 
\end{align*} 

Our first goal is to compute 
$ LZ_{\mathfrak{q}_{i+1}}^{\mathfrak{s}_i+\mathfrak{q}_{i+1}} \epsilon_{\mathfrak{p}_{i+1}}^{\mathfrak{q}_{i+1}} M^{\mathfrak{p}_{i+1}}(a_i, k-1, \ldots, 0). $

\begin{lemma}
\label{lemma10}
Let $ a_i = l. $  
\begin{enumerate}
\item Suppose $ s=0 $ or $ s=2. $ Then
$$ 
L_{s}Z_{\mathfrak{q}_{i+1}}^{\mathfrak{s}_i+\mathfrak{q}_{i+1}} \epsilon_{\mathfrak{p}_{i+1}}^{\mathfrak{q}_{i+1}} M^{\mathfrak{p}_{i+1}}(a_i, k-1, \ldots, 0) = 0 $$
\item Otherwise we have isomorphisms:
\begin{align*}
&L_{1}Z_{\mathfrak{q}_{i+1}}^{\mathfrak{s}_i+\mathfrak{q}_{i+1}} \epsilon_{\mathfrak{p}_{i+1}}^{\mathfrak{q}_{i+1}} M^{\mathfrak{p}_{i+1}}(a_i, k-1, \ldots, 0) \cong
L_1Z_{\mathfrak{q}_{i+1}}^{\mathfrak{s}_i+\mathfrak{q}_{i+1}}M^{\mathfrak{q}_{i+1}}(e.\alpha)/L_1Z_{\mathfrak{q}_{i+1}}^{\mathfrak{s}_i+\mathfrak{q}_{i+1}} K_0\\
&L_{1}Z_{\mathfrak{q}_{i+1}}^{\mathfrak{s}_i+\mathfrak{q}_{i+1}}K_0 \cong L_{1}Z_{\mathfrak{q}_{i+1}}^{\mathfrak{s}_i+\mathfrak{q}_{i+1}} M^{\mathfrak{q}_{i+1}}(\sigma_1.\alpha)/L_1Z_{\mathfrak{q}_{i+1}}^{\mathfrak{s}_i+\mathfrak{q}_{i+1}} K_1\\
&\cdots\\
&L_{1}Z_{\mathfrak{q}_{i+1}}^{\mathfrak{s}_i+\mathfrak{q}_{i+1}}K_{k-l-3} \cong L_{1}Z_{\mathfrak{q}_{i+1}}^{\mathfrak{s}_i+\mathfrak{q}_{i+1}} M^{\mathfrak{q}_{i+1}}(\sigma_1 \cdots \sigma_{k-l-2}.\alpha)/L_1Z_{\mathfrak{q}_{i+1}}^{\mathfrak{s}_i+\mathfrak{q}_{i+1}} 
K_{k-l-2}\\
&L_{s}Z_{\mathfrak{q}_{i+1}}^{\mathfrak{s}_i+\mathfrak{q}_{i+1}}K_{k-l-2} \cong L_{s-1}Z_{\mathfrak{q}_{i+1}}^{\mathfrak{s}_i+\mathfrak{q}_{i+1}} K_{k-l-1}\\
&L_{0}Z_{\mathfrak{q}_{i+1}}^{\mathfrak{s}_i+\mathfrak{q}_{i+1}}K_{k-l-1} \cong L_{0}Z_{\mathfrak{q}_{i+1}}^{\mathfrak{s}_i+\mathfrak{q}_{i+1}}M^{\mathfrak{q}_{i+1}}(\sigma_1 \cdots \sigma_{k-l}.\alpha)/L_{0}Z_{\mathfrak{q}_{i+1}}^{\mathfrak{s}_i+\mathfrak{q}_{i+1}}K_{k-l}\\
&L_{0}Z_{\mathfrak{q}_{i+1}}^{\mathfrak{s}_i+\mathfrak{q}_{i+1}}K_{k-l} \cong L_{0}Z_{\mathfrak{q}_{i+1}}^{\mathfrak{s}_i+\mathfrak{q}_{i+1}}M^{\mathfrak{q}_{i+1}}(\sigma_1 \cdots \sigma_{k-l+1}.\alpha)/L_{0}Z_{\mathfrak{q}_{i+1}}^{\mathfrak{s}_i+\mathfrak{q}_{i+1}}
K_{k-l+1}\\
& \cdots\\
&L_{0}Z_{\mathfrak{q}_{i+1}}^{\mathfrak{s}_i+\mathfrak{q}_{i+1}}K_{k-4} \cong L_{0}Z_{\mathfrak{q}_{i+1}}^{\mathfrak{s}_i+\mathfrak{q}_{i+1}}M^{\mathfrak{q}_{i+1}}(\sigma_1 \cdots \sigma_{k-3}.\alpha)/L_{0}Z_{\mathfrak{q}_{i+1}}^{\mathfrak{s}_i+\mathfrak{q}_{i+1}}K_{k-3}\\
&L_{0}Z_{\mathfrak{q}_{i+1}}^{\mathfrak{s}_i+\mathfrak{q}_{i+1}}K_{k-3} \cong L_{0}Z_{\mathfrak{q}_{i+1}}^{\mathfrak{s}_i+\mathfrak{q}_{i+1}}M^{\mathfrak{q}_{i+1}}(\sigma_1 \cdots \sigma_{k-2}.\alpha)/L_{0}Z_{\mathfrak{q}_{i+1}}^{\mathfrak{s}_i+\mathfrak{q}_{i+1}}
M^{\mathfrak{q}_{i+1}}(\sigma_1 \cdots \sigma_{k-1}.\alpha). 
\end{align*} 
\end{enumerate}
\end{lemma}

\begin{proof}
By proposition 5.5 of [ES], $ L_{s}Z_{\mathfrak{q}_{i+1}}^{\mathfrak{s}_i+\mathfrak{q}_{i+1}}M^{\mathfrak{q}_{i+1}}(\sigma_1 \cdots \sigma_t.\alpha) \cong $
\begin{eqnarray*} 
M^{\mathfrak{s}_i+\mathfrak{q}_{i+1}}(a_i, k-1-t, k-1, \ldots, \hat{k-1-t}, \ldots, 0) 	& \textrm{ if } & s=0, l>k-1-t \\
M^{\mathfrak{s}_i+\mathfrak{q}_{i+1}}(k-1-t, a_i, k-1, \ldots, \hat{k-1-t}, \ldots, 0) 	& \textrm{ if } & s=1, l<k-1-t\\
0		& \textrm{ if }		& \textrm{ otherwise }
\end{eqnarray*}

Now the functor $ LZ_{\mathfrak{q}_{i+1}}^{\mathfrak{s}_i+\mathfrak{q}_{i+1}} $ induces long exact sequences for all the short exact sequences stated before the lemma.
$$ \aligned 
0 \rightarrow &L_2Z_{\mathfrak{q}_{i+1}}^{\mathfrak{s}_i+\mathfrak{q}_{i+1}}K_0 \rightarrow L_2Z_{\mathfrak{q}_{i+1}}^{\mathfrak{s}_i+\mathfrak{q}_{i+1}}M^{\mathfrak{q}_{i+1}}(e.\alpha) \rightarrow L_2Z_{\mathfrak{q}_{i+1}}^{\mathfrak{s}_i+\mathfrak{q}_{i+1}}\epsilon_{\mathfrak{p}_{i+1}}^{\mathfrak{q}_{i+1}} M^{\mathfrak{p}_{i+1}}(\alpha) \rightarrow\\
&L_1Z_{\mathfrak{q}_{i+1}}^{\mathfrak{s}_i+\mathfrak{q}_{i+1}}K_0 \rightarrow L_1Z_{\mathfrak{q}_{i+1}}^{\mathfrak{s}_i+\mathfrak{q}_{i+1}}M^{\mathfrak{q}_{i+1}}(e) \rightarrow L_1Z_{\mathfrak{q}_{i+1}}^{\mathfrak{s}_i+\mathfrak{q}_{i+1}}\epsilon_{\mathfrak{p}_{i+1}}^{\mathfrak{q}_{i+1}}M^{\mathfrak{p}_{i+1}}(\alpha) \rightarrow\\
&L_0Z_{\mathfrak{q}_{i+1}}^{\mathfrak{s}_i+\mathfrak{q}_{i+1}}K_0 \rightarrow L_0Z_{\mathfrak{q}_{i+1}}^{\mathfrak{s}_i+\mathfrak{q}_{i+1}}M^{\mathfrak{q}_{i+1}}(e) \rightarrow L_0Z_{\mathfrak{q}_{i+1}}^{\mathfrak{s}_i+\mathfrak{q}_{i+1}}\epsilon_{\mathfrak{p}_{i+1}}^{\mathfrak{q}_{i+1}}M^{\mathfrak{p}_{i+1}}(\alpha) \rightarrow 0 \endaligned $$

$$ \aligned
0 \rightarrow &L_2Z_{\mathfrak{q}_{i+1}}^{\mathfrak{s}_i+\mathfrak{q}_{i+1}}K_1 \rightarrow L_2Z_{\mathfrak{q}_{i+1}}^{\mathfrak{s}_i+\mathfrak{q}_{i+1}}M^{\mathfrak{q}_{i+1}}(\sigma_1.\alpha) \rightarrow L_2Z_{\mathfrak{q}_{i+1}}^{\mathfrak{s}_i+\mathfrak{q}_{i+1}}K_0 \rightarrow\\
&L_1Z_{\mathfrak{q}_{i+1}}^{\mathfrak{s}_i+\mathfrak{q}_{i+1}}K_1 \rightarrow L_1Z_{\mathfrak{q}_{i+1}}^{\mathfrak{s}_i+\mathfrak{q}_{i+1}}M^{\mathfrak{q}_{i+1}}(\sigma_1.\alpha) \rightarrow L_1Z_{\mathfrak{q}_{i+1}}^{\mathfrak{s}_i+\mathfrak{q}_{i+1}}K_0 \rightarrow\\
&L_0Z_{\mathfrak{q}_{i+1}}^{\mathfrak{s}_i+\mathfrak{q}_{i+1}}K_1 \rightarrow L_0Z_{\mathfrak{q}_{i+1}}^{\mathfrak{s}_i+\mathfrak{q}_{i+1}}M^{\mathfrak{q}_{i+1}}(\sigma_1.\alpha) \rightarrow L_0Z_{\mathfrak{q}_{i+1}}^{\mathfrak{s}_i+\mathfrak{q}_{i+1}}K_0 \rightarrow 0 \endaligned $$

$$ \aligned
0 \rightarrow &L_2Z_{\mathfrak{q}_{i+1}}^{\mathfrak{s}_i+\mathfrak{q}_{i+1}}K_2 \rightarrow L_2Z_{\mathfrak{q}_{i+1}}^{\mathfrak{s}_i+\mathfrak{q}_{i+1}}M^{\mathfrak{q}_{i+1}}(\sigma_1 \sigma_2.\alpha) \rightarrow L_2Z_{\mathfrak{q}_{i+1}}^{\mathfrak{s}_i+\mathfrak{q}_{i+1}}K_1 \rightarrow\\
&L_1Z_{\mathfrak{q}_{i+1}}^{\mathfrak{s}_i+\mathfrak{q}_{i+1}}K_2 \rightarrow L_1Z_{\mathfrak{q}_{i+1}}^{\mathfrak{s}_i+\mathfrak{q}_{i+1}}M^{\mathfrak{q}_{i+1}}(\sigma_1 \sigma_2.\alpha) \rightarrow L_1Z_{\mathfrak{q}_{i+1}}^{\mathfrak{s}_i+\mathfrak{q}_{i+1}}K_1 \rightarrow\\
&L_0Z_{\mathfrak{q}_{i+1}}^{\mathfrak{s}_i+\mathfrak{q}_{i+1}}K_2 \rightarrow L_0Z_{\mathfrak{q}_{i+1}}^{\mathfrak{s}_i+\mathfrak{q}_{i+1}}M^{\mathfrak{q}_{i+1}}(\sigma_1 \sigma_2.\alpha) \rightarrow L_0Z_{\mathfrak{q}_{i+1}}^{\mathfrak{s}_i+\mathfrak{q}_{i+1}}K_1 \rightarrow 0  \endaligned $$

$$ \cdots $$

$$ \aligned 
0 \rightarrow &L_2Z_{\mathfrak{q}_{i+1}}^{\mathfrak{s}_i+\mathfrak{q}_{i+1}}K_{k-l-2} \rightarrow L_2Z_{\mathfrak{q}_{i+1}}^{\mathfrak{s}_i+\mathfrak{q}_{i+1}}M^{\mathfrak{q}_{i+1}}(\sigma_1 \cdots \sigma_{k-l-2}.\alpha) \rightarrow L_2Z_{\mathfrak{q}_{i+1}}^{\mathfrak{s}_i+\mathfrak{q}_{i+1}}K_{k-l-3} \rightarrow\\
&L_1Z_{\mathfrak{q}_{i+1}}^{\mathfrak{s}_i+\mathfrak{q}_{i+1}}K_{k-l-2} \rightarrow L_1Z_{\mathfrak{q}_{i+1}}^{\mathfrak{s}_i+\mathfrak{q}_{i+1}}M^{\mathfrak{q}_{i+1}}(\sigma_1 \cdots \sigma_{k-l-2}.\alpha) \rightarrow L_1Z_{\mathfrak{q}_{i+1}}^{\mathfrak{s}_i+\mathfrak{q}_{i+1}}K_{k-l-3} 
\rightarrow\\
&L_0Z_{\mathfrak{q}_{i+1}}^{\mathfrak{s}_i+\mathfrak{q}_{i+1}}K_{k-l-2} \rightarrow L_0Z_{\mathfrak{q}_{i+1}}^{\mathfrak{s}_i+\mathfrak{q}_{i+1}}M^{\mathfrak{q}_{i+1}}(\sigma_1 \cdots \sigma_{k-l-2}.\alpha) \rightarrow L_0Z_{\mathfrak{q}_{i+1}}^{\mathfrak{s}_i+\mathfrak{q}_{i+1}}K_{k-l-3} 
\rightarrow 0. \endaligned $$

This concludes the first set of exact sequences.  Then we have the following isomorphisms.
$$ \aligned
L_0Z_{\mathfrak{q}_{i+1}}^{\mathfrak{s}_i+\mathfrak{q}_{i+1}}M^{\mathfrak{q}_{i+1}}(\sigma_1 \cdots \sigma_{k-l-1}.\alpha) &\cong L_1Z_{\mathfrak{q}_{i+1}}^{\mathfrak{s}_i+\mathfrak{q}_{i+1}}M^{\mathfrak{q}_{i+1}}(\sigma_1 \cdots \sigma_{k-l-1}.\alpha)\\
&\cong L_2Z_{\mathfrak{q}_{i+1}}^{\mathfrak{s}_i+\mathfrak{q}_{i+1}}M^{\mathfrak{q}_{i+1}}(\sigma_1 \cdots \sigma_{k-l-1}.\alpha)\\
&\cong 0. 
\endaligned $$

Now we have this second set of exact sequences:
$$ \aligned 
0 \rightarrow &L_2Z_{\mathfrak{q}_{i+1}}^{\mathfrak{s}_i+\mathfrak{q}_{i+1}}K_{k-l} \rightarrow L_2Z_{\mathfrak{q}_{i+1}}^{\mathfrak{s}_i+\mathfrak{q}_{i+1}}M^{\mathfrak{q}_{i+1}}(\sigma_1 \cdots \sigma_{k-l}.\alpha) \rightarrow L_2Z_{\mathfrak{q}_{i+1}}^{\mathfrak{s}_i+\mathfrak{q}_{i+1}}K_{k-l-1} \rightarrow\\
&L_1Z_{\mathfrak{q}_{i+1}}^{\mathfrak{s}_i+\mathfrak{q}_{i+1}}K_{k-l} \rightarrow L_1Z_{\mathfrak{q}_{i+1}}^{\mathfrak{s}_i+\mathfrak{q}_{i+1}}M^{\mathfrak{q}_{i+1}}(\sigma_1 \cdots \sigma_{k-l}.\alpha) \rightarrow L_1Z_{\mathfrak{q}_{i+1}}^{\mathfrak{s}_i+\mathfrak{q}_{i+1}}K_{k-l-1} \rightarrow\\
&L_0Z_{\mathfrak{q}_{i+1}}^{\mathfrak{s}_i+\mathfrak{q}_{i+1}}K_{k-l} \rightarrow L_0Z_{\mathfrak{q}_{i+1}}^{\mathfrak{s}_i+\mathfrak{q}_{i+1}}M^{\mathfrak{q}_{i+1}}(\sigma_1 \cdots \sigma_{k-l}.\alpha) \rightarrow L_0Z_{\mathfrak{q}_{i+1}}^{\mathfrak{s}_i+\mathfrak{q}_{i+1}}K_{k-l-1}
\rightarrow 0 \endaligned $$

$$ \aligned 
0 \rightarrow &L_2Z_{\mathfrak{q}_{i+1}}^{\mathfrak{s}_i+\mathfrak{q}_{i+1}}K_{k-l+1} \rightarrow L_2Z_{\mathfrak{q}_{i+1}}^{\mathfrak{s}_i+\mathfrak{q}_{i+1}}M^{\mathfrak{q}_{i+1}}(\sigma_1 \cdots \sigma_{k-l+1}.\alpha) \rightarrow L_2Z_{\mathfrak{q}_{i+1}}^{\mathfrak{s}_i+\mathfrak{q}_{i+1}}K_{k-l} \rightarrow\\
&L_1Z_{\mathfrak{q}_{i+1}}^{\mathfrak{s}_i+\mathfrak{q}_{i+1}}K_{k-l+1} \rightarrow L_1Z_{\mathfrak{q}_{i+1}}^{\mathfrak{s}_i+\mathfrak{q}_{i+1}}M^{\mathfrak{q}_{i+1}}(\sigma_1 \cdots \sigma_{k-l+1}.\alpha) \rightarrow L_1Z_{\mathfrak{q}_{i+1}}^{\mathfrak{s}_i+\mathfrak{q}_{i+1}}K_{k-l} \rightarrow\\
&L_0Z_{\mathfrak{q}_{i+1}}^{\mathfrak{s}_i+\mathfrak{q}_{i+1}}K_{k-l+1} \rightarrow L_0Z_{\mathfrak{q}_{i+1}}^{\mathfrak{s}_i+\mathfrak{q}_{i+1}}M^{\mathfrak{q}_{i+1}}(\sigma_1 \cdots \sigma_{k-l+1}.\alpha) \rightarrow L_0Z_{\mathfrak{q}_{i+1}}^{\mathfrak{s}_i+\mathfrak{q}_{I+1}}K_{k-l}
\rightarrow 0 \endaligned $$

$$ \aligned
0 \rightarrow &L_2Z_{\mathfrak{q}_{i+1}}^{\mathfrak{s}_i+\mathfrak{q}_{i+1}}K_{k-l+2} \rightarrow L_2Z_{\mathfrak{q}_{i+1}}^{\mathfrak{s}_i+\mathfrak{q}_{i+1}}M^{\mathfrak{q}_{i+1}}(\sigma_1 \cdots \sigma_{k-l+2}.\alpha) \rightarrow L_2Z_{\mathfrak{q}_{i+1}}^{\mathfrak{s}_i+\mathfrak{q}_{i+1}}K_{k-l+1} \rightarrow\\
&L_1Z_{\mathfrak{q}_{i+1}}^{\mathfrak{s}_i+\mathfrak{q}_{i+1}}K_{k-l+2} \rightarrow L_1Z_{\mathfrak{q}_{i+1}}^{\mathfrak{s}_i+\mathfrak{q}_{i+1}}M^{\mathfrak{q}_{i+1}}(\sigma_1 \cdots \sigma_{k-l+2}.\alpha) \rightarrow L_1Z_{\mathfrak{q}_{i+1}}^{\mathfrak{s}_i+\mathfrak{q}_{i+1}}K_{k-l+1} \rightarrow\\
&L_0Z_{\mathfrak{q}_{i+1}}^{\mathfrak{s}_i+\mathfrak{q}_{i+1}}K_{k-l+2} \rightarrow L_0Z_{\mathfrak{q}_{i+1}}^{\mathfrak{s}_i+\mathfrak{q}_{i+1}}M^{\mathfrak{q}_{i+1}}(\sigma_1 \cdots \sigma_{k-l+2}.\alpha) \rightarrow L_0Z_{\mathfrak{q}_{i+1}}^{\mathfrak{s}_i+\mathfrak{q}_{i+1}}K_{k-l+1}
\rightarrow 0 \endaligned $$

$$ \cdots $$

$$ \aligned
0 \rightarrow &L_2Z_{\mathfrak{q}_{i+1}}^{\mathfrak{s}_i+\mathfrak{q}_{i+1}}K_{k-2} \rightarrow L_2Z_{\mathfrak{q}_{i+1}}^{\mathfrak{s}_i+\mathfrak{q}_{i+1}}M^{\mathfrak{q}_{i+1}}(\sigma_1 \cdots \sigma_{k-2}.\alpha) \rightarrow L_2Z_{\mathfrak{q}_{i+1}}^{\mathfrak{s}_i+\mathfrak{q}_{i+1}}K_{k-3} \rightarrow\\
&L_1Z_{\mathfrak{q}_{i+1}}^{\mathfrak{s}_i+\mathfrak{q}_{i+1}}K_{k-2} \rightarrow L_1Z_{\mathfrak{q}_{i+1}}^{\mathfrak{s}_i+\mathfrak{q}_{i+1}}M^{\mathfrak{q}_{i+1}}(\sigma_1 \cdots \sigma_{k-2}.\alpha) \rightarrow L_1Z_{\mathfrak{q}_{i+1}}^{\mathfrak{s}_i+\mathfrak{q}_{i+1}}K_{k-3} \rightarrow\\
&L_0Z_{\mathfrak{q}_{i+1}}^{\mathfrak{s}_i+\mathfrak{q}_{i+1}}K_{k-2} \rightarrow L_0Z_{\mathfrak{q}_{i+1}}^{\mathfrak{s}_i+\mathfrak{q}_{i+1}}M^{\mathfrak{q}_{i+1}}(\sigma_1 \cdots \sigma_{k-2}.\alpha) \rightarrow L_0Z_{\mathfrak{q}_{i+1}}^{\mathfrak{s}_i+\mathfrak{q}_{i+1}}K_{k-3}
\rightarrow 0. \endaligned $$

Now we know
$$ L_2Z_{\mathfrak{q}_{i+1}}^{\mathfrak{s}_i+\mathfrak{q}_{i+1}}K_{k-2} \cong L_1Z_{\mathfrak{q}_{i+1}}^{\mathfrak{s}_i+\mathfrak{q}_{i+1}}K_{k-2} \cong 0, $$ 
and that
$$ L_0Z_{\mathfrak{q}_{i+1}}^{\mathfrak{s}_i+\mathfrak{q}_{i+1}}K_{k-2} \cong M^{\mathfrak{s}_i+\mathfrak{q}_{i+1}}(a_i, 0, k-1, \ldots, 1). $$

In the first set of sequences, $ L_0Z_{\mathfrak{q}_{i+1}}^{\mathfrak{s}_i+\mathfrak{q}_{i+1}}M^{\mathfrak{q}_{i+1}}(\sigma_1 \cdots \sigma_t.\alpha) = 0. $
Thus $ L_0Z_{\mathfrak{q}_{i+1}}^{\mathfrak{s}_i+\mathfrak{q}_{i+1}}K_{t-1} \cong 0 $ and $ L_0Z_{\mathfrak{q}_{i+1}}^{\mathfrak{s}_i+\mathfrak{q}_{i+1}}M^{\mathfrak{p}_{i+1}}(\alpha) \cong 0. $
Also $ L_2Z_{\mathfrak{q}_{i+1}}^{\mathfrak{s}_i+\mathfrak{q}_{i+1}}M^{\mathfrak{q}_{i+1}}(\sigma_1 \cdots \sigma_t.\alpha) \cong 0 $ so $ L_2 Z_{\mathfrak{q}_{i+1}}^{\mathfrak{s}_i+\mathfrak{q}_{i+1}}K_t \cong 0. $
Therefore each sequence in the first set produces a short exact sequence
$$ 0 \rightarrow L_1Z_{\mathfrak{q}_{i+1}}^{\mathfrak{s}_i+\mathfrak{q}_{i+1}}K_t \rightarrow L_1Z_{\mathfrak{q}_{i+1}}^{\mathfrak{s}_i+\mathfrak{q}_{i+1}}M^{\mathfrak{q}_{i+1}}(\sigma_1 \cdots \sigma_t.\alpha) \rightarrow 
L_1Z_{\mathfrak{q}_{i+1}}^{\mathfrak{s}_i+\mathfrak{q}_{i+1}}K_{t-1} \rightarrow 0. $$

Thus we have part of the lemma:
\begin{align*}
L_1Z_{\mathfrak{q}_{i+1}}^{\mathfrak{s}_i+\mathfrak{q}_{i+1}}M^{\mathfrak{p}_{i+1}}(\alpha) &\cong L_1Z_{\mathfrak{q}_{i+1}}^{\mathfrak{s}_i+\mathfrak{q}_{i+1}}M^{\mathfrak{q}_{i+1}}(e.\alpha)/L_1Z_{\mathfrak{q}_{i+1}}^{\mathfrak{s}_i+\mathfrak{q}_{i+1}}K_0\\
L_1Z_{\mathfrak{q}_{i+1}}^{\mathfrak{s}_i+\mathfrak{q}_{i+1}}K_0 &\cong L_1Z_{\mathfrak{q}_{i+1}}^{\mathfrak{s}_i+\mathfrak{q}_{i+1}}M^{\mathfrak{q}_{i+1}}(\sigma_1.\alpha)/L_1Z_{\mathfrak{q}_{i+1}}^{\mathfrak{s}_i+\mathfrak{q}_{i+1}}K_1\\
&\cdots\\
L_1Z_{\mathfrak{q}_{i+1}}^{\mathfrak{s}_i+\mathfrak{q}_{i+1}}K_{k-l-3} &\cong L_1Z_{\mathfrak{q}_{i+1}}^{\mathfrak{s}_i+\mathfrak{q}_{i+1}}M^{\mathfrak{q}_{i+1}}(\sigma_1 \cdots \sigma_{k-l-2}.\alpha)/L_1Z_{\mathfrak{q}_{i+1}}^{\mathfrak{s}_i+\mathfrak{q}_{i+1}}K_{k-l-2}. 
\end{align*} 

In the second set of exact sequences
\begin{align*}
L_2Z_{\mathfrak{q}_{i+1}}^{\mathfrak{s}_i+\mathfrak{q}_{i+1}}M^{\mathfrak{q}_{i+1}}(\sigma_1 \cdots \sigma_t.\alpha) &\cong L_1Z_{\mathfrak{q}_{i+1}}^{\mathfrak{s}_i+\mathfrak{q}_{i+1}}M^{\mathfrak{q}_{i+1}}(\sigma_1 \cdots \sigma_t.\alpha)\\ 
&\cong 0. 
\end{align*} 
Thus 
\begin{align*}
L_2Z_{\mathfrak{q}_{i+1}}^{\mathfrak{s}_i+\mathfrak{q}_{i+1}}K_t &\cong L_1Z_{\mathfrak{q}_{i+1}}^{\mathfrak{s}_i+\mathfrak{q}_{i+1}}K_{t+1}\\ 
&\cong 0. 
\end{align*} 
Therefore each sequence in the second set produces a short exact sequence
$$ 0 \rightarrow L_0Z_{\mathfrak{q}_{i+1}}^{\mathfrak{s}_i+\mathfrak{q}_{i+1}}K_t \rightarrow L_0Z_{\mathfrak{q}_{i+1}}^{\mathfrak{s}_i+\mathfrak{q}_{i+1}}M^{\mathfrak{q}_{i+1}}(\sigma_1 \cdots \sigma_t.\alpha) \rightarrow 
L_0Z_{\mathfrak{q}_{i+1}}^{\mathfrak{s}_i+\mathfrak{q}_{i+1}}K_{t-1} \rightarrow 0. $$
These short exact sequences give the rest of the lemma.
\end{proof}

Next we would like to include this object back into $ \mathcal{O}^{\mathfrak{q}_{i+1}} $ and then apply $ LZ_{\mathfrak{q}_{i+1}}^{\mathfrak{p}_{i+1}}. $
That is, we would like to compute 
$$ LZ_{\mathfrak{q}_{i+1}}^{\mathfrak{p}_{i+1}} \epsilon_{\mathfrak{s}_i+\mathfrak{q}_{i+1}}^{\mathfrak{q}_{i+1}} LZ_{\mathfrak{q}_{i+1}}^{\mathfrak{s}_i+\mathfrak{q}_{i+1}} M^{\mathfrak{p}_{i+1}}(a_i, k-1, \ldots, 0). $$
Therefore we want to compute
$$ L_sZ_{\mathfrak{q}_{i+1}}^{\mathfrak{p}_{i+1}} \epsilon_{\mathfrak{s}_i+\mathfrak{q}_{i+1}}^{\mathfrak{q}_{i+1}} L_1Z_{\mathfrak{q}_{i+1}}^{\mathfrak{s}_i+\mathfrak{q}_{i+1}} M^{\mathfrak{p}_{i+1}}(a_i, k-1, \ldots, 0). $$
We begin with something simpler.

\begin{lemma}
\label{lemma11}
There is an isomorphism:
$$ L_sZ_{\mathfrak{q}_{i+1}}^{\mathfrak{p}_{i+1}} \epsilon_{\mathfrak{s}_i+\mathfrak{q}_{i+1}}^{\mathfrak{q}_{i+1}}(L_0Z_{\mathfrak{q}_{i+1}}^{\mathfrak{s}_i+\mathfrak{q}_{i+1}}M^{\mathfrak{q}_{i+1}}(\sigma_1 \cdots \sigma_{k-2}.\alpha)/
L_0Z_{\mathfrak{q}_{i+1}}^{\mathfrak{s}_i+\mathfrak{q}_{i+1}}M^{\mathfrak{q}_{i+1}}(\sigma_1 \cdots \sigma_{k-1}.\alpha)) \cong $$
\begin{eqnarray*}
M^{\mathfrak{p}_{i+1}}(a_i, k-1, \ldots, 0)	&\textrm{ if }	&s=k-2,k\\
0								&\textrm{ if } 	&s \neq k-2, k.
\end{eqnarray*}
\end{lemma}  

\begin{proof}
Consider the short exact sequence
$$ 
0 \rightarrow \epsilon_{\mathfrak{s}_i+\mathfrak{q}_{i+1}}^{\mathfrak{q}_{i+1}}L_0Z_{\mathfrak{q}_{i+1}}^{\mathfrak{s}_i+\mathfrak{q}_{i+1}}M^{\mathfrak{q}_{i+1}}(\sigma_1 \cdots \sigma_{k-1}.\alpha) \rightarrow
\epsilon_{\mathfrak{s}_i+\mathfrak{q}_{i+1}}^{\mathfrak{q}_{i+1}}L_0Z_{\mathfrak{q}_{i+1}}^{\mathfrak{s}_i+\mathfrak{q}_{i+1}}M^{\mathfrak{q}_{i+1}}(\sigma_1 \cdots \sigma_{k-2}.\alpha)
\rightarrow \epsilon_{\mathfrak{s}_i+\mathfrak{q}_{i+1}}^{\mathfrak{q}_{i+1}}L_0Z_{\mathfrak{q}_{i+1}}^{\mathfrak{s}_i+\mathfrak{q}_{i+1}}K_{k-3} \rightarrow 0.  $$
Once we understand the first two terms, we will be able to understand the functor applied to $ K_{k-3}. $

So first we compute $ L_sZ_{\mathfrak{q}_{i+1}}^{\mathfrak{p}_{i+1}} \epsilon_{\mathfrak{s}_i+\mathfrak{q}_{i+1}}^{\mathfrak{q}_{i+1}} L_0Z_{\mathfrak{q}_{i+1}}^{\mathfrak{s}_i+\mathfrak{q}_{i+1}}M^{\mathfrak{q}_{i+1}}(\sigma_1 \cdots \sigma_{k-1}.\alpha). $
This is isomorphic to
$$ L_sZ_{\mathfrak{q}_{i+1}}^{\mathfrak{p}_{i+1}} \epsilon_{\mathfrak{s}_i+\mathfrak{q}_{i+1}}^{\mathfrak{q}_{i+1}} M^{\mathfrak{s}_i+\mathfrak{q}_{i+1}}(a_i, 0, k-1, \ldots, 1). $$
For this, consider the short exact sequence
$$ 0 \rightarrow M^{\mathfrak{q}_{i+1}}(0, a_i, k-1, \ldots, 1) \rightarrow 
M^{\mathfrak{q}_{i+1}}(a_{i}, 0, k-1, \ldots, 1) \rightarrow
M^{\mathfrak{s}_i+\mathfrak{q}_{i+1}}(a_i, 0, k-1, \ldots, 1) \rightarrow 0.  $$
Now $ L_sZ_{\mathfrak{q}_{i+1}}^{\mathfrak{p}_{i+1}} M^{\mathfrak{q}_{i+1}}(0, a_i, k-1, \ldots, 1) \cong 0 $ for all $ s. $
Also, 
$$ L_sZ_{\mathfrak{q}_{i+1}}^{\mathfrak{p}_{i+1}} M^{\mathfrak{q}_{i+1}}(a_i, 0, k-1, \ldots, 1) \cong $$
\begin{eqnarray*}
M^{\mathfrak{p}_{i+1}}(a_i, k-1, \ldots, 0) &\textrm{ if } &s=k-1\\
0	&\textrm{ if } 	&s\neq k-1.
\end{eqnarray*}
Thus $ L_sZ_{\mathfrak{q}_{i+1}}^{\mathfrak{p}_{i+1}} \epsilon_{\mathfrak{s}_i+\mathfrak{q}_{i+1}}^{\mathfrak{q}_{i+1}} L_0Z_{\mathfrak{q}_{i+1}}^{\mathfrak{s}_i+\mathfrak{q}_{i+1}}M^{\mathfrak{q}_{i+1}}(\sigma_1 \cdots \sigma_{k-1}.\alpha) \cong $
\begin{eqnarray*}
M^{\mathfrak{p}_{i+1}}(a_i, k-1, \ldots, 0)	&\textrm{ if }	&s=k-1\\
0 	&\textrm{ if } 	&s \neq k-1.
\end{eqnarray*}

Next we compute $ L_sZ_{\mathfrak{q}_{i+1}}^{\mathfrak{p}_{i+1}} \epsilon_{\mathfrak{s}_i+\mathfrak{q}_{i+1}}^{\mathfrak{q}_{i+1}} L_0Z_{\mathfrak{q}_{i+1}}^{\mathfrak{s}_i+\mathfrak{q}_{i+1}}M^{\mathfrak{q}_{i+1}}(\sigma_1 \cdots \sigma_{k-2}.\alpha). $
This is isomorphic to
$$ L_sZ_{\mathfrak{q}_{i+1}}^{\mathfrak{p}_{i+1}} \epsilon_{\mathfrak{s}_i+\mathfrak{q}_{i+1}}^{\mathfrak{q}_{i+1}} M^{\mathfrak{s}_i+\mathfrak{q}_{i+1}}(a_i, 1, k-1, \ldots, 2, 0). $$
For this, consider the short exact sequence
$$ 0 \rightarrow M^{\mathfrak{q}_{i+1}}(1, a_i, k-1, \ldots, 2, 0) \rightarrow 
M^{\mathfrak{q}_{i+1}}(a_{i}, 1, k-1, \ldots, 2, 0) \rightarrow
M^{\mathfrak{s}_i+\mathfrak{q}_{i+1}}(a_i, 1, k-1, \ldots, 2, 0) \rightarrow 0. $$

Now $ L_sZ_{\mathfrak{q}_{i+1}}^{\mathfrak{p}_{i+1}} M^{\mathfrak{q}_{i+1}}(1, a_i, k-1, \ldots, 2, 0) \cong 0 $ for all $ s. $
Also, $$ L_sZ_{\mathfrak{q}_{i+1}}^{\mathfrak{p}_{i+1}} M^{\mathfrak{q}_{i+1}}(a_i, 1, k-1, \ldots, 2, 0) \cong $$
\begin{eqnarray*}
M^{\mathfrak{p}_{i+1}}(a_i, k-1, \ldots, 0)	&\textrm{ if } 	&s=k-2\\ 
0  & \textrm{ if } 	& s\neq k-2.
\end{eqnarray*}
Thus  
$$ L_sZ_{\mathfrak{q}_{i+1}}^{\mathfrak{p}_{i+1}} \epsilon_{\mathfrak{s}_i+\mathfrak{q}_{i+1}}^{\mathfrak{q}_{i+1}} L_0Z_{\mathfrak{q}_{i+1}}^{\mathfrak{s}_i+\mathfrak{q}_{i+1}}M^{\mathfrak{q}_{i+1}}(\sigma_1 \cdots \sigma_{k-2}.\alpha) \cong $$
\begin{eqnarray*}
M^{\mathfrak{p}_{i+1}}(a_i, k-1, \ldots, 0)	&\textrm{ if }	& s=k-2\\
0 	&\textrm{ if }	& s\neq k-2.
\end{eqnarray*}

Finally, we consider the long exact sequence for the functor $ LZ_{\mathfrak{q}_{i+1}}^{\mathfrak{p}_{i+1}} $ for the short exact sequence considered in the beginning of the proof.
It is
$$ \aligned
 \cdots \rightarrow &L_sZ_{\mathfrak{q}_{i+1}}^{\mathfrak{p}_{i+1}} \epsilon_{\mathfrak{s}_i+\mathfrak{q}_{i+1}}^{\mathfrak{q}_{i+1}} L_0Z_{\mathfrak{q}_{i+1}}^{\mathfrak{s}_i+\mathfrak{q}_{i+1}}M^{\mathfrak{q}_{i+1}}(\sigma_1 \cdots \sigma_{k-1}.\alpha) \rightarrow
L_sZ_{\mathfrak{q}_{i+1}}^{\mathfrak{p}_{i+1}} \epsilon_{\mathfrak{s}_i+\mathfrak{q}_{i+1}}^{\mathfrak{q}_{i+1}} L_0Z_{\mathfrak{q}_{i+1}}^{\mathfrak{s}_i+\mathfrak{q}_{i+1}}M^{\mathfrak{q}_{i+1}}(\sigma_1 \cdots \sigma_{k-2}.\alpha) \rightarrow\\
&L_sZ_{\mathfrak{q}_{i+1}}^{\mathfrak{p}_{i+1}} \epsilon_{\mathfrak{s}_i+\mathfrak{q}_{i+1}}^{\mathfrak{q}_{i+1}} \epsilon_{\mathfrak{s}_i+\mathfrak{q}_{i+1}}^{\mathfrak{q}_{i+1}}(L_0Z_{\mathfrak{q}_{i+1}}^{\mathfrak{s}_i+\mathfrak{q}_{i+1}} K_{k-3} \rightarrow
L_{s-1}Z_{\mathfrak{q}_{i+1}}^{\mathfrak{p}_{i+1}} \epsilon_{\mathfrak{s}_i+\mathfrak{q}_{i+1}}^{\mathfrak{q}_{i+1}} L_0Z_{\mathfrak{q}_{i+1}}^{\mathfrak{s}_i+\mathfrak{q}_{i+1}}M^{\mathfrak{q}_{i+1}}(\sigma_1 \cdots \sigma_{k-1}.\alpha) \rightarrow\\
&L_{s-1} Z_{\mathfrak{q}_{i+1}}^{\mathfrak{p}_{i+1}} \epsilon_{\mathfrak{s}_i+\mathfrak{q}_{i+1}}^{\mathfrak{q}_{i+1}} L_0Z_{\mathfrak{q}_{i+1}}^{\mathfrak{s}_i+\mathfrak{q}_{i+1}}M^{\mathfrak{q}_{i+1}}(\sigma_1 \cdots \sigma_{k-2}.\alpha) \rightarrow \cdots.  \endaligned $$

If $ s=k-2 $ a portion of this long exact sequence reduces to
$$ 0 \rightarrow M^{\mathfrak{p}_{i+1}}(a_i, k-1, \ldots, 0) \rightarrow 
L_{k-2} Z_{\mathfrak{q}_{i+1}}^{\mathfrak{p}_{i+1}} \epsilon_{\mathfrak{s}_i+\mathfrak{q}_{i+1}}^{\mathfrak{q}_{i+1}} (L_0Z_{\mathfrak{q}_{i+1}}^{\mathfrak{s}_i+\mathfrak{q}_{i+1}}K_{k-3})
\rightarrow 0. $$
Thus 
$$ L_{k-2} Z_{\mathfrak{q}_{i+1}}^{\mathfrak{p}_{i+1}} \epsilon_{\mathfrak{r}_i}^{\mathfrak{q}_{i+1}} K_{k-3} \cong M^{\mathfrak{p}_{i+1}}(a_i, k-1, \ldots, 0). $$
If $ s-1=k-1, $ then $ s=k $ and the sequence reduces to
$$ 0 \rightarrow L_{k} Z_{\mathfrak{q}_{i+1}}^{\mathfrak{p}_{i+1}} \epsilon_{\mathfrak{s}_i+\mathfrak{q}_{i+1}}^{\mathfrak{q}_{i+1}} (L_0Z_{\mathfrak{q}_{i+1}}^{\mathfrak{s}_i+\mathfrak{q}_{i+1}}K_{k-3}) \rightarrow M^{\mathfrak{p}_{i+1}}(a_i, k-1, \ldots, 0) 
\rightarrow 0. $$
Thus 
$$ L_{k} Z_{\mathfrak{q}_{i+1}}^{\mathfrak{p}_{i+1}} \epsilon_{\mathfrak{s}_i+\mathfrak{q}_{i+1}}^{\mathfrak{q}_{i+1}} K_{k-3} \cong M^{\mathfrak{p}_{i+1}}(a_i, k-1, \ldots, 0). $$
For all other values it is 0.
\end{proof}

\begin{lemma}
\label{lemma12}
Suppose $ 3 \leq t \leq l+1 $ and $ a_i=l. $ 
If $ s=k-t+1, k-t+3, \ldots, k+t-3, $ then
$$ L_sZ_{\mathfrak{q}_{i+1}}^{\mathfrak{p}_{i+1}} \epsilon_{\mathfrak{s}_i+\mathfrak{q}_{i+1}}^{\mathfrak{q}_{i+1}}(L_0Z_{\mathfrak{q}_{i+1}}^{\mathfrak{s}_i+\mathfrak{q}_{i+1}}M^{\mathfrak{q}_{i+1}}(\sigma_1 \cdots \sigma_{k-t+1}.\alpha)/
L_0Z_{\mathfrak{q}_{i+1}}^{\mathfrak{s}_i+\mathfrak{q}_{i+1}} K_{k-t+1}) $$
$$ \cong M^{\mathfrak{p}_{i+1}}(a_i, k-1, \ldots, 0). $$
Otherwise it is zero.
\end{lemma}

\begin{proof}
We proceed by induction on $ t $ .  The base case is the previous lemma.  

As in the proof of the previous lemma, we find that
$$ L_sZ_{\mathfrak{q}_{i+1}}^{\mathfrak{p}_{i+1}} \epsilon_{\mathfrak{s}_i+\mathfrak{q}_{i+1}}^{\mathfrak{q}_{i+1}} L_0Z_{\mathfrak{q}_{i+1}}^{\mathfrak{s}_i+\mathfrak{q}_{i+1}}M^{\mathfrak{q}_{i+1}}(\sigma_1 \cdots \sigma_{k-t+1}.\alpha) \cong $$
\begin{eqnarray*}
M^{\mathfrak{p}_{i+1}}(a_i, k-1, \ldots, 0) 	&\textrm{ if }		&s=k-t+1\\
0		&\textrm{ if }		&s \neq k-t+1.
\end{eqnarray*} 

We also have analogously a long exact sequence
$$ \aligned
\cdots \rightarrow &L_sZ_{\mathfrak{q}_{i+1}}^{\mathfrak{p}_{i+1}} \epsilon_{\mathfrak{s}_i+\mathfrak{q}_{i+1}}^{\mathfrak{q}_{i+1}} L_0Z_{\mathfrak{q}_{i+1}}^{\mathfrak{s}_i+\mathfrak{q}_{i+1}}K_{k-t+1} \rightarrow
L_sZ_{\mathfrak{q}_{i+1}}^{\mathfrak{p}_{i+1}} \epsilon_{\mathfrak{s}_i+\mathfrak{q}_{i+1}}^{\mathfrak{q}_{i+1}} L_0Z_{\mathfrak{q}_{i+1}}^{\mathfrak{s}_i+\mathfrak{q}_{i+1}}M^{\mathfrak{q}_{i+1}}(\sigma_1 \cdots \sigma_{k-t+1}.\alpha) \rightarrow\\
&L_sZ_{\mathfrak{q}_{i+1}}^{\mathfrak{p}_{i+1}} \epsilon_{\mathfrak{s}_i+\mathfrak{q}_{i+1}}^{\mathfrak{q}_{i+1}} L_0Z_{\mathfrak{q}_{i+1}}^{\mathfrak{s}_i+\mathfrak{q}_{i+1}}K_{k-t} \rightarrow
L_{s-1} Z_{\mathfrak{q}_{i+1}}^{\mathfrak{p}_{i+1}} \epsilon_{\mathfrak{s}_i+\mathfrak{q}_{i+1}}^{\mathfrak{q}_{i+1}} L_0Z_{\mathfrak{q}_{i+1}}^{\mathfrak{s}_i+\mathfrak{q}_{i+1}}K_{k-t+1} \rightarrow\\
&L_{s-1} Z_{\mathfrak{q}_{i+1}}^{\mathfrak{p}_{i+1}} \epsilon_{\mathfrak{s}_i+\mathfrak{q}_{i+1}}^{\mathfrak{q}_{i+1}} L_0Z_{\mathfrak{q}_{i+1}}^{\mathfrak{s}_i+\mathfrak{q}_{i+1}}M^{\mathfrak{q}_{i+1}}(\sigma_1 \cdots \sigma_{k-t+1}.\alpha) \rightarrow \cdots. 
\endaligned $$

Now note that by the induction hypothesis
$ L_{s} Z_{\mathfrak{q}_{i+1}}^{\mathfrak{p}_{i+1}} \epsilon_{\mathfrak{s}_i+\mathfrak{q}_{i+1}}^{\mathfrak{q}_{i+1}} L_0Z_{\mathfrak{q}_{i+1}}^{\mathfrak{s}_i+\mathfrak{q}_{i+1}}K_{k-t+1} \cong 0 $ unless
$ s = k-(t-1)+1, \ldots, k+(t-1)-3. $  

Suppose $ s \leq k-t. $  Then
$$ L_{s} Z_{\mathfrak{q}_{i+1}}^{\mathfrak{p}_{i+1}} \epsilon_{\mathfrak{s}_i+\mathfrak{q}_{i+1}}^{\mathfrak{q}_{i+1}} L_0Z_{\mathfrak{q}_{i+1}}^{\mathfrak{s}_i+\mathfrak{q}_{i+1}}M^{\mathfrak{q}_{i+1}}(\sigma_1 \cdots \sigma_{k-t+1}.\alpha) \cong
L_{s-1} Z_{\mathfrak{q}_{i+1}}^{\mathfrak{p}_{i+1}} \epsilon_{\mathfrak{s}_i+\mathfrak{q}_{i+1}}^{\mathfrak{q}_{i+1}} L_0Z_{\mathfrak{q}_{i+1}}^{\mathfrak{s}_i+\mathfrak{q}_{i+1}}M^{\mathfrak{q}_{i+1}}(\sigma_1 \cdots \sigma_{k-t+1}.\alpha) \cong 0. $$
Now just apply the induction hypothesis.

Suppose $ s \geq k-t+3. $  Then
$$ L_{s} Z_{\mathfrak{q}_{i+1}}^{\mathfrak{p}_{i+1}} \epsilon_{\mathfrak{s}_i+\mathfrak{q}_{i+1}}^{\mathfrak{q}_{i+1}} L_0Z_{\mathfrak{q}_{i+1}}^{\mathfrak{s}_i+\mathfrak{q}_{i+1}}M^{\mathfrak{q}_{i+1}}(\sigma_1 \cdots \sigma_{k-t+1}.\alpha) \cong
L_{s-1} Z_{\mathfrak{q}_{i+1}}^{\mathfrak{p}_{i+1}} \epsilon_{\mathfrak{s}_i+\mathfrak{q}_{i+1}}^{\mathfrak{q}_{i+1}} L_0Z_{\mathfrak{q}_{i+1}}^{\mathfrak{s}_i+\mathfrak{q}_{i+1}}M^{\mathfrak{q}_{i+1}}(\sigma_1 \cdots \sigma_{k-t+1}.\alpha) \cong 0. $$
Now just apply the induction hypothesis.

Suppose $ s=k-t+1.  $ Then the long exact sequence and the induction hypothesis gives
$$ 
0 \rightarrow M^{\mathfrak{p}_{i+1}}(a_i, k-1, \ldots, 0) \rightarrow 
L_{k-t+1}Z_{\mathfrak{q}_{i+1}}^{\mathfrak{p}_{i+1}} \epsilon_{\mathfrak{s}_i+\mathfrak{q}_{i+1}}^{\mathfrak{q}_{i+1}} L_0Z_{\mathfrak{q}_{i+1}}^{\mathfrak{s}_i+\mathfrak{q}_{i+1}} K_{k-t}
\rightarrow L_{k-t}Z_{\mathfrak{q}_{i+1}}^{\mathfrak{p}_{i+1}} \epsilon_{\mathfrak{s}_i+\mathfrak{q}_{i+1}}^{\mathfrak{q}_{i+1}} L_0Z_{\mathfrak{q}_{i+1}}^{\mathfrak{s}_i+\mathfrak{q}_{i+1}}K_{k-t+1} \cong 0. $$

Suppose $ s=k-t+2.  $ Then the long exact sequence and the induction hypothesis gives
$$ 0 \rightarrow L_{k-t+2}Z_{\mathfrak{q}_{i+1}}^{\mathfrak{p}_{i+1}} \epsilon_{\mathfrak{s}_i+\mathfrak{q}_{i+1}}^{\mathfrak{q}_{i+1}} L_0Z_{\mathfrak{q}_{i+1}}^{\mathfrak{s}_i+\mathfrak{q}_{i+1}} K_{k-t} \rightarrow
L_{k-t+1}Z_{\mathfrak{q}_{i+1}}^{\mathfrak{p}_{i+1}} \epsilon_{\mathfrak{s}_i+\mathfrak{q}_{i+1}}^{\mathfrak{q}_{i+1}} L_0Z_{\mathfrak{q}_{i+1}}^{\mathfrak{s}_i+\mathfrak{q}_{i+1}}K_{k-t+1}
\cong 0. $$
\end{proof}

\begin{corollary}
\label{corollary5}
If $ s=k-l, k-l+2, \ldots, k+l-2 $,
$$ L_{s} Z_{\mathfrak{q}_{i+1}}^{\mathfrak{p}_{i+1}} \epsilon_{\mathfrak{s}_i+\mathfrak{q}_{i+1}}^{\mathfrak{q}_{i+1}} L_0Z_{\mathfrak{q}_{i+1}}^{\mathfrak{s}_i+\mathfrak{q}_{i+1}} K_{k-l-1} \cong M^{\mathfrak{p}_{i+1}}(a_i, k-1, \ldots, 0) $$ Otherwise it is zero.
\end{corollary}

\begin{lemma}
\label{lemma13}
If $ s = k-l-1, k-l+1, \ldots, k+l-1 $ then
$$ L_sZ_{\mathfrak{q}_{i+1}}^{\mathfrak{p}_{i+1}} \epsilon_{\mathfrak{s}_i+\mathfrak{q}_{i+1}}^{\mathfrak{q}_{i+1}}(L_1Z_{\mathfrak{q}_{i+1}}^{\mathfrak{s}_i+\mathfrak{q}_{i+1}}M^{\mathfrak{q}_{i+1}}(\sigma_1 \cdots \sigma_{k-l-2}.\alpha)/
L_1Z_{\mathfrak{q}_{i+1}}^{\mathfrak{s}_i+\mathfrak{q}_{i+1}} K_{k-l-2}) $$
$$ \cong M^{\mathfrak{p}_{i+1}}(a_i, k-1, \ldots, 0). $$ 
Otherwise it is zero.
\end{lemma}

\begin{proof}
From the previous corollary we know that 
$$ L_sZ_{\mathfrak{q}_{i+1}}^{\mathfrak{p}_{i+1}} \epsilon_{\mathfrak{s}_i+\mathfrak{q}_{i+1}}^{\mathfrak{q}_{i+1}} L_0 Z_{\mathfrak{q}_{i+1}}^{\mathfrak{s}_i+\mathfrak{q}_{i+1}} K_{k-l-1} \cong $$
\begin{eqnarray*}
M^{\mathfrak{p}_{i+1}}(a_i, k-1, \ldots 0) 	&\textrm{ if }	&s = k-l, k-l+2, \ldots, k+l-2\\ 
0	&\textrm{ if } 	&s\neq k-l, k-l+2, \ldots, k+l-2.
\end{eqnarray*}

By lemma ~\ref{lemma10},
$ L_1 Z_{\mathfrak{q}_{i+1}}^{\mathfrak{s}_i+\mathfrak{q}_{i+1}} K_{k-l-2} \cong L_0 Z_{\mathfrak{q}_{i+1}}^{\mathfrak{s}_i+\mathfrak{q}_{i+1}} K_{k-l-1}. $
Thus 
$$ L_sZ_{\mathfrak{q}_{i+1}}^{\mathfrak{p}_{i+1}} \epsilon_{\mathfrak{s}_i+\mathfrak{q}_{i+1}}^{\mathfrak{q}_{i+1}} L_1 Z_{\mathfrak{q}_{i+1}}^{\mathfrak{s}_i+\mathfrak{q}_{i+1}} K_{k-l-2} \cong $$
\begin{eqnarray*}
M^{\mathfrak{p}_{i+1}}(a_i, k-1, \ldots 0)	&\textrm{ if }	& s=k-l, k-l+2, \ldots, k+l-2\\
0	&\textrm{ if }	& s \neq k-l, k-l+2, \ldots, k+l-2.
\end{eqnarray*}

Next we compute $ L_sZ_{\mathfrak{q}_{i+1}}^{\mathfrak{p}_{i+1}} \epsilon_{\mathfrak{s}_i+\mathfrak{q}_{i+1}}^{\mathfrak{q}_{i+1}} L_1 Z_{\mathfrak{q}_{i+1}}^{\mathfrak{s}_i+\mathfrak{q}_{i+1}} M^{\mathfrak{q}_{i+1}}(\sigma_1 \cdots \sigma_{k-l-2}.\alpha). $ It is
isomorphic to
$$ L_sZ_{\mathfrak{q}_{i+1}}^{\mathfrak{p}_{i+1}} \epsilon_{\mathfrak{s}_i+\mathfrak{q}_{i+1}}^{\mathfrak{q}_{i+1}} M^{\mathfrak{s}_i+\mathfrak{q}_{i+1}}(l+1, a_i, k-1, \ldots l+2, l, \ldots, 0). $$
In order to handle this, consider the short exact sequence
$$ 
M^{\mathfrak{q}_{i+1}}(l+1, k-1, \ldots, l+2, l, \ldots, 0) \hookrightarrow
M^{\mathfrak{q}_{i+1}}(l+1, a_i, k-1, \ldots, l+2, l, \ldots, 0) \twoheadrightarrow
M^{\mathfrak{s}_i+\mathfrak{q}_{i+1}}(l+1, a_i, k-1, \ldots, l+2, l, \ldots, 0).  $$
The long exact sequence for the functor $ LZ_{\mathfrak{q}_{i+1}}^{\mathfrak{p}_{i+1}} $ gives
$$ L_sZ_{\mathfrak{q}_{i+1}}^{\mathfrak{p}_{i+1}} \epsilon_{\mathfrak{s}_i+\mathfrak{q}_{i+1}}^{\mathfrak{q}_{i+1}} L_1 Z_{\mathfrak{q}_{i+1}}^{\mathfrak{s}_i+\mathfrak{q}_{i+1}} M^{\mathfrak{q}_{i+1}}(\sigma_1 \cdots \sigma_{k-l-2}.\alpha) \cong $$
\begin{eqnarray*}
M^{\mathfrak{p}_{i+1}}(a_i, k-1, \ldots, 0) 	&\textrm{ if }	&s=k-l-1\\
0	&\textrm { if }	&s \neq k-l-1.
\end{eqnarray*}

Now consider the short exact sequence
$$ 
0 \rightarrow \epsilon_{\mathfrak{s}_i+\mathfrak{q}_{i+1}}^{\mathfrak{q}_{i+1}} L_1 Z_{\mathfrak{q}_{i+1}}^{\mathfrak{s}_i+\mathfrak{q}_{i+1}} K_{k-l-2} \rightarrow \epsilon_{\mathfrak{s}_i+\mathfrak{q}_{i+1}}^{\mathfrak{q}_{i+1}} L_1 Z_{\mathfrak{q}_{i+1}}^{\mathfrak{s}_i+\mathfrak{q}_{i+1}} 
M^{\mathfrak{q}_{i+1}}(\sigma_1 \cdots \sigma_{k-l-2}.\alpha) \rightarrow
\epsilon_{\mathfrak{s}_i+\mathfrak{q}_{i+1}}^{\mathfrak{q}_{i+1}} L_1 Z_{\mathfrak{q}_{i+1}}^{\mathfrak{s}_i+\mathfrak{q}_{i+1}} K_{k-l-3} \rightarrow 0. $$
This gives rise to a long exact sequence for $ L Z_{\mathfrak{q}_{i+1}}^{\mathfrak{p}_{i+1}}. $
The long exact sequence together with the first two paragraphs of the proof easily give for $ s \geq k-l+1 $
$$ L_sZ_{\mathfrak{q}_{i+1}}^{\mathfrak{p}_{i+1}} \epsilon_{\mathfrak{s}_i+\mathfrak{q}_{i+1}}^{\mathfrak{q}_{i+1}} L_1 Z_{\mathfrak{q}_{i+1}}^{\mathfrak{s}_i+\mathfrak{q}_{i+1}} K_{k-l-3} \cong 
M^{\mathfrak{p}_{i+1}}(a_i, k-1, \ldots, 0) $$ 
if $ s = k-l+1, \ldots, k+l-1 $ and zero otherwise.

If $ s = k-l $ the long exact sequence gives
$$ L_sZ_{\mathfrak{q}_{i+1}}^{\mathfrak{p}_{i+1}} \epsilon_{\mathfrak{s}_i+\mathfrak{q}_{i+1}}^{\mathfrak{q}_{i+1}} L_1 Z_{\mathfrak{q}_{i+1}}^{\mathfrak{s}_i+\mathfrak{q}_{i+1}} K_{k-l-3} \cong 0. $$

If $ s = k-l-1, $ the long exact sequence gives
$$ L_sZ_{\mathfrak{q}_{i+1}}^{\mathfrak{p}_{i+1}} \epsilon_{\mathfrak{s}_i+\mathfrak{q}_{i+1}}^{\mathfrak{q}_{i+1}} L_1 Z_{\mathfrak{q}_{i+1}}^{\mathfrak{s}_i+\mathfrak{q}_{i+1}} K_{k-l-3} \cong
M^{\mathfrak{p}_{i+1}}(a_i, k-1, \ldots, 0). $$

If $ s < k-l-1, $ then
$ L_sZ_{\mathfrak{q}_{i+1}}^{\mathfrak{p}_{i+1}} \epsilon_{\mathfrak{s}_i+\mathfrak{q}_{i+1}}^{\mathfrak{q}_{i+1}} L_1 Z_{\mathfrak{q}_{i+1}}^{\mathfrak{s}_i+\mathfrak{q}_{i+1}} K_{k-l-3} \cong 0. $
\end{proof}

\begin{lemma}
\label{lemma14}
Suppose $ l+3 \leq t \leq k+1. $  If $ s = k-t+2, k-t+4, \ldots, k+t-4 $ then
$$ L_sZ_{\mathfrak{q}_{i+1}}^{\mathfrak{p}_{i+1}} \epsilon_{\mathfrak{s}_i+\mathfrak{q}_{i+1}}^{\mathfrak{q}_{i+1}}(L_1Z_{\mathfrak{q}_{i+1}}^{\mathfrak{s}_i+\mathfrak{q}_{i+1}}M^{\mathfrak{q}_{i+1}}(\sigma_1 \cdots \sigma_{k-t+1}.\alpha)/
L_1Z_{\mathfrak{q}_{i+1}}^{\mathfrak{s}_i+\mathfrak{q}_{i+1}} K_{k-t+1}) \cong M^{\mathfrak{p}_{i+1}}(a_i, k-1, \ldots, 0). $$ 
Otherwise it is zero.
\end{lemma}

\begin{proof}
We proceed by induction.  The base case is lemma ~\ref{lemma13}.

As in lemma ~\ref{lemma13}, we compute
$$ \aligned
&L_sZ_{\mathfrak{q}_{i+1}}^{\mathfrak{p}_{i+1}} \epsilon_{\mathfrak{s}_i+\mathfrak{q}_{i+1}}^{\mathfrak{q}_{i+1}} L_1 Z_{\mathfrak{q}_{i+1}}^{\mathfrak{s}_i+\mathfrak{q}_{i+1}} M^{\mathfrak{q}_{i+1}}(\sigma_1 \cdots \sigma_{k-t+1}.\alpha) \cong\\
&L_sZ_{\mathfrak{q}_{i+1}}^{\mathfrak{p}_{i+1}} \epsilon_{\mathfrak{s}_i+\mathfrak{q}_{i+1}}^{\mathfrak{q}_{i+1}} L_1 Z_{\mathfrak{q}_{i+1}}^{\mathfrak{s}_i+\mathfrak{q}_{i+1}} M^{\mathfrak{q}_{i+1}}(a_1, \ldots, a_i, t-2, k-1, \ldots, t-1, t-3, \ldots, 0, \ldots, a_n) \cong\\
&L_sZ_{\mathfrak{q}_{i+1}}^{\mathfrak{p}_{i+1}} \epsilon_{\mathfrak{s}_i+\mathfrak{q}_{i+1}}^{\mathfrak{q}_{i+1}} M^{\mathfrak{s}_i+\mathfrak{q}_{i+1}}(t-2, a_i, k-1, \ldots, t-1, t-3, \ldots 0). \endaligned $$

Now consider the short exact sequence
$$ \aligned
0 \rightarrow &M^{\mathfrak{q}_{i+1}}(a_i, t-2, k-1, \ldots, t-1, t-3, \ldots, 0) \rightarrow\\
&M^{\mathfrak{q}_{i+1}}(t-2, a_i, k-1, \ldots t-1, t-3, \ldots, 0) \rightarrow\\
&M^{\mathfrak{s}_i+\mathfrak{q}_{i+1}}(t-2, a_i, k-1, \ldots, t-1, t-3, \ldots, 0) \rightarrow 0.  
\endaligned $$

The long exact sequence for $ LZ_{\mathfrak{q}_{i+1}}^{\mathfrak{p}_{i+1}} $ gives
$$ L_sZ_{\mathfrak{q}_{i+1}}^{\mathfrak{p}_{i+1}} \epsilon_{\mathfrak{s}_i+\mathfrak{q}_{i+1}}^{\mathfrak{q}_{i+1}} L_1 Z_{\mathfrak{q}_{i+1}}^{\mathfrak{s}_i+\mathfrak{q}_{i+1}} M^{\mathfrak{q}_{i+1}}(\sigma_1 \cdots \sigma_{k-t+1}.\alpha) \cong
M^{\mathfrak{p}_{i+1}}(a_1, \ldots, a_i, k-1, \ldots, 0, \ldots, a_n) $$ 
if $ s = k-t+2 $ and zero otherwise.

Now consider the short exact sequence
$$ 
0 \rightarrow \epsilon_{\mathfrak{s}_i+\mathfrak{q}_{i+1}}^{\mathfrak{q}_{i+1}} L_1 Z_{\mathfrak{q}_{i+1}}^{\mathfrak{s}_i+\mathfrak{q}_{i+1}} K_{k-t+1} \rightarrow \epsilon_{\mathfrak{s}_i+\mathfrak{q}_{i+1}}^{\mathfrak{q}_{i+1}} L_1 Z_{\mathfrak{q}_{i+1}}^{\mathfrak{s}_i+\mathfrak{q}_{i+1}} 
M^{\mathfrak{q}_{i+1}}(\sigma_1 \cdots \sigma_{k-t+1}.\alpha) \rightarrow
\epsilon_{\mathfrak{s}_i+\mathfrak{q}_{i+1}}^{\mathfrak{q}_{i+1}} L_1 Z_{\mathfrak{q}_{i+1}}^{\mathfrak{s}_i+\mathfrak{q}_{i+1}} K_{k-t} \rightarrow 0. $$
This gives rise to a long exact sequence for $ LZ_{\mathfrak{q}_{i+1}}^{\mathfrak{p}_{i+1}}. $

Now by the induction hypothesis, the long exact sequence and the first paragraph, the lemma is true for $ s \geq k-t+4. $

Suppose $ s = k-t+3. $  Then the long exact sequence becomes
$$ 0 \rightarrow L_{k-t+3} Z_{\mathfrak{q}_{i+1}}^{\mathfrak{p}_{i+1}} \epsilon_{\mathfrak{s}_i+\mathfrak{q}_{i+1}}^{\mathfrak{q}_{i+1}} L_1 Z_{\mathfrak{q}_{i+1}}^{\mathfrak{s}_i+\mathfrak{q}_{i+1}} K_{k-t} \rightarrow
L_{k-t+2} Z_{\mathfrak{q}_{i+1}}^{\mathfrak{p}_{i+1}} \epsilon_{\mathfrak{s}_i+\mathfrak{q}_{i+1}}^{\mathfrak{q}_{i+1}} L_1 Z_{\mathfrak{q}_{i+1}}^{\mathfrak{s}_i+\mathfrak{q}_{i+1}} K_{k-t+1} \rightarrow \cdots. $$
By the induction hypothesis 
$$ L_{k-t+2} Z_{\mathfrak{q}_{i+1}}^{\mathfrak{p}_{i+1}} \epsilon_{\mathfrak{s}_i+\mathfrak{q}_{i+1}}^{\mathfrak{q}_{i+1}} L_1 Z_{\mathfrak{q}_{i+1}}^{\mathfrak{s}_i+\mathfrak{q}_{i+1}} K_{k-t+1} \cong 0, $$
so 
$$ L_{k-t+3} Z_{\mathfrak{q}_{i+1}}^{\mathfrak{p}_{i+1}} \epsilon_{\mathfrak{s}_i+\mathfrak{q}_{i+1}}^{\mathfrak{q}_{i+1}} L_1 Z_{\mathfrak{q}_{i+1}}^{\mathfrak{s}_i+\mathfrak{q}_{i+1}} K_{k-t} \cong 0. $$

Suppose $ s = k-t+2. $ Then the long exact sequence becomes
$$ \aligned
0 \rightarrow &L_{k-t+2} Z_{\mathfrak{q}_{i+1}}^{\mathfrak{p}_{i+1}} \epsilon_{\mathfrak{s}_i+\mathfrak{q}_{i+1}}^{\mathfrak{q}_{i+1}} L_1 Z_{\mathfrak{q}_{i+1}}^{\mathfrak{s}_i+\mathfrak{q}_{i+1}} K_{k-t+1} \rightarrow
L_{k-t+2} Z_{\mathfrak{q}_{i+1}}^{\mathfrak{p}_{i+1}} \epsilon_{\mathfrak{s}_i+\mathfrak{q}_{i+1}}^{\mathfrak{q}_{i+1}} L_1 Z_{\mathfrak{q}_{i+1}}^{\mathfrak{s}_i+\mathfrak{q}_{i+1}} M^{\mathfrak{q}_{i+1}}(\sigma_1 \cdots \sigma_{k-t+1}.\alpha) \rightarrow\\
&L_{k-t+2} Z_{\mathfrak{q}_{i+1}}^{\mathfrak{p}_{i+1}} \epsilon_{\mathfrak{s}_i+\mathfrak{q}_{i+1}}^{\mathfrak{q}_{i+1}} L_1 Z_{\mathfrak{q}_{i+1}}^{\mathfrak{s}_i+\mathfrak{q}_{i+1}} K_{k-t} 
\rightarrow L_{k-t+1} Z_{\mathfrak{q}_{i+1}}^{\mathfrak{p}_{i+1}} \epsilon_{\mathfrak{s}_i+\mathfrak{q}_{i+1}}^{\mathfrak{q}_{i+1}} L_1 Z_{\mathfrak{q}_{i+1}}^{\mathfrak{s}_i+\mathfrak{q}_{i+1}} K_{k-t+1} \rightarrow \cdots.  
\endaligned $$
By the induction hypothesis
$$ L_{k-t+2} Z_{\mathfrak{q}_{i+1}}^{\mathfrak{p}_{i+1}} \epsilon_{\mathfrak{s}_i+\mathfrak{q}_{i+1}}^{\mathfrak{q}_{i+1}} L_1 Z_{\mathfrak{q}_{i+1}}^{\mathfrak{s}_i+\mathfrak{q}_{i+1}} K_{k-t+1} \cong
L_{k-t+1} Z_{\mathfrak{q}_{i+1}}^{\mathfrak{p}_{i+1}} \epsilon_{\mathfrak{s}_i+\mathfrak{q}_{i+1}}^{\mathfrak{q}_{i+1}} L_1 Z_{\mathfrak{q}_{i+1}}^{\mathfrak{s}_i+\mathfrak{q}_{i+1}} K_{k-t+1} \cong 0. $$
Thus 
$$ \aligned
L_{k-t+2} Z_{\mathfrak{q}_{i+1}}^{\mathfrak{p}_{i+1}} \epsilon_{\mathfrak{s}_i+\mathfrak{q}_{i+1}}^{\mathfrak{q}_{i+1}} L_1 Z_{\mathfrak{q}_{i+1}}^{\mathfrak{s}_i+\mathfrak{q}_{i+1}} K_{k-t} &\cong
L_{k-t+2} Z_{\mathfrak{q}_{i+1}}^{\mathfrak{p}_{i+1}} \epsilon_{\mathfrak{s}_i+\mathfrak{q}_{i+1}}^{\mathfrak{q}_{i+1}} L_1 Z_{\mathfrak{q}_{i+1}}^{\mathfrak{s}_i+\mathfrak{q}_{i+1}} M^{\mathfrak{q}_{i+1}}(\sigma_1 \cdots \sigma_{k-t+1}.\alpha)\\ 
&\cong M^{\mathfrak{p}_{i+1}}(a_i, k-1, \ldots, 0).  \endaligned $$
Thus the claim is verified for $ s = k-t+2. $

Finally suppose $ s \leq k-t+1. $  Then
$$ \aligned
&L_{s} Z_{\mathfrak{q}_{i+1}}^{\mathfrak{p}_{i+1}} \epsilon_{\mathfrak{s}_i+\mathfrak{q}_{i+1}}^{\mathfrak{q}_{i+1}} L_1 Z_{\mathfrak{q}_{i+1}}^{\mathfrak{s}_i+\mathfrak{q}_{i+1}} M^{\mathfrak{q}_{i+1}}(\sigma_1 \cdots \sigma_{k-t+1}.\alpha) \cong\\
&L_{s-1} Z_{\mathfrak{q}_{i+1}}^{\mathfrak{p}_{i+1}} \epsilon_{\mathfrak{s}_i+\mathfrak{q}_{i+1}}^{\mathfrak{q}_{i+1}} L_1 Z_{\mathfrak{q}_{i+1}}^{\mathfrak{s}_i+\mathfrak{q}_{i+1}} M^{\mathfrak{q}_{i+1}}(\sigma_1 \cdots \sigma_{k-t+1}.\alpha) \cong 0. 
\endaligned $$
This case is then true by the induction hypothesis.
\end{proof}

\begin{corollary}
\label{corollary6}
If $ s = 1, 3, \ldots, 2(k-1)-1 $ then
$$ 
L_sZ_{\mathfrak{q}_{i+1}}^{\mathfrak{p}_{i+1}} \epsilon_{\mathfrak{s}_i+\mathfrak{q}_{i+1}}^{\mathfrak{q}_{i+1}} L_1 Z_{\mathfrak{q}_{i+1}}^{\mathfrak{s}_i+\mathfrak{q}_{i+1}} \epsilon_{\mathfrak{p}_{i+1}}^{\mathfrak{q}_{i+1}} 
M^{\mathfrak{p}_{i+1}}(a_i, k-1, \ldots, 0) \cong
M^{\mathfrak{p}_{i+1}}(a_i, k-1, \ldots, 0) $$ 
and zero otherwise.
\end{corollary}

\begin{proof}
This is a direct consequence of the previous lemma.
\end{proof}

\begin{corollary}
For $ s = 0, 2, \ldots, 2(k-1)-2 $ 
$$ 
L_sZ_{\mathfrak{q}_{i+1}}^{\mathfrak{p}_{i+1}} \epsilon_{\mathfrak{s}_i+\mathfrak{q}_{i+1}}^{\mathfrak{q}_{i+1}} LZ_{\mathfrak{q}_{i+1}}^{\mathfrak{s}_i+\mathfrak{q}_{i+1}} \epsilon_{\mathfrak{p}_{i+1}}^{\mathfrak{q}_{i+1}} 
M^{\mathfrak{p}_{i+1}}(a_i, k-1, \ldots, 0) \cong 
M^{\mathfrak{p}_{i+1}}(a_i, k-1, \ldots 0). $$
It is zero otherwise.
\end{corollary}

\begin{proof}
This follows from the previous corollary and noting that 
$$ 
LZ_{\mathfrak{q}_{i+1}}^{\mathfrak{s}_i+\mathfrak{q}_{i+1}} \epsilon_{\mathfrak{p}_{i+1}}^{\mathfrak{q}_{i+1}} M^{\mathfrak{p}_{i+1}}(a_i, k-1, \ldots, 0) \cong
L_{1}Z_{\mathfrak{q}_{i+1}}^{\mathfrak{s}_i+\mathfrak{q}_{i+1}} \epsilon_{\mathfrak{p}_{i+1}}^{\mathfrak{q}_{i+1}} M^{\mathfrak{p}_{i+1}}(a_i, k-1, \ldots, 0)[1].  $$
\end{proof}

By adjointness, there is an isomorphism 
$$ \Hom(LZ_{\mathfrak{q}_{i+1}}^{\mathfrak{s}_i+\mathfrak{q}_{i+1}} \epsilon_{\mathfrak{p}_{i+1}}^{\mathfrak{q}_{i+1}}, LZ_{\mathfrak{q}_{i+1}}^{\mathfrak{s}_i+\mathfrak{q}_{i+1}} \epsilon_{\mathfrak{p}_{i+1}}^{\mathfrak{q}_{i+1}}) \cong
\Hom(\epsilon_{\mathfrak{p}_{i+1}}^{\mathfrak{q}_{i+1}}, \epsilon_{\mathfrak{s}_i+\mathfrak{q}_{i+1}}^{\mathfrak{q}_{i+1}} LZ_{\mathfrak{q}_{i+1}}^{\mathfrak{s}_i+\mathfrak{q}_{i+1}} \epsilon_{\mathfrak{p}_{i+1}}^{\mathfrak{q}_{i+1}}). $$

Call the image of the identity morphism under this isomorphism $ \phi. $  Then 
$$ LZ_{\mathfrak{q}_{i+1}}^{\mathfrak{p}_{i+1}} \phi \colon
LZ_{\mathfrak{q}_{i+1}}^{\mathfrak{p}_{i+1}} \epsilon_{\mathfrak{p}_{i+1}}^{\mathfrak{q}_{i+1}} \rightarrow LZ_{\mathfrak{q}_{i+1}}^{\mathfrak{p}_{i+1}} \epsilon_{\mathfrak{s}_i+\mathfrak{q}_{i+1}}^{\mathfrak{q}_{i+1}} LZ_{\mathfrak{q}_{i+1}}^{\mathfrak{s}_i+\mathfrak{q}_{i+1}} \epsilon_{\mathfrak{p}_{i+1}}^{\mathfrak{q}_{i+1}}. $$
By proposition ~\ref{prop8}, we know 
$ LZ_{\mathfrak{q}_{i+1}}^{\mathfrak{p}_{i+1}} \epsilon_{\mathfrak{p}_{i+1}}^{\mathfrak{q}_{i+1}} \cong \oplus^{k-1}_{r=0} \Id [2r]. $

Finally we have an isomorphism of functors lifting the third graphical relation.

\begin{prop}
\label{prop10}
There is an isomorphism
$$ \tau^{\leq k-2}LZ_{\mathfrak{q}_{i+1}}^{\mathfrak{p}_{i+1}}[-k] \phi \colon
\oplus^{k-1}_{r=1} \Id [2r-k] \rightarrow LZ_{\mathfrak{q}_{i+1}}^{\mathfrak{p}_{i+1}} \epsilon_{\mathfrak{s}_i+\mathfrak{q}_{i+1}}^{\mathfrak{q}_{i+1}}[-1] LZ_{\mathfrak{q}_{i+1}}^{\mathfrak{s}_i+\mathfrak{q}_{i+1}} \epsilon_{\mathfrak{p}_{i+1}}^{\mathfrak{q}_{i+1}}[-(k-1)]  . $$
\end{prop}

\begin{proof}
We need to show that the cohomology functors of $ \tau^{\leq k-2}LZ_{\mathfrak{q}_{i+1}}^{\mathfrak{p}_{i+1}} \phi $ when applied to the generalized Verma module are
isomorphisms.
In the previous lemmas we needed to compute
$$ LZ_{\mathfrak{q}_{i+1}}^{\mathfrak{p}_{i+1}} \epsilon_{\mathfrak{s}_{i}+\mathfrak{q}_{i+1}}^{\mathfrak{q}_{i+1}}
M^{\mathfrak{s}_{i}+\mathfrak{q}_{i+1}}(\underbrace{a_i, a_{i+1},} \underbrace{a_{i+2}, \ldots, a_{i+k}}). $$
The calculation was made by considering the short exact sequence
$$ 
0 \rightarrow M^{\mathfrak{q}_{i+1}}(a_{i+1}, a_i, \underbrace{a_{i+2}, \ldots, a_{i+k}}) \rightarrow 
M^{\mathfrak{q}_{i+1}}(a_{i}, a_{i+1}, \underbrace{a_{i+2}, \ldots, a_{i+k}}) \rightarrow
M^{\mathfrak{s}_{i}+\mathfrak{q}_{i+1}}(\underbrace{a_i, a_{i+1},} \underbrace{a_{i+2}, \ldots, a_{i+k}}) \rightarrow 0. 
$$

It is clear that if $ \lbrace a_{i+1}, \ldots, a_{i+k} \rbrace = \lbrace k-1, \ldots, 0 \rbrace, $
then the isomorphism
$$
LZ_{\mathfrak{q}_{i+1}}^{\mathfrak{p}_{i+1}} \epsilon_{\mathfrak{s}_{i}+\mathfrak{q}_{i+1}}^{\mathfrak{q}_{i+1}}
M^{\mathfrak{s}_{i}+\mathfrak{q}_{i+1}}(\underbrace{a_i, a_{i+1},} \underbrace{a_{i+2}, \ldots, a_{i+k}}) \cong
M^{\mathfrak{p}_{i+1}}(a_{i}, \underbrace{k-1, \ldots, 0})[k-1-a_{i+1}] 
$$
comes from the natural transformation $ \Id \rightarrow \epsilon_{\mathfrak{s}_i+\mathfrak{q}_{i+1}}^{\mathfrak{q}_{i+1}} LZ_{\mathfrak{q}_{i+1}}^{\mathfrak{s}_i+\mathfrak{q}_{i+1}}. $
If instead it is the case that $ \lbrace a_{i}, a_{i+2}, \ldots, a_{i+k} \rbrace = \lbrace k-1, \ldots, 0 \rbrace, $ we need to demonstrate the isomorphism in a different way so that it is natural.

Consider the short exact sequence in the category $ \mathcal{O}(\mathfrak{gl}_{k+1}) $
$$ 0 \rightarrow M^{\mathfrak{q}_{2}}(a_{i+1}, a_i, \underbrace{a_{i+2}, \ldots, a_{i+k},}) \rightarrow 
M^{\mathfrak{q}_{2}}(a_{i}, a_{i+1}, \underbrace{a_{i+2}, \ldots, a_{i+k},}) \rightarrow
M^{\mathfrak{s}_{1}+\mathfrak{q}_{2}}(\underbrace{a_i, a_{i+1},} \underbrace{a_{i+2}, \ldots, a_{i+k},}) \rightarrow 0. $$
In this category, the objects at the beginning and end of the sequence are self dual with respect to the duality functor $ d $ because they are simple.
Thus there is an exact sequence
$$ M^{\mathfrak{s}_{1}+\mathfrak{q}_{2}}(\underbrace{a_i, a_{i+1},} \underbrace{a_{i+2}, \ldots, a_{i+k},}) \hookrightarrow
d M^{\mathfrak{q}_{2}}(a_{i}, a_{i+1}, \underbrace{a_{i+2}, \ldots, a_{i+k},}) \twoheadrightarrow
M^{\mathfrak{q}_{2}}(a_{i+1}, a_i, \underbrace{a_{i+2}, \ldots, a_{i+k},}). $$
Now we take an external tensor product on the left with the Verma module
$ M(a_1, \ldots, a_{i-1}) $ and on the right with $ M(a_{i+k+1}, \ldots, a_n). $
This will give an exact sequence in $ \mathcal{O}(\mathfrak{gl}_{i-1} \oplus \mathfrak{gl}_{k+1} \oplus \mathfrak{gl}_{n-k-i}). $  Then we may induce it to the following exact sequence in
$ \mathcal{O}(\mathfrak{gl}_{n}): $
$$ M^{\mathfrak{s}_{i}+\mathfrak{q}_{i+1}}(a_1, \ldots, a_{i-1}, \underbrace{a_i, a_{i+1},} \underbrace{a_{i+2}, \ldots, a_{i+k},} \ldots, a_n)  \hookrightarrow X \twoheadrightarrow
M^{\mathfrak{q}_{i+1}}(a_1, \ldots, a_{i-1}, a_{i+1}, a_i, \underbrace{a_{i+2}, \ldots, a_{i+k},} \ldots, a_n). $$
The object $ X $ is obtained by taking an external tensor product of Verma modules with a dual generalized Verma module followed by parabolic induction.  Since 
$ LZ_{\mathfrak{q}_{i+1}}^{\mathfrak{p}_{i+1}} $ naturally commutes with this composition of functors, it suffice to compute 
$$ LZ_{\mathfrak{q}_{2}}^{\mathfrak{p}_{1}} d M^{\mathfrak{q}_{2}}(a_{i}, a_{i+1}, \underbrace{a_{i+2}, \ldots, a_{i+k},}). $$
This is easily seen to be zero.  
Thus the isomorphism
$$
LZ_{\mathfrak{q}_{i+1}}^{\mathfrak{p}_{i+1}} \epsilon_{\mathfrak{s}_{i}+\mathfrak{q}_{i+1}}^{\mathfrak{q}_{i+1}}
M^{\mathfrak{s}_{i}+\mathfrak{q}_{i+1}}(\underbrace{a_i, a_{i+1},} \underbrace{a_{i+2}, \ldots, a_{i+k}}) \cong
M^{\mathfrak{p}_{i+1}}(a_{i}, \underbrace{k-1, \ldots, 0})[k-1-a_{i}] 
$$
also comes from the natural transformation 
$ \Id \rightarrow \epsilon_{\mathfrak{s}_i+\mathfrak{q}_{i+1}}^{\mathfrak{q}_{i+1}} LZ_{\mathfrak{q}_{i+1}}^{\mathfrak{s}_i+\mathfrak{q}_{i+1}}. $

Therefore the natural transformation stated in the proposition is an isomorphism for generalized Verma modules.  Then it must be true for any projective object by induction on the length of its generalized Verma flag.  The theorem is then true for an arbitrary objects by considering a projective resolution of it.
\end{proof}

\subsection{Diagram 4}

For ease of notation, we abbreviate the algebras defined before for this subsection as follows.
\begin{define}
Let
\begin{enumerate}
\item $ \alpha = \mathfrak{p}_i+\mathfrak{q}_{i+k} $
\item $ \beta = \mathfrak{q}_{i-1}+\mathfrak{q}_{i+k} $
\item $ \gamma = \mathfrak{q}_{i-1}+\mathfrak{s}_{i+k-1}+\mathfrak{q}_{i+k} $
\item $ \delta = \mathfrak{q}_{i-1}+\mathfrak{p}_{i+k}. $
\end{enumerate}
\end{define}

We begin with a generalized Verma module
$$ M^{\mathfrak{\alpha}}(a_1, \ldots a_{i-1}, \underbrace{k-1, \ldots 0}, a_{i+k}, \underbrace{a_{i+k+1}, \ldots a_{i+2k-1}}, a_{i+2k}, \ldots a_n). $$
Since the coordinates $ a_1, \ldots, a_{i-1}, a_{i+2k}, \ldots, a_n $ play no role in the computation, we omit them when labeling generalized Verma modules.
Suppose $ a_{i+k} = l. $
Let us denote the module above by $ M^{\mathfrak{\alpha}}(l). $

There are short exact sequences arising from the generalized BGG resolution:

$$ 0 \rightarrow K_0 \rightarrow M^{\beta}(\underbrace{k-1, \ldots, \hat{0}}, 0, a_{i+k}, \underbrace{a_{i+k+1}, \ldots, a_{i+2k-1}}) \rightarrow M^{\mathfrak{\alpha}}(l) \rightarrow 0 $$
$$ 0 \rightarrow K_1 \rightarrow M^{\beta}(\underbrace{k-1, \ldots, \hat{1}, 0}, 1, a_{i+k}, \underbrace{a_{i+k+1}, \ldots, a_{i+2k-1}}) \rightarrow K_0 \rightarrow 0 $$
$$ \cdots $$
$$ 0 \rightarrow K_{l-1} \rightarrow M^{\beta}(\underbrace{k-1, \ldots, \widehat{l-1}, 0}, l-1, a_{i+k}, \underbrace{a_{i+k+1}, \ldots, a_{i+2k-1}}) \rightarrow K_{l-2} \rightarrow 0 $$
$$ 0 \rightarrow K_{l} \rightarrow M^{\beta}(\underbrace{k-1, \ldots, \hat{l}, 0}, l, a_{i+k}, \underbrace{a_{i+k+1}, \ldots, a_{i+2k-1}}) \rightarrow K_{l-1} \rightarrow 0 $$
$$ 0 \rightarrow K_{l+1} \rightarrow M^{\beta}(\underbrace{k-1, \ldots, \widehat{l+1}, 0}, l+1, a_{i+k}, \underbrace{a_{i+k+1}, \ldots, a_{i+2k-1}}) \rightarrow K_{l} \rightarrow 0 $$
$$ \cdots $$
$$ 0 \rightarrow K_{k-1} \rightarrow M^{\beta}(\underbrace{\widehat{k-1}, \ldots, 0}, k-1, a_{i+k}, \underbrace{a_{i+k+1}, \ldots, a_{i+2k-1}}) \rightarrow K_{k-2} \rightarrow 0. $$

\begin{lemma}
There are exact sequences:
$$ 0 \rightarrow L_1Z_{\mathfrak{\beta}}^{\gamma} K_0 \rightarrow M^{\gamma}(\underbrace{k-1, \ldots, \hat{0}}, \underbrace{a_{i+k},0}, \underbrace{a_{i+k+1}, \ldots, a_{i+2k-1}}) \rightarrow 
L_1Z_{\mathfrak{\beta}}^{\mathfrak{\alpha}} \epsilon_{\alpha}^{\beta} M^{\mathfrak{\alpha}}(l) \rightarrow 0 $$
$$ 0 \rightarrow L_1Z_{\mathfrak{\beta}}^{\gamma} K_1 \rightarrow M^{\gamma}(\underbrace{k-1, \ldots, \hat{1}, 0}, \underbrace{a_{i+k},1}, \underbrace{a_{i+k+1}, \ldots, a_{i+2k-1}}) \rightarrow 
L_1 Z_{\mathfrak{\beta}}^{\gamma} K_0 \rightarrow 0 $$
$$ \cdots $$
$$ 0 \rightarrow L_1Z_{\mathfrak{\beta}}^{\gamma} K_{l-1} \rightarrow M^{\gamma}(\underbrace{k-1, \ldots, \widehat{l-1}, 0}, \underbrace{a_{i+k},l-1}, \underbrace{a_{i+k+1}, \ldots, a_{i+2k-1}}) \rightarrow L_1 Z_{\mathfrak{\beta}}^{\gamma} K_{l-2} \rightarrow 0 $$
$$ L_i Z_{\mathfrak{\beta}}^{\gamma} K_l = L_{i+1}Z_{\mathfrak{\beta}}^{\gamma} K_{l-1}, \forall i $$
$$ 0 \rightarrow L_0 Z_{\mathfrak{\beta}}^{\gamma} K_{l+1} \rightarrow M^{\gamma}(\underbrace{k-1, \ldots, \widehat{l+1}, 0}, \underbrace{l+1, a_{i+k}}, \underbrace{a_{i+k+1}, \ldots, a_{i+2k-1}}) \rightarrow L_0 Z_{\mathfrak{\beta}}^{\gamma} K_l \rightarrow 0 $$
$$ \cdots $$
$$ 0 \rightarrow L_0 Z_{\mathfrak{\beta}}^{\gamma} K_{k-1} \rightarrow M^{\gamma}(\underbrace{\widehat{k-1}, \ldots, 0}, \underbrace{k-1, a_{i+k}}, \underbrace{a_{i+k+1}, \ldots, a_{i+2k-1}}) \rightarrow L_0 Z_{\mathfrak{\beta}}^{\gamma} K_{k-2} \rightarrow 0 $$
\end{lemma}

\begin{proof}
In the first set of sequences, $ L_0 Z_{\mathfrak{\beta}}^{\gamma} $ must be zero.  
Also, $ L_1 Z_{\mathfrak{\beta}}^{\gamma} K_l = L_2 Z_{\mathfrak{\beta}}^{\gamma} K_{l-1} = 0. $
Then $ L_1 Z_{\mathfrak{\beta}}^{\gamma} K_{l+1} = L_2 Z_{\mathfrak{\beta}}^{\gamma} K_l = L_3 Z_{\mathfrak{\beta}}^{\gamma} K_{l-1} = 0. $
We also have $ L_1 Z_{\mathfrak{\beta}}^{\gamma} K_{l+2} = L_2 Z_{\mathfrak{\beta}}^{\gamma} K_{l+1} = 0 $ and so on.
This gives the second set of exact sequences.
\end{proof}

\begin{corollary}
\label{corollary8}
There are exact sequences:
$$ 0 \rightarrow \epsilon_{\gamma}^{\beta} L_1 Z_{\mathfrak{\beta}}^{\gamma} K_0 \rightarrow  \epsilon_{\gamma}^{\beta} M^{\gamma}(\underbrace{k-1, \ldots, \hat{0}}, \underbrace{a_{i+k},0} \underbrace{a_{i+k+1}, \ldots, a_{i+2k-1}}) \rightarrow 
 \epsilon_{\gamma}^{\beta} L_1Z_{\mathfrak{\beta}}^{\mathfrak{\alpha}} \epsilon_{\alpha}^{\beta} M^{\mathfrak{\alpha}}(l) \rightarrow 0 $$
$$ 0 \rightarrow  \epsilon_{\gamma}^{\beta} L_1Z_{\mathfrak{\beta}}^{\gamma} K_1 \rightarrow  \epsilon_{\gamma}^{\beta} M^{\gamma}(\underbrace{k-1, \ldots, \hat{1}, 0}, \underbrace{a_{i+k},1} \underbrace{a_{i+k+1}, \ldots, a_{i+2k-1}}) \rightarrow 
 \epsilon_{\gamma}^{\beta} L_1 Z_{\mathfrak{\beta}}^{\gamma} K_0 \rightarrow 0 $$
$$ \cdots $$
$$ 0 \rightarrow \epsilon_{\gamma}^{\beta} L_1Z_{\mathfrak{\beta}}^{\gamma} K_{l-1} \rightarrow  \epsilon_{\gamma}^{\beta} M^{\gamma}(\underbrace{k-1, \ldots, \widehat{l-1}, 0}, \underbrace{a_{i+k},l-1} \underbrace{a_{i+k+1}, \ldots, a_{i+2k-1}}) \rightarrow  \epsilon_{\gamma}^{\beta} L_1 Z_{\mathfrak{\beta}}^{\gamma} K_{l-2} \rightarrow 0 $$
$$  \epsilon_{\gamma}^{\beta} L_i Z_{\mathfrak{\beta}}^{\gamma} K_l =  \epsilon_{\gamma}^{\beta} L_{i+1}Z_{\mathfrak{\beta}}^{\gamma} K_{l-1}, \forall i $$
$$ 0 \rightarrow  \epsilon_{\gamma}^{\beta} L_0 Z_{\mathfrak{\beta}}^{\gamma} K_{l+1} \rightarrow  \epsilon_{\gamma}^{\beta} M^{\gamma}(\underbrace{k-1, \ldots, \widehat{l+1}, 0}, \underbrace{l+1, a_{i+k}} \underbrace{a_{i+k+1}, \ldots, a_{i+2k-1}}) \rightarrow  \epsilon_{\gamma}^{\beta} L_0 Z_{\mathfrak{\beta}}^{\gamma} K_l \rightarrow 0 $$
$$ \cdots $$
$$ 0 \rightarrow  \epsilon_{\gamma}^{\beta} L_0 Z_{\mathfrak{\beta}}^{\gamma} K_{k-1} \rightarrow  \epsilon_{\gamma}^{\beta} M^{\gamma}(\underbrace{\widehat{k-1}, \ldots, 0}, \underbrace{k-1, a_{i+k}} \underbrace{a_{i+k+1}, \ldots, a_{i+2k-1}}) \rightarrow  \epsilon_{\gamma}^{\beta} L_0 Z_{\mathfrak{\beta}}^{\gamma} K_{k-2} \rightarrow 0 $$
\end{corollary}

\begin{lemma}
Let $ \lbrace a_{i+k+1}, \ldots, a_{i+2k-1} \rbrace = \lbrace k-1, \ldots, \hat{m}, \ldots 0 \rbrace. $
Let $ m \neq l. $ If $ s=k-1, $ then
$$ L_s Z_{\mathfrak{\beta}}^{\delta} \epsilon_{\gamma}^{\beta} L_1 Z_{\mathfrak{\beta}}^{\gamma} \epsilon_{\alpha}^{\beta} M^{\mathfrak{\alpha}} =
M^{\delta}(\underbrace{k-1, \ldots, \hat{m}, \ldots, 0,} a_{i+k}, \underbrace{k-1, \ldots, 0}). $$
Otherwise it is zero.
\end{lemma}

\begin{proof}
This follows easily from the previous corollary.
\end{proof}

Now consider 
$$ M^{\delta} = M^{\delta}(\underbrace{a_i, \ldots, a_{i+k-2}}, a_{i+k-1}, \underbrace{k-1, \ldots, 0}) $$ and let $ a_{i+k-1} = l. $
There are short exact sequences

$$ 0 \rightarrow J_0 \rightarrow M^{\beta}(\underbrace{a_i, \ldots, a_{i+k-2}}, a_{i+k-1}, k-1, \underbrace{\widehat{k-1}, \ldots, 0}) \rightarrow M^{\delta}(l) \rightarrow 0 $$
$$ 0 \rightarrow J_1 \rightarrow M^{\beta}(\underbrace{a_i, \ldots, a_{i+k-2}}, a_{i+k-1}, k-2, \underbrace{k-1, \widehat{k-2}, \ldots, 0}) \rightarrow J_0 \rightarrow 0 $$
$$ \cdots $$
$$ 0 \rightarrow J_{k-1} \rightarrow M^{\beta}(\underbrace{a_i, \ldots, a_{i+k-2}}, a_{i+k-1}, 0, \underbrace{k-1, \ldots, \hat{0}}) \rightarrow J_{k-2} \rightarrow 0. $$

\begin{corollary}
\label{corollary9}
There are exact sequences:
$$ 0 \rightarrow \epsilon_{\gamma}^{\beta} L_1 Z_{\mathfrak{\beta}}^{\gamma} J_0 \rightarrow \epsilon_{\gamma}^{\beta} M^{\gamma}(\underbrace{a_i, \ldots, a_{i+k-2}}, \underbrace{k-1, a_{i+k-1}}, \underbrace{\widehat{k-1}, \ldots, 0}) \rightarrow \epsilon_{\gamma}^{\beta} L_1 Z_{\mathfrak{\beta}}^{\gamma} M^{\delta}(l) \rightarrow 0 $$
$$ 0 \rightarrow \epsilon_{\gamma}^{\beta} L_1 Z_{\mathfrak{\beta}}^{\gamma} J_1 \rightarrow \epsilon_{\gamma}^{\beta} M^{\gamma}(\underbrace{a_i, \ldots, a_{i+k-2}}, \underbrace{k-2, a_{i+k-1}}, \underbrace{k-1, \widehat{k-2}, \ldots, 0}) \rightarrow \epsilon_{\gamma}^{\beta} L_1 Z_{\mathfrak{\beta}}^{\gamma} J_0 \rightarrow 0 $$
$$ \cdots $$
$$ 0 \rightarrow \epsilon_{\gamma}^{\beta} L_1 Z_{\mathfrak{\beta}}^{\gamma} J_{k-l-2} \rightarrow \epsilon_{\gamma}^{\beta} M^{\gamma}(\underbrace{a_i, \ldots, a_{i+k-2}}, \underbrace{l+1, a_{i+k-1}}, \underbrace{k-1, \widehat{l+1}, \ldots, 0}) \rightarrow \epsilon_{\gamma}^{\beta} L_1 Z_{\mathfrak{\beta}}^{\gamma} J_{k-l-3} \rightarrow 0 $$
$$ \epsilon_{\gamma}^{\beta} L_i Z_{\mathfrak{\beta}}^{\gamma} J_{k-l-1} \cong \epsilon_{\gamma}^{\beta} L_{i+1} Z_{\mathfrak{\beta}}^{\gamma} L_{k-l-2}, \forall i $$
$$ 0 \rightarrow \epsilon_{\gamma}^{\beta} L_0 Z_{\mathfrak{\beta}}^{\gamma} J_{k-l} \rightarrow \epsilon_{\gamma}^{\beta} M^{\gamma}(\underbrace{a_i, \ldots, a_{i+k-2}}, \underbrace{a_{i+k-1}, l-1}, \underbrace{k-1, \widehat{l-1}, \ldots, 0}) \rightarrow \epsilon_{\gamma}^{\beta} L_0 Z_{\mathfrak{\beta}}^{\gamma} J_{k-l-1} \rightarrow 0 $$
$$ \cdots $$
$$ 0 \rightarrow \epsilon_{\gamma}^{\beta} L_0 Z_{\mathfrak{\beta}}^{\gamma} J_{k-1} \rightarrow \epsilon_{\gamma}^{\beta} M^{\gamma}(\underbrace{a_i, \ldots, a_{i+k-2}}, \underbrace{a_{i+k-1}, 0}, \underbrace{k-1, \ldots, \hat{0}}) \rightarrow \epsilon_{\gamma}^{\beta} L_0 Z_{\mathfrak{\beta}}^{\gamma} J_{k-2} \rightarrow 0 $$
\end{corollary}

\begin{proof}
This is similar to corollary ~\ref{corollary8}.
\end{proof}

\begin{lemma}
Let $ \lbrace a_{i}, \ldots, a_{i+k-2} \rbrace = \lbrace k-1, \ldots, \hat{m}, \ldots, 0 \rbrace. $
Let $ m \neq l. $ If $ s=k-1, $ then
$$ L_s Z_{\mathfrak{\beta}}^{\mathfrak{\alpha}} \epsilon_{\gamma}^{\beta} L_1 Z_{\mathfrak{\beta}}^{\gamma} \epsilon_{\delta}^{\beta} M^{\delta}(l) =
M^{\mathfrak{\alpha}}(\underbrace{k-1, \ldots, 0,} a_{i+k-1}, \underbrace{k-1, \ldots, \hat{m} \ldots, 0}). $$
Otherwise it is zero.
\end{lemma}

\begin{proof}
This follows from corollary ~\ref{corollary9}.
\end{proof}

We must now study the case when $ l = m. $

\begin{lemma}
\label{lemma19}
Let $ m=l. $  Then
$ L_s Z_{\mathfrak{\beta}}^{\delta} \epsilon_{\gamma}^{\beta} L_1 Z_{\beta}^{\gamma} \epsilon_{\alpha}^{\beta} M^{\mathfrak{\alpha}}(l) \cong $
\begin{eqnarray*} 
M^{\delta}(0)/\ldots /M^{\delta}(k-1) & \textrm{ if } & s=k-l \\
L_{k-l+1}Z_{\mathfrak{\beta}}^{\gamma} \epsilon_{\gamma}^{\beta}L_1 Z_{\mathfrak{\beta}}^{\gamma} K_{l-3}  & \textrm{ if } & s=k-1
\end{eqnarray*}
and there is an exact sequence
$$ 0 \rightarrow M^{\delta}(l)/\ldots/ M^{\delta}(k-1) \rightarrow L_{k-l+1}Z_{\mathfrak{\beta}}^{\gamma} \epsilon_{\gamma}^{\beta}L_1 Z_{\mathfrak{\beta}}^{\gamma} K_{l-3}
\rightarrow M^{\delta}(l+1)/\ldots/ M^{\delta}(k-1) \rightarrow 0. $$
\end{lemma}

\begin{proof}
From corollary ~\ref{corollary8} we get $ L_s Z_{\mathfrak{\beta}}^{\delta} \epsilon_{\gamma}^{\beta} L_0 Z_{\mathfrak{\beta}}^{\gamma} K_{k-2} \cong $
\begin{eqnarray*} 
M^{\delta}(k-1) & \textrm{ if } & s=k-l-1 \\
0  & \textrm{ if } & s \neq k-l-1.
\end{eqnarray*}
Continuing in this way and using the corollary we get
$ L_s Z_{\mathfrak{\beta}}^{\delta} \epsilon_{\gamma}^{\beta} L_0 Z_{\mathfrak{\beta}}^{\gamma} K_{l} \cong $
\begin{eqnarray*} 
M^{\delta}(l+1)/\ldots /M^{\delta}(k-1) & \textrm{ if } & s=k-l-1 \\
0  & \textrm{ if } & s \neq k-l-1
\end{eqnarray*}
and then 
$ L_s Z_{\mathfrak{\beta}}^{\delta} \epsilon_{\gamma}^{\beta} L_1 Z_{\mathfrak{\beta}}^{\gamma} K_{l-1} \cong $
\begin{eqnarray*} 
M^{\delta}(l+1)/\ldots/ M^{\delta}(k-1) & \textrm{ if } & s=k-l-1 \\
0  & \textrm{ if } & s \neq k-l-1.
\end{eqnarray*}
Corollary ~\ref{corollary8} then gives an exact sequence
$$ 0 \rightarrow M^{\delta}(l-1) \rightarrow L_{k-l} Z_{\mathfrak{\beta}}^{\delta} \epsilon_{\gamma}^{\beta} L_1 Z_{\mathfrak{\beta}}^{\gamma} K_{l-2} \rightarrow
M^{\delta}(l+1)/\ldots/ M^{\delta}(k-1) \rightarrow 0. $$
Next the corollary gives the following diagram:

\begin{tiny}
\xymatrix
{
&	&	&0\ar[d]	&	&	&\\
&	&	&M^{\delta}(l-1)\ar[d]		&	&	&\\
&0\ar[r]	&L_{k-l+1}Z_{\mathfrak{\beta}}^{\delta} \epsilon_{\gamma}^{\beta} L_1Z_{\mathfrak{\beta}}^{\gamma} K_{l-3}\ar[r]	&L_{k-l}Z_{\mathfrak{\beta}}^{\delta} \epsilon_{\gamma}^{\beta} L_1Z_{\mathfrak{\beta}}^{\gamma}K_{l-2}\ar[r]\ar[d]	&M^{\delta}(l-2)\ar[r] &L_{k-l}Z_{\mathfrak{\beta}}^{\delta}\epsilon_{\gamma}^{\beta} L_1Z_{\mathfrak{\beta}}^{\gamma} K_{l-3}\ar[r] &0\\
&	&	&M^{\delta}(l+1)/ \ldots /M^{\delta}(k-1)\ar[d]	&	&	&\\
&	&	&0	&	&	&
}
\end{tiny}

Since that the map $ M^{\delta}(l-1) \rightarrow M^{\delta}(l-2) $ is standard, the above diagram can be extended to the following commutative diagram.

\begin{tiny}
\xymatrix
{
&	&0\ar[d]	&0\ar[d]	&0\ar[d]	&	&\\
&	&M^{\delta}(l)/\ldots/M^{\delta}(k-1)\ar[r]\ar[d]	&M^{\delta}(l-1)\ar[d]	\ar[r]	&M^{\delta}(l-2)\ar[d]	&	&\\
&0\ar[r]	&L_{k-l+1}Z_{\mathfrak{\beta}}^{\delta} \epsilon_{\gamma}^{\beta} L_1Z_{\mathfrak{\beta}}^{\gamma} K_{l-3}\ar[r]\ar[d]	&L_{k-l}Z_{\mathfrak{\beta}}^{\delta} \epsilon_{\gamma}^{\beta} L_1Z_{\mathfrak{\beta}}^{\gamma}K_{l-2}\ar[r]\ar[d]	&M^{\delta}(l-2)\ar[r]\ar[d] &L_{k-l}Z_{\mathfrak{\beta}}^{\delta}\epsilon_{\gamma}^{\beta} L_1Z_{\mathfrak{\beta}}^{\gamma} K_{l-3}\ar[r] &0\\
&0\ar[r]	&M^{\delta}(l+1)/ \ldots/ M^{\delta}(k-1)\ar[r]\ar[d]	&M^{\delta}(l+1)/ \ldots/ M^{\delta}(k-1)\ar[d]\ar[r]	&0	&	&\\
&	&0	&0	&	&	&
}
\end{tiny}

The surjection in the first column arises from the snake lemma and formulas for dimensions of spaces of homomorphisms between generalized Verma modules given in [Sh].
It also implies that $ L_{k-l} Z_{\mathfrak{\beta}}^{\delta} \epsilon_{\gamma}^{\beta} L_1 Z_{\mathfrak{\beta}}^{\gamma} K_{l-3} \cong M^{\delta}(l-2)/ \ldots/ M^{\delta}(k-1) $ and 
there is an exact sequence
$$ 0 \rightarrow M^{\delta}(l)/ \ldots /M^{\delta}(k-1) \rightarrow L_{k-l+1} Z_{\mathfrak{\beta}}^{\delta} \epsilon_{\gamma}^{\beta} L_1 Z_{\mathfrak{\beta}}^{\gamma} K_{l-3} \rightarrow 
M^{\delta}(l+1)/ \ldots /M^{\delta}(k-1) \rightarrow 0. $$

Continuing in this manner and using corollary ~\ref{corollary8} gives us the lemma.
\end{proof}

\begin{lemma}
\label{lemma19a}
Let $ m=l. $  Then
$ L_s Z_{\mathfrak{\beta}}^{\mathfrak{\alpha}} \epsilon_{\gamma}^{\beta} L_1 Z_{\beta}^{\gamma} \epsilon_{\delta}^{\beta} M^{\delta}(l) \cong $
\begin{eqnarray*} 
M^{\mathfrak{\alpha}}(k-1)/\ldots /M^{\alpha}(0) & \textrm{ if } & s=l+1 \\
L_{l+2}Z_{\mathfrak{\beta}}^{\mathfrak{\alpha}} \epsilon_{\gamma}^{\beta}L_1 Z_{\mathfrak{\beta}}^{\gamma} J_{k-l-4}  & \textrm{ if } & s=k-1
\end{eqnarray*}
and there is an exact sequence
$$ 0 \rightarrow M^{\mathfrak{\alpha}}(l)/\ldots/ M^{\mathfrak{\alpha}}(0) \rightarrow L_{l+2}Z_{\mathfrak{\beta}}^{\mathfrak{\alpha}} \epsilon_{\gamma}^{\beta}L_1 Z_{\mathfrak{\beta}}^{\gamma} J_{k-l-4}
\rightarrow M^{\mathfrak{\alpha}}(l-1)/\ldots /M^{\mathfrak{\alpha}}(0) \rightarrow 0. $$
\end{lemma}

\begin{proof}
We use the exact sequences of corollary ~\ref{corollary9}.

If $ s=l, $ $ L_s Z_{\mathfrak{\beta}}^{\mathfrak{\alpha}} \epsilon_{\gamma}^{\beta} L_0 Z_{\mathfrak{\beta}}^{\gamma} J_{k-2} \cong M^{\mathfrak{\alpha}}(0) $ and is zero otherwise.

If $ s=l, $ $ L_s Z_{\mathfrak{\beta}}^{\mathfrak{\alpha}} \epsilon_{\gamma}^{\beta} L_0 Z_{\mathfrak{\beta}}^{\gamma} J_{k-3} \cong M^{\mathfrak{\alpha}}(1)/M^{\mathfrak{\alpha}}(0) $ and is zero otherwise.

Continuing in this way,

If $ s=l, $ $ L_s Z_{\mathfrak{\beta}}^{\mathfrak{\alpha}} \epsilon_{\gamma}^{\beta} L_0 Z_{\mathfrak{\beta}}^{\gamma} J_{k-l-1} \cong M^{\mathfrak{\alpha}}(l-1)/ \ldots /M^{\mathfrak{\alpha}}(0) $ and is zero otherwise.

Next we get an exact sequence
$$ 0 \rightarrow M^{\mathfrak{\alpha}}(l+1) \rightarrow L_{l+1} Z_{\mathfrak{\beta}}^{\mathfrak{\alpha}} \epsilon_{\gamma}^{\beta} L_1 Z_{\mathfrak{\beta}}^{\gamma} J_{k-l-3} \rightarrow
M^{\mathfrak{\alpha}}(l-1)/ \ldots /M^{\mathfrak{\alpha}}(0) \rightarrow 0. $$

Corollary ~\ref{corollary9} then produces the following diagram:

\begin{tiny}
\xymatrix
{
&	&	&0\ar[d]	&	&	&\\
&	&	&M^{\mathfrak{\alpha}}(l+1)\ar[d]		&	&	&\\
&0\ar[r]	&L_{l+2}Z_{\mathfrak{\beta}}^{\mathfrak{\alpha}} \epsilon_{\gamma}^{\beta} L_1Z_{\mathfrak{\beta}}^{\gamma} J_{k-l-4}\ar[r]	&L_{l+1}Z_{\mathfrak{\beta}}^{\mathfrak{\alpha}} \epsilon_{\gamma}^{\beta} L_1Z_{\mathfrak{\beta}}^{\gamma}J_{k-l-3}\ar[r]\ar[d]	&M^{\mathfrak{\alpha}}(l+2)\ar[r] &L_{l+1}Z_{\mathfrak{\beta}}^{\mathfrak{\alpha}}\epsilon_{\gamma}^{\beta} L_1Z_{\mathfrak{\beta}}^{\gamma} J_{k-l-4}\ar[r] &0\\
&	&	&M^{\mathfrak{\alpha}}(l-1)/ \ldots/ M^{\mathfrak{\alpha}}(0)\ar[d]	&	&	&\\
&	&	&0	&	&	&
}
\end{tiny}

Thus $ L_{l+1} Z_{\mathfrak{\beta}}^{\mathfrak{\alpha}} \epsilon_{\gamma}^{\beta} L_1 Z_{\mathfrak{\beta}}^{\gamma} J_{k-l-4} \cong M^{\mathfrak{\alpha}}(l+2)/ \ldots /M^{\mathfrak{\alpha}}(0) $ and there is an exact sequence
$$ 0 \rightarrow M^{\mathfrak{\alpha}}(l)/\ldots/ M^{\mathfrak{\alpha}}(0) \rightarrow L_{l+2}Z_{\mathfrak{\beta}}^{\mathfrak{\alpha}} \epsilon_{\gamma}^{\beta}L_1 Z_{\mathfrak{\beta}}^{\gamma} J_{k-l-4}
\rightarrow M^{\mathfrak{\alpha}}(l-1)/\ldots /M^{\mathfrak{\alpha}}(0) \rightarrow 0. $$

Using the rest of the exact sequences of corollary ~\ref{corollary9}, we easily obtain this lemma. (The details are the same as lemma ~\ref{lemma19}.)
\end{proof}

For the rest of this subsection, let $ F = LZ_{\mathfrak{\beta}}^{\mathfrak{\alpha}} \epsilon_{\gamma}^{\beta} L_1 Z_{\mathfrak{\beta}}^{\gamma} \epsilon_{\delta}^{\beta}. $
Let $ G = LZ_{\mathfrak{\beta}}^{\mathfrak{\delta}} \epsilon_{\gamma}^{\beta} L_1 Z_{\mathfrak{\beta}}^{\gamma} \epsilon_{\alpha}^{\beta}. $

\begin{lemma}
\label{lemma20}
Let $ l<k-1. $
\begin{enumerate}
\item $ H^{-s} F M^{\delta}(l)/ \ldots /M^{\delta}(k-1) \cong $
\begin{eqnarray*} 
M^{\mathfrak{\alpha}}(k-1)/\ldots /M^{\alpha}(0) & \textrm{ if } & s=l+1, l+3, \ldots, 2(k-1)-l-1 \\
M^{\mathfrak{\alpha}}(l-1)/ \ldots /M^{\mathfrak{\alpha}}(0)  & \textrm{ if } & s=k-1.
\end{eqnarray*}
\item $ H^{-s} G M^{\alpha}(l)/ \ldots/ M^{\alpha}(0) \cong $
\begin{eqnarray*} 
M^{\mathfrak{\delta}}(0)/\ldots /M^{\delta}(k-1) & \textrm{ if } & s=k-l, k-l+2, \ldots, k+l-2 \\
M^{\mathfrak{\delta}}(l+1)/ \ldots /M^{\mathfrak{\delta}}(k-1)  & \textrm{ if } & s=k-1.
\end{eqnarray*}
\end{enumerate}

\end{lemma}

\begin{proof}
\begin{enumerate}
\item Consider the exact sequence
$$ 0 \rightarrow M^{\delta}(k-1) \rightarrow \cdots \rightarrow M^{\delta}(l+1) \rightarrow
M^{\delta}(l) \rightarrow M^{\delta}(l)/ \ldots /M^{\delta}(k-1) \rightarrow 0. $$
By lemma ~\ref{lemma19a}, applying $ F $ gives the following complex with vertical short exact sequences:

\begin{tiny}
\xymatrix
{
&0\ar[d]	
&	
&0\ar[d]		
&0\ar[d]		
&\\
&M^{\mathfrak{\alpha}}(k-1)/\ldots /M^{\mathfrak{\alpha}}(0)\ar[d]
&
&M^{\mathfrak{\alpha}}(l+1)/\ldots /M^{\mathfrak{\alpha}}(0)\ar[d]
&M^{\mathfrak{\alpha}}(l)/\ldots /M^{\mathfrak{\alpha}}(0)\ar[d]
&\\
&L_{k-1} FM^{\mathfrak{\delta}}(k-1) \ar@{^{(}->}[r]\ar[d] 
&\cdots\ar[r] 
&L_{k-1} FM^{\mathfrak{\delta}}(l+1) \ar[r]\ar[d] 
&L_{k-1} FM^{\mathfrak{\delta}}(l) \ar@{->>}[r]\ar[d] 
&L_{k-1} FM^{\mathfrak{\delta}}(l)/\ldots/M^{\mathfrak{\delta}}(k-1) \\
&M^{\mathfrak{\alpha}}(k-2)/\ldots /M^{\mathfrak{\alpha}}(0)\ar[d]
&
&M^{\mathfrak{\alpha}}(l)/\ldots /M^{\mathfrak{\alpha}}(0)\ar[d]
&M^{\mathfrak{\alpha}}(l-1)/\ldots /M^{\mathfrak{\alpha}}(0)\ar[d]
&\\
&0
&
&0
&0
&
}
\end{tiny}

It suffices to show that the map $ L_{k-1} F M^{\mathfrak{\delta}}(l+1) \rightarrow 
L_{k-1} F M^{\mathfrak{\delta}}(l) $ factors through 
$ M^{\mathfrak{\alpha}}(l)/ \ldots/M^{\mathfrak{\alpha}}(0). $
Suppose that the generalized BGG resolutions for 
$ M^{\mathfrak{\delta}}(l+1) $ and $ M^{\mathfrak{\delta}}(l) $ give rise to short exact sequences involving modules $ B_i^{'} $ and $ B_i $ respectively.
Let $ M^{\mathfrak{\beta}}(r,s) = 
M^{\mathfrak{\beta}}(\underbrace{a_i, \ldots, a_{i+k-2},} r, s, \underbrace{k-1, \ldots \hat{s}, \ldots, 0}). $
There is a commutative diagram:

\xymatrix
{
&0\ar[r] 
&B_{k-l-3}^{'}\ar[r]\ar[d] 
&M^{\mathfrak{\beta}}(l+1, l+2)\ar[r]\ar[d]
&B_{k-l-4}^{'}\ar[r]\ar[d]
&0\\
&0\ar[r] 
&B_{k-l-3}\ar[r]
&M^{\mathfrak{\beta}}(l, l+2)\ar[r]
&B_{k-l-4}\ar[r]
&0\\
}

Now apply $ F $ to get the following commutative diagram:

\xymatrix
{
&L_{l+2}F B_{k-l-4}^{'}\ar[r]\ar[d]^{f} 
&L_{l+1}F B_{k-l-3}^{'} \cong M^{\mathfrak{\alpha}}(l)/ \ldots /M^{\mathfrak{\alpha}}(0) \ar[r]\ar[d]^{i}
&0\ar[d]\\
&L_{l+2}F B_{k-l-4}\ar[r]^{h}
&L_{l+1}F B_{k-l-3}\ar[r]^{j}
&M^{\mathfrak{\alpha}}(l+2)\\
}

Thus $ \text{im}(i) \subset \text{ker}(j) = \text{im}(h) = L_{l+2}F B_{k-l-4}. $
Thus the map 
$ L_{k-1} F M^{\mathfrak{\delta}}(l+1) \rightarrow 
L_{k-1} F M^{\mathfrak{\delta}}(l) $ factors through 
$ M^{\mathfrak{\alpha}}(l)/ \cdots/M^{\mathfrak{\alpha}}(0) $
as desired.

Now the long exact sequence for $ F $ gives
$ H^{-s} F M^{\mathfrak{\delta}}(l)/ \cdots/ M^{\mathfrak{\delta}}(k-1) \cong $

\begin{eqnarray*} 
M^{\mathfrak{\alpha}}(k-1)/\ldots /M^{\delta}(0) & \textrm{ if } & s=l+1, l+3, \ldots, 2(k-1)-l-1 \\
M^{\mathfrak{\alpha}}(l-1)/ \ldots M^{\mathfrak{\alpha}}(0)  & \textrm{ if } & s=k-1.
\end{eqnarray*}
\item This follows exactly as the first part.
\end{enumerate}
\end{proof}

\begin{lemma}
$ H^{-s} F LZ_{\mathfrak{\beta}}^{\delta} \epsilon_{\gamma}^{\beta} L_1 Z_{\mathfrak{\beta}}^{\gamma} \epsilon_{\alpha}^{\beta} M^{\mathfrak{\alpha}}(l) \cong $
\begin{eqnarray*} 
M^{\mathfrak{\alpha}}(k-1)/\ldots /M^{\delta}(0) & \textrm{ if } & s=k-l+1, k-l+3, \ldots k-l+2(k-1)-3 \\
M^{\mathfrak{\alpha}}(l)  & \textrm{ if } & s=2(k-1).
\end{eqnarray*}
\end{lemma}

\begin{proof}
By lemma ~\ref{lemma19}, there is a distinguished triangle
$$ L_{k-l+1} Z_{\mathfrak{\beta}}^{\delta} \epsilon_{\gamma}^{\beta} L_1 Z_{\mathfrak{\beta}}^{\gamma} K_{l-3}[k-1] \rightarrow 
LZ_{\mathfrak{\beta}}^{\delta} \epsilon_{\gamma}^{\beta} L_1 Z_{\mathfrak{\beta}}^{\gamma} \epsilon_{\alpha}^{\beta} M^{\mathfrak{\alpha}}(l) \rightarrow
M^{\delta}(0)/ \ldots /M^{\delta}(k-1)[k-l]. $$
By proposition ~\ref{prop10},
$ H^{-s}F M^{\delta}(0)/ \ldots / M^{\delta}(k-1)[k-l] \cong M^{\alpha}(k-1)/ \ldots /M^{\alpha}(0) $ for 
$ s=k-l+1, k-l+3, \ldots, k-l+2(k-1)-1. $

From the commutative diagram of lemma ~\ref{lemma19}, we have the short exact sequences:
\begin{align*}
&0 \rightarrow L_{k-l+1} Z_{\beta}^{\delta} \epsilon_{\gamma}^{\beta} L_1 Z_{\beta}^{\gamma} K_{l-3} \rightarrow 
L_{k-l} Z_{\beta}^{\delta} \epsilon_{\gamma}^{\beta} L_1 Z_{\beta}^{\gamma} K_{l-2} \rightarrow M^{\delta}(l-1)/ \ldots / M^{\delta}(k-1) \rightarrow 0\\
&0 \rightarrow M^{\delta}(l-1)/ \ldots / M^{\delta}(k-1) \rightarrow M^{\delta}(l-2) \rightarrow M^{\delta}(l-2)/ \ldots / M^{\delta}(k-1) \rightarrow 0. 
\end{align*} 
Therefore we get the following morphisms:
$$ M^{\delta}(0)/ \ldots / M^{\delta}(k-1) \rightarrow M^{\delta}(l-1)/ \ldots / M^{\delta}(k-1)[l-1] \rightarrow L_{k-l+1} Z_{\beta}^{\gamma} \epsilon_{\gamma}^{\beta} L_1 Z_{\beta}^{\gamma} K_{l-3}[l]. $$
This is the morphism $ M^{\delta}(0)/ \ldots /M^{\delta}(k-1)[k-l] \rightarrow L_{k-l+1} Z_{\mathfrak{\beta}}^{\delta} \epsilon_{\gamma}^{\beta} L_1 Z_{\mathfrak{\beta}}^{\gamma} K_{l-3}[k-1] $ in the distinguished triangle.
Assume $ s \leq k-l+2(k-1)-3.  $ Suppose 
$$ H^{-s}F(M^{\delta}(l-1)/ \ldots /M^{\delta}(k-1)[k-1]) \rightarrow H^{-s+1}F(L_{k-l+1} Z_{\beta}^{\delta} \epsilon_{\gamma}^{\beta} L_1 Z_{\beta}^{\gamma} K_{l-3}[k-1]) $$
is an isomorphism.  Then the long exact sequence for the triangle above and $ F $ gives a map:
$$ H^{-s}F M^{\delta}(0)/ \ldots / M^{\delta}(k-1)[k-l]  \rightarrow H^{-s+1}F(L_{k-l+1} Z_{\beta}^{\delta} \epsilon_{\gamma}^{\beta} L_1 Z_{\beta}^{\gamma} K_{l-3}[k-1]) $$
which is an isomorphism.
Then we also have 
$$ H^{-s}F(L_{k-l+1} Z_{\beta}^{\delta} \epsilon_{\gamma}^{\beta} L_1 Z_{\beta}^{\gamma} K_{l-3}[k-1]) \cong M^{\alpha}(k-1)/ \ldots / M^{\alpha}(0). $$ 
If the map is zero, then
$$ H^{-s}F(L_{k-l+1} Z_{\beta}^{\delta} \epsilon_{\gamma}^{\beta} L_1 Z_{\beta}^{\gamma} K_{l-3}[k-1]) \cong
H^{-s+1}F(L_{k-l+1} Z_{\beta}^{\delta} \epsilon_{\gamma}^{\beta} L_1 Z_{\beta}^{\gamma} K_{l-3}[k-1]) = 0. $$
Therefore in either case,
$$ H^{-s} F LZ_{\mathfrak{\beta}}^{\delta} \epsilon_{\gamma}^{\beta} L_1 Z_{\mathfrak{\beta}}^{\gamma} \epsilon_{\alpha}^{\beta} M^{\mathfrak{\alpha}}(l) \cong 
M^{\alpha}(k-1)/\ldots /M^{\delta}(0) $$ 
for $ s=k-l+1, k-l+3, \ldots, k-l+2(k-1)-3. $

For $ s = 2(k-1)-l-1, $ consider again the short exact sequence
$$ 0 \rightarrow L_{k-l+1} Z_{\beta}^{\delta} \epsilon_{\gamma}^{\beta} L_1 Z_{\beta}^{\gamma} K_{l-3} \rightarrow 
L_{k-l} Z_{\beta}^{\delta} \epsilon_{\gamma}^{\beta} L_1 Z_{\beta}^{\gamma} K_{l-2} \rightarrow M^{\delta}(l-1)/ \ldots / M^{\delta}(k-1) \rightarrow 0. $$
It is easy to calculate 
$$ H^{-(2(k-1)-l-1)} F L_{k-l+1} Z_{\beta}^{\delta} \epsilon_{\gamma}^{\beta} L_1 Z_{\beta}^{\gamma} K_{l-3} \cong M^{\alpha}(k-1)/ \ldots /M^{\alpha}(0) $$
and
$$ H^{-(2(k-1)-l-1)} F L_{k-l+1} Z_{\beta}^{\delta} \epsilon_{\gamma}^{\beta} L_1 Z_{\beta}^{\gamma} K_{l-2} \cong
H^{-(2(k-1)-l)} F L_{k-l+1} Z_{\beta}^{\delta} \epsilon_{\gamma}^{\beta} L_1 Z_{\beta}^{\gamma} K_{l-2} = 0. $$
Thus
$$ H^{-(2(k-1)-l)} F M^{\delta}(l-1)/ \ldots /M^{\delta}(k-1) \rightarrow H^{-(2(k-1)-l-1)} F L_{k-l+1} Z_{\beta}^{\delta} \epsilon_{\gamma}^{\beta} L_1 Z_{\beta}^{\gamma} K_{l-3} $$
is an isomorphism.
Therefore 
$$ H^{-(k-l+2(k-1))} F M^{\delta}(0)/ \ldots /M^{\delta}(k-1)[k-l] \rightarrow 
H^{-(k-l+2(k-1)-1)} F L_{k-l+1} Z_{\beta}^{\delta} \epsilon_{\gamma}^{\beta} L_1 Z_{\beta}^{\gamma} K_{l-3}[k-1] $$
is an isomorphism so
$$ H^{-(k-l+2(k-1)-1)} F LZ_{\beta}^{\delta} \epsilon_{\gamma}^{\beta} L_1 Z_{\mathfrak{\beta}}^{\gamma} \epsilon_{\alpha}^{\beta} M^{\mathfrak{\alpha}}(l) = 0. $$

In order to compute
$ H^{-s}F(L_{k-l+1} Z_{\beta}^{\delta} \epsilon_{\gamma}^{\beta} L_1 Z_{\beta}^{\gamma} K_{l-3}[k-1]) $
for $ s=2(k-1), $ we use the short exact sequence from lemma ~\ref{lemma19}:
$$ 0 \rightarrow M^{\delta}(l)/ \ldots / M^{\delta}(k-1) \rightarrow
L_{k-l+1} Z_{\beta}^{\delta} \epsilon_{\gamma}^{\beta} L_1 Z_{\beta}^{\gamma} K_{l-3} \rightarrow
M^{\delta}(l+1)/ \ldots / M^{\delta}(k-1) \rightarrow 0.  $$

The first goal is to prove that this sequence is non-split.
Let $ \bar{F} = LZ_{\beta}^{\alpha} \epsilon_{\gamma}^{\beta}[-1] LZ_{\beta}^{\gamma} \epsilon_{\delta}^{\beta}[-(k-1)]. $
Let the adjoint of $ \bar{F} $ be 
$ \bar{G} = LZ_{\beta}^{\delta} \epsilon_{\gamma}^{\beta}[-1] LZ_{\beta}^{\gamma} \epsilon_{\alpha}^{\beta}[-(k-1)]. $
Then the exact sequence becomes
$$ 0 \rightarrow M^{\delta}(l)/ \ldots /M^{\delta}(k-1) \rightarrow H^0 \bar{G} M^{\alpha}(l) \rightarrow M^{\delta}(l+1)/ \ldots / M^{\delta}(k-1) \rightarrow 0. $$
In order to show that it is non-split, it suffices to show that $ \Hom(M^{\delta}(l+1)/ \ldots / M^{\delta}(k-1), H^0 \bar{G} M^{\alpha}(l)) = 0. $
Consider the distinguished triangle
$$ H^0 \bar{G} M^{\alpha}(l) \rightarrow \bar{G} M^{\alpha}(l) \rightarrow M^{\delta}(0)/ \ldots / M^{\delta}(k-1) [-l+1]. $$
It then suffices to prove
$$ \Hom(M^{\delta}(l+1)/ \ldots / M^{\delta}(k-1), \bar{G} M^{\alpha}(l)) \cong \Hom(\bar{F} M^{\delta}(l+1)/ \ldots / M^{\delta}(k-1), M^{\alpha}(l)) = 0. $$
For this consider the distinguished triangle
$$ \tau^{\leq 0} \bar{F} M^{\delta}(l+1)/ \ldots / M^{\delta}(k-1) \rightarrow \bar{F} M^{\delta}(l+1)/ \ldots / M^{\delta}(k-1) \rightarrow \tau^{\geq 1} \bar{F} M^{\delta}(l+1)/ \ldots / M^{\delta}(k-1). $$
We will show that this triangle splits.
Note that $ M^{\alpha}(k-1)/ \ldots / M^{\alpha}(0) \cong \epsilon_{\alpha+\delta}^{\alpha} M^{\alpha+\delta}(\underbrace{k-1, \ldots, 0,} \underbrace{k-1,\ldots, 0}). $ 
Denote it simply by $ M^{\alpha+\delta}. $  Assume $ l-k $ is even.  The other case follows in the same way as below.
Clearly,
$$ \tau^{\geq 1} \bar{F} M^{\delta}(l+1)/ \ldots / M^{\delta}(k-1) \cong M^{\alpha+\delta}[-k+l+3] \oplus
M^{\alpha+\delta}[-k+l+5] \oplus \cdots \oplus M^{\alpha+\delta}[-1]. $$
Next consider $ \tau^{\leq 0} \bar{F} M^{\delta}(l+1)/ \ldots /M^{\delta}(k-1). $  There is a distinguished triangle:
\begin{align*}
&\tau^{\leq -1} \bar{F} M^{\delta}(l+1)/ \ldots /M^{\delta}(k-1) \rightarrow \tau^{\leq 0} \bar{F} M^{\delta}(l+1)/ \ldots /M^{\delta}(k-1) \rightarrow\\
&H^0 \bar{F} M^{\delta}(l+1)/ \ldots /M^{\delta}(k-1) \cong M^{\alpha}(l)/\ldots / M^{\alpha}(0) . 
\end{align*} 
By lemma ~\ref{lemma20},
$$ \tau^{\leq -1} \bar{F} M^{\delta}(l+1)/ \ldots /M^{\delta}(k-1) \cong
M^{\alpha+\delta}[1] \oplus M^{\alpha+\delta}[3] \oplus \cdots \oplus M^{\alpha+\delta}[k-l-3]. $$
We easily compute
$$ \Hom(M^{\alpha}(l)/\ldots / M^{\alpha}(0), \tau^{\leq -1} \bar{F} M^{\delta}(l+1)/ \ldots /M^{\delta}(k-1) [1]) = 0. $$
Thus
$$ \tau^{\leq 0} \bar{F} M^{\delta}(l+1)/ \ldots /M^{\delta}(k-1) \cong M^{\alpha+\delta}[1] \oplus \cdots \oplus M^{\alpha+\delta}[k-l-3] \oplus
M^{\alpha}(l)/\ldots / M^{\alpha}(0). $$
Next we check that
$$ \Hom(\tau^{\geq 1} \bar{F} M^{\delta}(l+1)/ \ldots / M^{\delta}(k-1),  \tau^{\leq 0} \bar{F} M^{\delta}(l+1)/ \ldots / M^{\delta}(k-1)[1]) =0. $$
This reduces to showing
$$ \Hom(M^{\alpha+\delta}[-k+l+3] \oplus M^{\alpha+\delta}[-k+l+5] \oplus \cdots \oplus M^{\alpha+\delta}[-1],
M^{\alpha}(l)/ \ldots / M^{\alpha}(0)[1]) = 0. $$

Then the above space of homomorphisms becomes by adjointness:
$$ \Hom(M^{\alpha+\delta}[-k+l+3] \oplus \cdots \oplus M^{\alpha+\delta}[-1], LZ_{\alpha}^{\alpha+\delta} M^{\alpha}(l)/\cdots/M^{\alpha}(0)[-2k+3]). $$
It is easy to calculate
$$ LZ_{\alpha}^{\alpha+\delta}M^{\alpha}(l)/\ldots/M^{\alpha}(0) \cong M^{\alpha+\delta}[k-l-1] \oplus \cdots \oplus M^{\alpha+\delta}[3k-l-3]. $$
Therefore the space of homomorphisms is zero so
$$ \bar{F} M^{\delta}(l+1)/ \ldots / M^{\delta}(k-1) \cong M^{\alpha}(l)/ \ldots /M^{\alpha}(0) \oplus M^{\alpha+\delta}[-k+l+3] \oplus \cdots \oplus M^{\alpha+\delta}[k-l-3]. $$
It is now easy to see that
$$ \Hom(\bar{F} M^{\delta}(l+1)/ \ldots / M^{\delta}(k-1), M^{\alpha}(l)) = 0. $$
Therefore 
$$ 0 \rightarrow M^{\delta}(l)/ \ldots / M^{\delta}(k-1) \rightarrow
L_{k-l+1} Z_{\beta}^{\delta} \epsilon_{\gamma}^{\beta} L_1 Z_{\beta}^{\gamma} K_{l-3} \rightarrow
M^{\delta}(l+1)/ \ldots / M^{\delta}(k-1) \rightarrow 0 $$
is non-split.

We want to show that when $ \bar{F} $ is applied to the above sequence, it remains non-split.
Since $ \bar{F} $ and $ \bar{G} $ are biadjoint, the composition $ \bar{G} \rightarrow \bar{G} \bar{F} \bar{G} \rightarrow \bar{G} $ is the identity.
This gives that the adjunction morphism induces an isomorphism on cohomology:
$$ M^{\delta}(l+1)/ \ldots / M^{\delta}(k-1) \rightarrow H^0 \bar{G} \bar{F} M^{\delta}(l+1)/ \ldots / M^{\delta}(k-1). $$
Therefore, by the five lemma, the adjunction morphism induces an isomorphism:
$$ L_{k-l+1} Z_{\beta}^{\delta} \epsilon_{\gamma}^{\beta} L_1 Z_{\beta}^{\gamma} K_{l-3} \rightarrow H^0 \bar{G} \bar{F} L_{k-l+1} Z_{\beta}^{\delta} \epsilon_{\gamma}^{\beta} L_1 Z_{\beta}^{\gamma} K_{l-3}. $$
This implies that
$$ 0 \rightarrow H^0 \bar{F} M^{\delta}(l)/ \ldots / M^{\delta}(k-1) \rightarrow
H^0 \bar{F} L_{k-l+1} Z_{\beta}^{\delta} \epsilon_{\gamma}^{\beta} L_1 Z_{\beta}^{\gamma} K_{l-3} \rightarrow
H^0 \bar{F} M^{\delta}(l+1)/ \ldots / M^{\delta}(k-1) \rightarrow 0 $$
is non-split.
This is the exact sequence 
$$ 0 \rightarrow M^{\alpha}(l-1)/ \ldots / M^{\alpha}(0) \rightarrow H^0 \bar{F} L_{k-l+1} Z_{\beta}^{\delta} \epsilon_{\gamma}^{\beta} L_1 Z_{\beta}^{\gamma} K_{l-3} \rightarrow 
M^{\alpha}(l)/ \ldots / M^{\alpha}(0) \rightarrow 0. $$
If the extension is unique up to isomorphism, the middle term is isomorphic to $ M^{\alpha}(l). $  This would complete the lemma.

Consider the truncated full subcategory where $ M^{\alpha}(l) $ is a projective object.
Apply the functor $ \Hom(M^{\alpha}(l), \bullet) $ to the short exact sequence
$$ 0 \rightarrow M^{\alpha}(l-1)/ \ldots /M^{\alpha}(0) \rightarrow M^{\alpha}(l) \rightarrow M^{\alpha}(l)/ \ldots / M^{\alpha}(0) \rightarrow 0 $$ 
in the truncated subcategory to get the exact sequence
$$ 0 \rightarrow \Hom(M^{\alpha}(l), M^{\alpha}(l-1)/ \ldots /M^{\alpha}(0)) \rightarrow \Hom(M^{\alpha}(l), M^{\alpha}(l)) \rightarrow \Hom(M^{\alpha}(l), M^{\alpha}(l)/ \ldots / M^{\alpha}(0)) \rightarrow 0. $$ 
Thus $ \Hom(M^{\alpha}(l), M^{\alpha}(l)/ \ldots / M^{\alpha}(0)) \cong \mathbb{C} $ and
$ \Ext^1(M^{\alpha}(l), M^{\alpha}(l-1)/ \ldots / M^{\alpha}(0)) = 0. $
If we apply the functor $ \Hom(\bullet, M^{\alpha}(l)/ \ldots / M^{\alpha}(0)) $ to the short exact sequence
$$ 0 \rightarrow M^{\alpha}(l-1)/ \ldots /M^{\alpha}(0) \rightarrow M^{\alpha}(l) \rightarrow M^{\alpha}(l)/ \ldots / M^{\alpha}(0) \rightarrow 0, $$ 
we get the exact sequence
$$ 0 \rightarrow \End(M^{\alpha}(l)/ \ldots / M^{\alpha}(0)) \rightarrow \Hom(M^{\alpha}(l),  M^{\alpha}(l)/ \ldots / M^{\alpha}(0)) \rightarrow \cdots. $$
This implies
$$ \End(M^{\alpha}(l)/ \ldots / M^{\alpha}(0)) \cong \mathbb{C}. $$
Now apply the functor $ \Hom(\bullet, M^{\alpha}(l-1)/ \ldots / M^{\alpha}(0)) $ to the short exact sequence
$$ 0 \rightarrow M^{\alpha}(l-1)/ \ldots /M^{\alpha}(0) \rightarrow M^{\alpha}(l) \rightarrow M^{\alpha}(l)/ \ldots / M^{\alpha}(0) \rightarrow 0. $$
This gives a long exact sequence
$$ 0 \rightarrow \mathbb{C} \rightarrow \Ext^1(M^{\alpha}(l)/ \ldots / M^{\alpha}(0), M^{\alpha}(l-1)/ \ldots / M^{\alpha}(0)) \rightarrow
\Ext^1(M^{\alpha}(l), M^{\alpha}(l-1)/ \ldots / M^{\alpha}(0)) = 0. $$
Therefore $ \Ext^1(M^{\alpha}(l)/ \ldots / M^{\alpha}(0), M^{\alpha}(l-1)/ \ldots / M^{\alpha}(0)) \cong \mathbb{C} $ and we have shown that the extension is unique.

\end{proof}

For the rest of this subsection we make the following definition.

\begin{define}
Let $ S = LZ_{\mathfrak{\beta}}^{\mathfrak{\alpha}} \epsilon_{\gamma}^{\beta}[-1] LZ_{\mathfrak{\beta}}^{\gamma} \epsilon_{\delta}^{\beta}[-(k-1)]LZ_{\mathfrak{\beta}}^{\delta} \epsilon_{\gamma}^{\beta}[-1] LZ_{\mathfrak{\beta}}^{\gamma} \epsilon_{\alpha}^{\beta}[-(k-1)]. $
\end{define}

The identity morphism in $ \End(LZ_{\mathfrak{\beta}}^{\delta} \epsilon_{\gamma}^{\beta}[-1] LZ_{\mathfrak{\beta}}^{\gamma} \epsilon_{\alpha}^{\beta}[-(k-1)]) $ gives rise to a morphism
$ S \rightarrow \Id. $
There is an adjunction morphism $ \epsilon_{\gamma}^{\beta}[-1] LZ_{\mathfrak{\beta}}^{\gamma} \rightarrow \Id[1] $ that gives rise to a morphism
$$ \theta \colon S \rightarrow LZ_{\mathfrak{\beta}}^{\mathfrak{\alpha}}  \epsilon_{\delta}^{\beta}[-(k-1)]LZ_{\mathfrak{\beta}}^{\delta} \epsilon_{\gamma}^{\beta}[-1] LZ_{\mathfrak{\beta}}^{\gamma} \epsilon_{\alpha}^{\beta}[-(k-1)][1]. $$
By proposition ~\ref{prop10}, this is a map
$$ \theta \colon S \rightarrow \oplus_{r=1}^{k-1} \epsilon_{\alpha+\delta}^{\alpha}[-(k-1)] LZ_{\mathfrak{\alpha}}^{\alpha+\delta}[1][2r-k]. $$
Composing this map with projection gives a map:
$$ \mu_2 \colon S \rightarrow \oplus_{r=1}^{k-2} \epsilon_{\alpha+\delta}^{\alpha}[-(k-1)] LZ_{\alpha}^{\alpha+\delta}[2r-(k-1)]. $$

\begin{prop}
\label{propdiagram4}
The map $ \mu_1 \oplus \mu_2 \colon S \rightarrow \Id \oplus \oplus_{r=1}^{k-2} \epsilon_{\alpha+\delta}^{\alpha+\beta}[-(k-1)] LZ_{\alpha+\mathfrak{\beta}}^{\alpha+\delta}[2r-(k-1)] $ is an isomorphism.
\end{prop}

\begin{proof}
This morphism is an isomorphism for generalized Verma modules as in the proof of proposition ~\ref{prop10}.  By considering a generalized Verma flag for a projective object, we see that $ \mu_1 \oplus \mu_2 $ is an isomorphism for projectives as well.  Since any bounded complex is quasi-isomorphic to a complex of projective objects, $ \mu_1 \oplus \mu_2 $ must be an isomorphism.
\end{proof}

\subsection{Diagram 5}
Recall the definitions of the subalgebras $ \mathfrak{s}_i $ and $ \mathfrak{t}_i. $
We would like to understand the functor $ LZ^{\mathfrak{s}_i} \epsilon_{\mathfrak{s}_{i+1}}[-1] LZ^{\mathfrak{s}_{i+1}} \epsilon_{\mathfrak{s}_i}[-1]. $  In [BFK] this functor turned out to be the identity.  When working on maximally singular categories, $ LZ^{\mathfrak{s}_{i+1}} \epsilon_{\mathfrak{s}_i}[-1] $ was an equivalence of categories.  On these less singular blocks this will turn out not to be the case.  It will be the identity plus some correction 
term which will depend on the two minimal parabolics involved.
The plan is to compute this functor on generalized Verma modules and then to construct some natural transformations which are isomorphisms on the generalized Verma modules.
These objects are parametrized by highest weights which satisfy the conditions of lying in the block and locally finite with respect to the appropriate subalgebra.
Since in these computations only $ \mathfrak{s}_i, \mathfrak{s}_{i+1}, 
$ and $ \mathfrak{t}_i $ are involved, for convenience of notation we will omit all the coordinates of the highest weight except the components in positions $ i, i+1, $ and $ i+2. $

\begin{lemma}
Let $ a_i > a_{i+1} = a_{i+2}. $  Then
$$ LZ^{\mathfrak{s}_i} \epsilon_{\mathfrak{s}_{i+1}}[-1] LZ^{\mathfrak{s}_{i+1}} \epsilon_{\mathfrak{s}_i}[-1] M^{\mathfrak{s}_i}(\underbrace{a_i, a_{i+1}}, a_{i+2}) \cong
M^{\mathfrak{s}_i}(\underbrace{a_i, a_{i+1}}, a_{i+2}). $$
\end{lemma}

\begin{proof}
There is a short exact sequence
$$ 0 \rightarrow M(a_{i+1}, a_i, a_{i+2}) \rightarrow M(a_i, a_{i+1}, a_{i+2}) \rightarrow
M^{\mathfrak{s}_i}(\underbrace{a_i, a_{i+1}}, a_{i+2}) \rightarrow 0. $$
Thus $ L_s Z^{\mathfrak{s}_{i+1}} \epsilon_{\mathfrak{s}_i} M^{\mathfrak{s}_i}(\underbrace{a_i, a_{i+1}}, a_{i+2}) \cong 0 $ if $ s = 0, 2 $ and is isomorphic to
$ M^{\mathfrak{s}_{i+1}}(a_{i+1}, \underbrace{a_i, a_{i+2}}) $ if $ s = 1. $

Now consider the short exact sequence
$$ 0 \rightarrow M(a_{i+1}, a_{i+2}, a_i) \rightarrow M(a_{i+1}, a_i, a_{i+2}) \rightarrow
M^{\mathfrak{s}_{i+1}}(a_{i+1}, \underbrace{a_i, a_{i+2}}) \rightarrow 0. $$ 
This gives
$ L_s Z^{\mathfrak{s}_i} \epsilon_{\mathfrak{s}_{i+1}} M^{\mathfrak{s}_{i+1}}(a_{i+1}, \underbrace{a_i, a_{i+2}}) \cong 0 $ if $ s = 0, 2 $ and is isomorphic to
$ M^{\mathfrak{s}_i}(\underbrace{a_i, a_{i+1}}, a_{i+2}) $ if $ s=1. $  The lemma now follows.
\end{proof}

\begin{lemma}
Let $ a_{i+2} = a_i > a_{i+1}. $  Then
$$ LZ^{\mathfrak{s}_i} \epsilon_{\mathfrak{s}_{i+1}}[-1] LZ^{\mathfrak{s}_{i+1}} \epsilon_{\mathfrak{s}_i}[-1] M^{\mathfrak{s}_i}(\underbrace{a_i, a_{i+1}}, a_{i+2}) \cong
M^{\mathfrak{s}_i}(\underbrace{a_i, a_{i+1}}, a_{i+2}). $$
\end{lemma}

\begin{proof}
Consider the short exact sequence
$$ 0 \rightarrow M(a_{i+1}, a_i, a_{i+2}) \rightarrow M(a_i, a_{i+1}, a_{i+2}) \rightarrow
M^{\mathfrak{s}_i}(\underbrace{a_i, a_{i+1}}, a_{i+2}) \rightarrow 0. $$
Thus $ L_s Z^{\mathfrak{s}_{i+1}} \epsilon_{\mathfrak{s}_i} M^{\mathfrak{s}_i}(\underbrace{a_i, a_{i+1}}, a_{i+2}) \cong 0 $ if $ s = 0, 2 $ and is isomorphic to
$ M^{s_{i+1}}(a_{i}, \underbrace{a_{i+2}, a_{i+1}}) $ if $ s = 1. $

Now consider the short exact sequence
$$ 0 \rightarrow M(a_{i}, a_{i+1}, a_{i+2}) \rightarrow M(a_{i}, a_{i+2}, a_{i+1}) \rightarrow
M^{\mathfrak{s}_{i+1}}(a_{i}, \underbrace{a_{i+2}, a_{i+1}}) \rightarrow 0. $$ 
This gives
$ L_s Z^{\mathfrak{s}_i} \epsilon_{\mathfrak{s}_{i+1}} M^{\mathfrak{s}_{i+1}}(a_{i}, \underbrace{a_{i+2}, a_{i+1}}) \cong 0 $ if $ s = 0, 2 $ and is isomorphic to
$ M^{\mathfrak{s}_i}(\underbrace{a_i, a_{i+1}}, a_{i+2}) $ if $ s=1. $  The lemma now follows.
\end{proof}

These were the two easiest cases and the computations were made in [BFK].  Now we come to the three new cases.

\begin{lemma}
Suppose $ a_{i+2} > a_i > a_{i+1}. $  Then
$$ LZ^{\mathfrak{s}_i} \epsilon_{\mathfrak{s}_{i+1}}[-1] LZ^{\mathfrak{s}_{i+1}} \epsilon_{\mathfrak{s}_i}[-1] M^{\mathfrak{s}_i}(\underbrace{a_i, a_{i+1}}, a_{i+2}) \cong
M^{\mathfrak{s}_i}(\underbrace{a_i, a_{i+1}}, a_{i+2}) \oplus
\epsilon_{\mathfrak{t}_i}^{\mathfrak{s}_i}M^{\mathfrak{t}_i}(\underbrace{a_{i+2}, a_{i}, a_{i+1}}). $$
\end{lemma}

\begin{proof}
Consider the short exact sequence
$$ 0 \rightarrow M(a_{i+1}, a_i, a_{i+2}) \rightarrow M(a_{i}, a_{i+1}, a_{i+2}) \rightarrow M^{\mathfrak{s}_i}(\underbrace{a_{i}, a_{i+1}}, a_{i+2}). $$ 
The long exact sequence for $ LZ^{\mathfrak{s}_{i+1}} $ reduces to the short exact sequence
$$ 0 \rightarrow M^{\mathfrak{s}_{i+1}}(a_{i+1}, \underbrace{a_{i+2}, a_{i}}) \rightarrow M^{\mathfrak{s}_{i+1}}(a_{i}, \underbrace{a_{i+2}, a_{i+1}}) \rightarrow
L_1 Z^{\mathfrak{s}_{i+1}} \epsilon_{\mathfrak{s}_i} M^{\mathfrak{s}_i}(\underbrace{a_{i}, a_{i+1}}, a_{i+2}) \rightarrow 0. $$

Thus if $ s=1 $ 
$$ L_s Z^{\mathfrak{s}_{i+1}} \epsilon_{\mathfrak{s}_i} M^{\mathfrak{s}_i}(\underbrace{a_{i}, a_{i+1}}, a_{i+2}) \cong
M^{\mathfrak{s}_{i+1}}(a_{i}, \underbrace{a_{i+2}, a_{i+1}})/M^{\mathfrak{s}_{i+1}}(a_{i+1}, \underbrace{a_{i+2}, a_{i}}) $$ 
and is 0 otherwise.

Now consider
$$ 0 \rightarrow M(a_{i+1}, a_i, a_{i+2}) \rightarrow M(a_{i+1}, a_{i+2}, a_{i}) \rightarrow M^{\mathfrak{s}_{i+1}}(a_{i+1}, \underbrace{a_{i+2}, a_{i}}) \rightarrow 0. $$
This gives a long exact sequence for $ LZ^{\mathfrak{s}_i} \epsilon_{\mathfrak{s}_{i+1}} $ reducing to
$$ 0 \rightarrow M^{\mathfrak{s}_i}(\underbrace{a_{i}, a_{i+1}}, a_{i+2}) \rightarrow
M^{\mathfrak{s}_i}(\underbrace{a_{i+2}, a_{i+1}}, a_{i}) \rightarrow
L_1 Z^{\mathfrak{s}_i} \epsilon_{\mathfrak{s}_{i+1}} M^{\mathfrak{s}_{i+1}}(a_{i+1}, \underbrace{a_{i+2}, a_{i}}) \rightarrow 0. $$
Thus if $ s=1, $
$$ L_s Z^{\mathfrak{s}_i} \epsilon_{\mathfrak{s}_{i+1}} M^{\mathfrak{s}_{i+1}}(a_{i+1}, \underbrace{a_{i+2}, a_{i}}) \cong
M^{\mathfrak{s}_i}(\underbrace{a_{i+2}, a_{i+1}}, a_{i})/M^{\mathfrak{s}_i}(\underbrace{a_{i}, a_{i+1}}, a_{i+2}) $$ 
and it is zero otherwise.

Now consider
$$ 0 \rightarrow M(a_{i}, a_{i+1}, a_{i+2}) \rightarrow M(a_{i}, a_{i+2}, a_{i+1}) \rightarrow
M^{\mathfrak{s}_{i+1}}(a_{i}, \underbrace{a_{i+2}, a_{i+1}},) \rightarrow 0. $$
This gives a long exact sequence for $ LZ^{\mathfrak{s}_i} \epsilon_{\mathfrak{s}_{i+1}} $ reducing to
$$ 0 \rightarrow M^{\mathfrak{s}_i}(\underbrace{a_{i+2}, a_{i}}, a_{i+1}) \rightarrow L_1 Z^{\mathfrak{s}_i} \epsilon_{\mathfrak{s}_{i+1}} M^{\mathfrak{s}_{i+1}}(a_{i}, \underbrace{a_{i+2}, a_{i+1}}) \rightarrow
M^{\mathfrak{s}_i}(\underbrace{a_{i}, a_{i+1}}, a_{i+2}) \rightarrow 0. $$

Thus the long exact sequence for 
$$ 0 \rightarrow M^{\mathfrak{s}_{i+1}}(a_{i+1}, \underbrace{a_{i+2}, a_{i}}) \rightarrow M^{\mathfrak{s}_{i+1}}(a_{i}, \underbrace{a_{i+2}, a_{i+1}})
\rightarrow
L_1 Z^{\mathfrak{s}_{i+1}} \epsilon_{\mathfrak{s}_i} M^{\mathfrak{s}_i}(\underbrace{a_{i}, a_{i+1}}, a_{i+2}) \rightarrow 0 $$ 
becomes
$$ 
M^{\mathfrak{s}_i}(\underbrace{a_{i+2}, a_{i+1}}, a_{i})/M^{\mathfrak{s}_i}(\underbrace{a_{i}, a_{i+1}}, a_{i+2})
\hookrightarrow L_1 Z^{\mathfrak{s}_i} \epsilon_{\mathfrak{s}_{i+1}} M^{\mathfrak{s}_{i+1}}(a_{i}, \underbrace{a_{i+2}, a_{i+1}}) \twoheadrightarrow
L_1 Z^{\mathfrak{s}_i} \epsilon_{\mathfrak{s}_{i+1}} L_1 Z^{\mathfrak{s}_{i+1}} \epsilon_{\mathfrak{s}_i} M^{\mathfrak{s}_i}(\underbrace{a_{i}, a_{i+1}}, a_{i+2}). $$

Thus we get the following diagram:

\begin{tiny}
\xymatrix{
&0 \ar[d]	&		&		\\
&M^{\mathfrak{s}_i}(\underbrace{a_{i}, a_{i+1}}, a_{i+2}) \ar[d]^f	&0 \ar[d]^h	&	\\
&M^{\mathfrak{s}_i}(\underbrace{a_{i+2}, a_{i+1}}, a_{i})\ar[d]^g	&M^{\mathfrak{s}_i}(\underbrace{a_{i+2}, a_{i}}, a_{i+1})\ar[d]^i	&	\\
&M^{\mathfrak{s}_i}(\underbrace{a_{i+2}, a_{i+1}}, a_{i})/M^{\mathfrak{s}_i}(\underbrace{a_{i}, a_{i+1}}, a_{i+2})\ar@{^{(}->}^j[r]	&L_1 Z^{\mathfrak{s}_i} \epsilon_{\mathfrak{s}_{i+1}} M^{\mathfrak{s}_{i+1}}(a_{i}, \underbrace{a_{i+2}, a_{i+1}})	\ar@{->>}[r]^k \ar[d]^l	&L_1 Z^{\mathfrak{s}_i} \epsilon_{\mathfrak{s}_{i+1}} L_1 Z^{\mathfrak{s}_{i+1}} \epsilon_{\mathfrak{s}_i} M^{\mathfrak{s}_i}(\underbrace{a_{i}, a_{i+1}}, a_{i+2})	\\
&	&M^{\mathfrak{s}_i}(\underbrace{a_{i}, a_{i+1}}, a_{i+2})\ar[d]	&	\\
&	&0	&	\\
}
\end{tiny}

Assume that $ l \circ j \neq 0. $  Then the following composition of maps is non-zero:
$$ M^{\mathfrak{s}_i}(\underbrace{a_{i+2}, a_{i+1}}, a_i) \rightarrow M^{\mathfrak{s}_i}(\underbrace{a_{i+2}, a_{i+1}}, a_i)/M^{\mathfrak{s}_i}(\underbrace{a_{i}, a_{i+1}}, a_{i+2}) \rightarrow M^{\mathfrak{s}_i}(\underbrace{a_{i}, a_{i+1}}, a_{i+2}) \rightarrow M^{\mathfrak{s}_i}(\underbrace{a_{i+2}, a_{i+1}}, a_i). $$
A non-zero endomorphism of a generalized Verma module must be an isomorphism so the first map above must be injective as well.  This is impossible so $ l \circ j =0. $
Thus 
$$ M^{\mathfrak{s}_i}(\underbrace{a_{i+2}, a_{i+1}}, a_{i})/M^{\mathfrak{s}_i}(\underbrace{a_{i}, a_{i+1}}, a_{i+2})
\hookrightarrow M^{\mathfrak{s}_i}(\underbrace{a_{i+2}, a_{i}}, a_{i+1}). $$

We claim that there is a well defined map from 
$ L_1 Z^{\mathfrak{s}_i} \epsilon_{\mathfrak{s}_{i+1}} L_1 Z^{\mathfrak{s}_{i+1}} \epsilon_{\mathfrak{s}_i} M^{\mathfrak{s}_i}(\underbrace{a_{i}, a_{i+1}}, a_{i+2}) $ to 
$ M^{\mathfrak{s}_i}(\underbrace{a_{i}, a_{i+1}}, a_{i+2}). $
Via $ k, $ pull an element back and then apply $ l. $
Suppose $ k(x)=k(y). $  Then $ k(x-y)=0 $ so $ x-y = j(z) $ for some $ z \in 
M^{\mathfrak{s}_i}(\underbrace{a_{i+2}, a_{i+1}}, a_{i})/M^{\mathfrak{s}_i}(\underbrace{a_{i}, a_{i+1}}, a_{i+2}). $
Then $ l(x-y) = lj(z). $ This is zero by the observation above.  Thus the map is well defined and call it $ \phi. $  
It is easy to see that $ ki $ is a surjection of $ M^{\mathfrak{s}_i}(\underbrace{a_{i+2}, a_i}, a_{i+1}) $ onto $ \text{ker} \phi. $ The kernel of this surjection is obviously
$ M^{\mathfrak{s}_i}(\underbrace{a_{i+2}, a_{i+1}}, a_{i})/M^{\mathfrak{s}_i}(\underbrace{a_{i}, a_{i+1}}, a_{i+2}). $
Thus the kernel of $ \phi $ is
$$ M^{\mathfrak{s}_i}(\underbrace{a_{i+2}, a_{i}}, a_{i+1})/(M^{\mathfrak{s}_i}(\underbrace{a_{i+2}, a_{i+1}}, a_{i})/M^{\mathfrak{s}_i}(\underbrace{a_{i}, a_{i+1}}, a_{i+2})). $$
Since this module is isomorphic to 
$ \epsilon_{\mathfrak{t}_i}^{\mathfrak{s}_i} M^{\mathfrak{t}_i}(\underbrace{a_{i+2}, a_{i}, a_{i+1}}), $ there is a short exact sequence
$$ 0 \rightarrow \epsilon_{\mathfrak{t}_i}^{\mathfrak{s}_i} M^{\mathfrak{t}_i}(\underbrace{a_{i+2}, a_{i}, a_{i+1}}) \rightarrow
L_1 Z^{\mathfrak{s}_i} \epsilon_{\mathfrak{s}_{i+1}} L_1 Z^{\mathfrak{s}_{i+1}} \epsilon_{\mathfrak{s}_i} M^{\mathfrak{s}_i}(\underbrace{a_{i}, a_{i+1}}, a_{i+2}) \rightarrow 
M^{\mathfrak{s}_i}(\underbrace{a_{i}, a_{i+1}}, a_{i+2}) \rightarrow 0. $$

Finally we must show that this exact sequence splits.
\begin{align*}
\Hom(M^{\mathfrak{s}_i}(\underbrace{a_{i}, a_{i+1}}, a_{i+2}), \epsilon_{\mathfrak{t}_i}^{\mathfrak{s}_i} M^{\mathfrak{t}_i}(\underbrace{a_{i+2}, a_{i}, a_{i+1}})[1]) 
&\cong \Hom(LZ_{\mathfrak{s}_i}^{\mathfrak{t}_i} M^{\mathfrak{s}_i}(\underbrace{a_{i}, a_{i+1}}, a_{i+2}), M^{\mathfrak{t}_i}(\underbrace{a_{i+2}, a_{i}, a_{i+1}})[1])\\
&\cong \Hom(M^{\mathfrak{t}_i}(\underbrace{a_{i+2}, a_{i}, a_{i+1}})[2], M^{\mathfrak{t}_i}(\underbrace{a_{i+2}, a_{i}, a_{i+1}})[1]). 
\end{align*} 
This last space is zero so the sequence must split.
\end{proof}

\begin{lemma}
\label{lemma26}
Suppose $ a_i > a_{i+2} > a_{i+1}. $  Then 
$$ LZ^{\mathfrak{s}_i} \epsilon_{\mathfrak{s}_{i+1}}[-1] LZ^{\mathfrak{s}_{i+1}} \epsilon_{\mathfrak{s}_i}[-1] M^{\mathfrak{s}_i}(\underbrace{a_i, a_{i+1}}, a_{i+2}) \cong \\
M^{\mathfrak{s}_i}(\underbrace{a_i, a_{i+1}}, a_{i+2}) \oplus
\epsilon_{\mathfrak{t}_i}^{\mathfrak{s}_i}M^{\mathfrak{t}_i}(\underbrace{a_{i}, a_{i+2}, a_{i+1}})[-1]. $$ 
\end{lemma}

\begin{proof}
Consider the short exact sequence
$$ 0 \rightarrow M(a_{i+1}, a_i, a_{i+2}) \rightarrow M(a_{i}, a_{i+1}, a_{i+2}) \rightarrow\
M^{\mathfrak{s}_i}(\underbrace{a_{i}, a_{i+1}}, a_{i+2}) \rightarrow 0.  $$  
The long exact sequence for $ LZ^{\mathfrak{s}_{i+1}} $ reduces to the short exact sequence
$$ 0 \rightarrow M^{\mathfrak{s}_{i+1}}(a_i, \underbrace{a_{i+2}, a_{i+1}}) \rightarrow
L_1 Z^{\mathfrak{s}_{i+1}} \epsilon_{\mathfrak{s}_i} M^{\mathfrak{s}_i}(\underbrace{a_{i}, a_{i+1}}, a_{i+2}) \rightarrow
M^{\mathfrak{s}_{i+1}}(a_{i+1}, \underbrace{a_{i}, a_{i+2}}) \rightarrow 0. $$
We compute $ LZ^{\mathfrak{s}_i} \epsilon_{\mathfrak{s}_{i+1}} $ on the first and third terms of this short exact sequence.

Consider 
$$ 0 \rightarrow M(a_{i+1}, a_{i+2}, a_{i}) \rightarrow M(a_{i+1}, a_{i}, a_{i+2}) \rightarrow
M^{\mathfrak{s}_{i+1}}(a_{i+1}, \underbrace{a_{i}, a_{i+2}}) \rightarrow 0. $$  
The long exact sequence gives
$$ L_s Z^{\mathfrak{s}_i} \epsilon_{\mathfrak{s}_{i+1}} M^{\mathfrak{s}_{i+1}}(a_{i+1}, \underbrace{a_{i}, a_{i+2}}) \cong\
M^{\mathfrak{s}_i}(\underbrace{a_i, a_{i+1}}, a_{i+2})/M^{\mathfrak{s}_i}(\underbrace{a_{i+2}, a_{i+1}}, a_{i})  $$
if $ s = 1 $ and is zero otherwise.

Next consider 
$$ 0 \rightarrow M(a_{i}, a_{i+1}, a_{i+2}) \rightarrow M(a_{i}, a_{i+2}, a_{i+1}) \rightarrow
M^{\mathfrak{s}_{i+1}}(a_{i}, \underbrace{a_{i+2}, a_{i+1}}) \rightarrow 0. $$  
The long exact sequence for $ LZ^{\mathfrak{s}_i} $ becomes
$$ 
0 \rightarrow L_1 Z^{\mathfrak{s}_i} \epsilon_{\mathfrak{s}_{i+1}} M^{\mathfrak{s}_{i+1}}(a_{i}, \underbrace{a_{i+2}, a_{i+1}}) \rightarrow
M^{\mathfrak{s}_i}(\underbrace{a_i, a_{i+1}}, a_{i+2}) \rightarrow
M^{\mathfrak{s}_i}(\underbrace{a_i, a_{i+2}}, a_{i+1}) \rightarrow
L_0 Z^{\mathfrak{s}_i} \epsilon_{\mathfrak{s}_{i+1}} M^{\mathfrak{s}_{i+1}}(a_{i}, \underbrace{a_{i+2}, a_{i+1}}) \rightarrow 0. $$
The middle map is clearly the standard map.  
Therefore
$$ L_0 Z^{\mathfrak{s}_i} \epsilon_{\mathfrak{s}_{i+1}} M^{\mathfrak{s}_{i+1}}(a_{i}, \underbrace{a_{i+2}, a_{i+1}}) \cong 
\epsilon_{\mathfrak{t}_i}^{\mathfrak{s}_i} M^{\mathfrak{t}_i}(\underbrace{a_{i}, a_{i+2}, a_{i+1}}) $$ 
and
$$ L_1 Z^{\mathfrak{s}_i} \epsilon_{\mathfrak{s}_{i+1}} M^{\mathfrak{s}_{i+1}}(a_{i}, \underbrace{a_{i+2}, a_{i+1}}) \cong
M^{\mathfrak{s}_i}(\underbrace{a_{i+2}, a_{i+1}}, a_{i}). $$

Now the long exact sequence for the functor $ LZ^{\mathfrak{s}_i} $ and
$$ 
0 \rightarrow M^{\mathfrak{s}_{i+1}}(a_{i}, \underbrace{a_{i+2}, a_{i+1}}) \rightarrow
L_1 Z^{\mathfrak{s}_{i+1}} \epsilon_{\mathfrak{s}_i} M^{\mathfrak{s}_i}(\underbrace{a_{i}, a_{i+1},} a_{i+2}) \rightarrow
M^{\mathfrak{s}_{i+1}}(a_{i+1}, \underbrace{a_{i}, a_{i+2}}) \rightarrow 0. $$ 
is
\begin{align*}
0 \rightarrow &M^{\mathfrak{s}_i}(\underbrace{a_{i+2}, a_{i+1}}, a_{i}) \rightarrow
L_1 Z^{\mathfrak{s}_i} \epsilon_{\mathfrak{s}_{i+1}} L_1 Z^{\mathfrak{s}_{i+1}}  M^{\mathfrak{s}_i}(\underbrace{a_{i}, a_{i+1}}, a_{i+2}) \rightarrow\\
&M^{\mathfrak{s}_i}(\underbrace{a_i, a_{i+1}}, a_{i+2})/M^{\mathfrak{s}_i}(\underbrace{a_{i+2}, a_{i+1}}, a_{i}) \rightarrow
\epsilon_{\mathfrak{t}_i}^{\mathfrak{s}_i} M^{\mathfrak{t}_i}(\underbrace{a_{i}, a_{i+2}, a_{i+1}}) \rightarrow
L_0 Z^{\mathfrak{s}_i} \epsilon_{\mathfrak{s}_{i+1}} L_1 Z^{\mathfrak{s}_{i+1}} M^{\mathfrak{s}_i}(\underbrace{a_{i}, a_{i+1}}, a_{i+2}) \rightarrow 0.  
\end{align*} 

Suppose that
$$ \Hom(M^{\mathfrak{s}_i}(\underbrace{a_i, a_{i+1}}, a_{i+2})/M^{s_i}(\underbrace{a_{i+2}, a_{i+1}}, a_{i}),
\epsilon_{\mathfrak{t}_i}^{\mathfrak{s}_i} M^{\mathfrak{t}_i}(\underbrace{a_{i}, a_{i+2}, a_{i+1}})) \neq 0.  $$
Then 
\begin{align*}
0 &\neq \Hom(M^{\mathfrak{s}_i}(\underbrace{a_i, a_{i+1}}, a_{i+2}), \epsilon_{\mathfrak{t}_i}^{\mathfrak{s}_i} M^{\mathfrak{t}_i}(\underbrace{a_{i}, a_{i+2}, a_{i+1}}))\\
&\cong \Hom(LZ_{\mathfrak{s}_i}^{\mathfrak{t}_i} M^{\mathfrak{s}_i}(\underbrace{a_{i}, a_{i+1}}, a_{i+2}),
M^{\mathfrak{t}_i}(\underbrace{a_{i}, a_{i+2}, a_{i+1}}))\\
&\cong \Hom(M^{\mathfrak{t}_i}(\underbrace{a_{i}, a_{i+2}, a_{i+1}})[1], 
M^{\mathfrak{t}_i}(\underbrace{a_{i}, a_{i+2}, a_{i+1}}))\\ 
&= 0. 
\end{align*} 

Thus
$ L_0 Z^{\mathfrak{s}_i} \epsilon_{\mathfrak{s}_{i+1}} L_1 Z^{\mathfrak{s}_{i+1}} M^{\mathfrak{s}_i}(\underbrace{a_{i}, a_{i+1}}, a_{i+2}) \cong
\epsilon_{\mathfrak{t}_i}^{\mathfrak{s}_i} M^{\mathfrak{t}_i}(\underbrace{a_{i}, a_{i+2}, a_{i+1}}) $ and there is the short exact sequence
$$ 0 \rightarrow M^{\mathfrak{s}_i}(\underbrace{a_{i+2}, a_{i+1}}, a_{i}) \rightarrow
L_1 Z^{\mathfrak{s}_i} \epsilon_{\mathfrak{s}_{i+1}} L_1 Z^{\mathfrak{s}_{i+1}} M^{\mathfrak{s}_i}(\underbrace{a_{i}, a_{i+1}}, a_{i+2}) \rightarrow
M^{\mathfrak{s}_i}(\underbrace{a_i, a_{i+1}}, a_{i+2})/M^{\mathfrak{s}_i}(\underbrace{a_{i+2}, a_{i+1}}, a_{i}) \rightarrow 0. $$

Next we show that 
$$ \Hom(L_1 Z^{\mathfrak{s}_i} \epsilon_{\mathfrak{s}_{i+1}} L_1 Z^{\mathfrak{s}_{i+1}} \epsilon_{s_i} M^{\mathfrak{s}_i}(\underbrace{a_{i}, a_{i+1}}, a_{i+2}), 
M^{\mathfrak{s}_i}(\underbrace{a_{i+2}, a_{i+1}}, a_{i})) \cong 0 $$ 
so that the sequence does not split.
The distinguished triangle
\begin{align*}
&H^0 LZ^{\mathfrak{s}_i} \epsilon_{\mathfrak{s}_{i+1}}[-1] LZ^{\mathfrak{s}_{i+1}} \epsilon_{s_i}[-1] M^{\mathfrak{s}_i}(\underbrace{a_{i}, a_{i+1}}, a_{i+2}) \rightarrow
\tau^{\geq 0} LZ^{\mathfrak{s}_i} \epsilon_{\mathfrak{s}_{i+1}}[-1] LZ^{\mathfrak{s}_{i+1}} \epsilon_{s_i}[-1] M^{\mathfrak{s}_i}(\underbrace{a_{i}, a_{i+1}}, a_{i+2}) \rightarrow\\
&\tau^{\geq 1} LZ^{\mathfrak{s}_i} \epsilon_{\mathfrak{s}_{i+1}}[-1] LZ^{\mathfrak{s}_{i+1}} \epsilon_{s_i}[-1] M^{\mathfrak{s}_i}(\underbrace{a_{i}, a_{i+1}}, a_{i+2}) 
\end{align*} 
is isomorphic to
$$ L_1 Z^{\mathfrak{s}_i} \epsilon_{\mathfrak{s}_{i+1}} L_1 Z^{\mathfrak{s}_{i+1}} \epsilon_{s_i} M^{\mathfrak{s}_i}(\underbrace{a_{i}, a_{i+1}}, a_{i+2}) \rightarrow
LZ^{\mathfrak{s}_i} \epsilon_{\mathfrak{s}_{i+1}}[-1] LZ^{\mathfrak{s}_{i+1}} \epsilon_{s_i}[-1] M^{\mathfrak{s}_i}(\underbrace{a_{i}, a_{i+1}}, a_{i+2}) \rightarrow
\epsilon_{t_i}^{s_i} M^{t_i}(\underbrace{a_i, a_{i+2}, a_{i+1}})[-1]. $$
Now apply the functor $ \Hom(\bullet, M^{\mathfrak{s}_i}(\underbrace{a_{i+2}, a_{i+1}}, a_{i})) $ to get a long exact sequence
\begin{align*}
\cdots \rightarrow 
&\Hom(LZ^{\mathfrak{s}_i} \epsilon_{\mathfrak{s}_{i+1}}[-1] LZ^{\mathfrak{s}_{i+1}} \epsilon_{s_i}[-1] M^{\mathfrak{s}_i}(\underbrace{a_{i}, a_{i+1}}, a_{i+2}), M^{\mathfrak{s}_i}(\underbrace{a_{i+2}, a_{i+1}}, a_{i})) \rightarrow\\
&\Hom(L_1 Z^{\mathfrak{s}_i} \epsilon_{\mathfrak{s}_{i+1}} L_1 Z^{\mathfrak{s}_{i+1}} \epsilon_{s_i} M^{\mathfrak{s}_i}(\underbrace{a_{i}, a_{i+1}}, a_{i+2}), 
M^{\mathfrak{s}_i}(\underbrace{a_{i+2}, a_{i+1}}, a_{i})) \rightarrow 0. 
\end{align*} 

We calculate
\begin{align*}
&\Hom(LZ^{\mathfrak{s}_i} \epsilon_{\mathfrak{s}_{i+1}}[-1] LZ^{\mathfrak{s}_{i+1}} \epsilon_{s_i}[-1] M^{\mathfrak{s}_i}(\underbrace{a_{i}, a_{i+1}}, a_{i+2}), M^{\mathfrak{s}_i}(\underbrace{a_{i+2}, a_{i+1}}, a_{i})) \cong\\
&\Hom(M^{\mathfrak{s}_i}(\underbrace{a_{i}, a_{i+1}}, a_{i+2}), LZ^{\mathfrak{s}_i} \epsilon_{\mathfrak{s}_{i+1}}[-1] LZ^{\mathfrak{s}_{i+1}} \epsilon_{s_i}[-1] M^{\mathfrak{s}_i}(\underbrace{a_{i+2}, a_{i+1}}, a_{i})) \cong\\
&\Hom(M^{\mathfrak{s}_i}(\underbrace{a_{i}, a_{i+1}}, a_{i+2}), M^{\mathfrak{s}_i}(\underbrace{a_{i+2}, a_{i+1}}, a_{i}) \oplus \epsilon_{t_i}^{s_i} M^{t_i}(\underbrace{a_i, a_{i+2}, a_{i+1}})) \cong\\
&\Hom(M^{\mathfrak{s}_i}(\underbrace{a_{i}, a_{i+1}}, a_{i+2}), \epsilon_{t_i}^{s_i} M^{t_i}(\underbrace{a_i, a_{i+2}, a_{i+1}})) \cong\\
&\Hom(LZ_{s_i}^{t_i} M^{\mathfrak{s}_i}(\underbrace{a_{i}, a_{i+1}}, a_{i+2}), M^{t_i}(\underbrace{a_i, a_{i+2}, a_{i+1}})) \cong\\
&\Hom(M^{t_i}(\underbrace{a_i, a_{i+2}, a_{i+1}})[1], M^{t_i}(\underbrace{a_i, a_{i+2}, a_{i+1}})) \cong 0. 
\end{align*} 
The second isomorphism is due to the previous lemma.  Therefore the sequence cannot split.

It remains to be shown that the extension is unique up to isomorphism.
Consider the full subcategory of objects which have composition factors factors with weights less than or equal to
$ (\underbrace{a_i, a_{i+1}}, a_{i+2}). $ In this full subcategory the generalized Verma module $ M^{\mathfrak{s}_i}(\underbrace{a_i, a_{i+1}}, a_{i+2}) $ is projective.
Thus for any exact sequence
$$ 0 \rightarrow M^{\mathfrak{s}_i}(\underbrace{a_{i+2}, a_{i+1}}, a_i) \rightarrow X \rightarrow M^{\mathfrak{s}_i}(\underbrace{a_{i}, a_{i+1}}, a_{i+2})/M^{\mathfrak{s}_i}(\underbrace{a_{i+2}, a_{i+1}}, a_i) \rightarrow 0 $$
there is a nontrivial map $ g \colon M^{\mathfrak{s}_i}(\underbrace{a_{i}, a_{i+1}}, a_{i+2}) \rightarrow X. $  Note that all of the objects in the above sequence are in this subcategory.
Now we can set up the following diagram.  

\begin{tiny}
$$
\CD
0 @>>> M^{\mathfrak{s}_i}(\underbrace{a_{i+2}, a_{i+1}}, a_i) @>j>> M^{\mathfrak{s}_i}(\underbrace{a_{i}, a_{i+1}}, a_{i+2}) @>f>> M^{\mathfrak{s}_i}(\underbrace{a_{i}, a_{i+1}}, a_{i+2})/
M^{\mathfrak{s}_i}(\underbrace{a_{i+2}, a_{i+1}}, a_i) @>>> 0\\
@.		@VVIdV		@VV g V		@VVIdV		@.\\
0 @>>> M^{\mathfrak{s}_i}(\underbrace{a_{i+2}, a_{i+1}}, a_i)  @>i>>	X @>h>> M^{\mathfrak{s}_i}(\underbrace{a_{i}, a_{i+1}}, a_{i+2})/M^{\mathfrak{s}_i}(\underbrace{a_{i+2}, a_{i+1}}, a_i) @>>> 0\\
\endCD
$$
\end{tiny}

The square on the right commutes by construction.
Take an element in the top left.  Applying the map $ hgj $ must give zero.  
Therefore $ \text{im} (gj) \subset \text{ker}(h) \cong M^{\mathfrak{s}_i}(\underbrace{a_{i+2}, a_{i+1}}, a_i). $
But then $ i $ must be the same as $ gj $ up to a scalar since scalar multiplication are the only endomorphisms of
$ M^{\mathfrak{s}_i}(\underbrace{a_{i+2}, a_{i+1}}, a_i). $
If the scalar is zero then there is a map from
$ M^{\mathfrak{s}_i}(\underbrace{a_{i}, a_{i+1}}, a_{i+2})/M^{\mathfrak{s}_i}(\underbrace{a_{i+2}, a_{i+1}}, a_i) $ to $ X $ giving a splitting, contrary to the above where it was shown that the extension
must be non-split.
Thus the scalar must be non-zero and so the middle map is an isomorphism by the five lemma.

Therefore,
$$ L_1 Z^{\mathfrak{s}_i} \epsilon_{\mathfrak{s}_{i+1}} L_1 Z^{\mathfrak{s}_{i+1}} M^{\mathfrak{s}_i}(\underbrace{a_{i}, a_{i+1}}, a_{i+2}) \cong
M^{\mathfrak{s}_i}(\underbrace{a_i, a_{i+1}}, a_{i+2}). $$  
\end{proof}

\begin{lemma}
Suppose $ a_{i} > a_{i+1} > a_{i+2}. $  Then
$$ LZ^{\mathfrak{s}_i} \epsilon_{\mathfrak{s}_{i+1}}[-1] LZ^{\mathfrak{s}_{i+1}} \epsilon_{\mathfrak{s}_i}[-1] M^{\mathfrak{s}_i}(\underbrace{a_i, a_{i+1}}, a_{i+2}) \cong
M^{\mathfrak{s}_i}(\underbrace{a_i, a_{i+1}}, a_{i+2}) \oplus
\epsilon_{\mathfrak{t}_i}^{\mathfrak{s}_i}M^{\mathfrak{t}_i}(\underbrace{a_{i}, a_{i+1}, a_{i+2}})[-2]. $$ 
\end{lemma}

\begin{proof}
Let $ X = LZ^{\mathfrak{s}_{i+1}} \epsilon_{\mathfrak{s}_i}[-1] M^{\mathfrak{s}_i}(\underbrace{a_i, a_{i+1}}, a_{i+2}). $

Consider the short exact sequence
$$ 0 \rightarrow M(a_{i+1}, a_{i}, a_{i+2}) \rightarrow
M(a_i, a_{i+1}, a_{i+2}) \rightarrow M^{\mathfrak{s}_i}(\underbrace{a_i, a_{i+1}}, a_{i+2}) \rightarrow 0. $$
This induces a long exact sequence for $ LZ^{\mathfrak{s}_{i+1}} $
$$ 
L_1 Z^{\mathfrak{s}_{i+1}} \epsilon_{\mathfrak{s}_i} M^{\mathfrak{s}_i}(\underbrace{a_i, a_{i+1}}, a_{i+2}) \hookrightarrow
M^{\mathfrak{s}_{i+1}}(a_{i+1}, \underbrace{a_{i}, a_{i+2}}) \rightarrow
M^{\mathfrak{s}_{i+1}}(a_{i}, \underbrace{a_{i+1}, a_{i+2}}) \twoheadrightarrow
L_0 Z^{\mathfrak{s}_{i+1}} \epsilon_{\mathfrak{s}_i} M^{\mathfrak{s}_i}(\underbrace{a_i, a_{i+1}}, a_{i+2}). $$
The middle map is clearly the standard map so
$$ L_1 Z^{\mathfrak{s}_{i+1}} \epsilon_{\mathfrak{s}_i} M^{\mathfrak{s}_i}(\underbrace{a_i, a_{i+1}}, a_{i+2}) \cong 
M^{\mathfrak{s}_{i+1}}(a_{i+2}, \underbrace{a_{i}, a_{i+1}}) $$ 
and
$$ L_0 Z^{\mathfrak{s}_{i+1}} \epsilon_{\mathfrak{s}_i} M^{\mathfrak{s}_i}(\underbrace{a_i, a_{i+1}}, a_{i+2}) \cong
\epsilon_{\mathfrak{t}_i}^{\mathfrak{s}_i}M^{\mathfrak{t}_i}(\underbrace{a_{i}, a_{i+1}, a_{i+2}}). $$
Thus
\begin{align*}
H^0(X) &\cong M^{\mathfrak{s}_{i+1}}(a_{i+2}, \underbrace{a_{i}, a_{i+1}})\\
H^1(X) &\cong \epsilon_{\mathfrak{t}_i}^{\mathfrak{s}_i}M^{\mathfrak{t}_i}(\underbrace{a_{i}, a_{i+1}, a_{i+2}})\\
H^s(X) &\cong 0, \forall s\neq 0, 1. 
\end{align*} 

Therefore we get the distinguished triangle
$$ H^0(X)[0] \rightarrow \tau^{\geq 0}X \rightarrow \tau^{\geq 1}X $$
which is
\begin{equation}
M^{\mathfrak{s}_{i+1}}(a_{i+2}, \underbrace{a_{i}, a_{i+1}}) \rightarrow X \rightarrow
\epsilon_{\mathfrak{t}_i}^{\mathfrak{s}_{i+1}}M^{\mathfrak{t}_i}(\underbrace{a_{i}, a_{i+1}, a_{i+2}})[-1]. 
\end{equation}

Next we compute $ L_s Z^{\mathfrak{s}_i} \epsilon_{\mathfrak{s}_{i+1}} M^{\mathfrak{s}_{i+1}}(a_{i+2}, \underbrace{a_{i}, a_{i+1}}). $
Consider the short exact sequence
$$ 0 \rightarrow M(a_{i+2}, a_{i+1}, a_{i}) \rightarrow M(a_{i+2}, a_{i}, a_{i+1}) \rightarrow M^{\mathfrak{s}_{i+1}}(a_{i+2}, \underbrace{a_{i}, a_{i+1}}) \rightarrow 0. $$
This gives a long exact sequence for $ LZ^{\mathfrak{s}_i}: $
$$ 0 \rightarrow M^{\mathfrak{s}_i}(\underbrace{a_{i+1}, a_{i+2}}, a_{i}) \rightarrow M^{\mathfrak{s}_i}(\underbrace{a_i, a_{i+2}}, a_{i+1}) \rightarrow 
L_1 Z^{\mathfrak{s}_i} \epsilon_{\mathfrak{s}_{i+1}} M^{\mathfrak{s}_{i+1}}(a_{i+2}, \underbrace{a_{i}, a_{i+1}}) \rightarrow 0. $$
Thus if $ s=1 $
$$ L_s Z^{\mathfrak{s}_i} \epsilon_{\mathfrak{s}_{i+1}} M^{\mathfrak{s}_{i+1}}(a_{i+2}, \underbrace{a_{i}, a_{i+1}}, \ldots a_n) \cong
M^{\mathfrak{s}_i}(\underbrace{a_i, a_{i+2}}, a_{i+1})/M^{\mathfrak{s}_i}(\underbrace{a_{i+1}, a_{i+2}}) $$
and it is zero otherwise.

Next to compute $ L_s Z^{\mathfrak{s}_i} \epsilon_{\mathfrak{s}_{i+1}} \epsilon_{\mathfrak{t}_i}^{\mathfrak{s}_{i+1}}M^{\mathfrak{t}_i}(\underbrace{a_{i}, a_{i+1}, a_{i+2}}): $
\begin{align*}
LZ^{\mathfrak{s}_i} \epsilon_{\mathfrak{s}_{i+1}} \epsilon_{\mathfrak{t}_i}^{\mathfrak{s}_{i+1}}M^{\mathfrak{t}_i}(\underbrace{a_{i}, a_{i+1}, a_{i+2}}) &\cong
LZ^{\mathfrak{s}_i} \epsilon_{\mathfrak{s}_i} \epsilon_{\mathfrak{t}_i}^{\mathfrak{s}_i}M^{\mathfrak{t}_i}(\underbrace{a_{i}, a_{i+1}, a_{i+2}})\\ 
&\cong 
\epsilon_{\mathfrak{t}_i}^{\mathfrak{s}_i}M^{\mathfrak{t}_i}(\underbrace{a_{i}, a_{i+1}, a_{i+2}}) \oplus
\epsilon_{\mathfrak{t}_i}^{\mathfrak{s}_i}M^{\mathfrak{t}_i}(\underbrace{a_{i}, a_{i+1}, a_{i+2}})[2]. 
\end{align*} 
The last isomorphism follows from a previous section.

Consider (1) again.  It gives rise to a distinguished triangle
$$ LZ^{\mathfrak{s}_i} \epsilon_{\mathfrak{s}_{i+1}}[-1] M^{\mathfrak{s}_{i+1}}(a_{i+2}, \underbrace{a_{i}, a_{i+1}}) \rightarrow LZ^{\mathfrak{s}_i} \epsilon_{\mathfrak{s}_{i+1}}[-1] X \rightarrow
LZ^{\mathfrak{s}_i} \epsilon_{\mathfrak{s}_{i+1}} [-1]\epsilon_{\mathfrak{t}_i}^{\mathfrak{s}_i}M^{\mathfrak{t}_i}(\underbrace{a_{i}, a_{i+1}, a_{i+2}})[-1]. $$
This gives rise to a long exact sequence
\begin{align*}
\rightarrow &H^s(LZ^{\mathfrak{s}_i} \epsilon_{\mathfrak{s}_{i+1}}[-1] M^{\mathfrak{s}_{i+1}}(a_{i+2}, \underbrace{a_{i}, a_{i+1}})) \rightarrow
H^s(LZ^{\mathfrak{s}_i} \epsilon_{\mathfrak{s}_{i+1}}[-1] X) \rightarrow\\
&H^s(LZ^{\mathfrak{s}_i} \epsilon_{\mathfrak{s}_{i+1}}[-1] \epsilon_{\mathfrak{t}_i}^{\mathfrak{s}_i}M^{\mathfrak{t}_i}(\underbrace{a_{i}, a_{i+1}, a_{i+2}})[-1]) \rightarrow 
\end{align*} 
This becomes
\begin{align*}
0 \rightarrow &M^{\mathfrak{s}_i}(\underbrace{a_i, a_{i+2}}, a_{i+1})/M^{\mathfrak{s}_i}(\underbrace{a_{i+1}, a_{i+2}}, a_{i}) \rightarrow
H^0(LZ^{\mathfrak{s}_i} \epsilon_{\mathfrak{s}_{i+1}}[-1] X) \rightarrow \epsilon_{\mathfrak{t}_i}^{\mathfrak{s}_i}M^{\mathfrak{t}_i}(\underbrace{a_{i}, a_{i+1}, a_{i+2}}) \rightarrow 0\rightarrow\\ 
&H^1(LZ^{\mathfrak{s}_i} \epsilon_{\mathfrak{s}_{i+1}}[-1] X) \rightarrow
L_1 Z^{\mathfrak{s}_i} \epsilon_{\mathfrak{s}_{i+1}} \epsilon_{\mathfrak{t}_i}^{\mathfrak{s}_i}M^{\mathfrak{t}_i}(\underbrace{a_{i}, a_{i+1}, a_{i+2}}) \rightarrow 0 \rightarrow\\
&H^2(LZ^{\mathfrak{s}_i} \epsilon_{\mathfrak{s}_{i+1}}[-1] X) \rightarrow \epsilon_{\mathfrak{t}_i}^{\mathfrak{s}_i}M^{\mathfrak{t}_i}(\underbrace{a_{i}, a_{i+1}, a_{i+2}}) \rightarrow 0.  \end{align*} 

Thus 
$$ H^2(LZ^{\mathfrak{s}_i} \epsilon_{\mathfrak{s}_{i+1}}[-1]X) \cong \epsilon_{\mathfrak{t}_i}^{\mathfrak{s}_i}M^{\mathfrak{t}_i}(\underbrace{a_{i}, a_{i+1}, a_{i+2}}) $$ 
and
$$ H^1 LZ^{\mathfrak{s}_i} \epsilon_{\mathfrak{s}_{i+1}}[-1] X \cong 0. $$

We must prove $ H^0(LZ^{\mathfrak{s}_i} \epsilon_{\mathfrak{s}_{i+1}}[-1]X) \cong M^{\mathfrak{s}_i}(\underbrace{a_{i}, a_{i+1}}, a_{i+2}). $
This reduces to showing that
$$ 0 \rightarrow M^{\mathfrak{s}_i}(\underbrace{a_i, a_{i+2}}, a_{i+1})/M^{\mathfrak{s}_i}(\underbrace{a_{i+1}, a_{i+2}}, a_{i})
\rightarrow H^0(LZ^{\mathfrak{s}_i} \epsilon_{\mathfrak{s}_{i+1}}[-1] X) \rightarrow \epsilon_{\mathfrak{t}_i}^{\mathfrak{s}_i}M^{\mathfrak{t}_i} (\underbrace{a_{i}, a_{i+1}, a_{i+2}}) \rightarrow 0 $$ 
is the unique non-trivial extension up to isomorphism.

The distinguished triangle
\begin{align*}
&H^0 (LZ^{s_{i}} \epsilon_{s_{i+1}}[-1] LZ^{s_{i+1}} \epsilon_{s_i}[-1] M^{s_i}(\underbrace{a_i, a_{i+1},} a_{i+2})) \rightarrow 
\tau^{\geq 0} (LZ^{s_{i}} \epsilon_{s_{i+1}}[-1] LZ^{s_{i+1}} \epsilon_{s_i}[-1] M^{s_i}(\underbrace{a_i, a_{i+1},} a_{i+2})) \rightarrow\\ 
&\tau^{\geq 1} (LZ^{s_{i}} \epsilon_{s_{i+1}}[-1] LZ^{s_{i+1}} \epsilon_{s_i}[-1] M^{s_i}(\underbrace{a_i, a_{i+1},} a_{i+2}))
\end{align*} 
is isomorphic to
\begin{align*}
&H^0 (LZ^{s_{i}} \epsilon_{s_{i+1}}[-1] LZ^{s_{i+1}} \epsilon_{s_i}[-1] M^{s_i}(\underbrace{a_i, a_{i+1},} a_{i+2})) \rightarrow 
LZ^{s_{i}} \epsilon_{s_{i+1}}[-1] LZ^{s_{i+1}} \epsilon_{s_i}[-1] M^{s_i}(\underbrace{a_i, a_{i+1},} a_{i+2}) \rightarrow\\
&\epsilon_{t_i}^{s_i} M^{t_i}(\underbrace{a_i, a_{i+1}, a_{i+2}})[-2]. 
\end{align*} 
Apply the functor  $ \Hom(\epsilon_{t_i}^{s_i} M^{t_i}(\underbrace{a_i, a_{i+1}, a_{i+2}}), \bullet) $ to get a long exact sequence
\begin{align*}
0 \rightarrow 
&\Hom(\epsilon_{t_i}^{s_i} M^{t_i}(\underbrace{a_i, a_{i+1}, a_{i+2}}), H^0 LZ^{s_{i}} \epsilon_{s_{i+1}}[-1] LZ^{s_{i+1}} \epsilon_{s_i}[-1] M^{s_i}(\underbrace{a_i, a_{i+1},} a_{i+2})) \rightarrow\\
&\Hom(\epsilon_{t_i}^{s_i} M^{t_i}(\underbrace{a_i, a_{i+1}, a_{i+2}}), LZ^{s_{i}} \epsilon_{s_{i+1}}[-1] LZ^{s_{i+1}} \epsilon_{s_i}[-1] M^{s_i}(\underbrace{a_i, a_{i+1},} a_{i+2})) \rightarrow \cdots. 
\end{align*} 
We calculate
\begin{align*}
&\Hom(\epsilon_{t_i}^{s_i} M^{t_i}(\underbrace{a_i, a_{i+1}, a_{i+2}}), LZ^{s_{i}} \epsilon_{s_{i+1}}[-1] LZ^{s_{i+1}} \epsilon_{s_i}[-1] M^{s_i}(\underbrace{a_i, a_{i+1},} a_{i+2})) \cong\\
&\Hom(M^{t_i}(\underbrace{a_i, a_{i+1}, a_{i+2}}), LZ_{s_i}^{t_i}[-4] LZ^{s_{i}} \epsilon_{s_{i+1}}[-1] LZ^{s_{i+1}} \epsilon_{s_i}[-1] M^{s_i}(\underbrace{a_i, a_{i+1},} a_{i+2})) \cong\\
&\Hom(M^{t_i}(\underbrace{a_i, a_{i+1}, a_{i+2}}), M^{t_i}(\underbrace{a_i, a_{i+1}, a_{i+2}})[-6] \oplus (M^{t_i}(\underbrace{a_i, a_{i+1}, a_{i+2}})[-4])^{\oplus 2} 
\oplus M^{t_i}(\underbrace{a_i, a_{i+1}, a_{i+2}})[-2])\\ 
& = 0. 
\end{align*} 
The second to last isomorphism comes from using proposition 8 twice.  This shows that the extension is non-trivial.  Now we show that it is unique.

Consider the short exact sequence
$$ 0 \rightarrow M^{\mathfrak{s}_i}(\underbrace{a_{i+1}, a_{i+2}}, a_{i}) \rightarrow
M^{\mathfrak{s}_i}(\underbrace{a_i, a_{i+2}}, a_{i+1}) \rightarrow
M^{\mathfrak{s}_i}(\underbrace{a_i, a_{i+2}}, a_{i+1})/M^{\mathfrak{s}_i}(\underbrace{a_{i+1}, a_{i+2}}, a_{i}) \rightarrow 0. $$
Applying $ LZ_{\mathfrak{s}_i}^{\mathfrak{t}_i} $ to get a long exact sequence we easily see that
$$ L_s Z_{\mathfrak{s}_i}^{\mathfrak{t}_i}(M^{\mathfrak{s}_i}(\underbrace{a_i, a_{i+2}}, a_{i+1})/M^{\mathfrak{s}_i}(\underbrace{a_{i+1}, a_{i+2}}, a_{i})) 
\cong M^{\mathfrak{t}_i}(\underbrace{a_{i}, a_{i+1}, a_{i+2}}) $$ 
if $ s = 1, 3 $ and 0 otherwise.

Now,
\begin{align*}
&\Ext^1(\epsilon_{\mathfrak{t}_i}^{\mathfrak{s}_i}M^{\mathfrak{t}_i}(\underbrace{a_{i}, a_{i+1}, a_{i+2}}),
M^{\mathfrak{s}_i}(\underbrace{a_i, a_{i+2}}, a_{i+1})/M^{\mathfrak{s}_i}(\underbrace{a_{i+1}, a_{i+2}}, a_{i})) \cong\\
&\Hom(\epsilon_{\mathfrak{t}_i}^{\mathfrak{s}_i}M^{\mathfrak{t}_i}(\underbrace{a_{i}, a_{i+1}, a_{i+2}}), M^{\mathfrak{s}_i}(\underbrace{a_i, a_{i+2}}, a_{i+1})/M^{\mathfrak{s}_i}(\underbrace{a_{i+1}, a_{i+2}}, a_{i})[1]) \cong\\
&\Hom(M^{\mathfrak{t}_i}(\underbrace{a_{i}, a_{i+1}, a_{i+2}})[0], M^{\mathfrak{t}_i}(\underbrace{a_{i}, a_{i+1}, a_{i+2}})[0] \oplus
M^{\mathfrak{t}_i}(\underbrace{a_{i}, a_{i+1}, a_{i+2}})[-2]) \cong \mathbb{C}.
\end{align*} 
\end{proof}

We would like to construct a natural transformation from $ LZ^{\mathfrak{s}_i} \epsilon_{\mathfrak{s}_{i+1}}[-1] LZ^{\mathfrak{s}_{i+1}} \epsilon_{\mathfrak{s}_i}[-1] $ to 
$ \Id \oplus \epsilon_{\mathfrak{t}_i}^{\mathfrak{s}_i}[-2] LZ_{\mathfrak{s}_i}^{\mathfrak{t}_i}. $

The identity map gives the adjunction map $ \alpha $ under the isomorphism
$ \Hom(LZ^{\mathfrak{s}_{i+1}} \epsilon_{\mathfrak{s}_i}, LZ^{\mathfrak{s}_{i+1}} \epsilon_{\mathfrak{s}_i}) \cong
\Hom(LZ^{\mathfrak{s}_i} \epsilon_{\mathfrak{s}_{i+1}}[-1] LZ^{\mathfrak{s}_{i+1}} \epsilon_{\mathfrak{s}_i}[-1], \Id). $

The identity map in $ \Hom(\epsilon_{\mathfrak{s}_i}, \epsilon_{\mathfrak{s}_i}) $ gives an adjunction map in $ \Hom(LZ^{\mathfrak{s}_i} \epsilon_{\mathfrak{s}_i}, \Id). $
Apply $ LZ_{\mathfrak{s}_i}^{\mathfrak{t}_i} $ to get a map in
\begin{align*}
\Hom(LZ_{\mathfrak{s}_i}^{\mathfrak{t}_i} LZ^{\mathfrak{s}_i} \epsilon_{\mathfrak{s}_i}, LZ_{\mathfrak{s}_i}^{\mathfrak{t}_i}) &\cong
\Hom(LZ^{\mathfrak{t}_i} \epsilon_{\mathfrak{s}_i}, LZ_{\mathfrak{s}_i}^{\mathfrak{t}_i})\\
&\cong \Hom(LZ_{\mathfrak{s}_{i+1}}^{\mathfrak{t}_i} LZ^{\mathfrak{s}_{i+1}} \epsilon_{\mathfrak{s}_i}, LZ_{\mathfrak{s}_i}^{\mathfrak{t}_i})\\
&\cong \Hom(LZ^{\mathfrak{s}_{i+1}} \epsilon_{\mathfrak{s}_i}, \epsilon_{\mathfrak{t}_i}^{\mathfrak{s}_{i+1}} LZ_{\mathfrak{s}_i}^{\mathfrak{t}_i}). 
\end{align*}

Now apply the functor $ \epsilon_{\mathfrak{s}_{i+1}} $ to get an element of
\begin{align*}
\Hom(\epsilon_{\mathfrak{s}_{i+1}} LZ^{\mathfrak{s}_{i+1}} \epsilon_{\mathfrak{s}_i}, \epsilon_{\mathfrak{s}_{i+1}} \epsilon_{\mathfrak{t}_i}^{\mathfrak{s}_{i+1}} LZ_{\mathfrak{s}_i}^{\mathfrak{t}_i}) &\cong
\Hom(\epsilon_{\mathfrak{s}_{i+1}} LZ^{\mathfrak{s}_{i+1}} \epsilon_{\mathfrak{s}_i}, \epsilon_{\mathfrak{t}_i}LZ_{\mathfrak{s}_i}^{\mathfrak{t}_i})\\
&\cong \Hom(\epsilon_{\mathfrak{s}_{i+1}} LZ^{\mathfrak{s}_{i+1}} \epsilon_{\mathfrak{s}_i}, \epsilon_{\mathfrak{s}_i} \epsilon_{\mathfrak{t}_i}^{\mathfrak{s}_i} LZ_{\mathfrak{s}_i}^{\mathfrak{t}_i})\\
&\cong \Hom(LZ^{\mathfrak{s}_i} \epsilon_{\mathfrak{s}_{i+1}} LZ^{\mathfrak{s}_{i+1}} \epsilon_{\mathfrak{s}_i}, \epsilon_{\mathfrak{t}_i}^{\mathfrak{s}_i} LZ_{\mathfrak{s}_i}^{\mathfrak{t}_i}). 
\end{align*} 
Now apply the shift functor $ [-2] $ to get a morphism $ \beta $ in
$$ \Hom(LZ^{\mathfrak{s}_i} \epsilon_{\mathfrak{s}_{i+1}}[-1] LZ^{\mathfrak{s}_{i+1}} \epsilon_{\mathfrak{s}_i}[-1], \epsilon_{\mathfrak{t}_i}^{\mathfrak{s}_i}[-2] LZ_{\mathfrak{s}_i}^{\mathfrak{t}_i}). $$

Therefore $ \alpha \oplus \beta \colon LZ^{\mathfrak{s}_i} \epsilon_{\mathfrak{s}_{i+1}}[-1] LZ^{\mathfrak{s}_{i+1}} \epsilon_{\mathfrak{s}_i}[-1] \rightarrow \Id \oplus \epsilon_{\mathfrak{t}_i}^{\mathfrak{s}_i}[-2] LZ_{\mathfrak{s}_i}^{\mathfrak{t}_i}. $

Now we come to the main result of the section.

\begin{prop}
\label{diagram5}
The map
$$ \alpha \oplus \beta \colon LZ^{\mathfrak{s}_i} \epsilon_{\mathfrak{s}_{i+1}}[-1] LZ^{\mathfrak{s}_{i+1}} \epsilon_{\mathfrak{s}_i}[-1] \rightarrow \Id \oplus \epsilon_{\mathfrak{t}_i}^{\mathfrak{s}_i}[-2] LZ_{\mathfrak{s}_i}^{\mathfrak{t}_i} $$
is an isomorphism.
\end{prop}

\begin{proof}
Lemmas 23-27 show that the cohomology functors of these morphisms are isomorphisms when applied to generalized Verma modules.
Thus the morphisms are isomorphisms when applied to projective modules by induction on the length of the generalized Verma flag.
Finally, they are isomorphisms on arbitrary modules by considering projective resolutions.
\end{proof}

The above theorem is the Koszul dual to lemma 4 of [BFK]. The correction term there vanishes when restricted to the appropriate parabolic subcategory just as in the case here where the correction terms vanishes when restricted to the appropriate singular block.  This result is not a surprise in light of [Rh] or [MOS] where it was shown that projective functors and Zuckerman functors are Koszul dual.  We end this section with an immediate corollary which is the functorial version of diagram 5.

\begin{corollary}
$ \epsilon_{\mathfrak{s}_i}[-1] LZ^{\mathfrak{s}_i} \epsilon_{\mathfrak{s}_{i+1}}[-1] LZ^{\mathfrak{s}_{i+1}} \epsilon_{\mathfrak{s}_i}[-1] LZ^{\mathfrak{s}_i} \oplus \epsilon_{\mathfrak{s}_{i+1}}[-1] LZ^{\mathfrak{s}_{i+1}} \cong $
$$ \epsilon_{\mathfrak{s}_{i+1}}[-1] LZ^{\mathfrak{s}_{i+1}} \epsilon_{\mathfrak{s}_i}[-1] LZ^{\mathfrak{s}_i} \epsilon_{\mathfrak{s}_{i+1}}[-1] LZ^{\mathfrak{s}_{i+1}} \oplus \epsilon_{\mathfrak{s}_i}[-1] LZ^{\mathfrak{s}_i}. $$
\end{corollary}

\section{Graded Category $ \mathcal{O} $}
In the previous sections, various aspects of the representation theory of $ \mathfrak{sl}_k $ were categorified.  Topological invariants arise from the representation theory of $ \mathcal{U}_{q}(\mathfrak{sl}_k). $  The categorification of quantum groups is accomplished through graded representation theory.  We will treat category $ \mathcal{O} $ as a category of graded modules.  Then a shift in this grading descends to multiplication by $ q $ in the Grothendieck group.  The idea of graded category $ \mathcal{O} $ originates from [Soe2].  It was continued in [Str1] where among other things, it was shown how to construct graded lifts of translation functors.  

\begin{define}
Let $ P_{\bf d} = \oplus_{x \in S_n/S_{\bf d}} P(x. \omega_{\bf d}) $ where $ {\bf d} = (d_{k-1}, \ldots, d_0), $ $ x $ is a minimal coset representative, and
$ \omega_{\bf d} = (\underbrace{k-1, \ldots, k-1}_{d_{k-1}}, \ldots, \underbrace{0, \ldots, 0}_{d_{0}}). $
\end{define}

Then $  P_{\bf d} $ is a minimal projective generator of $ \mathcal{O}_{\bf d}(\mathfrak{gl}_n). $
There is an equivalence of categories $ \mathcal{O}_{\bf d}(\mathfrak{gl}_n) \cong \text{mod}-A_{\bf d} $ where
$ A_{\bf d} = \End_g(P_{\bf d}) $ and $ \text{mod}-A_{\bf d} $ is the category of finitely generated right $ A_{\bf d}- $ modules.  We will interpret $ A_{\bf d} $ as a graded algebra.

\begin{lemma}
\label{lemma28}
Let $ R $ and $ S $ be any rings.  There is an equivalence of categories:
$$ \lbrace \text{right exact functors compatible with direct sums}: (\text{mod-}R \rightarrow \text{mod-}S) \rbrace \rightarrow R\text{-mod-}S. $$
Under this equivalence a functor $ F $ gets mapped to $ F(R). $ In the other direction,
a bimodule $ X $ gets mapped to $ \bullet \otimes_{R} X. $
\end{lemma}

\begin{proof}
This is in [Bass] as well as lemma 3.4 of [Str1].
\end{proof}

Let $ \alpha $ and $ \beta $ be parabolic subalgebras of $ \mathfrak{gl}_n $ such that $ \beta \subset \alpha. $
Let $ P(x) \in \mathcal{O}_{\bf d}(\mathfrak{gl}_n) $ be an indecomposable projective object.  Then $ Z^{\alpha} P(x) $ is either 0 or an indecomposable projective object in
$ \mathcal{O}_{\bf d}^{\alpha} (\mathfrak{gl}_n) $. 

Denote by $ P_{\bf d}^{\alpha} $ a minimal projective generator of $ \mathcal{O}_{\bf d}^{\alpha} (\mathfrak{gl}_n) $. 
Let $ A_{\bf d}^{\alpha} = \End(P_{\bf d}^{\alpha}) $ be its endomorphism algebra.
Then there is an equivalence of categories
$ \mathcal{O}_{\bf d}^{\alpha} (\mathfrak{sl}_n)  \cong \text{mod}-A_{\bf d}^{\alpha} $.
Let $ S = S(\mathfrak{h}) $ be the symmetric algebra associated to the Cartan subalgebra.  Let $ C = S/S_{+}^{W} $ be the associated coinvariant algebra.
Let $ C^{\lambda} $ be the subalgebra invariant under $ W_{\lambda}. $  Let $ w_0^{\lambda} $ be the longest element in the set of shortest coset representatives of $ W/W_{\lambda}. $  There is the is a well known isomorphism due to Soergel [Soe] between the endomorphism algebra of the indecomposable projective-injective module and this subalgebra of invariants: $ \End(P(w_{0}^{\lambda}.\lambda)) \cong C^{\lambda}. $

\begin{define}
Let $ \mathbb{V}_{\lambda} \colon \mathcal{O}_{\lambda} \rightarrow \text{mod}-C^{\lambda} $ be the Soergel functor defined by
$ M \mapsto \Hom_{\mathfrak{g}}(P(w_{0}^{\lambda}.\lambda), M). $
\end{define}

Soergel showed that the functor $ \mathbb{V}_{\lambda} $ is fully faithful on projective objects [Soe].  
Let I and J be compositions of n.  
If $ I = i_1 + \cdots+ i_r = n, $ associate to I the Young subgroup of $ S_n, $ $ S_{i_1} \times \cdots \times S_{i_r}. $
Let $ \mu_I $ and $ \mu_J $ be integral dominant weights stabalized by the subgroups associated to I and J.
Suppose $ J \subset I, $ (there is a containment of the associated subgroups.)
Then we can define the translation functor $ \theta_{I}^{J} $ from $ \mathcal{O}_{\mu_I}(\mathfrak{gl}_n) $ to $ \mathcal{O}_{\mu_J}(\mathfrak{gl}_n). $  It is the projective functor given be tensoring with the finite dimensional, irreducible module with highest weight $ \mu_J - \mu_I $ and then projecting onto the block $ \mathcal{O}_{\mu_J}(\mathfrak{gl}_n). $ 
Translation functors and Soergel functors are related by the following lemma.

\begin{lemma}
\label{lemma29}
Let $ \text{Res}_{J}^{I} \colon \text{mod}-C^J \rightarrow \text{mod}-C^I $ denote the restriction functor.  Then
\begin{enumerate}
\item $ {\mathbb{V}}_J \theta_I^J \cong C^J \otimes_{C^I} {\mathbb{V}}_I. $
\item $ {\mathbb{V}}_I \theta_J^I \cong Res_J^I {\mathbb{V}}_J. $
\end{enumerate}
\end{lemma}

\begin{proof}
This is proposition 3.3 of [FKS].
\end{proof}

This lemma together with the fact that $ \mathbb{V}_{\lambda} $ is a faithful functor on projective objects allows us to consider the endomorphism ring of a minimal projective generator of $ \mathcal{O} $ as a graded ring [Str1].  A projective object is a direct summand of a sequence of projective functors applied to a dominant Verma module.  Thus by the previous lemma, a projective object $ P, $ $ \mathbb{V}_{\lambda} P $ becomes a graded $ C^{\lambda}- $ module.  Then $ \End(\mathbb{V}_{\lambda} P) $ becomes a graded ring so there is a grading on $ A_{\bf d}. $

In [Str] it was actually shown that it is a graded quotient of $ A_{\bf d}. $
Now we can consider $ \mathcal{O}_{\bf d}(\mathfrak{gl}_n) $ and  $ \mathcal{O}_{\bf d}^{\alpha} (\mathfrak{gl}_n) $ as graded categories by considering
$ \text{gmod}-A_{\bf d} $ and $ \text{gmod}-A_{\bf d}^{\alpha}  $ respectively.

By lemma ~\ref{lemma28}, the dual Zuckerman functor becomes
$$ Z_{\beta}^{\alpha} \colon \text{mod}-A_{{\bf d}}^{\beta} \rightarrow \text{mod}-A_{{\bf d}}^{\alpha} $$ 
given by
$$ \bullet \otimes_{A_{{\bf d}}^{\beta}} A_{{\bf d}}^{\alpha}. $$

\begin{define}
The graded functor $ \widetilde{Z_{\beta}^{\alpha}} \colon \text{gmod}-A_{{\bf d}}^{\beta} \rightarrow \text{gmod}-A_{{\bf d}}^{\alpha} $ is given by
$ \bullet \otimes_{A_{{\bf d}}^{\beta}} A_{{\bf d}}^{\alpha}\langle -(\text{dim}(\alpha) - \text{dim}(\beta))/2 \rangle. $
\end{define}

We may now also consider
$$ L\widetilde{Z_{\beta}^{\alpha}} \colon D^b(\text{gmod}-A_{{\bf d}}^{\beta}) \rightarrow D^b(\text{gmod}-A_{{\bf d}}^{\alpha}). $$ 

We fix  a graded lift of a generalized Verma module $ \widetilde{M^{\mathfrak{p}}}(x.\lambda) $ so that its head is concentrated in degree zero.
The proof of proposition 5.2 of [FKS] gives a graded inclusion $ \widetilde{M}(x.\lambda)\langle l(x) \rangle \rightarrow \widetilde{M}(\lambda) $ of Verma modules.  This implies that the standard maps in a generalized BGG resolution are homogeneous of degree 1.
A Koszul module is an object that admits a projective resolution such that the object in degree $ i $ is generated by its degree $ i $ subspace.
We will need the following result later.

\begin{prop}
Generalized Verma modules are Koszul.
\end{prop}

\begin{proof}
This is theorem 3.11.4 of [BGS].
\end{proof}

\begin{prop}
Let $ \beta \subset \alpha $ be two parabolic subalgebras.  Assume $ \text{dim}(\alpha) - \text{dim}(\beta) = 2d. $  Then
\begin{enumerate}
\item The right adjoint of $ \widetilde{\epsilon_{\alpha}^{\beta}} $ is $ L\widetilde{Z_{\beta}^{\alpha}}[-2d] \langle -d \rangle. $ 
\item The right adjoint of $ L\widetilde{Z_{\beta}^{\alpha}} $ is $ \widetilde{\epsilon_{\alpha}^{\beta}} \langle d \rangle. $ 
\end{enumerate}
\end{prop}

\begin{proof}
\begin{enumerate}
\item In the ungraded case the right adjoint of $ \epsilon_{\alpha}^{\beta} $ is $ LZ_{\beta}^{\alpha}[-2d]. $
Therefore if the right adjoint of $ \widetilde{\epsilon}_{\beta}^{\alpha} $ exists, it must be a graded lift of $ LZ_{\beta}^{\alpha}[-2d]. $
Note that $ \widetilde{\epsilon}_{\alpha}^{\beta} M $ becomes a right $ A^{\beta}- $ module via the quotient map
$ A^{\beta} \rightarrow A^{\alpha}. $
Therefore
$$ \Hom_{\text{gmod}-A^{\beta}}(\widetilde{\epsilon_{\alpha}^{\beta}}M, N) \cong \Hom_{\text{gmod}-A^{\alpha}}(M, \Hom_{\text{gmod}-A^{\beta}}(A^{\alpha}, N)). $$
Now $ \Hom_{\text{gmod}-A^{\beta}}(A^{\alpha}, \bullet) $ is clearly a functor from $ \text{gmod}-A^{\beta} $ to $ \text{gmod}-A^{\alpha}. $  
The right adjoint is a direct sum of graded functors $ \oplus_i F_i $ which in the ungraded category is simply $ L_{2d} Z_{\beta}^{\alpha}. $  Let $ \bigtriangledown(x) $ denote a dual Verma module.  We have an isomorphism:
$$ \Hom_{\text{gmod}-A^{\beta}}(\widetilde{\epsilon}_{\alpha}^{\beta} {\widetilde{\bigtriangledown}^{\alpha}(x)}, {\widetilde{\bigtriangledown}^{\beta}(x)}) \cong
\Hom_{\text{gmod}-A^{\alpha}}({\widetilde{\bigtriangledown}^{\alpha}(x)}, \oplus_i F_i {\widetilde{\bigtriangledown}^{\beta}(x)}). $$
If $ F_i {\widetilde{\bigtriangledown}^{\beta}(x)} \neq 0, $ then 
$ F_i {\widetilde{\bigtriangledown}^{\beta}(x)} \cong L_{2d} \widetilde{Z}_{\beta}^{\alpha} {\widetilde{\bigtriangledown}^{\beta}(x)} \langle k_i \rangle. $
By the definition of $ \widetilde{Z}_{\beta}^{\alpha} $ and the fact that generalized Verma modules are Koszul, we get that the above is isomorphic to
$ {\widetilde{\bigtriangledown}^{\alpha}(x)} \langle d+k_i \rangle. $
Since the space of morphisms is one dimensional, $ k_i = -d. $
\item There is a degree zero homogenous map of $ A^{\beta}- $ modules $ \psi \colon M \rightarrow M \otimes_{A^{\beta}} A^{\alpha}  $ where
$ \psi(m) = m \otimes 1. $
Now suppose $ \phi \colon M \otimes_{A^{\beta}} A^{\alpha} \rightarrow N $ is a degree zero homogeneous $ A^{\beta}- $ module morphism. Then $ \phi \circ \psi \colon M \rightarrow N $ is a degree zero homogeneous $ A^{\beta}- $ module morphism.  Therefore we have a map
$ f \colon \Hom_{\text{gmod}-A^{\alpha}}(M \otimes_{A^{\beta}} A^{\alpha}, N) \rightarrow \Hom_{\text{gmod}-A^{\beta}}(M, N). $
It is clear that f is injective.  We must verify that it is an isomorphism.  Consider a two step free resolution of $ M, $ $ M_2 \rightarrow M_1 \rightarrow M \rightarrow 0. $
Then we get the following commutative diagram with exact rows:

\begin{tiny}
\xymatrix
{
&\Hom_{\text{gmod}-A^{\alpha}}(M_2 \otimes_{A^{\beta}} A^{\alpha}, N)\ar[d]	&\Hom_{\text{gmod}-A^{\alpha}}(M_1 \otimes_{A^{\beta}} A^{\alpha}, N)\ar[l]\ar[d]
&\Hom_{\text{gmod}-A^{\alpha}}(M \otimes_{A^{\beta}} A^{\alpha}, N)\ar[l]\ar[d]	& 0\ar[l]\\
&\Hom_{\text{gmod}-A^{\beta}}(M_2, N)	&\Hom_{\text{gmod}-A^{\beta}}(M_1, N)\ar[l]		&\Hom_{\text{gmod}-A^{\beta}}(M, N)\ar[l]	&0\ar[l].
}
\end{tiny}
Therefore it suffices to show that $ f $ is surjective for free modules $ M. $  It is clear for free modules that both spaces of morphisms have the same dimension so it must be a surjection.
\end{enumerate}
\end{proof}

\subsection{Graded Projective Functors and Quantum Serre Relations}
Define the k-tuple $ {\bf d} + t_0\epsilon_0 + \cdots + t_{k-1}\epsilon_{k-1} $ to be $ (d_{k-1}+t_{k-1}, \ldots, d_{0}+t_{0}) $ where the $ t_i $ are integers.
Assume $ t>0. $ 
Define $ ({\bf d}; {\bf d} + t \epsilon_i- t \epsilon_{i-1}) $ to be the (k+1)-tuple $ (d_{k-1}, \ldots, d_i, t, d_{i-1}-t, d_{i-2}, \ldots, d_0). $
Define $ ({\bf d}; {\bf d} - t \epsilon_i+ t \epsilon_{i-1}) $ to be the (k+1)-tuple $ (d_{k-1}, \ldots, d_i-t, t, d_{i-1}, d_{i-2}, \ldots, d_0). $

\begin{lemma}
\label{lemma30}
$ \mathcal{E}_{i}^{(t)} \colon \mathcal{O}_{\bf d} \rightarrow \mathcal{O}_{{\bf d} +t\epsilon_i - t\epsilon_{i-1}} $ is isomorphic to 
$ \theta_{{\bf d}; {\bf d} + t\epsilon_i- t\epsilon_{i-1}}^{{\bf d} + t\epsilon_i- t\epsilon_{i-1}} \theta^{{\bf d}; {\bf d} + t\epsilon_i- t\epsilon_{i-1}}_{\bf d}. $
\end{lemma}

\begin{proof}
See proposition 3.2b of [FKS].
\end{proof}
 
\begin{lemma}
\label{lemma31}
$ \mathcal{F}_{i}^{(t)} \colon \mathcal{O}_{\bf d} \rightarrow \mathcal{O}_{{\bf d} -t\epsilon_i + t\epsilon_{i-1}} $ is isomorphic to 
$ \theta_{{\bf d}; {\bf d} - t\epsilon_i+ t\epsilon_{i-1}}^{{\bf d} - t\epsilon_i+ t\epsilon_{i-1}} \theta^{{\bf d}; {\bf d} - t\epsilon_i+ t\epsilon_{i-1}}_{\bf d}. $
\end{lemma}

\begin{proof}
Same as the previous lemma.
\end{proof}

Now we are prepared to introduce the graded lifts $ \widetilde{\mathcal{E}}_i^{(t)}, $ $ \widetilde{\mathcal{F}}_i^{(t)}, $ $ \widetilde{\mathcal{H}}, $ and $ \widetilde{\mathcal{H}}^{-1}. $
By lemma ~\ref{lemma28}, 
$ \mathcal{E}_i^{(t)} \colon \text{mod}-A_{\bf d} \rightarrow \text{mod}-A_{{\bf d}+t\epsilon_i-t\epsilon_{i-1}} $ is given by
$$ \bullet \otimes_{\text{mod}-A_{\bf d}} \Hom_g(P_{{\bf d} + t\epsilon_i- t\epsilon_{i-1}}, \mathcal{E}_i^{(t)}P_{\bf d}). $$
By lemma ~\ref{lemma30} this is
$$ \bullet \otimes_{\text{mod}-A_{\bf d}} \Hom_g(P_{{\bf d} + t\epsilon_i- t\epsilon_{i-1}}, 
\theta_{{\bf d}; {\bf d} + t\epsilon_i- t\epsilon_{i-1}}^{{\bf d} + t\epsilon_i- t\epsilon_{i-1}} \theta^{{\bf d}; {\bf d} + t\epsilon_i- t\epsilon_{i-1}}_{\bf d}P_{\bf d}). $$
Then by lemma ~\ref{lemma29} this is isomorphic to
$$ \bullet \otimes_{\text{gmod}-A_{\bf d}} \Hom_{C^{{\bf d} + t\epsilon_i- t\epsilon_{i-1}}}({\mathbb{V}}_{{\bf d} + t\epsilon_i- t\epsilon_{i-1}}P_{{\bf d} + t\epsilon_i- t\epsilon_{i-1}},
\text{Res}^{{\bf d} + t\epsilon_i- t\epsilon_{i-1}}_{{\bf d}; {{\bf d} + t\epsilon_i- t\epsilon_{i-1}}} C^{{\bf d}; {{\bf d} + t\epsilon_i- t\epsilon_{i-1}}} \otimes_{C^{\bf d}} 
{\mathbb{V}}_{\bf d}P_{\bf d}). $$

\begin{define}
\begin{enumerate}
\item Let $ r_{i,t}=\sum_{r=1}^t d_{i-1}-r. $ Then we let
$ \widetilde{\mathcal{E}}_i^{(t)} = $
$$ \bullet \otimes_{\text{gmod}-A_{\bf d}} \Hom_{C^{{\bf d} + t\epsilon_i- t\epsilon_{i-1}}}({\mathbb{V}}_{{\bf d} + t\epsilon_i- t\epsilon_{i-1}}P_{{\bf d} + t\epsilon_i- t\epsilon_{i-1}},
\text{Res}^{{\bf d} + t\epsilon_i- t\epsilon_{i-1}}_{{\bf d}; {{\bf d} + t\epsilon_i- t\epsilon_{i-1}}} C^{{\bf d}; {{\bf d} + t\epsilon_i- t\epsilon_{i-1}}} \otimes_{C^{\bf d}} 
{\mathbb{V}}_{\bf d}P_{\bf d}\langle -r_{i,t}+(t-1) \rangle). $$

\item Let $ s_{i,t}=\sum_{s=1}^{t} d_i-s. $ Then we let
$ \widetilde{\mathcal{F}}_i^{(t)} = $
$$ \bullet \otimes_{\text{gmod}-A_{\bf d}} \Hom_{C^{{\bf d} - t\epsilon_i+ t\epsilon_{i-1}}}({\mathbb{V}}_{{\bf d} - t\epsilon_i+ t\epsilon_{i-1}}P_{{\bf d} - t\epsilon_i+ t\epsilon_{i-1}},
\text{Res}^{{\bf d} - t\epsilon_i+ t\epsilon_{i-1}}_{{\bf d}; {{\bf d} - t\epsilon_i+ t\epsilon_{i-1}}} C^{{\bf d}; {{\bf d} - t\epsilon_i+ t\epsilon_{i-1}}} \otimes_{C^{\bf d}} 
{\mathbb{V}}_{\bf d}P_{\bf d}\langle -s_{i,t}+(t-1) \rangle). $$

\item $ \widetilde{\mathcal{H}} = \Id\langle d_i-d_{i-1} \rangle. $

\item $ \widetilde{\mathcal{H}}^{-1} = \Id\langle -(d_i-d_{i-1}) \rangle. $
\end{enumerate}
\end{define}

\begin{lemma}
\label{lemma32}
$ \widetilde{\mathcal{E}_i} \widetilde{\mathcal{E}_i} \cong \widetilde{\mathcal{E}}_i^{(2)}\langle 1 \rangle \oplus \widetilde{\mathcal{E}}_i^{(2)}\langle -1 \rangle. $
\end{lemma}

\begin{proof}
First we consider the ungraded case and compute in the Grothendieck group.
\begin{align*} 
&[\mathcal{E}_i][\mathcal{E}_i][M(k-1, \ldots, k-1, \ldots, i+1, \ldots, i+1, i, \ldots, i, i-1, \ldots, i-1, \ldots, 0, \ldots, 0)] =\\
&2[P(k-1, \ldots, k-1, \ldots, i, \ldots, i, \underbrace{i-1, \ldots, i-1, i, i,}_{d_{i-1}} \ldots, 0, \ldots, 0)]. 
\end{align*} 
The indecomposable projective functor that takes the dominant Verma module to this projective object is $ \mathcal{E}_i^{(2)}. $
Thus there is an isomorphism of functors $ \mathcal{E}_i \mathcal{E}_i \cong \mathcal{E}_i^{(2)} \oplus \mathcal{E}_i^{(2)}. $
So there must be an isomorphism of graded functors 
$ \widetilde{\mathcal{E}}_i \widetilde{\mathcal{E}}_i \cong \widetilde{\mathcal{E}}_i^{(2)}\langle r \rangle \oplus \widetilde{\mathcal{E}}_i^{(2)}\langle s \rangle, $
for some integers r and s.

$ \widetilde{\mathcal{E}}_i^{(2)} $ is given by tensoring with the graded bimodule
$$ \Hom_{C^{{\bf d} + 2\epsilon_i- 2\epsilon_{i-1}}}({\mathbb{V}}_{{\bf d} + 2\epsilon_i- 2\epsilon_{i-1}}P_{{\bf d} + 2\epsilon_i- 2\epsilon_{i-1}},
C^{{\bf d}; {{\bf d} + 2\epsilon_i- 2\epsilon_{i-1}}} \otimes_{C^{\bf d}} 
{\mathbb{V}}_{\bf d}P_{\bf d}\langle -2d_{i-1}+4 \rangle). $$
As a $ C^{{\bf d} + 2\epsilon_i- 2\epsilon_{i-1}}- $ module, $ C^{{\bf d}; {{\bf d} + 2\epsilon_i- 2\epsilon_{i-1}}} $ is free of rank
$ |W^{{\bf d} + 2\epsilon_i- 2\epsilon_{i-1}} / W^{{\bf d}; {{\bf d} + 2\epsilon_i- 2\epsilon_{i-1}}}|. $
A basis can be chosen homogeneous in the degrees twice the length of minimal coset representatives of
\begin{align*}
&|W^{{\bf d} + 2\epsilon_i- 2\epsilon_{i-1}} / W^{{\bf d}; {{\bf d} + 2\epsilon_i- 2\epsilon_{i-1}}}| =\\
&{S_{d_{k-1}} \times \cdots \times S_{d_{i+1}} \times S_{d_i+2} \times S_{d_{i-1}-2} \times \cdots \times S_{d_0}}/
{S_{d_{k-1}} \times \cdots \times S_{d_{i}} \times S_2 \times S_{d_{i-1}-2} \times \cdots \times S_{d_0}}. 
\end{align*} 

So now we want to consider minimal coset representatives in $ S_{d_i+2}/(S_{d_i} \times S_2). $
Proposition A-2 in [Str] tells us the the length of these coset representatives.
It gives 
$$ C^{{\bf d}; {\bf d}+2\epsilon_i - 2\epsilon_{i-1}} \cong \oplus_{d_i+1 \geq r > s \geq 0} C^{{\bf d}+2\epsilon_i -2\epsilon_{i-1}}\langle 2r+2s-2 \rangle. $$

Now $ \widetilde{\mathcal{E}}_i \widetilde{\mathcal{E}}_i $ is given by tensoring with the bimodule
$$ \aligned
&\Hom_{C^{{\bf d} + \epsilon_i- \epsilon_{i-1}}}({\mathbb{V}}_{{\bf d} + \epsilon_i- \epsilon_{i-1}}P_{{\bf d} + \epsilon_i- \epsilon_{i-1}},
C^{{\bf d}; {{\bf d} + \epsilon_i- \epsilon_{i-1}}} \otimes_{C^{\bf d}} 
{\mathbb{V}}_{\bf d}P_{\bf d}\langle -d_{i-1}+1\rangle) \otimes_{C^{{\bf d}+\epsilon_i-\epsilon_{i-1}}}\\
&\Hom_{C^{{\bf d} + 2\epsilon_i- 2\epsilon_{i-1}}}({\mathbb{V}}_{{\bf d} + 2\epsilon_i- 2\epsilon_{i-1}}P_{{\bf d} + 2\epsilon_i- 2\epsilon_{i-1}},
C^{{\bf d+\epsilon_i -\epsilon_{i-1}}; {{\bf d} + 2\epsilon_i- 2\epsilon_{i-1}}} \otimes_{C^{{\bf d}+\epsilon_i-\epsilon_{i-1}}}
{\mathbb{V}}P_{{\bf d}+\epsilon_i-\epsilon_{i-1}}\langle -d_{i-1}+2 \rangle). \endaligned $$

Thus we need to compute 
$ C^{{\bf d+\epsilon_i -\epsilon_{i-1}}; {{\bf d} + 2\epsilon_i- 2\epsilon_{i-1}}} \otimes_{C^{{\bf d}+\epsilon_i-\epsilon_{i-1}}} C^{{\bf d}; {\bf d}+\epsilon_i-\epsilon_{i-1}} $
as a $ C^{{\bf d} + 2\epsilon_i- 2\epsilon_{i-1}}-$ module.
As above, $ C^{{\bf d}; {{\bf d}+\epsilon_i-\epsilon_{i-1}}} $ is a free  $ C^{{\bf d}+\epsilon_i-\epsilon_{i-1}}- $ module of rank
$ |W^{{\bf d} + \epsilon_i- \epsilon_{i-1}} / W^{{\bf d}; {{\bf d} + \epsilon_i- \epsilon_{i-1}}}| $ and a basis can be chosen homogeneous in degrees twice the length
of minimal coset representatives.  There is a similar description for
$ C^{{\bf d}+\epsilon_i -\epsilon_{i-1}; {{\bf d}+2\epsilon_i-2\epsilon_{i-1}}} $ as a $ C^{{\bf d}+2\epsilon_i-2\epsilon_{i-1}}- $ module.
Thus 
$$ C^{{\bf d}+\epsilon_i -\epsilon_{i-1}; {{\bf d} + 2\epsilon_i- 2\epsilon_{i-1}}} \otimes_{C^{{\bf d}+\epsilon_i-\epsilon_{i-1}}} C^{{\bf d}; {\bf d}+\epsilon_i-\epsilon_{i-1}} \cong
\oplus_{s=0}^{d_i+1} \oplus_{r=0}^{d_i} C^{{\bf d}+2\epsilon_i-2\epsilon_{i-1}}\langle 2r+2s \rangle. $$
Therefore the grading of $ \widetilde{\mathcal{E}}_i \widetilde{\mathcal{E}}_i $ is encoded by
$$ C^{{\bf d}+\epsilon_i -\epsilon_{i-1}; {{\bf d} + 2\epsilon_i- 2\epsilon_{i-1}}} \otimes_{C^{{\bf d}+\epsilon_i-\epsilon_{i-1}}} C^{{\bf d}; {\bf d}+\epsilon_i-\epsilon_{i-1}} \cong
\oplus_{s=0}^{d_i+1} \oplus_{r=0}^{d_i} C^{{\bf d}+2\epsilon_i-2\epsilon_{i-1}}\langle 2r+2s-2d_{i-3}+3 \rangle. $$
The grading of $ \widetilde{\mathcal{E}}_i^{(2)} $ is encoded by
$$ C^{{\bf d}; {\bf d}+2\epsilon_i - 2\epsilon_{i-1}} \cong \oplus_{d_i+1 \geq r > s \geq 0} C^{{\bf d}+2\epsilon_i -2\epsilon_{i-1}}\langle 2r+2s-2-2d_{i-1}+4 \rangle. $$
By examining the highest and lowest degrees,
$ \widetilde{\mathcal{E}}_i \widetilde{\mathcal{E}}_i \cong \widetilde{\mathcal{E}}_i^{(2)}\langle 1 \rangle \oplus \widetilde{\mathcal{E}}_i^{(2)}\langle -1 \rangle. $
\end{proof}

\begin{lemma}
\label{lemma33}
$ \widetilde{\mathcal{F}}_i \widetilde{\mathcal{F}}_i \cong \widetilde{\mathcal{F}}_i^{(2)}\langle 1 \rangle \oplus \widetilde{\mathcal{F}}_i^{(2)}\langle -1 \rangle. $
\end{lemma}

\begin{proof}
First we consider the ungraded case and compute in the Grothendieck group.
\begin{align*}
&[\mathcal{F}_i][\mathcal{F}_i][M(k-1, \ldots, k-1, \ldots i+1, \ldots, i+1, i, \ldots, i, i-1, \ldots, i-1, \ldots, 0, \ldots, 0)] =\\
&2[P(k-1, \ldots, k-1, \ldots, i+1, \ldots, i+1, \underbrace{i-1, i-1, i \ldots, i, i,}_{d_{i}} \ldots, 0, \ldots, 0)]. 
\end{align*} 
The indecomposable projective functor that takes the dominant Verma module to this projective object is $ \mathcal{F}_i^{(2)}. $
Thus there is an isomorphism of functors $ \mathcal{F}_i \mathcal{F}_i \cong \mathcal{F}_i^{(2)} \oplus \mathcal{F}_i^{(2)}. $
So there must be an isomorphism of graded functors 
$ \widetilde{\mathcal{F}}_i \widetilde{\mathcal{F}}_i \cong \widetilde{\mathcal{F}}_i^{(2)}\langle r \rangle \oplus \widetilde{\mathcal{F}}_i^{(2)}\langle s \rangle, $
for some integers r and s.

$ \widetilde{\mathcal{F}}_i^{(2)} $ is given by tensoring with the graded bimodule
$$ \Hom_{C^{{\bf d} - 2\epsilon_i+ 2\epsilon_{i-1}}}({\mathbb{V}}_{{\bf d} - 2\epsilon_i+ 2\epsilon_{i-1}}P_{{\bf d} - 2\epsilon_i+ 2\epsilon_{i-1}},
\text{Res}^{{\bf d} - 2\epsilon_i+ 2\epsilon_{i-1}}_{{\bf d}; {{\bf d} - 2\epsilon_i+ 2\epsilon_{i-1}}} C^{{\bf d}; {{\bf d} - 2\epsilon_i+ 2\epsilon_{i-1}}} \otimes_{C^{\bf d}} 
{\mathbb{V}}_{\bf d}P_{\bf d}\langle -2d_i+4 \rangle). $$
As a $ C^{{\bf d} - 2\epsilon_i+ 2\epsilon_{i-1}}- $ module, $ C^{{\bf d}; {{\bf d} - 2\epsilon_i+ 2\epsilon_{i-1}}} $ is free of rank
$ |W^{{\bf d} - 2\epsilon_i+ 2\epsilon_{i-1}} / W^{{\bf d}; {{\bf d} - 2\epsilon_i+ 2\epsilon_{i-1}}}|. $
A basis can be chosen homogenous in the degrees twice the length of minimal coset representatives of
\begin{align*}
&|W^{{\bf d} - 2\epsilon_i+ 2\epsilon_{i-1}} / W^{{\bf d}; {{\bf d} - 2\epsilon_i+ 2\epsilon_{i-1}}}| =\\
&{S_{d_{k-1}} \times \cdots \times S_{d_{i+1}} \times S_{d_i-2} \times S_{d_{i-1}+2} \times \cdots \times S_{d_0}}/
{S_{d_{k-1}} \times \cdots \times S_{d_{i}-2} \times S_2 \times S_{d_{i-1}} \times \cdots \times S_{d_0}}. 
\end{align*} 

So now we want to consider minimal coset representatives in $ S_{d_{i-1}+2}/(S_2 \times S_{d_{i-1}}). $
Proposition A-2 in [Str] tells us the the length of these coset representatives.
It gives 
$$ C^{{\bf d}; {\bf d}-2\epsilon_i + 2\epsilon_{i-1}} \cong \oplus_{d_{i-1}+1 \geq r > s \geq 0} C^{{\bf d}-2\epsilon_i +2\epsilon_{i-1}}\langle 2r+2s-2 \rangle. $$

Now $ \widetilde{\mathcal{F}}_i \widetilde{\mathcal{F}}_i $ is given by tensoring with the bimodule
$$ \aligned
&\Hom_{C^{{\bf d} - \epsilon_i+ \epsilon_{i-1}}}({\mathbb{V}}_{{\bf d} - \epsilon_i+ \epsilon_{i-1}}P_{{\bf d} - \epsilon_i  +\epsilon_{i-1}},
C^{{\bf d}; {{\bf d} - \epsilon_i+ \epsilon_{i-1}}} \otimes_{C^{\bf d}} 
{\mathbb{V}}_{\bf d}P_{\bf d}\langle -d_i+1\rangle) \otimes_{C^{{\bf d}-\epsilon_i+\epsilon_{i-1}}}\\
&\Hom_{C^{{\bf d} - 2\epsilon_i+ 2\epsilon_{i-1}}}({\mathbb{V}}_{{\bf d} - 2\epsilon_i+ 2\epsilon_{i-1}}P_{{\bf d} - 2\epsilon_i+ 2\epsilon_{i-1}},
C^{{\bf d-\epsilon_i +\epsilon_{i-1}}; {{\bf d} - 2\epsilon_i+ 2\epsilon_{i-1}}} \otimes_{C^{{\bf d}-\epsilon_i +\epsilon_{i-1}}}
{\mathbb{V}}P_{{\bf d}-\epsilon_i +\epsilon_{i-1}}\langle -d_i+2\rangle). \endaligned $$

Thus we need to compute 
$ C^{{{\bf d}-\epsilon_i +\epsilon_{i-1}}; {{\bf d} - 2\epsilon_i+ 2\epsilon_{i-1}}} \otimes_{C^{{\bf d}-\epsilon_i +\epsilon_{i-1}}} C^{{\bf d}; {\bf d}-\epsilon_i+\epsilon_{i-1}} $
as a $ C^{{\bf d} - 2\epsilon_i+ 2\epsilon_{i-1}}-$ module.
As above, $ C^{{\bf d}; {{\bf d}-\epsilon_i+\epsilon_{i-1}}} $ is a free  $ C^{{\bf d}-\epsilon_i+\epsilon_{i-1}}- $ module of rank
$ |W^{{\bf d} - \epsilon_i + \epsilon_{i-1}} / W^{{\bf d}; {{\bf d} - \epsilon_i + \epsilon_{i-1}}}| $ and a basis can be chosen homogeneous in degrees twice the length
of minimal coset representatives.  There is a similar description for
$ C^{{\bf d}-\epsilon_i +\epsilon_{i-1}; {{\bf d}-2\epsilon_i +2\epsilon_{i-1}}} $ as a $ C^{{\bf d}-2\epsilon_i+2\epsilon_{i-1}}- $ module.

Thus 
$$ C^{{{\bf d}-\epsilon_i +\epsilon_{i-1}}; {{\bf d} - 2\epsilon_i + 2\epsilon_{i-1}}} \otimes_{C^{{\bf d}-\epsilon_i +\epsilon_{i-1}}} C^{{\bf d}; {\bf d}-\epsilon_i+\epsilon_{i-1}} \cong
\oplus_{s=0}^{d_{i-1}} \oplus_{r=0}^{d_{i-1}+1} C^{{\bf d}-2\epsilon_i+2\epsilon_{i-1}}\langle 2r+2s \rangle. $$
Therefore the grading of $ \widetilde{\mathcal{F}}_i \widetilde{\mathcal{F}}_i $ is encoded by
$$ C^{{{\bf d}-\epsilon_i +\epsilon_{i-1}}; {{\bf d} - 2\epsilon_i+ 2\epsilon_{i-1}}} \otimes_{C^{{\bf d}-\epsilon_i +\epsilon_{i-1}}} C^{{\bf d}; {\bf d} -\epsilon_i +\epsilon_{i-1}} \cong
\oplus_{s=0}^{d_{i-1}} \oplus_{r=0}^{d_{i-1}+1} C^{{\bf d}+2\epsilon_i-2\epsilon_{i-1}}\langle 2r+2s-2d_i+3 \rangle. $$
The grading of $ \widetilde{\mathcal{F}}_i^{(2)} $ is encoded by
$$ C^{{\bf d}; {\bf d}-2\epsilon_i + 2\epsilon_{i-1}} \cong \oplus_{d_{i-1}+1 \geq r > s \geq 0} C^{{\bf d} -2\epsilon_i +2\epsilon_{i-1}}\langle 2r+2s-2d_i+2 \rangle. $$
Thus 
$ \widetilde{\mathcal{F}}_i \widetilde{\mathcal{F}}_i \cong \widetilde{\mathcal{F}}_i^{(2)}\langle 1 \rangle \oplus \widetilde{\mathcal{F}}_i^{(2)}\langle -1 \rangle. $
\end{proof}

\begin{lemma}
\label{lemma34}
For some indecomposable projective functor $ G_1, $
there are isomorphisms of graded functors:
\begin{enumerate}
\item $ \widetilde{\mathcal{E}}_i \widetilde{\mathcal{E}}_{i+1} \widetilde{\mathcal{E}}_i \cong \widetilde{\mathcal{E}}_i^{(2)} \widetilde{\mathcal{E}}_{i+1} \oplus
\widetilde{\mathcal{E}}_i^{(2)} \widetilde{\mathcal{E}}_{i+1} \oplus \widetilde{G}_1. $
\item $ \widetilde{\mathcal{E}}_{i+1} \widetilde{\mathcal{E}}_i^{(2)} \cong \widetilde{\mathcal{E}}_i^{(2)} \widetilde{\mathcal{E}}_{i+1} \oplus \widetilde{G}_1. $
\item $ \widetilde{\mathcal{E}}_i \widetilde{\mathcal{E}}_{i+1} \widetilde{\mathcal{E}}_i \cong 
\widetilde{\mathcal{E}}_i^{(2)} \widetilde{\mathcal{E}}_{i+1} \oplus \widetilde{\mathcal{E}}_{i+1} \widetilde{\mathcal{E}}_i^{(2)}. $
\end{enumerate}
\end{lemma}

\begin{proof}
First we consider the ungraded case and compute in the Grothendieck group.
\begin{align*}
&[\mathcal{E}_i][\mathcal{E}_{i+1}][\mathcal{E}_i][M(k-1, \ldots, k-1, \ldots, 0, \ldots, 0)] =\\ 
&2[P(k-1, \ldots, k-1, \ldots, i+1, \ldots, i+1, \underbrace{i, \ldots, i, i+1,}_{d_i} \underbrace{i-1, \ldots, i-1, i, i,}_{d_{i-1}} 0, \ldots, 0)]+\\
&[P(k-1, \ldots, k-1, \ldots, i, \ldots, i, \underbrace{i-1, \ldots, i-1, i, i+1,}_{d_{i-1}} \ldots, 0, \ldots, 0)]. 
\end{align*}

We also have
\begin{align*}
&[\mathcal{E}_i^{(2)}][\mathcal{E}_{i+1}][M(k-1, \ldots, k-1, \ldots, 0, \ldots, 0)] =\\
&[P(k-1, \ldots, k-1, \ldots, \underbrace{i, \ldots, i, i+1,}_{d_i} \underbrace{i-1, \ldots, i-1, i, i,}_{d_{i-1}} \ldots, 0, \ldots, 0)] 
\end{align*}  
and
\begin{align*}
&[\mathcal{E}_{i+1}][\mathcal{E}_i^{(2)}][M(k-1, \ldots, k-1, \ldots, 0, \ldots, 0)]=\\
&[P(k-1, \ldots, k-1, \ldots, i, \ldots, i, i-1, \ldots, i-1, i, i+1, \ldots, 0, \ldots, 0)]+\\
&[P(k-1, \ldots, k-1, \ldots, i+1, \ldots, i+1, i, \ldots, i, i+1, i-1, \ldots, i-1, i, i, \ldots, 0, \ldots, 0)]. 
\end{align*} 

Therefore we get the desired isomorphism of functors if we ignore the grading.
There is an isomorphism
$$ \widetilde{\mathcal{E}}_{i+1} \widetilde{\mathcal{E}}_i^{(2)} \cong \widetilde{\mathcal{E}}_i^{(2)} \widetilde{\mathcal{E}}_{i+1} \langle a \rangle  \oplus \widetilde{G}_1 $$
for some graded lift of the indecomposable projective functor $ G_1 $ and for some shift $ a. $

The functor $ \widetilde{\mathcal{E}}_{i+1} \widetilde{\mathcal{E}}_i^{(2)} $ is given by tensoring with the graded $ (A_{\bf d}, A_{{\bf d}+\epsilon_{i+1}+\epsilon_i-2\epsilon_{i-1}})- $ bimodule
\begin{align*}
&\Hom_{C^{{\bf d}+2\epsilon_i-2\epsilon_{i-1}}}(\mathbb{V}P_{{\bf d}+2\epsilon_i-2\epsilon_{i-1}}, C^{{\bf d}; {\bf d}+2\epsilon_i-2\epsilon_{i-1}}
\otimes_{C^{\bf d}} 
\mathbb{V}P_{\bf d} \langle -2d_{i-1}+4 \rangle) \otimes_{A_{{\bf d}+2\epsilon_i-2\epsilon_{i-1}}}\\
&\Hom_{C^{{\bf d}+\epsilon_{i+1}+\epsilon_i-2\epsilon_{i-1}}}(\mathbb{V}P_{{\bf d}+\epsilon_{i+1}+\epsilon_i-2\epsilon_{i-1}}, 
C^{{\bf d}+2\epsilon_i-2\epsilon_{i-1}; {\bf d}+\epsilon_{i+1}+\epsilon_i-2\epsilon_{i-1}}
\otimes_{C^{{\bf d}+2\epsilon_i-2\epsilon_{i-1}}}
\mathbb{V}P_{{\bf d}+2\epsilon_i-2\epsilon_{i-1}} \langle -d_i-1 \rangle).
\end{align*} 

Therefore we must consider
$$ X = C^{{\bf d}+2\epsilon_i-2\epsilon_{i-1}; {\bf d}+\epsilon_{i+1}+\epsilon_i-2\epsilon_{i-1}} 
\otimes_{C^{{\bf d}+2\epsilon_i-2\epsilon_{i-1}}} C^{{\bf d}; {\bf d}+2\epsilon_i-2\epsilon_{i-1}} \langle -2d_{i-1}-d_i+3 \rangle $$
as a $ C^{{\bf d}+\epsilon_{i+1}+\epsilon_i-2\epsilon_{i-1}}- $ module.

As a $ C^{{\bf d}+2\epsilon_i-2\epsilon_{i-1}}- $ module, $ C^{{\bf d}, {\bf d}+2\epsilon_i-2\epsilon_{i-1}} $ is free of rank
$ |W^{{\bf d}+2\epsilon_i-2\epsilon_{i-1}}/W^{{\bf d}; {\bf d}+2\epsilon_i-2\epsilon_{i-1}}| $ and a basis can be chosen in degrees twice the length of minimal coset representatives.  Therefore we would like to consider minimal coset representatives in $ S_{d_i+2}/S_{d_i} \times S_2. $
Proposition A-2 of [Str] tells us the length of these coset representatives.  We get
$ C^{{\bf d}; {\bf d}+2\epsilon_i-2\epsilon_{i-1}} \cong \oplus _{d_i+1 \geq r > s \geq 0} C^{{\bf d}+2\epsilon_i-2\epsilon_{i-1}} \langle 2r+2s-2 \rangle. $

Similarly, as a $ C^{{\bf d}+\epsilon_{i+1}+\epsilon_i-2\epsilon_{i-1}}- $ module, $ C^{{\bf d}+2\epsilon_i-2\epsilon_{i-1}; {\bf d}+\epsilon_{i+1}+\epsilon_i-2\epsilon_{i-1}} $ is free and a basis can be chosen in degrees twice the length of minimal coset representatives of 
$$ S_{d_k} \times \cdots \times S_{d_{i+1}+1} \times S_{d_i+1} \times S_{d_{i-1}-2} \times \cdots \times S_{d_0}/S_{d_k} \times \cdots \times S_{d_{i+1}} \times S_1 \times S_{d_i+1} \times S_{d_{i-1}-2}
\times \cdots \times S_{d_0}. $$
Therefore,
$$ X \cong \oplus_{t=0}^{d_{i+1}} \oplus_{d_i+1 \geq r>s \geq 0} C^{{\bf d}+\epsilon_{i+1}+\epsilon_i-2\epsilon_{i-1}} \langle 2r+2s+2t-2d_{i-1}-d_i+1 \rangle. $$

Next, $ \widetilde{\mathcal{E}}_i^{(2)} \widetilde{\mathcal{E}}_{i+1} $ is given by tensoring with the graded $ (A_{\bf d}, A_{{\bf d}+\epsilon_{i+1}+\epsilon_i-2\epsilon_{i-1}})- $ bimodule
\begin{align*}
&\Hom_{C^{{\bf d}+\epsilon_{i+1}-\epsilon_{i}}}(\mathbb{V}P_{{\bf d}+\epsilon_{i+1}-\epsilon_{i}}, C^{{\bf d}; {\bf d}+\epsilon_{i+1}-\epsilon_{i}}
\otimes_{C^{\bf d}} 
\mathbb{V}P_{\bf d} \langle -d_i+1 \rangle) \otimes_{A_{{\bf d}+\epsilon_{i+1}-\epsilon_{i-1}}}\\
&\Hom_{C^{{\bf d}+\epsilon_{i+1}+\epsilon_i-2\epsilon_{i-1}}}(\mathbb{V}P_{{\bf d}+\epsilon_{i+1}+\epsilon_i-2\epsilon_{i-1}}, 
C^{{\bf d}+\epsilon_{i+1}-\epsilon_{i}; {\bf d}+\epsilon_{i+1}+\epsilon_i-2\epsilon_{i-1}}
\otimes_{C^{{\bf d}+\epsilon_{i+1}-\epsilon_{i}}}
\mathbb{V}P_{{\bf d}+\epsilon_{i+1}-\epsilon_{i}} \langle -2d_{i-1}+4 \rangle). 
\end{align*} 
Thus we must consider
$$ Y = C^{{\bf d}+\epsilon_{i+1}-\epsilon_{i}; {\bf d}+\epsilon_{i+1}+\epsilon_i-2\epsilon_{i-1}} 
\otimes_{C^{{\bf d}+\epsilon_{i+1}-\epsilon_{i}}} C^{{\bf d}; {\bf d}+\epsilon_{i+1}-\epsilon_{i}} \langle -2d_{i-1}-d_i+5 \rangle $$
as a $ C^{{\bf d}+\epsilon_{i+1}+\epsilon_i-2\epsilon_{i-1}}- $ module.
As a $ C^{{\bf d}+\epsilon_{i+1}-\epsilon_i}- $ module, $ C^{{\bf d}; {\bf d}+\epsilon_{i+1}-\epsilon_i} $ is free and a basis can be chosen in degrees twice the length of minimal coset representatives from
$$ S_{d_{k-1}} \times \cdots \times S_{d_{i+1}+1} \times S_{d_i-1} \times \cdots \times S_{d_0}/
S_{d_{k-1}} \times \cdots \times S_{d_{i+1}} \times S_1 \times S_{d_i-1} \times \cdots \times S_{d_0}. $$
So now we want to consider minimal coset representatives in $ S_{d_{i+1}+1}/S_{d_{i+1}} \times S_1. $
Proposition A-2 of [Str] gives 
$$ C^{{\bf d}; {\bf d}+\epsilon_{i+1}-\epsilon_i}  \cong \oplus_{t=0}^{d_{i+1}} C^{{\bf d}+\epsilon_{i+1}-\epsilon_i} \langle 2t \rangle. $$
Similarly, as a $ C^{{\bf d}+\epsilon_{i+1}+\epsilon_i-2\epsilon_{i-1}}- $ module, $ C^{{\bf d}+\epsilon_{i+1}-\epsilon_i; {\bf d}+\epsilon_{i+1}+\epsilon_i-2\epsilon_{i-1}} $ is free and a basis can be chosen in degrees twice the length of minimal coset representatives of
$$ S_{d_{k-1}} \times \cdots \times S_{d_{i+1}+1} \times S_{d_i + 1} \times S_{d_{i-1}-2} \times \cdots \times S_{d_0}/
S_{d_{k-1}} \times \cdots \times S_{d_{i+1}+1} \times S_{d_i-1} \times S_2 \times S_{d_{i-1}-2} \times \cdots \times S_{d_0}. $$
Therefore,
$$ Y \cong \oplus_{t=0}^{d_{i+1}} \oplus_{d_i \geq r > s \geq 0} C^{{\bf d}+\epsilon_{i+1}+\epsilon_i-2\epsilon_{i-1}} \langle 2r+2s+2t-2d_{i-1}-d_i+3 \rangle. $$

Finally, $ \widetilde{\mathcal{E}}_i \widetilde{\mathcal{E}}_{i+1} \widetilde{\mathcal{E}}_i $ is given by tensoring with the graded 
$ (A_{{\bf d}}, A_{{\bf d}+\epsilon_{i+1}+\epsilon_i-2\epsilon_{i-1}})- $ bimodule
\begin{align*}
&\Hom_{C^{{\bf d}+\epsilon_{i}-\epsilon_{i-1}}}(\mathbb{V}P_{{\bf d}+\epsilon_{i}-\epsilon_{i-1}}, C^{{\bf d}; {\bf d}+\epsilon_{i}-\epsilon_{i-1}}
\otimes_{C^{\bf d}} 
\mathbb{V}P_{\bf d} \langle -d_{i-1}+1 \rangle) \otimes_{A_{{\bf d}+\epsilon_{i}-\epsilon_{i-1}}}\\
&\Hom_{C^{{\bf d}+\epsilon_{i+1}-\epsilon_{i-1}}}(\mathbb{V}P_{{\bf d}+\epsilon_{i+1}-\epsilon_{i-1}}, C^{{\bf d}; {\bf d}+\epsilon_{i}-\epsilon_{i-1}}
\otimes_{C^{{\bf d}+\epsilon_i-\epsilon_{i-1}}} 
\mathbb{V}P_{{\bf d}+\epsilon_i-\epsilon_{i-1}} \langle -d_i \rangle) \otimes_{A_{{\bf d}+\epsilon_{i+1}-\epsilon_{i-1}}}\\
&\Hom_{C^{{\bf d}+\epsilon_{i+1}+\epsilon_{i}-2\epsilon_{i-1}}}(\mathbb{V}P_{{\bf d}+\epsilon_{i+1}+\epsilon_{i}-2\epsilon_{i-1}}, C^{{\bf d}+\epsilon_{i+1}-\epsilon_{i-1}; {\bf d}+\epsilon_{i+1}-\epsilon_{i}-2\epsilon_{i-1}}
\otimes_{C^{{\bf d}+\epsilon_{i+1}-\epsilon_{i-1}}} 
\mathbb{V}P_{{\bf d}+\epsilon_{i+1}-\epsilon_{i-1}} \langle -d_{i-1}+2\rangle). 
\end{align*} 
Thus we must consider
\begin{align*}
W= &C^{{\bf d}+\epsilon_{i+1}-\epsilon_{i-1}; {\bf d}+\epsilon_{i+1}-\epsilon_{i}-2\epsilon_{i-1}}
\otimes_{C^{{\bf d}+\epsilon_{i+1}-\epsilon_{i-1}}}
C^{{\bf d}+\epsilon_{i}-\epsilon_{i-1}; {\bf d}+\epsilon_{i+1}-\epsilon_{i-1}}
\otimes_{C^{{\bf d}+\epsilon_{i}-\epsilon_{i-1}}}\\
&C^{{\bf d}+\epsilon_{i}-\epsilon_{i-1}; {\bf d}+\epsilon_{i}-\epsilon_{i-1}} \langle -2d_{i-1}-d_i+3 \rangle. 
\end{align*} 
As in all the other computations,
$$ W \cong \oplus_{t=0}^{d_{i+1}} \oplus_{r=0}^{d_i} \oplus_{s=0}^{d_i} C^{{\bf d}+\epsilon_{i+1}+\epsilon_i-2\epsilon_{i-1}} \langle 2r+2s+2t-2d_{i-1}-d_i+3 \rangle. $$

The grading for $ \widetilde{\mathcal{E}}_i \widetilde{\mathcal{E}}_{i+1} \widetilde{\mathcal{E}}_i $ is given by the summation
$$ \sum_{t=0}^{d_{i+1}} \sum_{r=0}^{d_i} \sum_{s=0}^{d_i} q^{2r+2s+2t-2d_{i-1}-d_i+3} = $$
$$ 2 \sum_{t=0}^{d_{i+1}} \sum_{d_i \geq r > s \geq 0} q^{2r+2s+2t-2d_{i-1}-d_i+3} +
\sum_{t=0}^{d_{i+1}} \sum_{r=0}^{d_i} q^{4r+2t-2d_{i-1}-d_i+3}. $$

The grading for $ \widetilde{\mathcal{E}}_{i+1} \widetilde{\mathcal{E}}_i^{(2)} $ is given by
$$ \sum_{t=0}^{d_{i+1}} \sum_{d_i+1 \geq r > s \geq 0} q^{2(r-1)+2s+2t-2d_{i-1}-d_i+3}. $$
Letting $ u=r-1, $ we get
$$ \sum_{t=0}^{d_{i+1}} \sum_{d_i \geq u \geq s \geq 0} q^{2u+2s+2t-2d_{i-1}-d_i+3} = $$
$$ \sum_{t=0}^{d_{i+1}} \sum_{d_i \geq u>s \geq 0} q^{2u+2s+2t-2d_{i-1}-d_i+3} + 
\sum_{t=0}^{d_{i+1}} \sum_{u=0}^{d_i} q^{4u+2t-2d_{i-1}-d_i+3}. $$

The grading for $ \widetilde{\mathcal{E}}_i^{(2)} \widetilde{\mathcal{E}}_{i+1} $ is given by
$$ \sum_{t=0}^{d_{i+1}} \sum_{d_i \geq u>s \geq 0} q^{2u+2s+2t-2d_{i-1}-d_i+3} . $$
Therefore we get the desired isomorphisms of graded functors.
\end{proof}

\begin{lemma}
\label{lemma35}
For some indecomposable projective functor $ G_2, $
there are isomorphisms of graded functors:
\begin{enumerate}
\item $ \widetilde{\mathcal{E}}_i \widetilde{\mathcal{E}}_{i-1} \widetilde{\mathcal{E}}_i \cong \widetilde{\mathcal{E}}_{i-1} \widetilde{\mathcal{E}}_i^{(2)}  \oplus
\widetilde{\mathcal{E}}_{i-1} \widetilde{\mathcal{E}}_i^{(2)}  \oplus \widetilde{G}_2. $
\item $ \widetilde{\mathcal{E}}_i^{(2)} \widetilde{\mathcal{E}}_{i-1} \cong \widetilde{\mathcal{E}}_{i-1} \widetilde{\mathcal{E}}_i^{(2)}  \oplus \widetilde{G}_2. $
\item $ \widetilde{\mathcal{E}}_i \widetilde{\mathcal{E}}_{i-1} \widetilde{\mathcal{E}}_i \cong 
\widetilde{\mathcal{E}}_i^{(2)} \widetilde{\mathcal{E}}_{i-1} \oplus \widetilde{\mathcal{E}}_{i-1} \widetilde{\mathcal{E}}_i^{(2)}. $
\end{enumerate}
\end{lemma}

\begin{proof}
First we consider the ungraded case and compute in the Grothendieck group.
\begin{align*}
&[\mathcal{E}_i][\mathcal{E}_{i-1}][\mathcal{E}_i][M(k-1, \ldots, k-1, \ldots, 0, \ldots, 0)] =\\
&2[P(k-1, \ldots, k-1, \ldots, i, \ldots, i, \underbrace{i-1, \ldots, i-1, i,i,}_{d_{i-1}} \underbrace{i-2, \ldots, i-2, i-1,}_{d_{i-2}} 0, \ldots, 0)]+\\
&[P(k-1, \ldots, k-1, \ldots, i, \ldots, i, \underbrace{i-1, \ldots, i-1, i, i,}_{d_{i-1}} \underbrace{i-2,\ldots, i-2, i,}_{d_{i-2}} \ldots, 0, \ldots, 0)]. 
\end{align*} 

We also have
\begin{align*}
&[\mathcal{E}_{i-1}][\mathcal{E}_i^{(2)}][M(k-1, \ldots, k-1, \ldots, 0, \ldots, 0)] =\\
&[P(k-1, \ldots, k-1, \ldots, \underbrace{i-1, \ldots, i-1, i,}_{d_{i-1}} \underbrace{i-2, \ldots, i-2, i,}_{d_{i-2}} \ldots, 0, \ldots, 0)] 
\end{align*} 

and
\begin{align*}
&[\mathcal{E}_i^{(2)}][\mathcal{E}_{i-1}][M(k-1, \ldots, k-1, \ldots, 0, \ldots, 0)]=\\
&[P(k-1, \ldots, k-1, \underbrace{i-1, \ldots, i-1, i, i,}_{d_{i-1}} \underbrace{i-2, \ldots, i-2, i-1,}_{d_{i-2}} \ldots, 0, \ldots, 0)]+\\
&[P(k-1, \ldots, k-1, \underbrace{i-1, \ldots, i-1,i,}_{d_{i-1}} \underbrace{i-2, \ldots, i-2, i,}_{d_{i-2}} \ldots, 0, \ldots, 0)]. 
\end{align*} 

Therefore we get the desired isomorphism of functors if we ignore the grading.
The functor $ \widetilde{\mathcal{E}}_{i}^{(2)} \widetilde{\mathcal{E}}_{i-1} $ is given by tensoring with the graded $ (A_{\bf d}, A_{{\bf d}-\epsilon_{i-2}-\epsilon_{i-1}+2\epsilon_{i}})- $ bimodule
\begin{align*}
&\Hom_{C^{{\bf d}+\epsilon_{i-1}-2\epsilon_{i-2}}}(\mathbb{V}P_{{\bf d}+\epsilon_{i-1}-\epsilon_{i-2}}, C^{{\bf d}; {\bf d}+\epsilon_{i-1}-\epsilon_{i-2}}
\otimes_{C^{\bf d}} 
\mathbb{V}P_{\bf d} \langle -d_{i-2}+1 \rangle) \otimes_{A_{{\bf d}+\epsilon_{i-1}-\epsilon_{i-2}}}\\
&\Hom_{C^{{\bf d}+2\epsilon_{i}-\epsilon_{i-1}-\epsilon_{i-2}}}(\mathbb{V}P_{{\bf d}+2\epsilon_{i}-\epsilon_{i-1}-\epsilon_{i-2}}, 
C^{{\bf d}+\epsilon_{i-1}-\epsilon_{i-2}; {\bf d}+2\epsilon_{i}-\epsilon_{i-1}-\epsilon_{i-2}}
\otimes_{C^{{\bf d}+\epsilon_{i-1}-\epsilon_{i-2}}}
\mathbb{V}P_{{\bf d}+\epsilon_{i-1}-\epsilon_{i-2}} \langle -2d_{i-1}+2 \rangle).
\end{align*} 

Therefore we must consider
$$ X = C^{{\bf d}+\epsilon_{i-1}-\epsilon_{i-2}; {\bf d}+2\epsilon_{i}-\epsilon_{i-1}-\epsilon_{i-2}} 
\otimes_{C^{{\bf d}+\epsilon_{i-1}-\epsilon_{i-2}}} C^{{\bf d}; {\bf d}+\epsilon_{i-1}-\epsilon_{i-2}} \langle -2d_{i-1}-d_{i-2}+3 \rangle $$
as a $ C^{{\bf d}+2\epsilon_{i}-\epsilon_{i-1}-\epsilon_{i-2}}- $ module.

As a $ C^{{\bf d}+\epsilon_{i-1}-\epsilon_{i-2}}- $ module, $ C^{{\bf d}; {\bf d}+\epsilon_{i-1}-\epsilon_{i-2}} $ is free of rank
$ |W^{{\bf d}+\epsilon_{i-1}-\epsilon_{i-2}}/W^{{\bf d}; {\bf d}+\epsilon_{i-1}-\epsilon_{i-2}}| $ and a basis can be chosen in degrees twice the length of minimal coset representatives.  Therefore we would like to consider minimal coset representatives in $ S_{d_{i-1}+1}/S_{d_{i-1}} \times S_1. $
Proposition A-2 of [Str] tells us the length of these coset representatives.  We get
$ C^{{\bf d}; {\bf d}+\epsilon_{i-1}-\epsilon_{i-2}} \cong \oplus _{t=0}^{d_{i-1}} C^{{\bf d}+\epsilon_{i-1}-\epsilon_{i-2}} \langle 2t \rangle. $

Similarly, as a $ C^{{\bf d}+2\epsilon_{i}-\epsilon_{i-1}-\epsilon_{i-2}}- $ module, $ C^{{\bf d}+\epsilon_{i-1}-\epsilon_{i-2}; {\bf d}+2\epsilon_{i}-\epsilon_{i-1}-\epsilon_{i-2}} $ is free and a basis can be chosen in degrees twice the length of minimal coset representatives of 
$$ S_{d_{k-1}} \times \cdots \times S_{d_{i}+2} \times S_{d_{i-1}-1} \times S_{d_{i-2}-1} \times \cdots \times S_{d_0}/
S_{d_{k-1}} \times \cdots \times S_{d_{i}} \times S_2 \times S_{d_{i-1}-1} \times S_{d_{i-2}-1}
\times \cdots \times S_{d_0}. $$
Therefore,
$$ X \cong \oplus_{t=0}^{d_{i-1}} \oplus_{d_i+1 \geq r>s \geq 0} C^{{\bf d}+2\epsilon_{i}-\epsilon_{i-1}-\epsilon_{i-2}} \langle 2r+2s+2t-2d_{i-1}-d_{i-2}+1 \rangle. $$

Next, $ \widetilde{\mathcal{E}}_{i-1} \widetilde{\mathcal{E}}_{i}^{(2)} $ is given by tensoring with the graded 
$ (A_{\bf d}, A_{{\bf d}+2\epsilon_{i}-\epsilon_{i-1}-\epsilon_{i-2}})- $ bimodule
\begin{align*}
&\Hom_{C^{{\bf d}+2\epsilon_{i}-2\epsilon_{i-1}}}(\mathbb{V}P_{{\bf d}+2\epsilon_{i}-2\epsilon_{i-1}}, C^{{\bf d}; {\bf d}+2\epsilon_{i}-2\epsilon_{i-1}}
\otimes_{C^{\bf d}} 
\mathbb{V}P_{\bf d} \langle -2d_{i-1}+4 \rangle) \otimes_{A_{{\bf d}+2\epsilon_{i}-2\epsilon_{i-1}}}\\
&\Hom_{C^{{\bf d}+2\epsilon_{i}-\epsilon_{i-1}-\epsilon_{i-2}}}(\mathbb{V}P_{{\bf d}+2\epsilon_{i}-\epsilon_{i-1}-\epsilon_{i-2}}, 
C^{{\bf d}+2\epsilon_{i}-2\epsilon_{i-1}; {\bf d}+2\epsilon_{i}-\epsilon_{i-1}-\epsilon_{i-2}}
\otimes_{C^{{\bf d}+2\epsilon_{i}-2\epsilon_{i-1}}}
\mathbb{V}P_{{\bf d}+2\epsilon_{i}-2\epsilon_{i-1}} \langle -d_{i-2}+1 \rangle). 
\end{align*} 
Thus we must consider
$$ Y = C^{{\bf d}+2\epsilon_{i}-2\epsilon_{i-1}; {\bf d}+2\epsilon_{i}-\epsilon_{i-1}-\epsilon_{i-2}} 
\otimes_{C^{{\bf d}+2\epsilon_{i}-2\epsilon_{i-1}}} C^{{\bf d}; {\bf d}+2\epsilon_{i}-2\epsilon_{i-1}} \langle -2d_{i-1}-d_{i-2}+5 \rangle. $$
As a $ C^{{\bf d}+2\epsilon_{i}-2\epsilon_{i-1}}- $ module, $ C^{{\bf d}; {\bf d}+2\epsilon_{i}-2\epsilon_{i-1}} $ is free and a basis can be chosen in degrees twice the length of minimal coset representatives from
$$ S_{d_{k-1}} \times \cdots \times S_{d_{i}+2} \times S_{d_{i-1}-2} \times \cdots \times S_{d_0}/
S_{d_{k-1}} \times \cdots \times S_{d_{i}} \times S_2 \times S_{d_{i-1}-2} \times \cdots \times S_{d_0}. $$
So now we want to consider minimal coset representatives in $ S_{d_{i}+2}/S_{d_{i}} \times S_2. $
Proposition A-2 of [Str] gives 
$$ C^{{\bf d}; {\bf d}+2\epsilon_{i}-2\epsilon_{i-1}}  \cong \oplus_{d_i+1 \geq r>s \geq 0} C^{{\bf d}+2\epsilon_{i}-2\epsilon_{i-1}} \langle 2r+2s-2 \rangle. $$
Similarly, as a $ C^{{\bf d}+2\epsilon_{i}-\epsilon_{i-1}-\epsilon_{i-2}}- $ module, $ C^{{\bf d}+2\epsilon_{i}-2\epsilon_{i-1}; {\bf d}+2\epsilon_{i}-\epsilon_{i-1}-\epsilon_{i-2}} $ is free and a basis can be chosen in degrees twice the length of minimal coset representatives of
$$ S_{d_{k-1}} \times \cdots \times S_{d_{i}+2} \times S_{d_{i-1}-1} \times S_{d_{i-2}-1} \times \cdots \times S_{d_0}/
S_{d_{k-1}} \times \cdots \times S_{d_{i}+2} \times S_{d_{i-1}-2} \times S_1 \times S_{d_{i-2}-1} \times \cdots \times S_{d_0}. $$
Therefore,
$$ Y \cong \oplus_{t=0}^{d_{i-1}-2} \oplus_{d_i+1 \geq r > s \geq 0} C^{{\bf d}+2\epsilon_{i}-\epsilon_{i-1}-\epsilon_{i-2}} \langle 2r+2s+2t-2d_{i-1}-d_{i-2}+3 \rangle. $$

Finally, $ \widetilde{\mathcal{E}}_i \widetilde{\mathcal{E}}_{i-1} \widetilde{\mathcal{E}}_i $ is given by tensoring with the graded 
$ (A_{{\bf d}}, A_{{\bf d}+2\epsilon_{i}-\epsilon_{i-1}-\epsilon_{i-2}})- $ bimodule
\begin{align*}
&\Hom_{C^{{\bf d}+\epsilon_{i}-\epsilon_{i-1}}}(\mathbb{V}P_{{\bf d}+\epsilon_{i}-\epsilon_{i-1}}, C^{{\bf d}; {\bf d}+\epsilon_{i}-\epsilon_{i-1}}
\otimes_{C^{\bf d}} 
\mathbb{V}P_{\bf d} \langle -d_{i-1}+1 \rangle) \otimes_{A_{{\bf d}+\epsilon_{i}-\epsilon_{i-1}}}\\
&\Hom_{C^{{\bf d}+\epsilon_{i}-\epsilon_{i-2}}}(\mathbb{V}P_{{\bf d}+\epsilon_{i}-\epsilon_{i-2}}, C^{{\bf d}+\epsilon_i-\epsilon_{i-1}; {\bf d}+\epsilon_{i}-\epsilon_{i-2}}
\otimes_{C^{{\bf d}+\epsilon_i-\epsilon_{i-1}}} 
\mathbb{V}P_{{\bf d}+\epsilon_i-\epsilon_{i-1}} \langle -d_{i-2}+1 \rangle) \otimes_{A_{{\bf d}+\epsilon_{i}-\epsilon_{i-2}}}\\
&\Hom_{C^{{\bf d}+2\epsilon_{i}-\epsilon_{i-1}-\epsilon_{i-2}}}(\mathbb{V}P_{{\bf d}+2\epsilon_{i}+\epsilon_{i-1}-\epsilon_{i-2}}, 
C^{{\bf d}+\epsilon_{i}-\epsilon_{i-2}; {\bf d}+2\epsilon_{i}-\epsilon_{i-1}-2\epsilon_{i-2}}
\otimes_{C^{{\bf d}+\epsilon_{i}-\epsilon_{i-2}}} 
\mathbb{V}P_{{\bf d}+\epsilon_{i}-\epsilon_{i-2}} \langle -d_{i-1}+1 \rangle). 
\end{align*} 
Thus we must consider
$$ W= C^{{\bf d}+\epsilon_{i}-\epsilon_{i-2}; {\bf d}+2\epsilon_{i}-\epsilon_{i-1}-\epsilon_{i-2}}
\otimes_{C^{{\bf d}+\epsilon_{i}-\epsilon_{i-2}}}
C^{{\bf d}+\epsilon_{i}-\epsilon_{i-1}; {\bf d}+\epsilon_{i}-\epsilon_{i-2}}
\otimes_{C^{{\bf d}+\epsilon_{i}-\epsilon_{i-1}}}
C^{{\bf d}; {\bf d}+\epsilon_{i}-\epsilon_{i-1}} \langle -2d_{i-1}-d_{i-2}+3 \rangle. $$
As in all the other computations,
$$ W \cong \oplus_{t=0}^{d_{i-1}-1} \oplus_{r=0}^{d_i} \oplus_{s=0}^{d_i+1} C^{{\bf d}+2\epsilon_{i}-\epsilon_{i-1}-\epsilon_{i-2}} \langle 2r+2s+2t-2d_{i-1}-d_{i-1}+3 \rangle. $$

The grading for $ \widetilde{\mathcal{E}}_{i}^{(2)} \widetilde{\mathcal{E}}_{i-1} \oplus \widetilde{\mathcal{E}}_{i-1} \widetilde{\mathcal{E}}_{i}^{(2)} $ is given by 
$$ \sum_{t=0}^{d_{i-1}} \sum_{d_i+1 \geq r > s \geq 0} q^{2r+2s+2t-2d_{i-1}-d_{i-2}+1} +
\sum_{t=0}^{d_{i-1}-2} \sum_{d_i+1 \geq r > s \geq 0} q^{2r+2s+2t-2d_{i-1}-d_{i-2}+3}= $$
$$ \sum_{t=0}^{d_{i-1}-1} \sum_{d_i+1 \geq r > s \geq 0} q^{2r+2s+2t-2d_{i-1}-d_{i-2}+1} +
\sum_{t=0}^{d_{i-1}-1} \sum_{d_i+1 \geq r > s \geq 0} q^{2r+2s+2t-2d_{i-1}-d_{i-2}+3} = $$
$$ (1+q^{-2}) \sum_{t=0}^{d_{i-1}-1} \sum_{d_i+1 \geq r > s \geq 0} q^{2r+2s+2t-2d_{i-1}-d_{i-2}+3} = $$
$$ \sum_{t=0}^{d_{i-1}-1} \sum_{r=0}^{d_i+1} \sum_{s=0}^{d_i} q^{2r+2s+2t-2d_{i-1}-d_{i-2}+3}. $$
This is the grading for $ \widetilde{\mathcal{E}}_i \widetilde{\mathcal{E}}_{i-1} \widetilde{\mathcal{E}}_i $ so we have the desired isomorphism of graded functors.
\end{proof}

\begin{lemma}
\label{lemma36}
For some indecomposable projective functor $ G_3, $
there are isomorphisms of graded functors:
\begin{enumerate}
\item $ \widetilde{\mathcal{F}}_i \widetilde{\mathcal{F}}_{i+1} \widetilde{\mathcal{F}}_i \cong \widetilde{\mathcal{F}}_{i+1} \widetilde{\mathcal{F}}_i^{(2)} \oplus
\widetilde{\mathcal{F}}_{i+1} \widetilde{\mathcal{F}}_i^{(2)}  \oplus \widetilde{G}_3. $
\item $ \widetilde{\mathcal{F}}_i^{(2)} \widetilde{\mathcal{F}}_{i+1}  \cong \widetilde{\mathcal{F}}_{i+1} \widetilde{\mathcal{F}}_i^{(2)}  \oplus \widetilde{G}_3. $
\item $ \widetilde{\mathcal{F}}_i \widetilde{\mathcal{F}}_{i+1} \widetilde{\mathcal{F}}_i \cong 
\widetilde{\mathcal{F}}_i^{(2)} \widetilde{\mathcal{F}}_{i+1} \oplus \widetilde{\mathcal{F}}_{i+1} \widetilde{\mathcal{F}}_i^{(2)}. $
\end{enumerate}
\end{lemma}

\begin{proof}
First we consider the ungraded case and compute in the Grothendieck group.
\begin{align*}
&[\mathcal{F}_i][\mathcal{F}_{i+1}][\mathcal{F}_i][M(k-1, \ldots, k-1, \ldots, 0, \ldots, 0)] =\\
&2[P(k-1, \ldots, k-1, \ldots, \underbrace{i, i+1, i+1, \ldots, i+1,}_{d_{i+1}} \underbrace{i-1, i-1, i, \ldots, i-1,}_{d_{i}} 0, \ldots, 0)]+\\
&[P(k-1, \ldots, k-1, \ldots, \underbrace{i-1,i+1, i+1, \ldots, i+1,}_{d_{i}} \ldots, 0, \ldots, 0)]. 
\end{align*} 

We also have
\begin{align*}
&[\mathcal{F}_{i+1}][\mathcal{F}_{i}^{(2)}][M(k-1, \ldots, k-1, \ldots, 0, \ldots, 0)] =\\
&[P(k-1, \ldots, k-1, \ldots, \underbrace{i, i+1,i+1, \ldots, i+1,}_{d_{i+1}} \underbrace{i-1,i-1,i, \ldots, i,}_{d_{i}} \ldots, 0, \ldots, 0)] 
\end{align*}  
and
\begin{align*}
&[\mathcal{F}_{i}^{(2)}][\mathcal{F}_{i+1}][M(k-1, \ldots, k-1, \ldots, 0, \ldots, 0)]=\\
&[P(k-1, \ldots, k-1, \ldots, \underbrace{i-1, i+1, i+1, \ldots, i+1,}_{d_{i+1}} \underbrace{i-1, i, \ldots, i}_{d_{i}}, \ldots, 0, \ldots, 0)]+\\
&[P(k-1, \ldots, k-1, \ldots, \underbrace{i, i+1, i+1, \ldots, i+1,}_{d_i} \underbrace{i-1, i-1, i, \ldots, i,}_{d_{i}} \ldots, 0, \ldots, 0)]. 
\end{align*} 

Therefore we get the desired isomorphism of functors if we ignore the grading.

The functor $ \widetilde{\mathcal{F}}_{i}^{(2)} \widetilde{\mathcal{F}}_{i+1} $ is given by tensoring with the graded 
$ (A_{\bf d}, A_{{\bf d}+2\epsilon_{i-1}-\epsilon_{i}-\epsilon_{i+1}})- $ bimodule
\begin{align*}
&\Hom_{C^{{\bf d}+\epsilon_i-2\epsilon_{i+1}}}(\mathbb{V}P_{{\bf d}+\epsilon_i-\epsilon_{i+1}}, C^{{\bf d}; {\bf d}+\epsilon_i-\epsilon_{i+1}}
\otimes_{C^{{\bf d}} }
\mathbb{V}P_{\bf d} \langle -d_{i+1}+1 \rangle) \otimes_{A_{{\bf d}+\epsilon_i-\epsilon_{i+1}}}\\
&\Hom_{C^{{\bf d}+2\epsilon_{i-1}-\epsilon_{i}-\epsilon_{i+1}}}(\mathbb{V}P_{{\bf d}+2\epsilon_{i-1}-\epsilon_{i}-\epsilon_{i+1}}, 
C^{{\bf d}+\epsilon_{i}-\epsilon_{i+1}; {\bf d}+2\epsilon_{i-1}-\epsilon_{i}-\epsilon_{i+1}}
\otimes_{C^{{\bf d}+\epsilon_i-\epsilon_{i+1}}}
\mathbb{V}P_{{\bf d}+\epsilon_i-\epsilon_{i+1}} \langle -2d_i+2 \rangle).
\end{align*} 

Therefore we must consider
$$ X = C^{{\bf d}+\epsilon_{i}-\epsilon_{i+1}; {\bf d}+2\epsilon_{i-1}-\epsilon_{i}-\epsilon_{i+1}} 
\otimes_{C^{{\bf d}+\epsilon_i-\epsilon_{i+1}}} C^{{\bf d}; {\bf d}+\epsilon_i-\epsilon_{i+1}} \langle -2d_i-d_{i+1}+3 \rangle $$
as a $ C^{{\bf d}+2\epsilon_{i-1}-\epsilon_{i}-\epsilon_{i+1}}- $ module.

As a $ C^{{\bf d}+\epsilon_i-\epsilon_{i+1}}- $ module, $ C^{{\bf d}; {\bf d}+\epsilon_i-\epsilon_{i+1}} $ is free 
and a basis can be chosen in degrees twice the length of minimal coset representatives.  
Therefore we would like to consider minimal coset representatives in $ S_{d_{i}+1}/S_1 \times S_{d_{i}}. $
Proposition A-2 of [Str] tells us the length of these coset representatives.  We get
$ C^{{\bf d}; {\bf d}+\epsilon_i-\epsilon_{i+1}} \cong \oplus _{t=0}^{d_i} C^{{\bf d}+\epsilon_i-\epsilon_{i+1}} \langle 2t \rangle. $

Similarly, as a $ C^{{\bf d}+2\epsilon_{i-1}-\epsilon_{i}-\epsilon_{i+1}}- $ module, 
$ C^{{\bf d}+\epsilon_i-\epsilon_{i+1}; {\bf d}+2\epsilon_{i-1}-\epsilon_{i}-\epsilon_{i+1}} $ is free and a basis can be chosen in degrees twice the length of minimal coset representatives of 
$$ S_{d_{k-1}} \times \cdots \times S_{d_{i+1}-1} \times S_{d_{i}-1} \times S_{d_{i-1}+2} \times \cdots \times S_{d_0}/
S_{d_{k-1}} \times \cdots \times S_{d_{i+1}-1} \times S_{d_{i}-1} \times S_2 \times S_{d_{i-1}} \times \cdots \times S_{d_0}. $$
Therefore,
$$ X \cong \oplus_{t=0}^{d_{i}} \oplus_{d_{i-1}+1 \geq r>s \geq 0} C^{{\bf d}+2\epsilon_{i-1}-\epsilon_{i}-\epsilon_{i+1}} \langle 2r+2s+2t-2d_i-d_{i+1}+3 \rangle. $$

Next, $ \widetilde{\mathcal{F}}_{i+1} \widetilde{\mathcal{F}}_{i}^{(2)} $ is given by tensoring with the graded 
$ (A_{\bf d}, A_{{\bf d}+2\epsilon_{i-1}-\epsilon_{i}-\epsilon_{i+1}})- $ bimodule
\begin{align*}
&\Hom_{C^{{\bf d}-2\epsilon_{i}+2\epsilon_{i-1}}}(\mathbb{V}P_{{\bf d}-2\epsilon_{i}+2\epsilon_{i-1}}, C^{{\bf d}; {\bf d}-2\epsilon_{i}+2\epsilon_{i-1}}
\otimes_{C^{\bf d}} 
\mathbb{V}P_{\bf d} \langle -2d_i+4 \rangle) \otimes_{A_{{\bf d}-2\epsilon_{i}+2\epsilon_{i-1}}}\\
&\Hom_{C^{{\bf d}+2\epsilon_{i-1}-\epsilon_{i}-\epsilon_{i+1}}}(\mathbb{V}P_{{\bf d}+2\epsilon_{i-1}-\epsilon_{i}-\epsilon_{i+1}}, 
C^{{\bf d}-2\epsilon_{i}+2\epsilon_{i-1}; {\bf d}+2\epsilon_{i-2}-\epsilon_{i}-\epsilon_{i+1}}
\otimes_{C^{{\bf d}-2\epsilon_{i}+2\epsilon_{i-1}}} 
\mathbb{V}P_{{\bf d}-2\epsilon_{i}+2\epsilon_{i-1}} \langle -d_{i+1}+1 \rangle). 
\end{align*} 
Thus we must consider
$$ Y = C^{{\bf d}-2\epsilon_{i}+2\epsilon_{i-1}; {\bf d}+2\epsilon_{i-2}-\epsilon_{i-1}-\epsilon_{i+1}} 
\otimes_{C^{{\bf d}-2\epsilon_{i}+2\epsilon_{i-1}}} C^{{\bf d}; {\bf d}-2\epsilon_{i}+2\epsilon_{i-1}} \langle -2d_i-d_{i+1}+5 \rangle. $$
As a $ C^{{\bf d}-2\epsilon_{i}+2\epsilon_{i-1}}- $ module, $ C^{{\bf d}; {\bf d}-2\epsilon_{i}+2\epsilon_{i-1}} $ is free and a basis can be chosen in degrees twice the length of minimal coset representatives of
$$ S_{d_{k-1}} \times \cdots \times S_{d_{i}-2} \times S_{d_{i-1}+2} \times \cdots \times S_{d_0}/
S_{d_{k-1}} \times \cdots \times S_{d_{i}-2} \times S_2 \times S_{d_{i-1}} \times \cdots \times S_{d_0}. $$
So now we want to consider minimal coset representatives in $ S_{d_{i-1}+2}/S_2 \times S_{d_{i-1}}. $
Proposition A-2 of [Str] gives 
$$ C^{{\bf d}; {\bf d}-2\epsilon_{i}+2\epsilon_{i-1}}  \cong \oplus_{d_{i-1}+1 \geq r > s \geq 0} C^{{\bf d}-2\epsilon_{i}+2\epsilon_{i-1}} \langle 2r+2s-2 \rangle. $$
Similarly, as a $ C^{{\bf d}+2\epsilon_{i-1}-\epsilon_{i}-\epsilon_{i+1}}- $ module, 
$ C^{{\bf d}-2\epsilon_{i}+2\epsilon_{i-1}; {\bf d}+2\epsilon_{i-1}-\epsilon_{i}-\epsilon_{i+1}} $ is free and a basis can be chosen in degrees twice the length of minimal coset representatives of
$$ S_{d_{k-1}} \times \cdots \times S_{d_{i+1}-1} \times S_{d_{i} - 1} \times S_{d_{i-1}+2} \times \cdots \times S_{d_0}/
S_{d_{k-1}} \times \cdots \times S_{d_{i+1}-1} \times S_1 \times S_{d_{i}-2} \times S_{d_{i-1}+2} \times \cdots \times S_{d_0}. $$
Therefore,
$$ Y \cong \oplus_{t=0}^{d_{i}-2} \oplus_{d_{i-1}+1 \geq r > s \geq 0} C^{{\bf d}+2\epsilon_{i-1}-\epsilon_{i}-\epsilon_{i+1}} \langle 2r+2s+2t-2d_i-d_{i+1}+3 \rangle. $$

Finally, $ \widetilde{\mathcal{F}}_i \widetilde{\mathcal{F}}_{i+1} \widetilde{\mathcal{F}}_i $ is given by tensoring with the graded 
$ (A_{{\bf d}}, A_{{\bf d}+2\epsilon_{i-1}-\epsilon_{i}-\epsilon_{i+1}})- $ bimodule
\begin{align*}
&\Hom_{C^{{\bf d}-\epsilon_{i}+\epsilon_{i-1}}}(\mathbb{V}P_{{\bf d}-\epsilon_{i}+\epsilon_{i-1}}, C^{{\bf d}; {\bf d}-\epsilon_{i}+\epsilon_{i-1}}
\otimes_{C^{\bf d}} 
\mathbb{V}P_{\bf d} \langle -d_i+1 \rangle) \otimes_{A_{{\bf d}-\epsilon_{i}+\epsilon_{i-1}}}\\
&\Hom_{C^{{\bf d}-\epsilon_{i+1}+\epsilon_{i-1}}}(\mathbb{V}P_{{\bf d}-\epsilon_{i+1}+\epsilon_{i-1}}, C^{{\bf d}-\epsilon_i+\epsilon_{i-1}; {\bf d}-\epsilon_{i+1}+\epsilon_{i-1}}
\otimes_{C^{{\bf d}-\epsilon_i+\epsilon_{i-1}}} 
\mathbb{V}P_{{\bf d}-\epsilon_i+\epsilon_{i-1}} \langle -d_{i+1}+1 \rangle) \otimes_{A_{{\bf d}-\epsilon_{i+1}+\epsilon_{i-1}}}\\
&\Hom_{C^{{\bf d}+2\epsilon_{i-1}-\epsilon_{i}-\epsilon_{i+1}}}(\mathbb{V}P_{{\bf d}+2\epsilon_{i-1}-\epsilon_{i}-\epsilon_{i+1}}, C^{{\bf d}-\epsilon_{i+1}+\epsilon_{i-1}; 
{\bf d}+2\epsilon_{i-1}-\epsilon_{i}-\epsilon_{i+1}}
\otimes_{C^{{\bf d}-\epsilon_{i+1}+\epsilon_{i-1}}} 
\mathbb{V}P_{{\bf d}-\epsilon_{i+1}+\epsilon_{i-1}} \langle -d_i+1 \rangle). 
\end{align*} 
Thus we must consider
$$ W= C^{{\bf d}-\epsilon_{i+1}+\epsilon_{i-1}; {\bf d}+2\epsilon_{i-1}-\epsilon_{i}-\epsilon_{i+1}}
\otimes_{C^{{\bf d}-\epsilon_{i+1}+\epsilon_{i-1}}}
C^{{\bf d}-\epsilon_{i}+\epsilon_{i-1}; {\bf d}-\epsilon_{i+1}+\epsilon_{i-1}}
\otimes_{C^{{\bf d}-\epsilon_{i}+\epsilon_{i-1}}}
C^{{\bf d}; {\bf d}-\epsilon_{i}+\epsilon_{i-1}} \langle -2d_i-d_{i+1}+3 \rangle. $$
As in all the other computations,
$$ W \cong \oplus_{t=0}^{d_{i-1}+1} \oplus_{r=0}^{d_{i-1}} \oplus_{s=0}^{d_{i}-1} C^{{\bf d}+\epsilon_{i-2}+\epsilon_{i-1}-2\epsilon_{i}} \langle 2r+2s+2t-2d_i-d_{i+1}+3 \rangle. $$

The grading for $ \widetilde{\mathcal{F}}_{i}^{(2)} \widetilde{\mathcal{F}}_{i+1} \oplus \widetilde{\mathcal{F}}_{i+1} \widetilde{\mathcal{F}}_{i}^{(2)} $ is given by the summation
$$ \sum_{t=0}^{d_i} \sum_{d_{i-1}+1 \geq r>s \geq 0} q^{2r+2s+2t-2d_i-d_{i+1}+1} +
\sum_{t=0}^{d_i-2} \sum_{d_{i-1}+1 \geq r>s \geq 0} q^{2r+2s+2t-2d_i-d_{i+1}+3} = $$
$$ \sum_{t=0}^{d_i-1} \sum_{d_{i-1}+1 \geq r>s \geq 0} q^{2r+2s+2t-2d_i-d_{i+1}+1} +
\sum_{t=0}^{d_i-1} \sum_{d_{i-1}+1 \geq r>s \geq 0} q^{2r+2s+2t-2d_i-d_{i+1}+3} = $$
$$ (q^{-2}+1) \sum_{t=0}^{d_i-1} \sum_{d_{i-1}+1 \geq r>s \geq 0} q^{2r+2s+2t-2d_i-d_{i+1}+3}. $$
This is the grading for $ \widetilde{\mathcal{F}}_i \widetilde{\mathcal{F}}_{i+1} \widetilde{\mathcal{F}}_i  $
so we have the desired isomorphism of graded functors.
\end{proof}

\begin{lemma}
\label{lemma37}
For some indecomposable projective functor $ G_4, $
there are isomorphisms of graded functors:
\begin{enumerate}
\item $ \widetilde{\mathcal{F}}_i \widetilde{\mathcal{F}}_{i-1} \widetilde{\mathcal{F}}_i \cong \widetilde{\mathcal{F}}_i^{(2)} \widetilde{\mathcal{F}}_{i-1} \oplus
\widetilde{\mathcal{F}}_i^{(2)} \widetilde{\mathcal{F}}_{i-1} \oplus \widetilde{G}_4. $
\item $ \widetilde{\mathcal{F}}_{i-1} \widetilde{\mathcal{F}}_i^{(2)} \cong \widetilde{\mathcal{F}}_i^{(2)} \widetilde{\mathcal{F}}_{i-1} \oplus \widetilde{G}_4. $
\item $ \widetilde{\mathcal{F}}_i \widetilde{\mathcal{F}}_{i-1} \widetilde{\mathcal{F}}_i \cong 
\widetilde{\mathcal{F}}_i^{(2)} \widetilde{\mathcal{F}}_{i-1} \oplus \widetilde{\mathcal{F}}_{i-1} \widetilde{\mathcal{F}}_i^{(2)}. $
\end{enumerate}
\end{lemma}

\begin{proof}
First we consider the ungraded case and compute in the Grothendieck group.
\begin{align*}
&[\mathcal{F}_i][\mathcal{F}_{i-1}][\mathcal{F}_i][M(k-1, \ldots, k-1, \ldots, 0, \ldots, 0)] =\\
&2[P(k-1, \ldots, k-1, \ldots, \underbrace{i-1, i-1, i, \ldots, i,}_{d_i} \underbrace{i-2, i-1, \ldots, i-1,}_{d_{i-1}} 0, \ldots, 0)]+\\
&[P(k-1, \ldots, k-1, \ldots, \underbrace{i-2,i-1, i, \ldots, i,}_{d_{i-1}} \ldots, 0, \ldots, 0)]. 
\end{align*} 

We also have
\begin{align*}
&[\mathcal{F}_i^{(2)}][\mathcal{F}_{i-1}][M(k-1, \ldots, k-1, \ldots, 0, \ldots, 0)] =\\
&[P(k-1, \ldots, k-1, \ldots, \underbrace{i-2, i-1,i, \ldots, i,}_{d_i} \underbrace{i-2,i-1, \ldots, i-1,}_{d_{i-1}} \ldots, 0, \ldots, 0)] 
\end{align*} 
and
\begin{align*}
&[\mathcal{F}_{i-1}][\mathcal{F}_i^{(2)}][M(k-1, \ldots, k-1, \ldots, 0, \ldots, 0)]=\\
&[P(k-1, \ldots, k-1, \ldots, \underbrace{i-1, i-1, i, \ldots, i,}_{d_i} \underbrace{i-2, i-1, \ldots, i-1}_{d_{i-1}}, \ldots, 0, \ldots, 0)]+\\
&[P(k-1, \ldots, k-1, \ldots, \underbrace{i-2, i-1, i, \ldots, i,}_{d_i} \underbrace{i-1, \ldots, i-1,}_{d_{i-1}} \ldots, 0, \ldots, 0)]. 
\end{align*} 

Therefore we get the desired isomorphism of functors if we ignore the grading.

The functor $ \widetilde{\mathcal{F}}_{i-1} \widetilde{\mathcal{F}}_i^{(2)} $ is given by tensoring with the graded $ (A_{\bf d}, A_{{\bf d}+\epsilon_{i-2}+\epsilon_{i-1}-2\epsilon_{i}})- $ bimodule
\begin{align*}
&\Hom_{C^{{\bf d}-2\epsilon_i+2\epsilon_{i-1}}}(\mathbb{V}P_{{\bf d}-2\epsilon_i+2\epsilon_{i-1}}, C^{{\bf d}; {\bf d}-2\epsilon_i+2\epsilon_{i-1}}
\otimes_{C^{\bf d}} 
\mathbb{V}P_{\bf d} \langle -2d_i+4 \rangle) \otimes_{A_{{\bf d}-2\epsilon_i+2\epsilon_{i-1}}}\\
&\Hom_{C^{{\bf d}+\epsilon_{i-2}+\epsilon_{i-1}-2\epsilon_{i}}}(\mathbb{V}P_{{\bf d}+\epsilon_{i-2}+\epsilon_{i-1}-2\epsilon_{i}}, 
C^{{\bf d}+2\epsilon_{i-1}-2\epsilon_{i}; {\bf d}+\epsilon_{i-2}+\epsilon_{i-1}-2\epsilon_{i}}
\otimes_{C^{{\bf d}-2\epsilon_i+2\epsilon_{i-1}}}
\mathbb{V}P_{{\bf d}-2\epsilon_i+2\epsilon_{i-1}} \langle -d_{i-1}-1 \rangle).
\end{align*} 

Therefore we must consider
$$ X = C^{{\bf d}+2\epsilon_{i-1}-2\epsilon_{i}; {\bf d}+\epsilon_{i-2}+\epsilon_{i-1}-2\epsilon_{i}} 
\otimes_{C^{{\bf d}-2\epsilon_i+2\epsilon_{i-1}}} C^{{\bf d}; {\bf d}-2\epsilon_i+2\epsilon_{i-1}} \langle -2d_i-d_{i-1}+3 \rangle $$
as a $ C^{{\bf d}+\epsilon_{i-2}+\epsilon_{i-1}-2\epsilon_{i}}- $ module.

As a $ C^{{\bf d}-2\epsilon_i+2\epsilon_{i-1}}- $ module, $ C^{{\bf d}, {\bf d}-2\epsilon_i+2\epsilon_{i-1}} $ is free of rank
$ |W^{{\bf d}-2\epsilon_i+2\epsilon_{i-1}}/W^{{\bf d}, {\bf d}-2\epsilon_i+2\epsilon_{i-1}}| $ and a basis can be chosen in degrees twice the length of minimal coset representatives.  Therefore we would like to consider minimal coset representatives in $ S_{d_{i-1}+2}/S_2 \times S_{d_{i-1}}. $
Proposition A-2 of [Str] tells us the length of these coset representatives.  We get
$ C^{{\bf d}, {\bf d}-2\epsilon_i+2\epsilon_{i-1}} \cong \oplus _{d_{i-1}+1 \geq r > s \geq 0} C^{{\bf d}-2\epsilon_i+2\epsilon_{i-1}} \langle 2r+2s-2 \rangle. $

Similarly, as a $ C^{{\bf d}+\epsilon_{i-2}+\epsilon_{i-1}-2\epsilon_{i}}- $ module, 
$ C^{{\bf d}-2\epsilon_i+2\epsilon_{i-1}, {\bf d}+\epsilon_{i-2}+\epsilon_{i-1}-2\epsilon_{i}} $ is free and a basis can be chosen in degrees twice the length of minimal coset representatives of 
$$ S_{d_k} \times \cdots S_{d_{i}-2} \times S_{d_{i-1}+1} \times S_{d_{i-2}+1} \times \cdots \times S_{d_0}/
S_{d_k} \times \cdots \times S_{d_{i}-2} \times S_{d_{i-1}+1} S_1 \times S_{d_{i-2}} \times \cdots \times S_{d_0}. $$
Therefore,
$$ X \cong \oplus_{t=0}^{d_{i-2}} \oplus_{d_{i-1}+1 \geq r>s \geq 0} C^{{\bf d}+\epsilon_{i-2}+\epsilon_{i-1}-2\epsilon_{i}} \langle 2r+2s+2t-2d_i-d_{i-1}+1 \rangle. $$

Next, $ \widetilde{\mathcal{F}}_i^{(2)} \widetilde{\mathcal{F}}_{i-1} $ is given by tensoring with the graded 
$ (A_{\bf d}, A_{{\bf d}+\epsilon_{i-2}+\epsilon_{i-1}-2\epsilon_{i}})- $ bimodule
\begin{align*}
&\Hom_{C^{{\bf d}-\epsilon_{i-1}+\epsilon_{i-2}}}(\mathbb{V}P_{{\bf d}-\epsilon_{i-1}+\epsilon_{i-2}}, C^{{\bf d}; {\bf d}-\epsilon_{i-1}+\epsilon_{i-2}}
\otimes_{C^{\bf d}} 
\mathbb{V}P_{\bf d} \langle -d_{i-1}+1 \rangle) \otimes_{A_{{\bf d}-\epsilon_{i-1}+\epsilon_{i-2}}}\\
&\Hom_{C^{{\bf d}+\epsilon_{i-2}+\epsilon_{i-1}-2\epsilon_{i}}}(\mathbb{V}P_{{\bf d}+\epsilon_{i-2}+\epsilon_{i-1}-2\epsilon_{i}}, 
C^{{\bf d}+\epsilon_{i-2}-\epsilon_{i-1}; {\bf d}+\epsilon_{i-2}+\epsilon_{i-1}-2\epsilon_{i}}
\otimes_{C^{{\bf d}+\epsilon_{i-2}-\epsilon_{i-1}}} 
\mathbb{V}P_{{\bf d}-\epsilon_{i-2}-\epsilon_{i-1}} \langle -2d_i+4 \rangle). 
\end{align*} 
Thus we must consider
$$ Y = C^{{\bf d}+\epsilon_{i-2}-\epsilon_{i-1}; {\bf d}+\epsilon_{i-2}+\epsilon_{i-1}+2\epsilon_{i}} 
\otimes_{C^{{\bf d}+\epsilon_{i-2}-\epsilon_{i-1}}} C^{{\bf d}; {\bf d}+\epsilon_{i-2}-\epsilon_{i-1}} \langle -2d_i-d_{i-1}+5 \rangle $$
as a $ C^{{\bf d}+\epsilon_{i-2}+\epsilon_{i-1}-2\epsilon_{i}}- $ module.
As a $ C^{{\bf d}+\epsilon_{i-2}-\epsilon_{i-1}}- $ module, $ C^{{\bf d}; {\bf d}+\epsilon_{i-2}-\epsilon_{i-1}} $ is free and a basis can be chosen in degrees twice the length of minimal coset representatives of
$$ S_{d_{k-1}} \times \cdots \times S_{d_{i-1}-1} \times S_{d_{i-2}+1} \times \cdots \times S_{d_0}/
S_{d_{k-1}} \times \cdots \times S_{d_{i-1}-1} \times S_1 \times S_{d_{i-2}} \times \cdots \times S_{d_0}. $$
So now we want to consider minimal coset representatives in $ S_{d_{i-2}+1}/S_1 \times S_{d_{i-2}}. $
Proposition A-2 of [Str] gives 
$$ C^{{\bf d}, {\bf d}+\epsilon_{i-2}-\epsilon_{i-1}}  \cong \oplus_{t=0}^{d_{i-2}} C^{{\bf d}+\epsilon_{i-2}-\epsilon_{i-1}} \langle 2t \rangle. $$
Similarly, as a $ C^{{\bf d}+\epsilon_{i-2}+\epsilon_{i-1}-2\epsilon_{i}}- $ module, 
$ C^{{\bf d}+\epsilon_{i-2}-\epsilon_{i-1}; {\bf d}+\epsilon_{i-2}+\epsilon_{i-1}-2\epsilon_{i}} $ is free and a basis can be chosen in degrees twice the length of minimal coset representatives of
$$ S_{d_{k-1}} \times \cdots \times S_{d_{i}-2} \times S_{d_{i-1} + 1} \times S_{d_{i-2}+1} \times \cdots \times S_{d_0}/
S_{d_{k-1}} \times \cdots \times S_{d_{i}-2} \times S_2 \times S_{d_{i-1}-1} \times S_{d_{i-2}+1} \times \cdots \times S_{d_0}. $$
Therefore,
$$ Y \cong \oplus_{t=0}^{d_{i-2}} \oplus_{d_{i-1} \geq r > s \geq 0} C^{{\bf d}+\epsilon_{i-2}+\epsilon_{i-1}-2\epsilon_{i}} \langle 2r+2s+2t-2d_i-d_{i-1}+3 \rangle. $$

Finally, $ \widetilde{\mathcal{F}}_i \widetilde{\mathcal{F}}_{i-1} \widetilde{\mathcal{F}}_i $ is given by tensoring with the graded 
$ (A_{{\bf d}}, A_{{\bf d}+\epsilon_{i-2}+\epsilon_{i-1}-2\epsilon_{i}})- $ bimodule
\begin{align*}
&\Hom_{C^{{\bf d}-\epsilon_{i}+\epsilon_{i-1}}}(\mathbb{V}P_{{\bf d}-\epsilon_{i}+\epsilon_{i-1}}, C^{{\bf d}; {\bf d}-\epsilon_{i}+\epsilon_{i-1}}
\otimes_{C^{\bf d}} 
\mathbb{V}P_{\bf d} \langle -d_i+1 \rangle) \otimes_{A_{{\bf d}-\epsilon_{i}+\epsilon_{i-1}}}\\
&\Hom_{C^{{\bf d}+\epsilon_{i-2}-\epsilon_{i}}}(\mathbb{V}P_{{\bf d}+\epsilon_{i-2}-\epsilon_{i}}, C^{{\bf d}-\epsilon_i+\epsilon_{i-1}; {\bf d}-\epsilon_{i}+\epsilon_{i-2}}
\otimes_{C^{{\bf d}-\epsilon_i+\epsilon_{i-1}}} 
\mathbb{V}P_{{\bf d}-\epsilon_i+\epsilon_{i-1}} \langle -d_{i-1} \rangle) \otimes_{A_{{\bf d}+\epsilon_{i-2}-\epsilon_{i}}}\\
&\Hom_{C^{{\bf d}+\epsilon_{i-2}+\epsilon_{i-1}-2\epsilon_{i}}}(\mathbb{V}P_{{\bf d}+\epsilon_{i-2}+\epsilon_{i-1}-2\epsilon_{i}}, C^{{\bf d}+\epsilon_{i-2}-\epsilon_{i};
{\bf d}+\epsilon_{i-2}+\epsilon_{i-1}-2\epsilon_{i}}
\otimes_{C^{{\bf d}+\epsilon_{i-2}-\epsilon_{i}}} 
\mathbb{V}P_{{\bf d}+\epsilon_{i-2}-\epsilon_{i}} \langle -d_i+2 \rangle). 
\end{align*} 
Thus we must consider
$$ W= C^{{\bf d}+\epsilon_{i-2}-\epsilon_{i}; {\bf d}+\epsilon_{i-2}+\epsilon_{i-1}-2\epsilon_{i}}
\otimes_{C^{{\bf d}+\epsilon_{i-2}-\epsilon_{i}}}
C^{{\bf d}+\epsilon_{i-1}-\epsilon_{i}; {\bf d}+\epsilon_{i-2}-\epsilon_{i}}
\otimes_{C^{{\bf d}+\epsilon_{i-1}-\epsilon_{i}}}
C^{{\bf d}; {\bf d}+\epsilon_{i-1}-\epsilon_{i}} \langle -2d_i-d_{i-1}+3 \rangle. $$
As in all the other computations,
$$ W \cong \oplus_{t=0}^{d_{i-2}} \oplus_{r=0}^{d_{i-1}} \oplus_{s=0}^{d_{i-1}} C^{{\bf d}+\epsilon_{i-2}+\epsilon_{i-1}-2\epsilon_{i}} \langle 2r+2s+2t-2d_i-d_{i-1}+3 \rangle. $$

The grading for $ \widetilde{\mathcal{F}}_i \widetilde{\mathcal{F}}_{i-1} \widetilde{\mathcal{F}}_i $ is given by 
$$ \sum_{t=0}^{d_{i-2}} \sum_{r=0}^{d_{i-1}} \sum_{s=0}^{d_{i-1}} q^{2r+2s+2t-2d_i-d_{i-1}+3} = $$
$$ 2 \sum_{t=0}^{d_{i-2}} \sum_{d_{i-1} \geq r>s \geq 0} q^{2r+2s+2t-2d_i-d_{i-1}+3} +
\sum_{t=0}^{d_{i-2}} \sum_{r=0}^{d_{i-1}} q^{4r+2t-2d_i-d_{i-1}+3}. $$

The grading for $ \widetilde{\mathcal{F}}_{i-1} \widetilde{\mathcal{F}}_i^{(2)} $ is given by
$$ \sum_{t=0}^{d_{i-2}} \sum_{d_{i-1}+1 \geq r>s \geq 0} q^{2(r-1)+2s+2t-2d_i-d_{i-1}+3}. $$
Letting $ u=r-1, $ this is equal to
$$ \sum_{t=0}^{d_{i-2}} \sum_{d_{i-1} \geq u \geq s \geq 0} q^{2u+2s+2t-2d_i-d_{i-1}+3} = $$
$$ \sum_{t=0}^{d_{i-2}} \sum_{d_{i-1} \geq u>s \geq 0} q^{2u+2s+2t-2d_i-d_{i-1}+3} +
\sum_{t=0}^{d_{i-2}} \sum_{u=0}^{d_{i-1}} q^{4u+2t-2d_i-d_{i-1}+3}. $$

Since the grading for $ \widetilde{\mathcal{F}}_i^{(2)} $ is given by
$ \sum_{t=0}^{d_{i-2}} \sum_{d_{i-1} \geq r>s \geq 0} q^{2r+2s+2t-2d_i-d_{i-1}+3}, $
we get the desired isomorphism of graded functors.
\end{proof}

Now we may state a categorification of the quantum Serre relations.

\begin{theorem}
There are isomorphisms of graded projective functors:
\begin{enumerate}
\item $ \widetilde{\mathcal{H}}_i \widetilde{\mathcal{H}}_i^{-1} \cong \Id \cong \widetilde{\mathcal{H}}^{-1}_i \widetilde{\mathcal{H}}_i. $
\item $ \widetilde{\mathcal{H}}_i \widetilde{\mathcal{H}}_j \cong \widetilde{\mathcal{H}}_j \widetilde{\mathcal{H}}_i. $
\item If $ |i-j|>1, $ $ \widetilde{\mathcal{H}}_i \widetilde{\mathcal{E}}_j \cong \widetilde{\mathcal{E}}_j \widetilde{\mathcal{H}}_i. $
\item  $ \widetilde{\mathcal{H}}_i \widetilde{\mathcal{E}}_i \cong \widetilde{\mathcal{E}}_i \widetilde{\mathcal{H}}_i \langle 2 \rangle. $
\item If $ |i-j|=1, $ $ \widetilde{\mathcal{H}}_i \widetilde{\mathcal{E}}_j \cong \widetilde{\mathcal{E}}_j \widetilde{\mathcal{H}}_i \langle -1 \rangle. $
\item If $ |i-j|>1, $ $ \widetilde{\mathcal{H}}_i \widetilde{\mathcal{F}}_j \cong \widetilde{\mathcal{F}}_j \widetilde{\mathcal{H}}_i. $
\item  $ \widetilde{\mathcal{H}}_i \widetilde{\mathcal{F}}_i \cong \widetilde{\mathcal{F}}_i \widetilde{\mathcal{H}}_i \langle -2 \rangle. $
\item If $ |i-j|=1, $ $ \widetilde{\mathcal{H}}_i \widetilde{\mathcal{F}}_j \cong \widetilde{\mathcal{F}}_j \widetilde{\mathcal{H}}_i \langle 1 \rangle. $
\item If $ i \neq j, $ $ \widetilde{\mathcal{E}}_i \widetilde{\mathcal{F}}_j \cong \widetilde{\mathcal{F}}_j \widetilde{\mathcal{E}}_i. $
\item If $ |i-j|>1, $ $ \widetilde{\mathcal{E}}_i \widetilde{\mathcal{E}}_j \cong \widetilde{\mathcal{E}}_j \widetilde{\mathcal{E}}_i. $
\item If $ |i-j|>1, $ $ \widetilde{\mathcal{F}}_i \widetilde{\mathcal{F}}_j \cong \widetilde{\mathcal{F}}_j \widetilde{\mathcal{F}}_i. $
\item $ \widetilde{\mathcal{E}}_i \widetilde{\mathcal{E}}_i \widetilde{\mathcal{E}}_{i+1} \oplus \widetilde{\mathcal{E}}_{i+1} \widetilde{\mathcal{E}}_i
\widetilde{\mathcal{E}}_i \cong \widetilde{\mathcal{E}}_i \widetilde{\mathcal{E}}_{i+1} \widetilde{\mathcal{E}}_i \langle 1 \rangle \oplus
\widetilde{\mathcal{E}}_i \widetilde{\mathcal{E}}_{i+1} \widetilde{\mathcal{E}}_i \langle -1 \rangle. $
\item $ \widetilde{\mathcal{E}}_i \widetilde{\mathcal{E}}_i \widetilde{\mathcal{E}}_{i-1} \oplus \widetilde{\mathcal{E}}_{i-1} \widetilde{\mathcal{E}}_i
\widetilde{\mathcal{E}}_i \cong \widetilde{\mathcal{E}}_i \widetilde{\mathcal{E}}_{i-1} \widetilde{\mathcal{E}}_i \langle 1 \rangle \oplus
\widetilde{\mathcal{E}}_i \widetilde{\mathcal{E}}_{i-1} \widetilde{\mathcal{E}}_i \langle -1 \rangle. $
\item $ \widetilde{\mathcal{F}}_i \widetilde{\mathcal{F}}_i \widetilde{\mathcal{F}}_{i+1} \oplus \widetilde{\mathcal{F}}_{i+1} \widetilde{\mathcal{F}}_i
\widetilde{\mathcal{F}}_i \cong \widetilde{\mathcal{F}}_i \widetilde{\mathcal{F}}_{i+1} \widetilde{\mathcal{F}}_i \langle 1 \rangle \oplus
\widetilde{\mathcal{F}}_i \widetilde{\mathcal{F}}_{i+1} \widetilde{\mathcal{F}}_i \langle -1 \rangle. $
\item $ \widetilde{\mathcal{F}}_i \widetilde{\mathcal{F}}_i \widetilde{\mathcal{F}}_{i-1} \oplus \widetilde{\mathcal{F}}_{i-1} \widetilde{\mathcal{F}}_i
\widetilde{\mathcal{F}}_i \cong \widetilde{\mathcal{F}}_i \widetilde{\mathcal{F}}_{i-1} \widetilde{\mathcal{F}}_i \langle 1 \rangle \oplus
\widetilde{\mathcal{F}}_i \widetilde{\mathcal{F}}_{i-1} \widetilde{\mathcal{F}}_i \langle -1 \rangle. $
\item If $ d_{i-1} > d_i $ then
$$ \widetilde{\mathcal{F}}_i \widetilde{\mathcal{E}}_i \cong \widetilde{\mathcal{E}}_i \widetilde{\mathcal{F}}_i 
\oplus_{r=0}^{d_{i-1}-d_i-1} \Id \langle d_{i-1}-d_i-1-2r \rangle. $$
If $ d_{i} > d_{i-1} $ then
$$ \widetilde{\mathcal{E}}_i \widetilde{\mathcal{F}}_i \cong \widetilde{\mathcal{F}}_i \widetilde{\mathcal{E}}_i 
\oplus_{r=0}^{d_{i}-d_{i-1}-1} \Id \langle d_{i}-d_{i-1}-1-2r \rangle. $$		
If $ d_i = d_{i-1}, $ then $ \widetilde{\mathcal{F}}_i \widetilde{\mathcal{E}}_i \cong \widetilde{\mathcal{E}}_i \widetilde{\mathcal{F}}_i. $
\end{enumerate}
\end{theorem}

\begin{proof}
\begin{enumerate}
\item This is obvious.
\item This is obvious.
\item If $ |i-j|>1, $ then the value of $ |d_i - d_{i-1}| $ remains unchanged after applying $ \widetilde{\mathcal{E}}_j. $
\item After applying $ \widetilde{\mathcal{E}}_i, $ $ d_i $ increases by 1 and $ d_{i-1} $ decreases by 1 so there is a shift of 2.
\item If $ j=i+1, $ then $ d_i $ decreases by 1 and $ d_{i-1} $ remains unchanged.
If $ j=i-1, $ then $ d_{i-1} $ increases by 1 and $ d_i $ remains unchanged so there is a shift of -1.
\item This is the same as the proof of 3.
\item After applying $ \widetilde{\mathcal{F}}_i, $ $ d_i $ decreases by 1 and $ d_{i-1} $ increases by 1 so there is a shift of -2.
\item If $ j=i+1, $ then $ d_i $ increases by 1 and $ d_{i-1} $ remains unchanged.
If $ j=i-1, $ then $ d_{i-1} $ decreases by 1 and $ d_{i} $ remains unchanged.
\item We would like to prove that $ \mathcal{E}_i \mathcal{F}_j $ and $ \mathcal{F}_j \mathcal{E}_i $ are indecomposable projective functors if $ i\neq j. $
We show this by computing $ [\mathcal{E}_i \mathcal{F}_j] $ and $ [\mathcal{F}_j \mathcal{E}_i] $ on 
$ [M(k-1, \ldots, k-1, \ldots, 0, \ldots, 0)]. $
If $ i-j > 1, $ both are equal to
$$ [P(k-1, \ldots, k-1, \ldots, i, \ldots, i, \underbrace{i-1, \ldots, i-1, i,}_{d_{i-1}} \ldots, \underbrace{j-1, j, j, \ldots, j,}_{d_j}, \ldots, 0, \ldots, 0)]. $$
The case of $ i-j<1 $ is similar.
If $ i-j=1, $ then both are equal to 
$$ [P(k-1, \ldots, k-1, \ldots, i, \ldots, i, \underbrace{i-2, i-1, i-1, \ldots, i-1, i,}_{d_{i-1}}, i-2, \ldots, i-2, \ldots, 0, \ldots, 0)]. $$
The case of $ i-j=-1 $ is similar.
Thus $ \mathcal{E}_i \mathcal{F}_j $ and $ \mathcal{F}_j \mathcal{E}_i $ are indecomposable projective functors that are equal on the Grothendieck group so they
are isomorphic.  
Since their lowest and highest gradings obviously coincide, their graded lifts are also isomorphic to each other with no shift.
\item Again, we first compute the ungraded functors on the dominant generalized Verma module in the Grothendieck group to get 
$$ [P(k-1, \ldots, k-1, \ldots, \underbrace{j-1, \ldots, j-1, j}_{d_{j-1}}, \ldots, \underbrace{i-1, i-1, \ldots, i-1, i,}_{d_{i-1}} \ldots, 0, \ldots, 0)]. $$
Thus $ \mathcal{E}_i \mathcal{E}_j $ and $ \mathcal{E}_j \mathcal{E}_i $ are indecomposable isomorphic projective functors.
Since their lowest and highest gradings coincide, their graded lifts are isomorphic with no shift.
\item The proof of this is the same as the previous one.  Applying these functors to the dominant generalized Verma module in the Grothendieck group gives
$$ [P(k-1, \ldots, k-1, \ldots, \underbrace{j-1, j, \ldots, j}_{d_j}, \ldots, \underbrace{i-1, i, \ldots, i}_{d_i}, \ldots, 0, \ldots, 0)]. $$  
The rest of the proof is the same.
\item Lemma ~\ref{lemma32} implies both
$$ \widetilde{\mathcal{E}}_{i+1} \widetilde{\mathcal{E}}_i \widetilde{\mathcal{E}}_i \cong \widetilde{\mathcal{E}}_{i+1} \widetilde{\mathcal{E}}_i^{(2)}\langle 1 \rangle \oplus 			\widetilde{\mathcal{E}}_{i+1} \widetilde{\mathcal{E}}_i^{(2)}\langle -1 \rangle $$
$$ \widetilde{\mathcal{E}}_{i} \widetilde{\mathcal{E}}_i \widetilde{\mathcal{E}}_{i+1} \cong \widetilde{\mathcal{E}}_i^{(2)}\langle 1 \rangle \widetilde{\mathcal{E}}_{i+1}  \oplus 		\widetilde{\mathcal{E}}_i^{(2)}\langle -1 \rangle \widetilde{\mathcal{E}}_{i+1}.  $$
Then lemma ~\ref{lemma34} gives
$$  \widetilde{\mathcal{E}}_i \widetilde{\mathcal{E}}_{i+1} \widetilde{\mathcal{E}}_i \cong 
\widetilde{\mathcal{E}}_i^{(2)} \widetilde{\mathcal{E}}_{i+1} \oplus \widetilde{\mathcal{E}}_{i+1} \widetilde{\mathcal{E}}_i^{(2)}. $$
These facts now easily prove the claim.
\item This is similar to the proof of 12 using lemmas ~\ref{lemma32} and ~\ref{lemma35}.
\item This is similar to the proof of 12 using lemmas ~\ref{lemma33} and ~\ref{lemma36}.
\item This is similar to the proof of 12 using lemmas ~\ref{lemma33} and ~\ref{lemma37}.
\item We compute $ [\mathcal{E}_i][\mathcal{F}_i]([M(k-1, \ldots, k-1, \ldots, 0, \ldots, 0]). $
In general, 
$$  [\mathcal{E}_i][\mathcal{F}_i]([M(a_1, \ldots, a_n)] = \sum_{r=1; a_r=i}^{n} \sum_{s=1; a_s=i-1}^n[M(a_1, \ldots, a_r', \ldots, a_s'', \ldots, a_n)] 
+ \sum_{r=1; a_r=i}^{n} [M(a_1, \ldots, a_n)] $$
where $ a_r' =i-1 $ and $ a_s''=i. $
Applied to the dominant Verma module it is equal to
$$ [P(k-1, \ldots, k-1, \ldots, \underbrace{i-1, i, \ldots, i,}_{d_i} \underbrace{i-1, i-1, \ldots, i-1, i,}_{d_{i-1}} 0, \ldots, 0)] + \sum_{r=1; a_r=i}^{n} [M(k-1, \ldots, 0)]. $$
		
On the other hand, 
$$ [\mathcal{F}_i][\mathcal{E}_i]([M(a_1, \ldots, a_n)] = \sum_{r=1; a_r=i}^{n} \sum_{s=1; a_s=i-1}^n[M(a_1, \ldots, a_r' \ldots, a_s'', \ldots, a_n)] 
+ \sum_{s=1; a_s=i-1}^{n} [M(a_1, \ldots, a_n)]. $$
So for the dominant Verma module it is equal to
$$ [P(k-1, \ldots, k-1, \ldots, \underbrace{i-1, i, \ldots, i,}_{d_i} \underbrace{i-1, i-1, \ldots, i-1, i,}_{d_{i-1}} 0, \ldots, 0)] + \sum_{s=1; a_s=i-1}^{n} [M(k-1, \ldots, 0)]. $$
Thus
$$ \mathcal{E}_i \mathcal{F}_i \cong G \oplus \oplus_{r=0}^{d_i -1} \Id $$
and
$$ \mathcal{F}_i \mathcal{E}_i \cong G \oplus \oplus_{s=0}^{d_{i-1}-1} \Id $$
where G is an indecomposable projective functor.
		
As in lemma 3.4 of [FKS], as a $ C^{\bf d}- $ module, $ C^{{\bf d}; {\bf d}+\epsilon_i-\epsilon_{i-1}} \cong \oplus_{r=0}^{d_{i-1}-1} C^{{\bf d}} \langle 2r \rangle $ and
$ C^{{\bf d}; {\bf d}-\epsilon_i+\epsilon_{i-1}} \cong \oplus_{r=0}^{d_{i}-1} C^{{\bf d}} \langle 2r \rangle. $ 
Thus
\begin{equation}
\label{eq3}
C^{{\bf d}; {\bf d}+\epsilon_i-\epsilon_{i-1}} \otimes_{C^{{\bf d}+\epsilon_i-\epsilon_{i-1}}} C^{{\bf d}; {\bf d}+\epsilon_i-\epsilon_{i-1}}  \otimes_{C^{{\bf d}}} \mathbb{C} \cong
\oplus_{l=0}^{d_i} \oplus_{k=0}^{d_{i-1}-1} \mathbb{C} \langle 2k+2l \rangle  
\end{equation} 

\begin{equation}
\label{eq4}
C^{{\bf d}; {\bf d}-\epsilon_i+\epsilon_{i-1}} \otimes_{C^{{\bf d}-\epsilon_i+\epsilon_{i-1}}} C^{{\bf d}; {\bf d}-\epsilon_i+\epsilon_{i-1}}  \otimes_{C^{{\bf d}}} \mathbb{C} \cong
\oplus_{l=0}^{d_{i-1}} \oplus_{k=0}^{d_{i}-1} \mathbb{C} \langle 2k+2l \rangle.   
\end{equation} 
		
Since their lowest degrees coincide and G is indecomposable, there is an isomorphism of graded functors
$$ \widetilde{\mathcal{E}}_i \widetilde{\mathcal{F}}_i \cong \widetilde{G} \oplus \oplus_{r=0}^{d_i -1} \Id \langle m_r \rangle $$
and 
$$ \widetilde{\mathcal{F}}_i \widetilde{\mathcal{E}}_i \cong \widetilde{G} \oplus \oplus_{s=0}^{d_{i-1}-1} \Id \langle n_r \rangle $$
		
Now $ \widetilde{\mathcal{F}}_i \widetilde{\mathcal{E}}_i $ is given by equation ~\ref{eq3} but shifted by $ \langle -(d_i -1) \rangle \langle -(d_{i-1}) \rangle $ and
$ \widetilde{\mathcal{E}}_i \widetilde{\mathcal{F}}_i $ is given by equation ~\ref{eq4} but shifted by 
$ \langle -(d_{i-1}) \rangle \langle -(d_{i}-1) \rangle. $
		
Hence $ [\widetilde{\mathcal{F}}_i][\widetilde{\mathcal{E}}_i]-[\widetilde{\mathcal{E}}_i][\widetilde{\mathcal{F}}_i] = $
$$ \sum_{l=0}^{d_i} \sum_{k=0}^{d_{i-1}-1}[\mathbb{C}\langle -2k-2l+d_i+d_{i-1}-1 \rangle] -
\sum_{k=0}^{d_{i-1}} \sum_{l=0}^{d_{i}-1}[\mathbb{C}\langle -2k-2l+d_i+d_{i-1}-1 \rangle]. $$
Suppose $ d_i < d_{i-1}. $ Then the above is equal to
\begin{align*}
&\sum_{l=0}^{d_i-1} \sum_{k=0}^{d_{i-1}-1}[\mathbb{C}\langle -2k-2l+d_i+d_{i-1}-1 \rangle] +
\sum_{k=0}^{d_{i-1}-1}[\mathbb{C}\langle -2k-2d_i+d_i+d_{i-1}-1 \rangle]-\\
&\sum_{k=0}^{d_{i-1}-1} \sum_{l=0}^{d_{i}-1}[\mathbb{C}\langle -2k-2l+d_i+d_{i-1}-1 \rangle]-
\sum_{l=0}^{d_{i}-1}[\mathbb{C}\langle -2l-2d_{i-1}+d_i+d_{i-1}-1 \rangle]  = \\
&\sum_{k=0}^{d_{i-1}-1} [\mathbb{C} \langle -2k - 2d_i + d_i + d_{i-1}-1 \rangle] - \sum_{l=0}^{d_{i}-1} [\mathbb{C} \langle -2l - 2d_{i-1} + d_i + d_{i-1}-1 \rangle] =\\
&\sum_{k=0}^{d_{i-1}-1} [\mathbb{C} \langle -2k - d_i + d_{i-1}-1 \rangle] - \sum_{l=0}^{d_{i}-1} [\mathbb{C} \langle -2l - d_{i-1} + d_i -1 \rangle] . 
\end{align*} 
Let $ m= l - d_i + d_{i-1}. $ Then the above is equal to 
\begin{align*}
&\sum_{k=0}^{d_{i-1}-1} [\mathbb{C} \langle -2k - d_i + d_{i-1}-1 \rangle] - 
\sum_{m=d_{i-1}-d_i}^{d_{i-1}-1} [\mathbb{C} \langle -2m - 2d_i +2d_{i-1} +d_i-d_{i-1} -1 \rangle]  =\\
&\sum_{k=0}^{d_{i-1}-1} [\mathbb{C} \langle -2k - d_i + d_{i-1}-1 \rangle] - 
\sum_{m=d_{i-1}-d_i}^{d_{i-1}-1} [\mathbb{C} \langle -2m +d_{i-1} - d_i -1 \rangle]  =\\
&\sum_{k=0}^{d_{i-1}-d_i-1} [\mathbb{C} \langle -2k + d_{i-1} -d_i-1 \rangle] = [d_{i-1}-d_i]. 
\end{align*} 
Thus
$$ \widetilde{\mathcal{F}}_i \widetilde{\mathcal{E}}_i \cong \widetilde{\mathcal{E}}_i \widetilde{\mathcal{F}}_i 
\oplus_{r=0}^{d_{i-1}-d_i-1} \Id \langle d_{i-1}-d_i-1-2r \rangle. $$
		 
Suppose $ d_i = d_{i-1}. $ Then clearly $ \widetilde{\mathcal{F}}_i \widetilde{\mathcal{E}}_i \cong \widetilde{\mathcal{E}}_i \widetilde{\mathcal{F}}_i. $
		 
Finally suppose $ d_i > d_{i-1}. $ As in the case $ d_{i-1} > d_i, $
$$ [\widetilde{\mathcal{E}}_i \widetilde{\mathcal{F}}_i] - [\widetilde{\mathcal{F}}_i \widetilde{\mathcal{E}}_i] = [d_i - d_{i-1}]. $$
Thus
$$ \widetilde{\mathcal{E}}_i \widetilde{\mathcal{F}}_i \cong \widetilde{\mathcal{F}}_i \widetilde{\mathcal{E}}_i 
\oplus_{r=0}^{d_{i}-d_{i-1}-1} \Id \langle d_{i}-d_{i-1}-1-2r \rangle. $$		
\end{enumerate}
\end{proof}

\subsection{Graded Equivalence of Categories}
We would like to prove a graded version of the equivalence given in section 3.
In the ungraded case there is an isomorphism of functors
$$ LZ_{\mathfrak{q}_j}^{\mathfrak{p}_j} \epsilon_{\mathfrak{p}_{j+1}}^{\mathfrak{q}_j}[-(k-1)] LZ_{\mathfrak{q}_j}^{\mathfrak{p}_{j+1}} \epsilon_{\mathfrak{p}_j}^{\mathfrak{q}_j}[-(k-1)] \cong \Id. $$
Thus in the graded case there is an isomorphism for some shift $ a $:
$$ L\widetilde{Z}_{\mathfrak{q}_j}^{\mathfrak{p}_j} \widetilde{\epsilon}_{\mathfrak{p}_{j+1}}^{\mathfrak{q}_j}[-(k-1)] L\widetilde{Z}_{\mathfrak{q}_j}^{\mathfrak{p}_{j+1}} \widetilde{\epsilon}_{\mathfrak{p}_j}^{\mathfrak{q}_j}[-(k-1)] \cong \Id\langle a \rangle. $$
To determine the shift, we compute on a generalized Verma module.

\begin{lemma}
In the notation of section ~\ref{sec3.2},
$ L_s\widetilde{Z}_{\mathfrak{q}_j}^{\mathfrak{p}_{j+1}} \widetilde{\epsilon}_{\mathfrak{p}_j}^{\mathfrak{q}_j} \widetilde{M}^{\mathfrak{p}_j}(\alpha) \cong \widetilde{M}^{\mathfrak{p}_{j+1}}(\beta) $ if $ s=k-1 $ and 0 otherwise.
\end{lemma}

\begin{proof}
We know that it is zero for $ s \neq k-1 $ so we only need to concentrate on the case that $ s=k-1. $
Consider the graded sequences from lemma 5:
\begin{align*}
&0 \rightarrow \widetilde{K}_0 \rightarrow \widetilde{M}^{\mathfrak{q}_j}(\alpha) \rightarrow \widetilde{M}^{\mathfrak{p}_j}(\alpha) \rightarrow 0.\\
&0 \rightarrow \widetilde{K}_1 \rightarrow \widetilde{M}^{\mathfrak{q}_j}(\sigma_1. \alpha) \langle 1 \rangle \rightarrow \widetilde{K}_0 \rightarrow 0.\\
&0 \rightarrow \widetilde{K}_2 \rightarrow \widetilde{M}^{\mathfrak{q}_j}(\sigma_1 \sigma_2 .\alpha) \langle 2 \rangle \rightarrow \widetilde{K}_1 \rightarrow 0.\\
&\cdots\\
&0 \rightarrow \widetilde{K}_{k-2} \rightarrow \widetilde{M}^{\mathfrak{q}_j}(\sigma_1 \cdots \sigma_{k-2}.\alpha) \langle k-2 \rangle \rightarrow \widetilde{K}_{k-3} \rightarrow 0.\\
&0 \rightarrow \widetilde{K}_{k-1} \rightarrow \widetilde{M}^{\mathfrak{q}_j}(\sigma_1 \cdots \sigma_{k-2} \sigma_{k-1}. \alpha) \langle k-1 \rangle \rightarrow \widetilde{K}_{k-2} \rightarrow 0. 
\end{align*} 

As in the ungraded we get
\begin{align*}
L_s \widetilde{Z}_{\mathfrak{q}_j}^{\mathfrak{p}_j} \widetilde{\epsilon}_{\mathfrak{p}_{j+1}}^{\mathfrak{q}_j} \widetilde{M}^{\mathfrak{q}_j}(\sigma_1 \cdots \sigma_l.\alpha) \langle l \rangle \cong
L_s \widetilde{Z}_{\mathfrak{q}_j}^{\mathfrak{p}_j} \widetilde{\epsilon}_{\mathfrak{p}_{j+1}}^{\mathfrak{q}_j} \widetilde{K}_{l-1} \cong 
L_{s+1} \widetilde{Z}_{\mathfrak{q}_j}^{\mathfrak{p}_j} \widetilde{\epsilon}_{\mathfrak{p}_{j+1}}^{\mathfrak{q}_j} \widetilde{K}_{l-2} \cong \cdots &\cong
L_{s+l-1} \widetilde{Z}_{\mathfrak{q}_j}^{\mathfrak{p}_j} \widetilde{\epsilon}_{\mathfrak{p}_{j+1}}^{\mathfrak{q}_j} \widetilde{K}_{0}\\ 
&\cong
L_{s+l} \widetilde{Z}_{\mathfrak{q}_j}^{\mathfrak{p}_j} \widetilde{\epsilon}_{\mathfrak{p}_{j+1}}^{\mathfrak{q}_j} \widetilde{M}^{\mathfrak{p}_j}(\alpha). 
\end{align*}

We know 
$ L_s \widetilde{Z}_{\mathfrak{q}_j}^{\mathfrak{p}_j} \widetilde{\epsilon}_{\mathfrak{p}_{j+1}}^{\mathfrak{q}_j} \widetilde{M}^{\mathfrak{q}_j}(\sigma_1 \cdots \sigma_l.\alpha) \langle l \rangle \cong 0 $ if $ s \neq k-l-1. $
If $ s=k-l-1, $ it is isomorphic to 
$$ \widetilde{M}^{\mathfrak{p}_{j+1}}(\beta) \langle l \rangle \langle k-l-1 \rangle \langle -(k-1) \rangle \cong \widetilde{M}^{\mathfrak{p}_{j+1}}(\beta). $$
The second shift comes from the fact the $ \widetilde{M}^{\mathfrak{p}_{j+1}}(\beta) $ is a Koszul module.  The third shift comes from the grading convention for the lift
of the dual Zuckerman functor.
Thus there is no shift overall.
\end{proof}

We let $ \widetilde{\Delta} $ denote a composition of graded equivalences.

\subsection{Graded Diagram Relation 1}
It suffices to repeat the calculation of lemma ~\ref{lemmadiagram1} with the grading.

\begin{lemma}
If $ i=2r, $ 
$$ L_i \widetilde{Z}_{\mathfrak{q}_j}^{\mathfrak{p}_j} \widetilde{\epsilon}_{\mathfrak{p}_{j+1}}^{\mathfrak{q}_j} \widetilde{M}^{\mathfrak{p}_j}(\alpha) \cong \widetilde{M}^{\mathfrak{p}_j}(\alpha) \langle 2r+1-k \rangle. $$
\end{lemma}

\begin{proof}
From the graded exact sequences of the previous subsection we easily see that
$$ L_i \widetilde{Z}_{\mathfrak{q}_j}^{\mathfrak{p}_j} \widetilde{\epsilon}_{\mathfrak{p}_{j+1}}^{\mathfrak{q}_j} \widetilde{M}^{\mathfrak{p}_j}(\alpha) \cong \cdots \cong
L_{i-r} \widetilde{Z}_{\mathfrak{q}_j}^{\mathfrak{p}_j} \widetilde{M}^{\mathfrak{q}_j}(\sigma_1 \cdots \sigma_{r}.\alpha) \langle r \rangle. $$
Now, $ \widetilde{M}^{\mathfrak{q}_j}(\sigma_1 \cdots \sigma_{r}.\alpha) $ is a Koszul module and has a projective resolution whose $r$th component is generated in degree $r.$
Thus 
$$ L_{i-r} \widetilde{Z}_{\mathfrak{q}_j}^{\mathfrak{p}_j} \widetilde{\epsilon}_{\mathfrak{p}_{j+1}}^{\mathfrak{q}_j} \widetilde{M}^{\mathfrak{q}_j}(\sigma_1 \cdots \sigma_{r}.\alpha) \langle r \rangle \cong
\widetilde{M}^{\mathfrak{p}_j}(\alpha) \langle r \rangle \langle r \rangle \langle -(k-1) \rangle \cong 
\widetilde{M}^{\mathfrak{p}_j}(\alpha) \langle 2r+1-k \rangle. $$
The second shift comes from the fact that the module is Koszul, and the third shift comes from the definition of the functor.
\end{proof}

\begin{corollary}
There is an isomorphism of graded functors:
$$ L\widetilde{Z}_{\mathfrak{q}_j}^{\mathfrak{p}_j} \widetilde{\epsilon}_{\mathfrak{p}_j}^{\mathfrak{q}_j}[-(k-1)] \cong \oplus_{r=0}^{k-1} \Id[2r-(k-1)] \langle 2r-(k-1) \rangle. $$
\end{corollary}

\subsection{Graded Diagram Relation 2}
In this subsection we include only the $ i $ and $ (i+1)\text{st} $ entries of the weights for the Verma modules.

\begin{lemma}
There is an isomorphism, of graded objects
$$ L\widetilde{Z}^{\mathfrak{s}_i} \widetilde{\epsilon}_{\mathfrak{s}_i} \widetilde{M}^{\mathfrak{s}_i}(\underbrace{a_i, a_{i+1}}) \cong
\widetilde{M}^{\mathfrak{s}_i}(\underbrace{a_i, a_{i+1}}) \langle -1 \rangle \oplus \widetilde{M}^{\mathfrak{s}_i}(\underbrace{a_i, a_{i+1}})[2] \langle 1 \rangle. $$
\end{lemma}

\begin{proof}
Consider the object $ \widetilde{M}^{\mathfrak{s}_i}(\underbrace{a_i, a_{i+1}}). $  There is an exact sequence
$$ 0 \rightarrow \widetilde{M}(a_{i+1}, a_{i}) \langle 1 \rangle \rightarrow \widetilde{M}(a_i, a_{i+1})
\rightarrow \widetilde{\epsilon}_{\mathfrak{s}_i} \widetilde{M}^{\mathfrak{s}_i}(\underbrace{a_i, a_{i+1}}) \rightarrow 0. $$
Now apply the functor $ L\widetilde{Z}^{\mathfrak{s}_i} $ to get a long exact sequence.
It follows that 
$$ L_0 \widetilde{Z}^{\mathfrak{s}_i} \widetilde{\epsilon}_{\mathfrak{s}_i} \widetilde{M}^{\mathfrak{s}_i}(\underbrace{a_i, a_{i+1}}) \cong 
\widetilde{M}^{\mathfrak{s}_i}(\underbrace{a_i, a_{i+1}}) \langle -1 \rangle. $$
The shift comes from the convention taken for the functor.
Thus,
$$ L_2 \widetilde{Z}^{\mathfrak{s}_i} \widetilde{\epsilon}_{\mathfrak{s}_i} \widetilde{M}^{\mathfrak{s}_i}(\underbrace{a_i, a_{i+1}}) \cong 
\widetilde{M}^{\mathfrak{s}_i}(\underbrace{a_i, a_{i+1}}) \langle -1 \rangle \langle 1 \rangle \langle 1 \rangle. $$
The second shift comes from the exact sequence and the third shift comes from that the module is Koszul..
Therefore
$$ L\widetilde{Z}^{\mathfrak{s}_i} \widetilde{\epsilon}_{\mathfrak{s}_i} \widetilde{M}^{\mathfrak{s}_i}(\underbrace{a_i, a_{i+1}}) \cong
\widetilde{M}^{\mathfrak{s}_i}(\underbrace{a_i, a_{i+1}}) \langle -1 \rangle \oplus \widetilde{M}^{\mathfrak{s}_i}(\underbrace{a_i, a_{i+1}}) [2] \langle 1 \rangle. $$
\end{proof}

\begin{corollary}
There is an isomorphism of graded functors:
$ L\widetilde{Z}^{\mathfrak{s}_i} \widetilde{\epsilon}_{\mathfrak{s}_i}[-1] \cong \Id[-1] \langle -1 \rangle \oplus \Id[1] \langle 1 \rangle. $
\end{corollary}

\subsection{Graded Diagram Relation 3}
As in section ~\ref{sec4.3}, we study the graded generalized Verma module of highest weight $ \alpha, $ $ \widetilde{M}^{\mathfrak{p}_{i+1}}(a_1, \ldots, a_i, k-1, \ldots, 0, a_{i+k}, \ldots, a_n). $

\begin{lemma}
Let $ a_i = l. $  
Suppose $ s=0 $ or $ s=2. $ Then
\begin{align*}
&L_{s}\widetilde{Z}_{\mathfrak{q}_{i+1}}^{\mathfrak{s}_i+\mathfrak{q}_{i+1}} \widetilde{\epsilon}_{\mathfrak{p}_{i+1}}^{\mathfrak{q}_{i+1}} \widetilde{M}^{\mathfrak{p}_{i+1}}(a_i, k-1, \ldots, 0) = 0.\\
&L_{1}\widetilde{Z}_{\mathfrak{q}_{i+1}}^{\mathfrak{s}_i+\mathfrak{q}_{i+1}} \widetilde{\epsilon}_{\mathfrak{p}_{i+1}}^{\mathfrak{q}_{i+1}} \widetilde{M}^{\mathfrak{p}_{i+1}}(a_i, k-1, \ldots, 0, a_{i+k}) \cong
L_1\widetilde{Z}_{\mathfrak{q}_{i+1}}^{\mathfrak{s}_i+\mathfrak{q}_{i+1}} \widetilde{M}^{\mathfrak{q}_{i+1}}(e.\alpha)/L_1\widetilde{Z}_{\mathfrak{q}_{i+1}}^{\mathfrak{s}_i+\mathfrak{q}_{i+1}} \widetilde{K}_0\\
&L_{1} \widetilde{Z}_{\mathfrak{q}_{i+1}}^{\mathfrak{s}_i+\mathfrak{q}_{i+1}} \widetilde{K}_0 \cong L_{1} \widetilde{Z}_{\mathfrak{q}_{i+1}}^{\mathfrak{s}_i+\mathfrak{q}_{i+1}} \widetilde{M}^{\mathfrak{q}_{i+1}}(\sigma_1.\alpha)\langle 1 \rangle/ 
L_1 \widetilde{Z}_{\mathfrak{q}_{i+1}}^{\mathfrak{s}_i+\mathfrak{q}_{i+1}} \widetilde{K}_1\\
&\cdots\\
&L_{1}\widetilde{Z}_{\mathfrak{q}_{i+1}}^{\mathfrak{s}_i+\mathfrak{q}_{i+1}} \widetilde{K}_{k-l-3} \cong L_{1} \widetilde{Z}_{\mathfrak{q}_{i+1}}^{\mathfrak{s}_i+\mathfrak{q}_{i+1}} \widetilde{M}^{\mathfrak{q}_{i+1}}(\sigma_1 \cdots \sigma_{k-l-2}.\alpha)\langle k-l-2 \rangle /
L_1 \widetilde{Z}_{\mathfrak{q}_{i+1}}^{\mathfrak{s}_i+\mathfrak{q}_{i+1}} \widetilde{K}_{k-l-2}\\
&L_{s} \widetilde{Z}_{\mathfrak{q}_{i+1}}^{\mathfrak{s}_i+\mathfrak{q}_{i+1}} \widetilde{K}_{k-l-2} \cong L_{s-1} \widetilde{Z}_{\mathfrak{q}_{i+1}}^{\mathfrak{s}_i+\mathfrak{q}_{i+1}} \widetilde{K}_{k-l-1}\\
&L_{0} \widetilde{Z}_{\mathfrak{q}_{i+1}}^{\mathfrak{s}_i+\mathfrak{q}_{i+1}} \widetilde{K}_{k-l-1} \cong L_{0} \widetilde{Z}_{\mathfrak{q}_{i+1}}^{\mathfrak{s}_i+\mathfrak{q}_{i+1}} \widetilde{M}^{\mathfrak{q}_{i+1}}(\sigma_1 \cdots \sigma_{k-l}.\alpha) \langle k-l \rangle /
L_{0} \widetilde{Z}_{\mathfrak{q}_{i+1}}^{\mathfrak{s}_i+\mathfrak{q}_{i+1}} \widetilde{K}_{k-l}\\
&L_{0} \widetilde{Z}_{\mathfrak{q}_{i+1}}^{\mathfrak{s}_i+\mathfrak{q}_{i+1}} \widetilde{K}_{k-l} \cong L_{0} \widetilde{Z}_{\mathfrak{q}_{i+1}}^{\mathfrak{s}_i+\mathfrak{q}_{i+1}} \widetilde{M}^{\mathfrak{q}_{i+1}}(\sigma_1 \cdots \sigma_{k-l+1}.\alpha)\langle k-l+1 \rangle/
L_{0} \widetilde{Z}_{\mathfrak{q}_{i+1}}^{\mathfrak{s}_i+\mathfrak{q}_{i+1}} \widetilde{K}_{k-l+1}\\
&\cdots\\
&L_{0} \widetilde{Z}_{\mathfrak{q}_{i+1}}^{\mathfrak{s}_i+\mathfrak{q}_{i+1}} \widetilde{K}_{k-4} \cong L_{0} \widetilde{Z}_{\mathfrak{q}_{i+1}}^{\mathfrak{s}_i+\mathfrak{q}_{i+1}} \widetilde{M}^{\mathfrak{q}_{i+1}}(\sigma_1 \cdots \sigma_{k-3}.\alpha)\langle k-3 \rangle /
L_{0} \widetilde{Z}_{\mathfrak{q}_{i+1}}^{\mathfrak{s}_i+\mathfrak{q}_{i+1}} \widetilde{K}_{k-3}\\
&L_{0} \widetilde{Z}_{\mathfrak{q}_{i+1}}^{\mathfrak{s}_i+\mathfrak{q}_{i+1}} \widetilde{K}_{k-3} \cong L_{0} \widetilde{Z}_{\mathfrak{q}_{i+1}}^{\mathfrak{s}_i+\mathfrak{q}_{i+1}} \widetilde{M}^{\mathfrak{q}_{i+1}}(\sigma_1 \cdots \sigma_{k-2}.\alpha)\langle k-2 \rangle /
L_{0} \widetilde{Z}_{\mathfrak{q}_{i+1}}^{\mathfrak{s}_i+\mathfrak{q}_{i+1}} \widetilde{M}^{\mathfrak{q}_{i+1}}(\sigma_1 \cdots \sigma_{k-1}.\alpha)\langle k-1 \rangle. 
\end{align*} 
\end{lemma}

\begin{proof}
This easily follows from the proof of lemma ~\ref{lemma10}.
\end{proof}

\begin{lemma}
Suppose $ s \neq k-2, k. $  Then
$$ L_s \widetilde{Z}_{\mathfrak{q}_{i+1}}^{\mathfrak{p}_{i+1}} \widetilde{\epsilon}_{\mathfrak{s}_i+\mathfrak{q}_{i+1}}^{\mathfrak{q}_{i+1}}(L_0 \widetilde{Z}_{\mathfrak{q}_{i+1}}^{\mathfrak{s}_i+\mathfrak{q}_{i+1}}
\widetilde{M}^{\mathfrak{q}_{i+1}}(\sigma_1 \cdots \sigma_{k-2}.\alpha)\langle k-2 \rangle / L_0 \widetilde{Z}_{\mathfrak{q}_{i+1}}^{\mathfrak{s}_i+\mathfrak{q}_{i+1}} \widetilde{M}^{\mathfrak{q}_{i+1}}(\sigma_1 \cdots \sigma_{k-1}.\alpha)\langle k-1 \rangle) \cong 0. $$
For $ s=k-2, $ it is isomorphic to $ \widetilde{M}^{\mathfrak{p}_{i+1}}(a_i, k-1, \ldots, 0) \langle k-4 \rangle. $
For $ s=k, $ it is isomorphic to $ \widetilde{M}^{\mathfrak{p}_{i+1}}(a_i, k-1, \ldots, 0) \langle k-2 \rangle. $
\end{lemma}  

\begin{proof}
Note that if $ s=k-2, $
\begin{align*}
&L_s \widetilde{Z}_{\mathfrak{q}_{i+1}}^{\mathfrak{p}_{i+1}} \widetilde{\epsilon}_{\mathfrak{s}_i+\mathfrak{q}_{i+1}}^{\mathfrak{q}_{i+1}}
L_0 \widetilde{Z}_{\mathfrak{q}_{i+1}}^{\mathfrak{s}_i+\mathfrak{q}_{i+1}} \widetilde{M}^{\mathfrak{q}_{i+1}}(\sigma_1 \cdots \sigma_{k-2}.\alpha)) \langle k-2 \rangle \cong\\
&\widetilde{M}^{\mathfrak{p}_{i+1}}(a_i, k-1, \ldots, 0) \langle (k-2) \rangle \langle -1 \rangle \langle (k-2) \rangle  \langle -(k-1) \rangle. 
\end{align*} 
The second shift comes from the innermost functor.  The third shift comes from Koszulness of the module.  The fourth shift comes from the outermost functor.
For all other s it is zero.

Also if $ s = k-1, $
\begin{align*}
&L_s \widetilde{Z}_{\mathfrak{q}_{i+1}}^{\mathfrak{p}_{i+1}} \widetilde{\epsilon}_{\mathfrak{s}_i+\mathfrak{q}_{i+1}}^{\mathfrak{q}_{i+1}}
L_0 \widetilde{Z}_{\mathfrak{q}_{i+1}}^{\mathfrak{s}_i+\mathfrak{q}_{i+1}} \widetilde{M}^{\mathfrak{q}_{i+1}}(\sigma_1 \cdots \sigma_{k-1}.\alpha)) \langle k-1 \rangle \cong\\
&\widetilde{M}^{\mathfrak{p}_{i+1}}(a_i, k-1, \ldots, 0) \langle (k-1) \rangle \langle -1 \rangle \langle (k-1) \rangle  \langle -(k-1) \rangle. 
\end{align*} 
For all other s it is zero.
Then the lemma follows from the of lemma ~\ref{lemma11}.
\end{proof}

\begin{lemma}
If $ s = 1, 3, \ldots, 2(k-1)-1, $ then 
$$ L_s \widetilde{Z}_{\mathfrak{q}_{i+1}}^{\mathfrak{p}_{i+1}} \widetilde{\epsilon}_{\mathfrak{s}_i+\mathfrak{q}_{i+1}}^{\mathfrak{q}_{i+1}} L_1 \widetilde{Z}_{\mathfrak{q}_{i+1}}^{\mathfrak{s}_i+\mathfrak{q}_{i+1}}
\widetilde{M}^{\mathfrak{p}_{i+1}}(a_i, k-1, \ldots, 0) \cong 
\widetilde{M}^{\mathfrak{p}_{i+1}}(a_i, k-1, \ldots, 0) \langle -k+1+s \rangle. $$
Otherwise it is zero.
\end{lemma}

\begin{proof}
This follows from the previous lemma and following the exact sequences of lemmas ~\ref{lemma12}, ~\ref{lemma13}, ~\ref{lemma14} and corollaries ~\ref{corollary5} and ~\ref{corollary6}.
\end{proof}

The graded version of proposition ~\ref{prop10} now holds.

\begin{prop}
\label{prop13}
There is an isomorphism of graded functors
$$ L\widetilde{Z}_{\mathfrak{q}_{i+1}}^{\mathfrak{p}_{i+1}} \widetilde{\epsilon}_{\mathfrak{s}_i+\mathfrak{q}_{i+1}}^{\mathfrak{q}_{i+1}}[-1]  L\widetilde{Z}_{\mathfrak{q}_{i+1}}^{\mathfrak{s}_i+\mathfrak{q}_{i+1}}
\widetilde{\epsilon}_{\mathfrak{p}_{i+1}}^{\mathfrak{q}_{i+1}}[-(k-1)] \cong \oplus^{k-1}_{r=1} \Id [2r-k] \langle 2r-k \rangle. $$
\end{prop}

\subsection{Graded Diagram Relation 4}

We begin with a generalized Verma module
$$ \widetilde{M}^{\alpha}(\underbrace{k-1, \ldots, 0}, a_{i+k}, \underbrace{a_{i+k+1}, \ldots, a_{i+2k-1}}). $$
Suppose $ a_{i+k} = l. $
Let us denote the module above by $ \widetilde{M}^{\alpha}(l). $

There are short exact sequences coming from the generalized BGG resolution:

$$ 0 \rightarrow \widetilde{K}_0 \rightarrow \widetilde{M}^{\beta}(\underbrace{k-1, \ldots, \hat{0}}, 0, a_{i+k}, \underbrace{a_{i+k+1}, \ldots, a_{i+2k-1}})\langle 0 \rangle \rightarrow 
\widetilde{M}^{\alpha}(l) \rightarrow 0 $$
$$ 0 \rightarrow \widetilde{K}_1 \rightarrow \widetilde{M}^{\beta}(\underbrace{k-1, \ldots, \hat{1}, 0}, 1, a_{i+k}, \underbrace{a_{i+k+1}, \ldots, a_{i+2k-1}})\langle 1 \rangle \rightarrow 
\widetilde{K}_0 \rightarrow 0 $$
$$ \cdots $$
$$ 0 \rightarrow \widetilde{K}_{l-1} \rightarrow \widetilde{M}^{\beta}(\underbrace{k-1, \ldots, \widehat{l-1}, 0}, l-1, a_{i+k}, \underbrace{a_{i+k+1}, \ldots, a_{i+2k-1}})\langle l-1 \rangle \rightarrow \widetilde{K}_{l-2} \rightarrow 0 $$
$$ 0 \rightarrow \widetilde{K}_{l} \rightarrow \widetilde{M}^{\beta}(\underbrace{k-1, \ldots, \hat{l}, 0}, l, a_{i+k}, \underbrace{a_{i+k+1}, \ldots, a_{i+2k-1}})\langle l \rangle \rightarrow 
\widetilde{K}_{l-1} \rightarrow 0 $$
$$ 0 \rightarrow \widetilde{K}_{l+1} \rightarrow \widetilde{M}^{\beta}(\underbrace{k-1, \ldots, \widehat{l+1}, 0}, l+1, a_{i+k}, \underbrace{a_{i+k+1}, \ldots, a_{i+2k-1}})\langle l+1 \rangle \rightarrow \widetilde{K}_{l} \rightarrow 0 $$
$$ \cdots $$
$$ 0 \rightarrow \widetilde{K}_{k-1} \rightarrow \widetilde{M}^{\beta}(\underbrace{\widehat{k-1}, \ldots, 0}, k-1, a_{i+k}, \underbrace{a_{i+k+1}, \ldots, a_{i+2k-1}})\langle k-1 \rangle \rightarrow \widetilde{K}_{k-2} \rightarrow 0. $$

As in the ungraded case, the following lemma follows immediately.

\begin{lemma}
There are exact sequences:
$$ 0 \rightarrow L_1\widetilde{Z}_{\beta}^{\gamma} \widetilde{K}_0 \rightarrow \widetilde{M}^{\gamma}(\underbrace{k-1, \ldots, \hat{0}}, \underbrace{a_{i+k},0,} \underbrace{a_{i+k+1}, \ldots, a_{i+2k-1}})\langle 0 \rangle \rightarrow 
L_1\widetilde{Z}_{\beta}^{\alpha} \widetilde{\epsilon}_{\alpha}^{\beta} \widetilde{M}^{\alpha}(l) \rightarrow 0 $$
$$ 0 \rightarrow L_1\widetilde{Z}_{\beta}^{\gamma} \widetilde{K}_1 \rightarrow \widetilde{M}^{\gamma}(\underbrace{k-1, \ldots, \hat{1}, 0}, \underbrace{a_{i+k},1,} \underbrace{a_{i+k+1}, \ldots, a_{i+2k-1}})\langle 1 \rangle \rightarrow 
L_1 \widetilde{Z}_{\beta}^{\gamma} \widetilde{K}_0 \rightarrow 0 $$
$$ \cdots $$
$$ 0 \rightarrow L_1\widetilde{Z}_{\beta}^{\gamma} \widetilde{K}_{l-1} \rightarrow \widetilde{M}^{\gamma}(\underbrace{k-1, \ldots, \widehat{l-1}, 0}, \underbrace{a_{i+k},l-1,} \underbrace{a_{i+k+1}, \ldots, a_{i+2k-1}})\langle l-1 \rangle \rightarrow L_1 \widetilde{Z}_{\beta}^{\gamma} \widetilde{K}_{l-2} \rightarrow 0 $$
$$ L_i \widetilde{Z}_{\beta}^{\gamma} \widetilde{K}_l = L_{i+1}\widetilde{Z}_{\beta}^{\gamma} \widetilde{K}_{l-1}, \forall i $$
$$ 0 \rightarrow L_0 \widetilde{Z}_{\beta}^{\gamma} \widetilde{K}_{l+1} \rightarrow \widetilde{M}^{\gamma}(\underbrace{k-1, \ldots, \widehat{l+1}, 0}, \underbrace{l+1, a_{i+k},} \underbrace{a_{i+k+1}, \ldots, a_{i+2k-1}})\langle l \rangle \rightarrow L_0 \widetilde{Z}_{\beta}^{\gamma} \widetilde{K}_l \rightarrow 0 $$
$$ \cdots $$
$$ 0 \rightarrow L_0 \widetilde{Z}_{\beta}^{\gamma} \widetilde{K}_{k-1} \rightarrow \widetilde{M}^{\gamma}(\underbrace{\widehat{k-1}, \ldots, 0}, \underbrace{k-1, a_{i+k},} \underbrace{a_{i+k+1}, \ldots, a_{i+2k-1}})\langle k-2 \rangle \rightarrow L_0 \widetilde{Z}_{\beta}^{\gamma} \widetilde{K}_{k-2} \rightarrow 0 $$
\end{lemma}

\begin{proof}
The shifts in the grading come from the definition of the functors and Koszulness of the modules.
\end{proof}

\begin{corollary}
\label{cor12}
There are exact sequences:
$$ 0 \rightarrow \widetilde{\epsilon}_{\gamma}^{\beta}L_1\widetilde{Z}_{\beta}^{\gamma} \widetilde{K}_0 \rightarrow \widetilde{\epsilon}_{\gamma}^{\beta}\widetilde{M}^{\gamma}(\underbrace{k-1, \ldots, \hat{0}}, \underbrace{a_{i+k},0,} \underbrace{a_{i+k+1}, \ldots, a_{i+2k-1}})\langle 0 \rangle \rightarrow 
\widetilde{\epsilon}_{\gamma}^{\beta} L_1\widetilde{Z}_{\beta}^{\alpha} \widetilde{\epsilon}_{\alpha}^{\beta} \widetilde{M}^{\alpha}(l) \rightarrow 0 $$
$$ 0 \rightarrow \widetilde{\epsilon}_{\gamma}^{\beta} L_1\widetilde{Z}_{\beta}^{\gamma} \widetilde{K}_1 \rightarrow \widetilde{\epsilon}_{\gamma}^{\beta} \widetilde{M}^{\gamma}(\underbrace{k-1, \ldots, \hat{1}, 0}, \underbrace{a_{i+k},1,} \underbrace{a_{i+k+1}, \ldots, a_{i+2k-1}})\langle 1 \rangle \rightarrow 
\widetilde{\epsilon}_{\gamma}^{\beta} L_1 \widetilde{Z}_{\beta}^{\gamma} \widetilde{K}_0 \rightarrow 0 $$
$$ \cdots $$
$$ 0 \rightarrow \widetilde{\epsilon}_{\gamma}^{\beta}L_1\widetilde{Z}_{\beta}^{\gamma} \widetilde{K}_{l-1} \rightarrow \widetilde{\epsilon}_{\gamma}^{\beta} \widetilde{M}^{\gamma}(\underbrace{k-1, \ldots, \widehat{l-1}, 0}, \underbrace{a_{i+k},l-1,} \underbrace{a_{i+k+1}, \ldots, a_{i+2k-1}})\langle l-1 \rangle \rightarrow \widetilde{\epsilon}_{\gamma}^{\beta} L_1 \widetilde{Z}_{\beta}^{\gamma} \widetilde{K}_{l-2} \rightarrow 0 $$
$$ \widetilde{\epsilon}_{\gamma}^{\beta} L_i \widetilde{Z}_{\beta}^{\gamma} \widetilde{K}_l \cong \widetilde{\epsilon}_{\gamma}^{\beta}L_{i+1}\widetilde{Z}_{\beta}^{\gamma} \widetilde{K}_{l-1}, \forall i $$
$$ 0 \rightarrow \widetilde{\epsilon}_{\gamma}^{\beta} L_0 \widetilde{Z}_{\beta}^{\gamma} \widetilde{K}_{l+1} \rightarrow \widetilde{\epsilon}_{\gamma}^{\beta} \widetilde{M}^{\gamma}(\underbrace{k-1, \ldots, \widehat{l+1}, 0}, \underbrace{l+1, a_{i+k},} \underbrace{a_{i+k+1}, \ldots, a_{i+2k-1}})\langle l \rangle \rightarrow \widetilde{\epsilon}_{\gamma}^{\beta} L_0 \widetilde{Z}_{\beta}^{\gamma} \widetilde{K}_l \rightarrow 0 $$
$$ \cdots $$
$$ 0 \rightarrow \widetilde{\epsilon}_{\gamma}^{\beta}L_0 \widetilde{Z}_{\beta}^{\gamma} \widetilde{K}_{k-1} \rightarrow \widetilde{\epsilon}_{\gamma}^{\beta} \widetilde{M}^{\gamma}(\underbrace{\widehat{k-1}, \ldots, 0}, \underbrace{k-1, a_{i+k},} \underbrace{a_{i+k+1}, \ldots, a_{i+2k-1}})\langle k-2 \rangle \rightarrow \widetilde{\epsilon}_{\gamma}^{\beta} L_0 \widetilde{Z}_{\beta}^{\gamma} \widetilde{K}_{k-2} \rightarrow 0 $$
\end{corollary}

\begin{lemma}
Let $ \lbrace a_{i+k+1}, \ldots, a_{i+2k-1} \rbrace = \lbrace k-1, \ldots, \hat{m}, \ldots, 0 \rbrace. $
Let $ m \neq l. $ If $ s=k-1, $ then
$$ L_s \widetilde{Z}_{\beta}^{\delta} \widetilde{\epsilon}_{\gamma}^{\beta} L_1 \widetilde{Z}_{\beta}^{\gamma} \widetilde{\epsilon}_{\alpha}^{\beta} \widetilde{M}^{\alpha}(l) =
\widetilde{M}^{\delta}(\underbrace{k-1, \ldots, \hat{m}, \ldots, 0,} a_{i+k}, \underbrace{k-1, \ldots, 0})\langle 0 \rangle. $$
Otherwise it is zero.
\end{lemma}

\begin{proof}
This follows easily from the previous corollary and the Koszul property of the modules.
\end{proof}

Now consider $ \widetilde{M}^{\delta}(\underbrace{a_i, \ldots, a_{i+k-2}}, l, \underbrace{k-1, \ldots, 0})=\widetilde{M}^{\delta}(l). $
There is a graded resolution for this module as there was for $ \widetilde{M}^{\alpha}(l). $ Similarly we get the following lemma.

\begin{lemma}
\label{lemma46}
There are exact sequences:
$$ 0 \rightarrow \widetilde{\epsilon}_{\gamma}^{\beta} L_1 \widetilde{Z}_{\beta}^{\gamma} \widetilde{J}_0 \rightarrow \widetilde{\epsilon}_{\gamma}^{\beta} \widetilde{M}^{\gamma}(\underbrace{a_i, \ldots, a_{i+k-2},} \underbrace{k-1, a_{i+k-1}}, \underbrace{\widehat{k-1}, \ldots, 0})\langle 0 \rangle \rightarrow \widetilde{\epsilon}_{\gamma}^{\beta} L_1 \widetilde{Z}_{\beta}^{\gamma} \widetilde{M}^{\delta}(l) \rightarrow 0 $$
$$ 0 \rightarrow \widetilde{\epsilon}_{\gamma}^{\beta} L_1 \widetilde{Z}_{\beta}^{\gamma} \widetilde{J}_1 \rightarrow \widetilde{\epsilon}_{\gamma}^{\beta} \widetilde{M}^{\gamma}(\underbrace{a_i, \ldots, a_{i+k-2},} \underbrace{k-2, a_{i+k-1}}, \underbrace{k-1, \widehat{k-2}, \ldots, 0})\langle 1 \rangle  \rightarrow \widetilde{\epsilon}_{\gamma}^{\beta} L_1 \widetilde{Z}_{\beta}^{\gamma} \widetilde{J}_0 \rightarrow 0 $$
$$ \cdots $$
$$ 0 \rightarrow \widetilde{\epsilon}_{\gamma}^{\beta} L_1 \widetilde{Z}_{\beta}^{\gamma} \widetilde{J}_{k-l-2} \rightarrow \widetilde{\epsilon}_{\gamma}^{\beta} \widetilde{M}^{\gamma}(\underbrace{a_i, \ldots, a_{i+k-2},} \underbrace{l+1, a_{i+k-1}}, \underbrace{k-1, \widehat{l+1}, \ldots, 0})\langle k-l-2 \rangle \rightarrow \widetilde{\epsilon}_{\gamma}^{\beta} L_1 \widetilde{Z}_{\beta}^{\gamma} \widetilde{J}_{k-l-3} \rightarrow 0 $$
$$ \widetilde{\epsilon}_{\gamma}^{\beta} L_i \widetilde{Z}_{\beta}^{\gamma} \widetilde{J}_{k-l-1} \cong \widetilde{\epsilon}_{\gamma}^{\beta} L_{i+1} \widetilde{Z}_{\beta}^{\gamma} \widetilde{J}_{k-l-2}, \forall i $$
$$ 0 \rightarrow \widetilde{\epsilon}_{\gamma}^{\beta} L_0 \widetilde{Z}_{\beta}^{\gamma} \widetilde{J}_{k-l} \rightarrow \widetilde{\epsilon}_{\gamma}^{\beta} \widetilde{M}^{\gamma}(\underbrace{a_i, \ldots, a_{i+k-2},} \underbrace{a_{i+k-1}, l-1}, \underbrace{k-1, \widehat{l-1}, \ldots, 0})\langle k-l-1 \rangle \rightarrow \widetilde{\epsilon}_{\gamma}^{\beta} L_1 \widetilde{Z}_{\beta}^{\gamma} \widetilde{J}_{k-l-1} \rightarrow 0 $$
$$ \cdots $$
$$ 0 \rightarrow \widetilde{\epsilon}_{\gamma}^{\beta} L_0 \widetilde{Z}_{\beta}^{\gamma} \widetilde{J}_{k-1} \rightarrow \widetilde{\epsilon}_{\gamma}^{\beta} \widetilde{M}^{\gamma}(\underbrace{a_i, \ldots, a_{i+k-2},} \underbrace{a_{i+k-1}, 0}, \underbrace{k-1, \ldots, \hat{0}})\langle k-2 \rangle \rightarrow \widetilde{\epsilon}_{\gamma}^{\beta} L_0 \widetilde{Z}_{\beta}^{\gamma} \widetilde{J}_{k-2} \rightarrow 0 $$
\end{lemma}

\begin{lemma}
Let $ \lbrace a_{i}, \ldots, a_{i+k-2} \rbrace = \lbrace k-1, \ldots, \hat{m}, \ldots, 0 \rbrace. $
Let $ m \neq l. $ If $ s=k-1, $ then
$$ L_s \widetilde{Z}_{\beta}^{\alpha} \widetilde{\epsilon}_{\gamma}^{\beta} L_1 \widetilde{Z}_{\beta}^{\gamma} \widetilde{\epsilon}_{\delta}^{\beta} \widetilde{M}^{\delta}(l) =
\widetilde{M}^{\alpha}(\underbrace{k-1, \ldots, 0,} a_{i+k-1}, \underbrace{k-1, \ldots, \hat{m}, \ldots,  0}) \langle 0 \rangle. $$
Otherwise it is zero.
\end{lemma}

\begin{proof}
This follows from the previous lemma.
\end{proof}

We must now study the case when $ l = m. $

\begin{lemma}
\label{lemma48}
Let $ m=l. $  Then
$ L_s \widetilde{Z}_{\beta}^{\delta} \widetilde{\epsilon}_{\gamma}^{\beta} L_1 \widetilde{Z}_{\beta}^{\gamma} \widetilde{\epsilon}_{\alpha}^{\beta} \widetilde{M}^{\alpha}(l) \cong $
\begin{eqnarray*} 
\widetilde{M}^{\delta}(0)\langle -l+1 \rangle /\ldots /\widetilde{M}^{\delta}(k-1)\langle k-l \rangle & \textrm{ if } & s=k-l \\
L_{k-l+1}\widetilde{Z}_{\beta}^{\gamma} \widetilde{\epsilon}_{\gamma}^{\beta}L_1 \widetilde{Z}_{\beta}^{\gamma} \widetilde{K}_{l-3}  & \textrm{ if } & s=k-1
\end{eqnarray*}
and there is an exact sequence
$$ 0 \rightarrow \widetilde{M}^{\delta}(l)\langle 1 \rangle /\ldots /\widetilde{M}^{\delta}(k-1)\langle k-l \rangle \rightarrow L_{k-l+1}\widetilde{Z}_{\beta}^{\gamma} \widetilde{\epsilon}_{\gamma}^{\beta}L_1 \widetilde{Z}_{\beta}^{\gamma} \widetilde{K}_{l-3}
\rightarrow \widetilde{M}^{\delta}(l+1) \langle 0 \rangle /\ldots /\widetilde{M}^{\delta}(k-1)\langle k-l-2 \rangle \rightarrow 0. $$
\end{lemma}

\begin{proof}
By corollary ~\ref{cor12}, we get $ L_s \widetilde{Z}_{\beta}^{\delta} \widetilde{\epsilon}_{\gamma}^{\beta} L_0 \widetilde{Z}_{\beta}^{\gamma} \widetilde{K}_{k-2} \cong $
\begin{eqnarray*} 
\widetilde{M}^{\delta}(k-1)\langle k-2-l \rangle & \textrm{ if } & s=k-l-1 \\
0  & \textrm{ if } & s \neq k-l-1.
\end{eqnarray*}
Continuing in this way and using corollary ~\ref{cor12} we get
$ L_s \widetilde{Z}_{\beta}^{\delta} \widetilde{\epsilon}_{\gamma}^{\beta} L_0 \widetilde{Z}_{\beta}^{\gamma} \widetilde{K}_{l} \cong $
\begin{eqnarray*} 
\widetilde{M}^{\delta}(l+1)\langle 0 \rangle/\ldots /\widetilde{M}^{\delta}(k-1)\langle k-2-l \rangle & \textrm{ if } & s=k-l-1 \\
0  & \textrm{ if } & s \neq k-l-1
\end{eqnarray*}
and then 
$ L_s \widetilde{Z}_{\beta}^{\delta} \widetilde{\epsilon}_{\gamma}^{\beta} L_1 \widetilde{Z}_{\beta}^{\gamma} \widetilde{K}_{l-1} \cong $
\begin{eqnarray*} 
\widetilde{M}^{\delta}(l+1)\langle 0 \rangle /\ldots /\widetilde{M}^{\delta}(k-1)\langle k-l-2 \rangle & \textrm{ if } & s=k-l-1 \\
0  & \textrm{ if } & s \neq k-l-1.
\end{eqnarray*}
Corollary ~\ref{cor12} then gives an exact sequence
$$ 0 \rightarrow \widetilde{M}^{\delta}(l-1)\langle 0 \rangle \rightarrow L_{k-l} \widetilde{Z}_{\beta}^{\delta} \widetilde{\epsilon}_{\gamma}^{\beta} L_1 \widetilde{Z}_{\beta}^{\gamma}\widetilde{K}_{l-2} \rightarrow
\widetilde{M}^{\delta}(l+1)\langle 0 \rangle/\ldots /\widetilde{M}^{\delta}(k-1)\langle k-l-2 \rangle \rightarrow 0. $$
Next, corollary ~\ref{cor12} gives the following diagram:

\begin{tiny}
\xymatrix
{
&	&0\ar[d]	&	&	\\
&	&\widetilde{M}^{\delta}(l-1)\langle 0 \rangle\ar[d]		&	&	\\
&L_{k-l+1}\widetilde{Z}_{\beta}^{\delta} \widetilde{\epsilon}_{\gamma}^{\beta} L_1\widetilde{Z}_{\beta}^{\gamma} \widetilde{K}_{l-3}\ar@{^{(}->}[r]	&L_{k-l}\widetilde{Z}_{\beta}^{\delta} \widetilde{\epsilon}_{\gamma}^{\beta} L_1\widetilde{Z}_{\beta}^{\gamma} \widetilde{K}_{l-2}\ar[r]\ar[d]	&\widetilde{M}^{\delta}(l-2)\langle -1 \rangle \ar@{->>}[r] &L_{k-l}
\widetilde{Z}_{\beta}^{\delta} \widetilde{\epsilon}_{\gamma}^{\beta} L_1\widetilde{Z}_{\beta}^{\gamma} \widetilde{K}_{l-3}\\
&	&\widetilde{M}^{\delta}(l+1)\langle 0 \rangle/ \ldots /\widetilde{M}^{\delta}(k-1)\langle k-l-2 \rangle \ar[d]	&	&\\
&	&0	&	&	
}
\end{tiny}

Thus $ L_{k-l} \widetilde{Z}_{\beta}^{\gamma} \widetilde{\epsilon}_{\gamma}^{\beta} L_1 \widetilde{Z}_{\beta}^{\gamma} \widetilde{K}_{l-3} \cong \widetilde{M}^{\delta}(l-2)\langle -1 \rangle/ \ldots / \widetilde{M}^{\delta}(k-1)\langle k-l \rangle $ and there is an exact sequence
$$ 0 \rightarrow \widetilde{M}^{\delta}(l)\langle 1 \rangle / \ldots /\widetilde{M}^{\delta}(k-1)\langle k-l \rangle \rightarrow L_{k-l+1} \widetilde{Z}_{\beta}^{\gamma} \widetilde{\epsilon}_{\gamma}^{\beta} L_1 \widetilde{Z}_{\beta}^{\gamma} \widetilde{K}_{l-3} \rightarrow 
\widetilde{M}^{\delta}(l+1)\langle 0 \rangle / \ldots /\widetilde{M}^{\delta}(k-1)\langle k-l-2 \rangle \rightarrow 0. $$

Continuing in this manner and using the corollary gives us the lemma.
\end{proof}

\begin{lemma}
\label{lemma49}
Let $ m=l. $  Then
$ L_s \widetilde{Z}_{\beta}^{\alpha} \widetilde{\epsilon}_{\gamma}^{\beta} L_1 \widetilde{Z}_{\beta}^{\gamma} \widetilde{\epsilon}_{\delta}^{\beta} \widetilde{M}^{\delta}(l) \cong $
\begin{eqnarray*} 
\widetilde{M}^{\alpha}(k-1)\langle l-k+2 \rangle /\ldots /\widetilde{M}^{\alpha}(0)\langle l+1 \rangle & \textrm{ if } & s=l+1 \\
L_{l+2}\widetilde{Z}_{\beta}^{\alpha} \widetilde{\epsilon}_{\gamma}^{\beta}L_1 \widetilde{Z}_{\beta}^{\gamma} \widetilde{J}_{k-l-4}  & \textrm{ if } & s=k-1
\end{eqnarray*}
and there is an exact sequence
$$ 0 \rightarrow \widetilde{M}^{\alpha}(l)\langle 1 \rangle /\ldots \widetilde{M}^{\alpha}(0)\langle l+1 \rangle \rightarrow L_{l+2}\widetilde{Z}_{\beta}^{\alpha} \widetilde{\epsilon}_{\gamma}^{\beta}L_1 \widetilde{Z}_{\beta}^{\gamma} \widetilde{J}_{k-l-4}
\rightarrow \widetilde{M}^{\alpha}(l-1)\langle 0 \rangle /\ldots /\widetilde{M}^{\alpha}(0)\langle l-1 \rangle \rightarrow 0. $$
\end{lemma}

\begin{proof}
We use the exact sequences of lemma ~\ref{lemma46}.

If $ s=l, $ $ L_s \widetilde{Z}_{\beta}^{\alpha} \widetilde{\epsilon}_{\gamma}^{\beta} L_0 \widetilde{Z}_{\beta}^{\gamma} \widetilde{J}_{k-2} \cong \widetilde{M}^{\alpha}(0)\langle l-1 \rangle $ and is zero otherwise.

If $ s=l, $ $ L_s \widetilde{Z}_{\beta}^{\alpha} \widetilde{\epsilon}_{\gamma}^{\beta} L_0 \widetilde{Z}_{\beta}^{\gamma} \widetilde{J}_{k-3} \cong \widetilde{M}^{\alpha}(1)\langle l-2 \rangle /\widetilde{M}^{\alpha}(0)\langle l-1 \rangle $ and is zero otherwise.

Continuing in this way,

If $ s=l, $ $ L_s \widetilde{Z}_{\beta}^{\alpha} \widetilde{\epsilon}_{\gamma}^{\beta} L_0 \widetilde{Z}_{\beta}^{\gamma} \widetilde{J}_{k-l-1} \cong \widetilde{M}^{\alpha}(l-1)\langle 0 \rangle/ \ldots /\widetilde{M}^{\alpha}(0)\langle l-1 \rangle $ and is zero otherwise.

Next we get an exact sequence
$$ 0 \rightarrow \widetilde{M}^{\alpha}(l+1)\langle 0 \rangle \rightarrow L_{l+1} \widetilde{Z}_{\beta}^{\alpha} \widetilde{\epsilon}_{\gamma}^{\beta} L_1 \widetilde{Z}_{\beta}^{\gamma} \widetilde{J}_{k-l-3} \rightarrow
\widetilde{M}^{\alpha}(l-1)\langle 0 \rangle/ \ldots /\widetilde{M}^{\alpha}(0)\langle l-1 \rangle \rightarrow 0. $$

Lemma ~\ref{lemma46} then produces the following diagram:

\begin{tiny}
\xymatrix
{
&	&0\ar[d]	&	&	&\\
&	&\widetilde{M}^{\alpha}(l+1)\langle 0 \rangle\ar[d]		&	&	&\\
&L_{l+2}\widetilde{Z}_{\beta}^{\alpha} \widetilde{\epsilon}_{\gamma}^{\beta} L_1\widetilde{Z}_{\beta}^{\gamma} \widetilde{J}_{k-l-4}\ar@{^{(}->}[r]	&L_{l+1}\widetilde{Z}_{\beta}^{\alpha} \widetilde{\epsilon}_{\gamma}^{\beta} L_1\widetilde{Z}_{\beta}^{\gamma} \widetilde{J}_{k-l-3}\ar[r]\ar[d]	&\widetilde{M}^{\alpha}(l+2)\langle -1 \rangle \ar[r] &L_{l+1}\widetilde{Z}_{\beta}^{\alpha} \widetilde{\epsilon}_{\gamma}^{\beta} L_1\widetilde{Z}_{\beta}^{\gamma} \widetilde{J}_{k-l-4}\ar[r] &0\\
&	&\widetilde{M}^{\alpha}(l-1)\langle 0 \rangle/ \ldots/ \widetilde{M}^{\alpha}(0)\langle l-1\rangle \ar[d]	&	&	&\\
&	&0	&	&	&
}
\end{tiny}

Thus $ L_{l+1} \widetilde{Z}_{\beta}^{\alpha} \widetilde{\epsilon}_{\gamma}^{\beta} L_1 \widetilde{Z}_{\beta}^{\gamma} \widetilde{J}_{k-l-4} \cong \widetilde{M}^{\alpha}(l+2)\langle -1 \rangle/ \ldots /\widetilde{M}^{\alpha}(0)\langle l+1 \rangle $ and there is an exact sequence
$$ 0 \rightarrow \widetilde{M}^{\alpha}(l)\langle 1 \rangle /\ldots /\widetilde{M}^{\alpha}(0)\langle l+1 \rangle \rightarrow L_{l+2}\widetilde{Z}_{\beta}^{\alpha} \widetilde{\epsilon}_{\gamma}^{\beta}L_1 \widetilde{Z}_{\beta}^{\gamma} \widetilde{J}_{k-l-4}
\rightarrow \widetilde{M}^{\alpha}(l-1)\langle 0 \rangle /\ldots/ \widetilde{M}^{\alpha}(0)\langle l-1 \rangle \rightarrow 0. $$

Using the rest of the exact sequences of lemma ~\ref{lemma46}, we easily obtain this lemma. (The details are the same as lemma ~\ref{lemma48}.)
\end{proof}

Let $ F=L\widetilde{Z}_{\beta}^{\alpha}
\widetilde{\epsilon}_{\beta}^{\alpha}L_1\widetilde{Z}_{\beta}^{\gamma} \widetilde{\epsilon}_{\delta}^{\beta}. $

\begin{lemma}
For $ s=0, 1, \ldots, k-l-2, $
$ H^{-r} F \widetilde{M}^{\delta}(l)\langle l \rangle / \ldots /\widetilde{M}^{\delta}(k-1)\langle k-1 \rangle \cong $
\begin{eqnarray*} 
(\widetilde{M}^{\alpha}(k-1)\langle 0 \rangle /\ldots /\widetilde{M}^{\delta}(0)\langle k-1 \rangle) \langle k-(2s+1)-l-1 \rangle & \textrm{ if } & r=2(k-1)-l-(2s+1) \\
\widetilde{M}^{\alpha}(l-1)\langle l \rangle / \ldots \widetilde{M}^{\alpha}(0)\langle 2l-1 \rangle  & \textrm{ if } & r=k-1.
\end{eqnarray*}
\end{lemma}

\begin{proof}
We will verify this by induction on $ l. $
By proposition ~\ref{prop13},
$$ F \widetilde{M}^{\delta}(0)\langle 0 \rangle / \ldots /\widetilde{M}^{\delta}(k-1) \langle k-1 \rangle \cong \oplus_{r=0}^{k-2} \widetilde{M}^{\alpha}(k-1) \langle 0 \rangle/ \ldots /\widetilde{M}^{\alpha}(0)\langle k-1 \rangle [2r+1] \langle 2r-k+2 \rangle. $$
This is the base case.

Consider the short exact sequence
$$ 0 \rightarrow \widetilde{M}^{\delta}(l)\langle l \rangle / \ldots / \widetilde{M}^{\delta}(k-1)\langle k-1 \rangle \rightarrow \widetilde{M}^{\delta}(l-1)\langle l-1 \rangle \rightarrow 
\widetilde{M}^{\delta}(l-1)\langle l-1 \rangle / \ldots \widetilde{M}^{\delta}(k-1)\langle k-1 \rangle \rightarrow 0. $$
By the induction hypothesis, for $ s=0, 1, \ldots, k-l-1, $ $ H^{-r} F \widetilde{M}^{\delta}(l-1)\langle l-1 \rangle/ \ldots \widetilde{M}^{\delta}(k-1) \langle k-1\rangle \cong $
\begin{eqnarray*} 
(\widetilde{M}^{\alpha}(k-1)\langle 0 \rangle /\ldots /\widetilde{M}^{\delta}(0)\langle k-1\rangle) \langle k-(2s+1)-l \rangle & \textrm{ if } & r=2(k-1)-l+1-(2s+1) \\
\widetilde{M}^{\alpha}(l-2)\langle 1 \rangle/ \ldots / \widetilde{M}^{\alpha}(0)\langle l-1 \rangle  & \textrm{ if } & r=k-1.
\end{eqnarray*}

By lemma ~\ref{lemma49}, $ H^{-s} F \widetilde{M}^{\delta}(l-1) \cong $
\begin{eqnarray*} 
\widetilde{M}^{\alpha}(k-1) \langle l-k+1 \rangle /\ldots /\widetilde{M}^{\delta}(0)\langle l \rangle & \textrm{ if } & s=l \\
X_{l-1}  & \textrm{ if } & s=k-1,
\end{eqnarray*}
such that $ X_{l-1} $ fits into the following diagram:

\begin{tiny}
\xymatrix
{
0\ar[r]	&\widetilde{M}^{\alpha}(l-1)\langle l \rangle / \ldots/ \widetilde{M}^{\alpha}(0)\langle 2l-1 \rangle \ar[r]^{f}	&X_{l-1}\ar[r]^{g}\ar[d]^{h}		&\widetilde{M}^{\alpha}(l-2)\langle l-1 \rangle / \ldots/ \widetilde{M}^{\alpha}(0)\langle 2l-3 \rangle \ar[r]		&0\\
	&	&\widetilde{M}^{\alpha}(l-2)\langle l-1 \rangle / \ldots / \widetilde{M}^{\alpha}(0)\langle 2l-3 \rangle.	&	&
}
\end{tiny}
The fact that the maps $ g $ and $ h $ are the same follows from the proof of the ungraded case.
This finishes the lemma as in the ungraded case as well.
\end{proof}

\begin{lemma}
$ H^{-s} F L\widetilde{Z}_{\beta}^{\delta} \widetilde{\epsilon}_{\gamma}^{\beta} L_1 \widetilde{Z}_{\beta}^{\gamma} \widetilde{\epsilon}_{\alpha}^{\beta} \widetilde{M}^{\alpha}(l) \cong $
\begin{eqnarray*} 
(\widetilde{M}^{\alpha}(k-1)\langle 0 \rangle /\ldots /\widetilde{M}^{\delta}(0)\langle k-1)\langle s-2k+l+1 \rangle & \textrm{ if } & s=k-l+1, k-l+3, \ldots, k-l+2(k-1)-3 \\
\widetilde{M^{\alpha}}(l)  & \textrm{ if } & s=2(k-1).
\end{eqnarray*}
\end{lemma}

\begin{proof}
By lemma ~\ref{lemma48}, there is a distinguished triangle
$$ L_{k-l+1} \widetilde{Z}_{\beta}^{\delta} \widetilde{\epsilon}_{\gamma}^{\beta} L_1 \widetilde{Z}_{\beta}^{\gamma} \widetilde{K}_{l-3}[k-1] \rightarrow 
L\widetilde{Z}_{\beta}^{\delta} \widetilde{\epsilon}_{\gamma}^{\beta} L_1 \widetilde{Z}_{\beta}^{\gamma} \widetilde{\epsilon}_{\alpha}^{\beta} \widetilde{M}^{\alpha}(l) \rightarrow
\widetilde{M}^{\delta}(0)\langle -l+1 \rangle / \ldots /\widetilde{M}^{\delta}(k-1)\langle k-l \rangle [k-l]. $$
By proposition ~\ref{prop13}, we know that $ H^{-s} F \widetilde{M}^{\delta}(0)\langle 0 \rangle/ \ldots /\widetilde{M}^{\delta}(k-1)\langle k-1 \rangle [k-l] \cong (\widetilde{M}^{\alpha}(k-1)\langle 0 \rangle/ \ldots /\widetilde{M}^{\alpha}(0)\langle k-1 \rangle)\langle s-2k+l+1 \rangle $ for
$ s = k-l+1, k-l+3, \ldots, k-l+2(k-1)-1. $

We easily compute that for $ s=2(k-1), $
$$ H^{-s} F L_{k-l+1} \widetilde{Z}_{\beta}^{\delta} \widetilde{\epsilon}_{\gamma}^{\beta} L_1 \widetilde{Z}_{\beta}^{\gamma} \widetilde{K}_{l-3}[k-1] \cong
\widetilde{M}^{\alpha}(l) \langle 0 \rangle. $$

Therefore the long exact sequence for the distinguished triangle and $ F $  gives
$ H^{-s} F L\widetilde{Z}_{\beta}^{\delta} \widetilde{\epsilon}_{\gamma}^{\beta} L_1 \widetilde{Z}_{\beta}^{\gamma} \widetilde{\epsilon}_{\alpha}^{\beta} \widetilde{M}^{\alpha}(l) \cong $
\begin{eqnarray*} 
(\widetilde{M}^{\alpha}(k-1)\langle 0 \rangle /\ldots /\widetilde{M}^{\delta}(0)\langle k-1 \rangle)\langle s-2k+2 \rangle & \textrm{ if } & s=k-l+1, k-l+3, \ldots, k-l+2(k-1)-3 \\
\widetilde{M}^{\alpha}(l)  & \textrm{ if } & s=2(k-1).
\end{eqnarray*}
\end{proof}

This computation along with proposition ~\ref{propdiagram4} gives the following corollary.

\begin{corollary}
There is an isomorphism of graded functors:
$$ L\widetilde{Z}_{\mathfrak{\beta}}^{\mathfrak{\alpha}} \widetilde{\epsilon}_{\gamma}^{\beta}[-1] L\widetilde{Z}_{\mathfrak{\beta}}^{\gamma} \widetilde{\epsilon}_{\delta}^{\beta}[-(k-1)]L\widetilde{Z}_{\mathfrak{\beta}}^{\delta} \widetilde{\epsilon}_{\gamma}^{\beta}[-1] L\widetilde{Z}_{\mathfrak{\beta}}^{\gamma} \widetilde{\epsilon}_{\alpha}^{\beta}[-(k-1)] \cong \Id \oplus
\oplus_{r=1}^{k-2} \widetilde{\epsilon}_{\alpha+\delta}^{\alpha}[-(k-1)] L\widetilde{Z}_{\alpha}^{\alpha+\delta}[2r-(k-1)] \langle 2r-(k-1) \rangle. $$
\end{corollary}

\subsection{Graded Diagram Relation 5}
We need only compute in the graded case on a convenient generalized Verma module.  We choose
$ \widetilde{M}^{\mathfrak{s}_i}(\underbrace{a_i, a_{i+1}}, a_{i+2}) $ with $ a_i > a_{i+2} > a_{i+1}. $  
Now we just repeat most of the arguments of lemma ~\ref{lemma26}.

\begin{prop}
There is an isomorphism of graded functors
$$ L\widetilde{Z}^{\mathfrak{s}_i} \widetilde{\epsilon}_{\mathfrak{s}_{i+1}}[-1] L\widetilde{Z}^{\mathfrak{s}_{i+1}} \widetilde{\epsilon}_{\mathfrak{s}_i}[-1] \cong
Id \oplus \widetilde{\epsilon}_{\mathfrak{t}_i}^{\mathfrak{s}_i}[-2] L \widetilde{Z}_{\mathfrak{s}_i}^{\mathfrak{t}_i}. $$
\end{prop}

\begin{proof}
Consider the short exact sequence
$$ 0 \rightarrow \widetilde{M}(a_{i+1}, a_i, a_{i+2}) \langle 1 \rangle \rightarrow \widetilde{M}(a_{i}, a_{i+1}, a_{i+2}) \rightarrow  \widetilde{M}^{\mathfrak{s}_i}(\underbrace{a_{i}, a_{i+1}}, a_{i+2}) 
\rightarrow 0. $$

This gives rise to the short exact sequence
\begin{equation}
\label{eq2}
0 \rightarrow \widetilde{M}^{\mathfrak{s}_{i+1}}(a_i, \underbrace{a_{i+2}, a_{i+1}}) \rightarrow
L_1 \widetilde{Z}^{\mathfrak{s}_{i+1}} \widetilde{\epsilon}_{\mathfrak{s}_i} \widetilde{M}^{\mathfrak{s}_i}(\underbrace{a_i, a_{i+1}}, a_{i+2}) \rightarrow 
\widetilde{M}^{\mathfrak{s}_{i+1}}(a_{i+1}, \underbrace{a_{i}, a_{i+2}}) \rightarrow 0. 
\end{equation}
There are two shifts in the first module of this sequence.  One comes from the Koszul property of the module.  The other is due to the definition of the functor, but they cancel each other out.

Now consider the exact sequence
$$ 0 \rightarrow \widetilde{M}(a_{i+1}, a_{i+2}, a_{i}) \langle 1 \rangle \rightarrow \widetilde{M}(a_{i+1}, a_{i}, a_{i+2}) \rightarrow
\widetilde{M}^{\mathfrak{s}_{i+1}}(a_{i+1}, \underbrace{a_{i}, a_{i+2}}) \rightarrow 0. $$
Applying the functor $ L\widetilde{Z}^{\mathfrak{s}_i} $ gives
\begin{align*}
L_1 \widetilde{Z}^{\mathfrak{s}_{i+1}} \widetilde{\epsilon}_{\mathfrak{s}_{i+1}} \widetilde{M}^{\mathfrak{s}_{i+1}}(a_{i+1}, a_{i}, a_{i+2}) &\cong
\widetilde{M}^{\mathfrak{s}_i}(\underbrace{a_i, a_{i+1}}, a_{i+2})/\widetilde{M}^{\mathfrak{s}_i}(\underbrace{a_{i+1}, a_{i+2}}, a_{i}) \langle 1 \rangle\\
L_0 \widetilde{Z}^{\mathfrak{s}_{i+1}} \widetilde{\epsilon}_{\mathfrak{s}_{i+1}} \widetilde{M}^{\mathfrak{s}_{i+1}}(a_{i+1}, a_{i}, a_{i+2}) &\cong
L_2 \widetilde{Z}^{\mathfrak{s}_{i+1}} \widetilde{\epsilon}_{\mathfrak{s}_{i+1}} \widetilde{M}^{\mathfrak{s}_{i+1}}(a_{i+1}, a_{i}, a_{i+2}) \cong 0. 
\end{align*} 

Next we look at
$$ 0 \rightarrow \widetilde{M}(a_{i}, a_{i+1}, a_{i+2}) \langle 1 \rangle \rightarrow \widetilde{M}(a_{i}, a_{i+2}, a_{i+1}) \rightarrow
\widetilde{M}^{\mathfrak{s}_{i+1}}(a_{i}, \underbrace{a_{i+2}, a_{i+1}}) \rightarrow 0. $$
This leads to the exact sequence
$$ \aligned
0 \rightarrow &L_1 \widetilde{Z}^{\mathfrak{s}_i} \widetilde{\epsilon}_{\mathfrak{s}_{i+1}} \widetilde{M}^{\mathfrak{s}_{i+1}}(a_{i}, \underbrace{a_{i+2}, a_{i+1}}) \rightarrow
\widetilde{M}^{\mathfrak{s}_i}(\underbrace{a_i, a_{i+1}}, a_{i+2}) \rightarrow \widetilde{M}^{\mathfrak{s}_i}(\underbrace{a_i, a_{i+1}}, a_{i+2}) \langle -1 \rangle \rightarrow\\
&L_0 \widetilde{Z}^{\mathfrak{s}_i} \widetilde{\epsilon}_{\mathfrak{s}_{i+1}} \widetilde{M}^{\mathfrak{s}_{i+1}}(a_{i}, \underbrace{a_{i+2}, a_{i+1}}) \rightarrow 0. \endaligned $$
As in the ungraded case, the middle map is the standard map so
\begin{align*}
L_0 \widetilde{Z}^{\mathfrak{s}_i} \widetilde{\epsilon}_{\mathfrak{s}_{i+1}} \widetilde{M}^{\mathfrak{s}_{i+1}}(a_{i}, a_{i+2}, a_{i+1}) &\cong 
\widetilde{\epsilon}_{\mathfrak{t}_i}^{\mathfrak{s}_i} \widetilde{M}^{\mathfrak{t}_i}(\underbrace{a_i, a_{i+2}, a_{i+1}}) \langle -1 \rangle\\
L_1 \widetilde{Z}^{\mathfrak{s}_i} \widetilde{\epsilon}_{\mathfrak{s}_{i+1}} \widetilde{M}^{\mathfrak{s}_{i+1}}(a_{i}, a_{i+2}, a_{i+1}) &\cong
\widetilde{M}^{\mathfrak{s}_i}(\underbrace{a_{i+2}, a_{i+1}}, a_{i}) \langle 1 \rangle. 
\end{align*} 

There is a long exact sequence coming from equation ~\ref{eq2}.  Substituting in the above into this exact sequence, we get
\begin{align*}
0 \rightarrow &\widetilde{M}^{\mathfrak{s}_i}(\underbrace{a_{i+2}, a_{i+1}}, a_{i}) \langle 1 \rangle \rightarrow
L_1 \widetilde{Z}^{\mathfrak{s}_i} \widetilde{\epsilon}_{\mathfrak{s}_{i+1}} L_1 \widetilde{Z}^{\mathfrak{s}_{i+1}} \widetilde{\epsilon}_{\mathfrak{s}_{i}}
\widetilde{M}^{\mathfrak{s}_i}(\underbrace{a_{i}, a_{i+1}}, a_{i+2}) \rightarrow\\
&\widetilde{M}^{\mathfrak{s}_i}(\underbrace{a_{i}, a_{i+1}}, a_{i+2})/\widetilde{M}^{\mathfrak{s}_i}(\underbrace{a_{i+1}, a_{i+2}}, a_{i}) \langle 1 \rangle \rightarrow
\widetilde{\epsilon}_{\mathfrak{t}_i}^{\mathfrak{s}_i} \widetilde{M}^{\mathfrak{t}_i}(\underbrace{a_i, a_{i+2}, a_{i+1}}) \langle -1 \rangle \rightarrow\\
&L_0 \widetilde{Z}^{\mathfrak{s}_i} \widetilde{\epsilon}_{\mathfrak{s}_{i+1}} L_1 \widetilde{Z}^{\mathfrak{s}_{i+1}} \widetilde{\epsilon}_{\mathfrak{s}_{i}}
\widetilde{M}^{\mathfrak{s}_i}(\underbrace{a_{i}, a_{i+1}}, a_{i+2}) \rightarrow 0. 
\end{align*}

Thus as in the ungraded case,
$$ L_0 \widetilde{Z}^{\mathfrak{s}_i} \widetilde{\epsilon}_{\mathfrak{s}_{i+1}} L_1 \widetilde{Z}^{\mathfrak{s}_{i+1}} \widetilde{\epsilon}_{\mathfrak{s}_{i}}
\widetilde{M}^{\mathfrak{s}_i}(\underbrace{a_{i}, a_{i+1}}, a_{i+2})  \cong 
\widetilde{\epsilon}_{\mathfrak{t}_i}^{\mathfrak{s}_i} \widetilde{M}^{\mathfrak{t}_i}(\underbrace{a_i, a_{i+2}, a_{i+1}}) \langle -1 \rangle $$
and
$$ L_1 \widetilde{Z}^{\mathfrak{s}_i} \widetilde{\epsilon}_{\mathfrak{s}_{i+1}} L_1 \widetilde{Z}^{\mathfrak{s}_{i+1}} \widetilde{\epsilon}_{\mathfrak{s}_{i}}
\widetilde{M}^{\mathfrak{s}_i}(\underbrace{a_{i}, a_{i+1}}, a_{i+2})  \cong \widetilde{M}^{\mathfrak{s}_i}(\underbrace{a_{i}, a_{i+1}}, a_{i+2}). $$

Now, 
$$ \widetilde{\epsilon}_{\mathfrak{t}_i}^{\mathfrak{s}_i}[-2] L\widetilde{Z}_{\mathfrak{s}_i}^{\mathfrak{t}_i} \widetilde{M}^{\mathfrak{s}_i}(\underbrace{a_{i}, a_{i+1}}, a_{i+2}) \cong
\widetilde{M}^{\mathfrak{t}_i}(\underbrace{a_i, a_{i+2}, a_{i+1}}) \langle -2 \rangle \langle 1 \rangle. $$
The first shift comes from the definition of the functor.  The second shift comes from the fact that module is Koszul.
Thus
$$ L\widetilde{Z}^{\mathfrak{s}_i} \widetilde{\epsilon}_{\mathfrak{s}_{i+1}}[-1] L\widetilde{Z}^{\mathfrak{s}_{i+1}} \widetilde{\epsilon}_{\mathfrak{s}_i}[-1] \cong
\Id \oplus \widetilde{\epsilon}_{\mathfrak{t}_i}^{\mathfrak{s}_i}[-2] L \widetilde{Z}_{\mathfrak{s}_i}^{\mathfrak{t}_i}. $$
\end{proof}

\section{Crossings}
The category $ \mathcal{T} $ of oriented tangles has finite sequences of $ +, - $ signs for objects and isotopy classes of oriented tangles as morphisms.  This is  a strict tensor category.  Oriented caps and cups are clearly morphisms in this category.  All eight different types of oriented crossings are morphisms as well.  It is shown in [Ka], theorem XII.2.2 that the morphisms are generated by oriented cups and caps and only two oriented crossings.  That theorem also gives a complete list of the relations that are satisfied, (along with planar isotopies for tangles without crossings.)  The goal of this section is to assign functors to the two crossings which serve as generators and show that the defining relations are satisfied on a functorial level.  That is the main result of this paper.  

First we state the functorial relations guaranteeing invariance under isotopy for tangles without crossings.
This follows from theorem 6.2 of [Str2] except for the base case in the induction hypothesis of proposition 6.4.  This follows from lemma 11.80 of [KV].

\begin{theorem}
Let $ j \geq i. $  Then there are isomorphisms of functors
\begin{enumerate}
\item $$ \cap_{i+1, \pm, r+2} \circ \cup_{i, \mp, r} \cong \Id $$
\item $$ \cap_{i, \pm, r+2} \circ \cup_{i+1, \mp, r} \cong \Id $$
\item $$ \cap_{j, \pm, r-2} \circ \cap_{i, \pm, r} \cong \cap_{i, \pm, r-2} \circ \cap_{j+2, \pm, r} $$
\item $$ \cup_{j, \pm, r-2} \circ \cap_{i, \pm, r} \cong \cap_{i, \pm, r+2} \circ \cup_{j+2, \pm, r} $$
\item $$ \cup_{i, \pm, r-2} \circ  \cap_{j, \pm, r} \cong \cap_{j+2, \pm, r+2} \circ \cup_{i, \pm, r} $$
\item $$ \cup_{i, \pm, r+2} \circ \cup_{j, \pm, r} \cong \cup_{j+2, \pm, r+2} \circ \cup_{i,\pm, r} $$
\item $$ \cap_{i, \pm, r+2} \circ \cup_{i, \pm, r} \cong \Id \oplus \Id. $$
\end{enumerate}
\end{theorem}

Now to the crossing 
$$ \xymatrix{ {}\ar@{<-}[dr]|\hole	&	{}\ar@{<-}[dl]\\
{}	&	{}} $$
we assign the functor
$$ \Omega_{i}^{} = \text{Cone}(\widetilde{\epsilon}_{s_i}[-1] L\widetilde{Z}_{s_i} \rightarrow \Id[1] \langle 1 \rangle)[-k] \langle -k \rangle. $$

To the opposite crossing 
$$ \xymatrix{ {}\ar@{<-}[dr]	&	{}\ar@{<-}[dl]|\hole \\
{}	&	{}} $$
we assign the functor
$$ \Pi_{i}^{} = \text{Cone}(\Id \langle -1 \rangle \rightarrow \widetilde{\epsilon}_{s_i} L\widetilde{Z}_{s_i})[k-2] \langle k \rangle. $$

For ease of notation, we will assume that all arcs not shown in the diagrams are straight arrows oriented up.  Now we must check the following eight relations.

\begin{theorem}
\label{main}
The functors assigned to the generators of the set of morphisms for the category $ \mathcal{T} $ satisfy the following relations:
\begin{enumerate}
\item 
$$ \cap_{i+1, +, r+2} \circ \cup_{i, -, r} \cong \Id \cong \cap_{i, -, r+2} \circ \cup_{i+1, +, r}. $$
\begin{figure}[htb]
  \centering
  \includegraphics{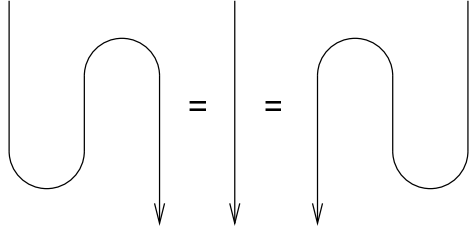}
  \caption{Relation 1}
  \label{Relation1}
\end{figure}

\item 
$$ \cap_{i+1, -, r+2} \circ \cup_{i, +, r} \cong \Id \cong \cap_{i, +, r+2} \circ \cup_{i+1, -, r}. $$
\begin{figure}[htb]
  \centering
  \includegraphics{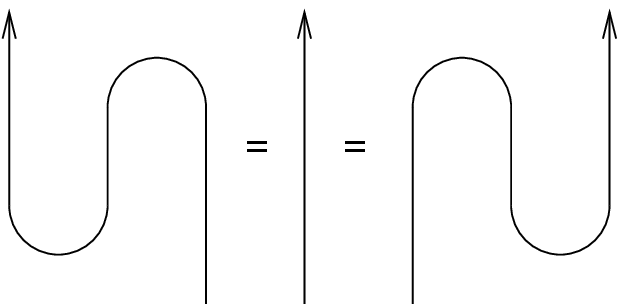}
  \caption{Relation 2}
  \label{Relation2}
\end{figure}

\item 
\begin{align*}
&\cap_{i-2,-,r+2} \circ \cap_{i-1,-,r+4} \circ \Omega_{i+2k-4} \circ \cup_{i+1,+,r+2} \circ \cup_{i,+,r} \cong\\
&\cap_{i,+,r+2} \circ \cap_{i+1,+,r+4} \circ \Omega_{i+2k-4} \circ \cup_{i-1, -,r+2} \circ \cup_{i-2,-,r} 
\end{align*} 

\begin{align*}
&\cap_{i-2,-,r+2} \circ \cap_{i-1,-,r+4} \circ \Omega_{i+2k-4} \circ \cup_{i+1,+,r+2} \circ \cup_{i,+,r} \cong\\
&\cap_{i,+,r+2} \circ \cap_{i+1,+,r+4} \circ \Omega_{i+2k-4} \circ \cup_{i-1, -,r+2} \circ \cup_{i-2,-,r} 
\end{align*} 

\begin{figure}[htb]
  \centering
  \includegraphics{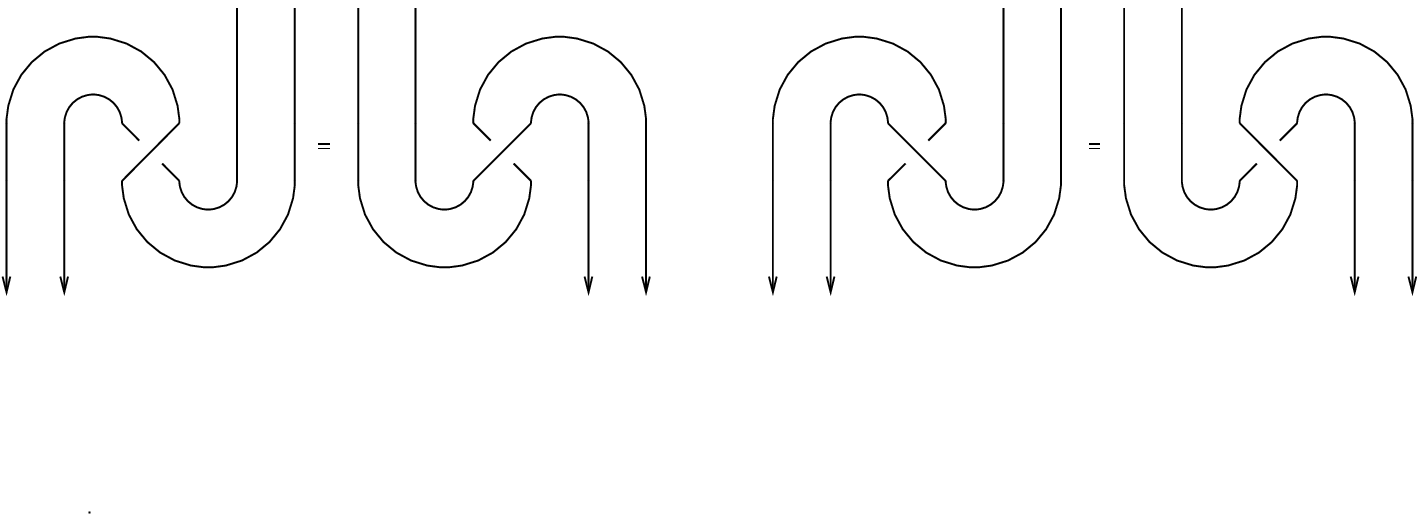}
  \caption{Relation 3}
  \label{Relation3}
\end{figure}

\item 
$$ \Omega_i \circ \Pi_i \cong \Id \cong \Pi_i \circ \Omega_i $$
\begin{figure}[htb]
  \centering
  \includegraphics{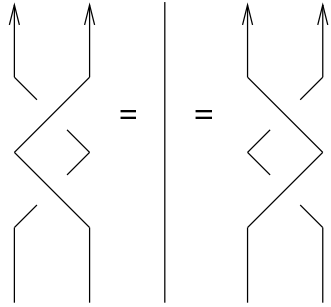}
  \caption{Relation 4}
  \label{Relation4}
\end{figure}

\item 
$$ \Omega_i \circ \Omega_{i+1} \circ \Omega_i \cong \Omega_{i+1} \circ \Omega_i \circ \Omega_{i+1} $$
\begin{figure}[htb]
  \centering
  \includegraphics{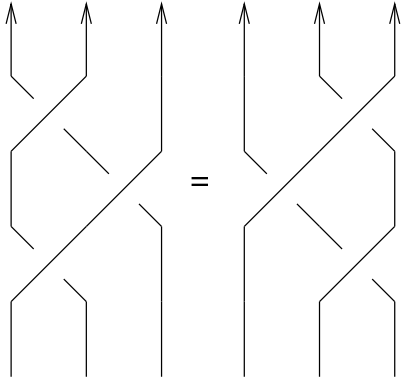}
  \caption{Relation 5}
  \label{Relation5}
\end{figure}

\item
$$ \cap_{i+1,+,r+2} \circ \Omega_{i} \circ \cup_{i+1,+,r} \cong \Id \cong \cap_{i+1,+,r+2} \circ \Pi_{i+k-2} \circ \cup_{i+1,+,r} $$
\begin{figure}[htb]
  \centering
  \includegraphics{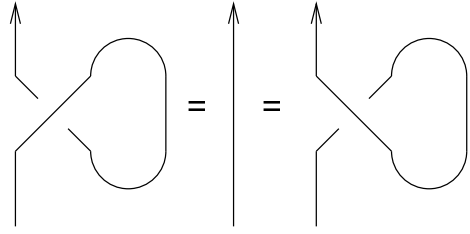}
  \caption{Relation 6}
  \label{Relation6}
\end{figure}

\item 
\begin{align*}
&\cap_{i,-,r+2} \circ \Omega_{i+k-1} \circ \cup_{i+2,+,r} \circ \cap_{i+2,+,r+2} \circ \Pi_{i+k-1} \circ \cup_{i, -,r} \cong
\Id \cong\\
&\cap_{i,-,r+2} \circ \Pi_{i+k-1} \circ \cup_{i+2,+,r} \circ \cap_{i+2,+,r+2} \circ \Omega_{i+k-1} \circ \cup_{i, -,r}
\end{align*} 
\begin{figure}[htb]
  \centering
  \includegraphics{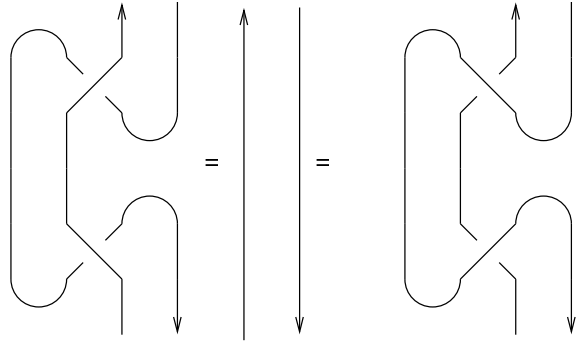}
  \caption{Relation 7}
  \label{Relation7}
\end{figure}

\item
\begin{align*}
&\cap_{i,+,r+2} \circ \Omega_{i+k-3} \circ \cup_{i-2,-,r} \circ \cap_{i-2,-,r+2} \circ \Pi_{i+k-3} \circ \cup_{i,+,r} 
\cong \Id \cong\\
&\cap_{i,+,r+2} \circ \Omega_{i+k-3} \circ \cup_{i-2,-,r} \circ \cap_{i-2,-,r+2} \circ \Pi_{i+k-3} \circ \cup_{i,+,r} 
\end{align*} 
\begin{figure}[htb]
  \centering
  \includegraphics{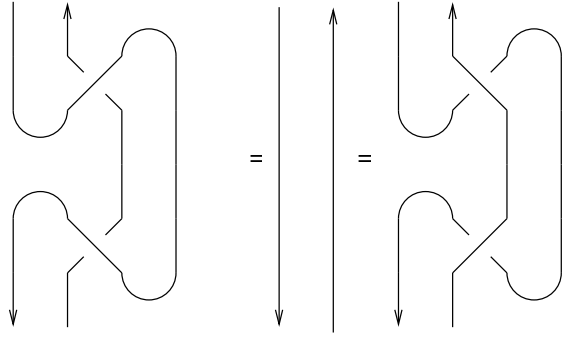}
  \caption{Relation 8}
  \label{Relation8}
\end{figure}
\end{enumerate}
\end{theorem}

\subsection{Relation 1}
The goal of this section is to prove that the functors assigned to the tangles in diagram ~\ref{Relation1} are isomorphic.  We interpret the integers $ a_i $ and $ -a_i, $ mod $ k. $

\begin{lemma}
$ \cap_{i+1, -a_i, r+2} \circ \cup_{i, a_i, r} \cong \Id. $
\end{lemma}

\begin{proof}
The left hand side is isomorphic to
$$ \nu \circ {\widetilde{\Delta}}_{(a_1, \ldots, a_i, k, \ldots, a_r)}^{(k, a_1, \ldots, a_r)} \circ L\widetilde{Z}_{(a_1, \ldots, a_{i-1}, a_i, -a_i, a_i, \ldots, a_r)}^{(a_1, \ldots, a_i, k, \ldots, a_r)} \circ
\widetilde{\epsilon}_{(a_1, \ldots, a_{i-1}, k, a_i, \ldots, a_r)}^{(a_1, \ldots, a_{i-1}, a_i, -a_i, a_i, \ldots, a_r)}[-(k-1)] \circ \widetilde{\Delta}_{(k, a_1, \ldots, a_r)}^{(a_1, \ldots, a_{i-1}, k, a_i, \ldots, a_r)} 
\circ \zeta. $$
This is easily seen to be isomorphic to
$$ \cong \nu \circ \widetilde{\Delta}_{(a_1, \ldots, a_i, k, \ldots, a_r)}^{(k, a_1, \ldots, a_r)} \circ
\widetilde{\Delta}_{(k, a_1, \ldots, a_r)}^{(a_1, \ldots, a_{i}, k, a_{i+1}, \ldots, a_r)} 
\circ \zeta \cong \Id.  $$

\end{proof}

\begin{lemma}
$ \cap_{i, a_i, r+2} \circ \cup_{i+1, -a_i, r} \cong \Id. $
\end{lemma}

\begin{proof}
The left hand side is isomorphic to
$$ \nu \circ \widetilde{\Delta}_{(a_1, \ldots, a_{i-1}, k, \ldots, a_r)}^{(k, a_1, \ldots, a_r)} \circ L\widetilde{Z}_{(a_1, \ldots, a_i, -a_i, a_i, a_{i+1},  \ldots, a_r)}^{(a_1, \ldots, a_{i-1}, k, \ldots, a_r)} \circ
\widetilde{\epsilon}_{(a_1, \ldots, a_{i}, k, \ldots, a_r)}^{(a_1, \ldots, a_{i}, -a_i, a_i, a_{i+1}, \ldots, a_r)}[-(k-1)] \circ \widetilde{\Delta}_{(k, a_1, \ldots, a_r)}^{(a_1, \ldots, a_{i}, k, \ldots, a_r)} 
\circ \zeta. $$
This is easily seen to be isomorphic to 
$$ \cong \nu \circ \widetilde{\Delta}_{(a_1, \ldots, a_{i-1}, k, \ldots, a_r)}^{(k, a_1, \ldots, a_r)} \circ
\widetilde{\Delta}_{(k, a_1, \ldots, a_r)}^{(a_1, \ldots, a_{i-1}, k, a_{i}, \ldots, a_r)} 
\circ \zeta \cong \Id.  $$
\end{proof}

Now we easily get relation 1.

\begin{corollary}
$ \cap_{i+1, +, r+2} \circ \cup_{i, -, r} \cong \Id \cong \cap_{i, -, r+2} \circ \cup_{i+1, +, r}. $
\end{corollary}

\begin{proof}
Take $ a_i = -1 $ and apply the two lemmas.
\end{proof}

\subsection{Relation 2}
There is the tangle relation for diagram 39.  
As in the previous subsection there is an analogous isomorphism of functors.

\begin{corollary}
$ \cap_{i+1, -, r+2} \circ \cup_{i, +, r} \cong \Id \cong \cap_{i, +, r+2} \circ \cup_{i+1, -, r}. $
\end{corollary}

\begin{proof}
Take $ a_i = 1 $ from the two lemmas in the previous subsection.
\end{proof}

\subsection{Relation 3}
In this section we define functors for crossings with both arrows pointing down.  It will be shown that these functors are consistent in a certain way with the functors defined for the other types of crossings.  First we must define functors categorifying intertwiners $ \Lambda^{k-1} V_{k-1} \otimes \Lambda^{k-1} V_{k-1} \rightarrow \Lambda^{k-2} \otimes \Lambda^k V_{k-1} $ and $ \Lambda^{k-2} \otimes \Lambda^k V_{k-1} \rightarrow \Lambda^{k-1} V_{k-1} \otimes \Lambda^{k-1} V_{k-1}. $ 
Let $ F_1 $ and $ F_2 $ be the functors for the diagrams on the left and right respectively of figure ~\ref{F_1,F_2}.
\begin{figure}[htb]
  \centering
  \includegraphics{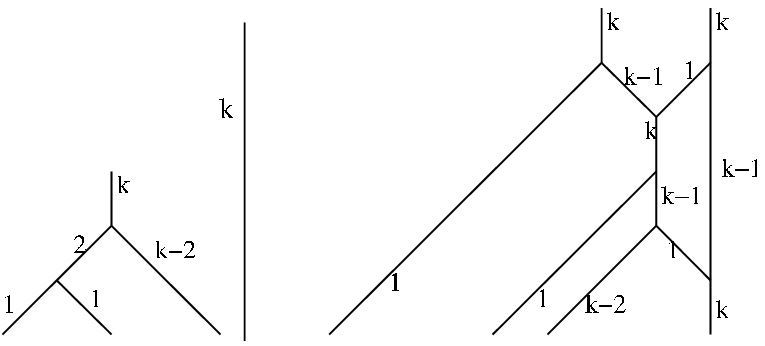}
  \caption{}
  \label{F_1,F_2}
\end{figure}

Explicitly,
$$ F_1 = L\widetilde{Z}_{\mathfrak{q}_i+\mathfrak{p}_{i+k}}^{\mathfrak{p}_i+\mathfrak{p}_{i+k}} L\widetilde{Z}_{\mathfrak{r}_{i+2}+\mathfrak{p}_{i+k}}^{\mathfrak{q}_i+\mathfrak{p}_{i+k}} $$
$$ F_2 = L\widetilde{Z}_{\mathfrak{p}_i+\mathfrak{q}_{i+k}}^{\mathfrak{p}_i+\mathfrak{p}_{i+k}} L\widetilde{Z}_{\mathfrak{q}_i+\mathfrak{q}_{i+k}}^{\mathfrak{p}_i+\mathfrak{q}_{i+k}} \widetilde{\epsilon}_{\mathfrak{p}_{i+1}+\mathfrak{q}_{i+k}}^{\mathfrak{q}_i+\mathfrak{q}_{i+k}}[-(k-1)] 
L\widetilde{Z}_{\mathfrak{q}_{i+1}+\mathfrak{q}_{i+k}}^{\mathfrak{p}_{i+1}+\mathfrak{q}_{i+k}} L\widetilde{Z}_{\mathfrak{r}_{i+2}+\mathfrak{q}_{i+k}}^{\mathfrak{q}_{i+1}+\mathfrak{q}_{i+k}} \widetilde{\epsilon}_{\mathfrak{r}_{i+2}+\mathfrak{p}_{i+k}}^{\mathfrak{r}_{i+2}+\mathfrak{q}_{i+k}}[-(k-1)]. $$

If $ \phi \colon R \rightarrow S $ is a natural transformation, then there is a natural transformation $ \phi^* \colon S^* \rightarrow R^* $ of adjoint functors.
There is a map constructed by taking the adjoint of an adjunction morphism:
$$ \alpha_1 \colon F_1 \rightarrow L\widetilde{Z}_{\mathfrak{r}_{i+2}+\mathfrak{q}_{i+k}}^{\mathfrak{p}_i+\mathfrak{p}_{i+k}} \widetilde{\epsilon}_{\mathfrak{p}_{i+1}+\mathfrak{q}_{i+k}}^{\mathfrak{r}_{i+2}+\mathfrak{q}_{i+k}}[-(2k-3)]
L\widetilde{Z}_{\mathfrak{r}_{i+2}+\mathfrak{q}_{i+k}}^{\mathfrak{p}_{i+1}+\mathfrak{q}_{i+k}} \widetilde{\epsilon}_{\mathfrak{r}_{i+2}+\mathfrak{p}_{i+k}}^{\mathfrak{r}_{i+2}+\mathfrak{q}_{i+k}}[-(k-1)][k-2]\langle k-2 \rangle. $$

We may rewrite the above line as
$$ \alpha_1 \colon F_1 \rightarrow L\widetilde{Z}_{\mathfrak{q}_{i}+\mathfrak{q}_{i+k}}^{\mathfrak{p}_i+\mathfrak{p}_{i+k}} L\widetilde{Z}_{\mathfrak{r}_{i+2}+\mathfrak{q}_{i+k}}^{\mathfrak{q}_i+\mathfrak{q}_{i+k}} 
\widetilde{\epsilon}_{\mathfrak{q}_{i}+\mathfrak{q}_{i+k}}^{\mathfrak{r}_{i+2}+\mathfrak{q}_{i+k}}[-(k-2)] \widetilde{\epsilon}_{\mathfrak{p}_{i+1}+\mathfrak{q}_{i+k}}^{\mathfrak{q}_{i}+\mathfrak{q}_{i+k}}[-(k-1)]
L\widetilde{Z}_{\mathfrak{r}_{i+2}+\mathfrak{q}_{i+k}}^{\mathfrak{p}_{i+1}+\mathfrak{q}_{i+k}} \widetilde{\epsilon}_{\mathfrak{r}_{i+2}+\mathfrak{p}_{i+k}}^{\mathfrak{r}_{i+2}+\mathfrak{q}_{i+k}}[-(k-1)][k-2]\langle k-2 \rangle =F_3. $$

Since 
$$ L\widetilde{Z}_{\mathfrak{r}_{i+2}+\mathfrak{q}_{i+k}}^{\mathfrak{q}_i+\mathfrak{q}_{i+k}} \widetilde{\epsilon}_{\mathfrak{q}_{i}+\mathfrak{q}_{i+k}}^{\mathfrak{r}_{i+2}+\mathfrak{q}_{i+k}} \cong \oplus_{j=0}^{k-2} \Id[2j]\langle -(k-2)+2j \rangle, $$
there is a projection $ \alpha_2 \colon F_3 \rightarrow F_2. $

\begin{lemma}
\label{lemma54}
The composite $ \alpha_2 \circ \alpha_1 \colon F_1 \rightarrow F_2 $ is an isomorphism.
\end{lemma}

\begin{proof}
The map coming from adjunction is clearly non-zero. Since $ \alpha_2 $ is just projection,
$ \alpha_2 \circ \alpha_1 $ is non-zero.
An easily calculation shows that $ F_1 $ an $ F_2 $ applied to a generalized Verma module 
$$ M^{\mathfrak{r}_{i+2}+\mathfrak{p}_{i+k}}(a_1, \ldots, a_{i-1}, a_i, a_{i+1}, \underbrace{a_{i+2}, \ldots, a_{i+k-1}} \underbrace{k-1, \ldots, 0,} a_{i+2k}, \ldots, a_n), $$
gives a generalized Verma module 
$$ M^{\mathfrak{p}_i+\mathfrak{p}_{i+k}}(a_1, \ldots, a_{i-1}, \underbrace{k-1, \ldots, 0,} \underbrace{k-1, \ldots, 0,} a_{i+2}, \ldots, a_n) $$
or zero.  Therefore $ \alpha_2 \circ \alpha_1 $ is a non-zero endomorphism of a generalized Verma module (when it is non-trivial) so it is an isomorphism.  This composition of maps is then an isomorphism when the functors are applied to a projective object by considering a filtration with subquotients of generalized Verma modules.  Then it it an isomorphism of functors because a bounded complex is quasi-isomorphic to a complex of projective objects.
\end{proof}

Let $ G_1 $ and $ G_2 $ be functors for the diagrams on the left and right respectively of figure ~\ref{G_1,G_2}.

\begin{figure}[htb]
  \centering
  \includegraphics{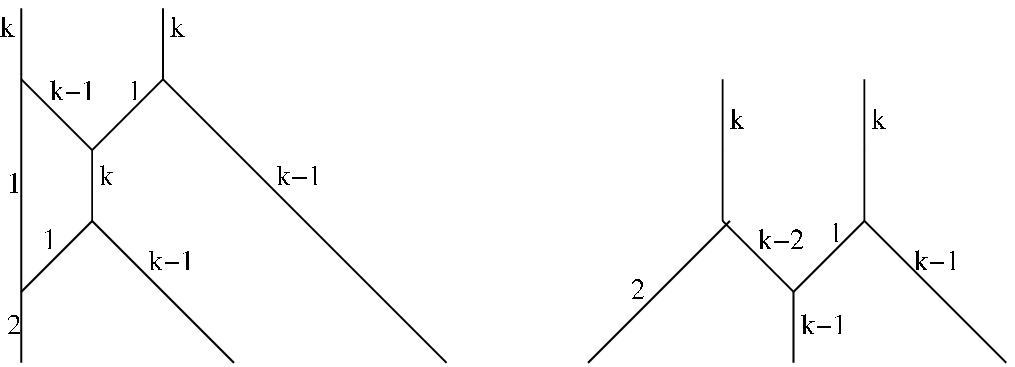}
  \caption{}
  \label{G_1,G_2}
\end{figure}

Explicity, 
\begin{align*}
&G_1 = L\widetilde{Z}_{\mathfrak{p}_i+\mathfrak{q}_{i+k}}^{\mathfrak{p}_i+\mathfrak{p}_{i+k}} L\widetilde{Z}_{\mathfrak{q}_i+\mathfrak{q}_{i+k}}^{\mathfrak{p}_i+\mathfrak{q}_{i+k}} 
\widetilde{\epsilon}_{\mathfrak{p}_{i+1}+\mathfrak{q}_{i+k}}^{\mathfrak{q}_i+\mathfrak{q}_{i+k}}[-(k-1)] L\widetilde{Z}_{\mathfrak{q}_{i+1}+\mathfrak{q}_{i+k}}^{\mathfrak{p}_{i+1}+\mathfrak{q}_{i+k}} \widetilde{\epsilon}_{\mathfrak{s}_i+\mathfrak{q}_{i+1}+\mathfrak{q}_{i+k}}^{\mathfrak{q}_{i+1}+\mathfrak{q}_{i+k}}[-1]\\
&G_2 = L\widetilde{Z}_{\mathfrak{p}_i+\mathfrak{q}_{i+k}}^{\mathfrak{p}_i+\mathfrak{p}_{i+k}} L\widetilde{Z}_{\mathfrak{s}_i+\mathfrak{r}_{i+2}+\mathfrak{q}_{i+k}}^{\mathfrak{p}_i+\mathfrak{q}_{i+k}} \widetilde{\epsilon}_{\mathfrak{s}_i+\mathfrak{q}_{i+1}+\mathfrak{q}_{i+k}}^{\mathfrak{s}_i+\mathfrak{r}_{i+2}+\mathfrak{q}_{i+k}}[-(k-2)]. 
\end{align*} 

There is a morphism constructed similarly to $ \alpha_1 $ above:
\begin{align*}
\beta_1 \colon G_2 \rightarrow 
&L\widetilde{Z}_{\mathfrak{q}_{i}+\mathfrak{q}_{i+k}}^{\mathfrak{p}_i+\mathfrak{p}_{i+k}} L\widetilde{Z}_{\mathfrak{r}_{i+2}+\mathfrak{q}_{i+k}}^{\mathfrak{q}_i+\mathfrak{q}_{i+k}} \widetilde{\epsilon}_{\mathfrak{q}_{i}+\mathfrak{q}_{i+k}}^{\mathfrak{r}_{i+2}+\mathfrak{q}_{i+k}}[-(k-2)]
\widetilde{\epsilon}_{\mathfrak{p}_{i+1}+\mathfrak{q}_{i+k}}^{\mathfrak{q}_{i}+\mathfrak{q}_{i+k}}[-(k-1)] \circ\\
&L\widetilde{Z}_{\mathfrak{q}_{i+1}+\mathfrak{q}_{i+k}}^{\mathfrak{p}_{i+1}+\mathfrak{q}_{i+k}} L\widetilde{Z}_{\mathfrak{r}_{i+2}+\mathfrak{q}_{i+k}}^{\mathfrak{q}_{i+1}+\mathfrak{q}_{i+k}} 
\widetilde{\epsilon}_{\mathfrak{q}_{i+1}+\mathfrak{q}_{i+k}}^{\mathfrak{r}_{i+2}+\mathfrak{q}_{i+k}}[-(k-2)]
\widetilde{\epsilon}_{\mathfrak{s}_i+\mathfrak{q}_{i+1}+\mathfrak{q}_{i+k}}^{\mathfrak{q}_{i+1}+\mathfrak{q}_{i+k}}[-1][2k-4]\langle -2k+4 \rangle. 
\end{align*} 

Then there is an obvious projection $ \beta_2 $ from the latter functor to $ G_1. $

\begin{lemma}
\label{lemma55}
The composite $ \beta_2 \circ \beta_1  \colon G_2 \rightarrow G_1 $ is an isomorphism.
\end{lemma}

\begin{proof}
The map $ \beta_1 $ is non-zero for an object when the functors applied to the object do not vanish because it comes from adjunction.
Since $ \beta_2 $ is projection, it too is non-zero when the functors are non-zero.
An easy calculation shows that $ G_1 $ and $ G_2 $ applied to the generalized Verma module
$$ M^{\mathfrak{s}_i \mathfrak{q}_{i+1}+\mathfrak{q}_{i+k}}(a_1, \ldots, \underbrace{a_i, a_{i+1},} \underbrace{a_{i+2}, \ldots, a_{i+k},} \underbrace{a_{i+k-1}, \ldots, a_{i+2k-1},} a_{i+2k}, \ldots, a_n) $$ 
is isomorphic to a shifted generalized Verma module
$$ M^{\mathfrak{p}_i+\mathfrak{p}_{i+k}}(a_1, \ldots, a_{i-1}, \underbrace{k-1, \ldots, 0,} \underbrace{k-1, \ldots, 0,} a_{i+2k}, \ldots, a_n) $$ 
or zero.
Therefore for generalized Verma modules $ \beta_2 \circ \beta_1  $ is an isomorphism.
It is then an isomorphism for projective objects by considering a Verma flag. It is then an isomorphism for all bounded complexed because there exists a complex of projective quasi-isomorphic to it.
\end{proof}

\begin{figure}[htb]
  \centering
  \includegraphics{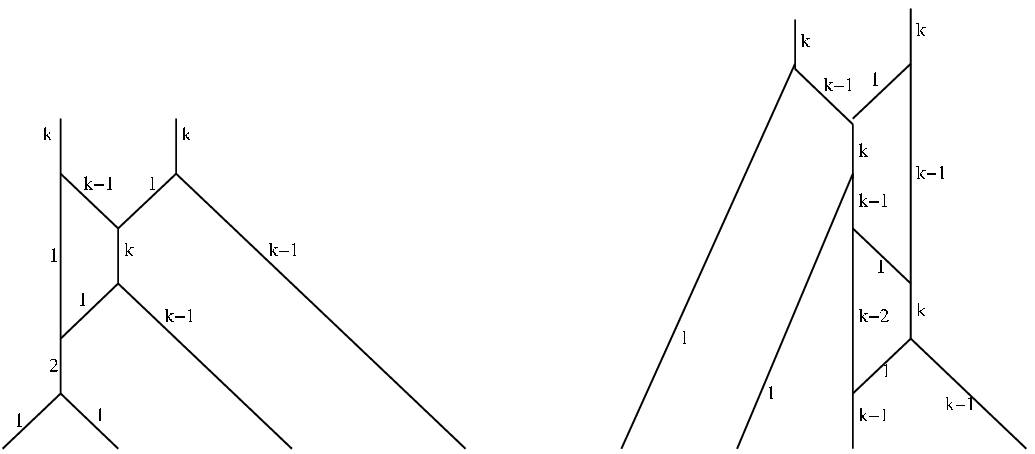}
  \caption{}
  \label{1,1,k-1,k-1->k,k}
\end{figure}

\begin{corollary}
There is an isomorphism of functors for the two diagrams in figure ~\ref{1,1,k-1,k-1->k,k}.  
$$ G_1 \circ L\widetilde{Z}_{\mathfrak{q}_{i+1}+\mathfrak{q}_{i+k}}^{\mathfrak{s}_i+\mathfrak{q}_{i+1}+\mathfrak{q}_{i+k}} \cong F_2 \circ L\widetilde{Z}_{\mathfrak{r}_{i+2}+\mathfrak{q}_{i+k}}^{\mathfrak{r}_{i+2}+\mathfrak{p}_{i+k}} \widetilde{\epsilon}_{\mathfrak{q}_{i+1}+\mathfrak{q}_{i+k}}^{\mathfrak{r}_{i+2}+\mathfrak{q}_{i+k}}[-(k-2)]. $$
\end{corollary}

\begin{proof}
By lemma ~\ref{lemma55}, the left hand side is isomorphic to
$$ L\widetilde{Z}_{\mathfrak{p}_i+\mathfrak{q}_{i+k}}^{\mathfrak{p}_i+\mathfrak{p}_{i+k}} L\widetilde{Z}_{\mathfrak{s}_i+\mathfrak{r}_{i+2}+\mathfrak{q}_{i+k}}^{\mathfrak{p}_i+\mathfrak{q}_{i+k}} \widetilde{\epsilon}_{\mathfrak{s}_i+\mathfrak{q}_{i+1}+\mathfrak{q}_{i+k}}^{\mathfrak{s}_i+\mathfrak{r}_{i+2}+\mathfrak{q}_{i+k}}[-(k-2)]
L\widetilde{Z}_{\mathfrak{q}_{i+1}+\mathfrak{q}_{i+k}}^{\mathfrak{s}_i+\mathfrak{q}_{i+1}+\mathfrak{q}_{i+k}}. $$
Commuting the last two functors in the expression above, we get that it is isomorphic to
$$ L\widetilde{Z}_{\mathfrak{p}_i+\mathfrak{q}_{i+k}}^{\mathfrak{p}_i+\mathfrak{p}_{i+k}} L\widetilde{Z}_{\mathfrak{s}_i+\mathfrak{r}_{i+2}+\mathfrak{q}_{i+k}}^{\mathfrak{p}_i+\mathfrak{q}_{i+k}}
L\widetilde{Z}_{\mathfrak{r}_{i+2}+\mathfrak{q}_{i+k}}^{\mathfrak{s}_i+\mathfrak{r}_{i+2}+\mathfrak{q}_{i+k}} \widetilde{\epsilon}_{\mathfrak{q}_{i+1}+\mathfrak{q}_{i+k}}^{\mathfrak{r}_{i+2}+\mathfrak{q}_{i+k}}[-(k-2)]. $$
Lemma ~\ref{lemma54} implies that this is isomorphic to the right hand side of the desired isomorphism.
\end{proof}

The next corollary gives an isomorphism between the two diagrams in figure ~\ref{k,k->1,1,k-1,k-1}.
\begin{figure}[htb]
  \centering
  \includegraphics{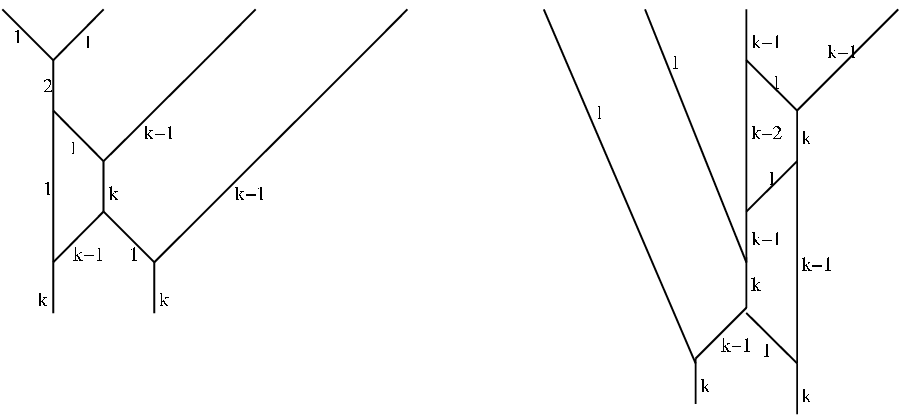}
  \caption{}
  \label{k,k->1,1,k-1,k-1}
\end{figure}

\begin{corollary}
Let $ G_1^* $ and $ F_2^* $ denote the adjoint functors.
Then
$$ \widetilde{\epsilon}_{\mathfrak{s}_i+\mathfrak{q}_{i+1}+\mathfrak{q}_{i+k}}^{\mathfrak{q}_{i+1}+\mathfrak{q}_{i+k}}[-1] \circ G_1^* \cong L\widetilde{Z}_{\mathfrak{r}_{i+2}+\mathfrak{q}_{i+k}}^{\mathfrak{q}_{i+1}+\mathfrak{q}_{i+k}} \circ \widetilde{\epsilon}_{\mathfrak{r}_{i+2}+\mathfrak{p}_{i+k}}^{\mathfrak{r}_{i+2}+\mathfrak{q}_{i+k}}[-(k-1)] \circ F_2^*. $$
\end{corollary}

\begin{proof}
This follows from the previous corollary by taking adjoints.
\end{proof}

Let $ H_1 $ be the functor for the diagram on the left and $ H_2 $ be the functor for the diagram on the right of figure ~\ref{H_1,H_2}.
\begin{figure}[htb]
  \centering
  \includegraphics{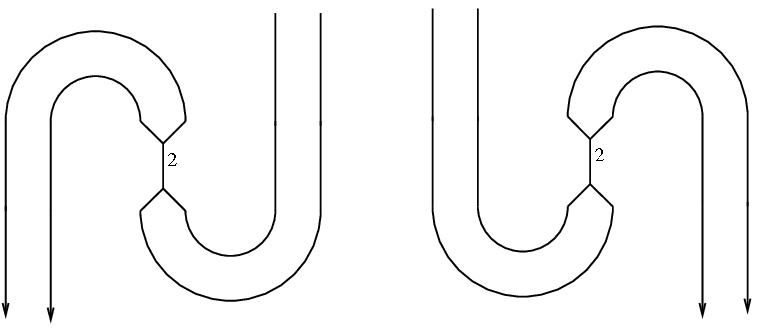}
  \caption{}
  \label{H_1,H_2}
\end{figure}
 
Let $ H_3 $ be the functor for the diagram on the left and $ H_4 $ be the functor for the diagram on the right of figure ~\ref{H_3,H_4}. 
\begin{figure}[htb]
  \centering
  \includegraphics{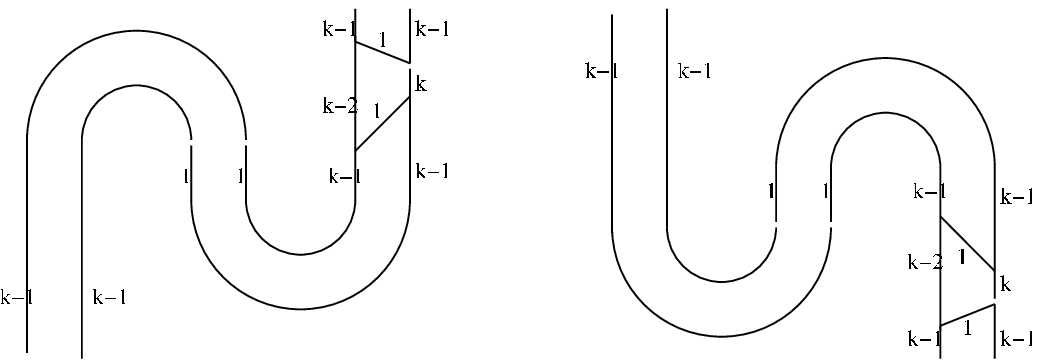}
  \caption{}
  \label{H_3,H_4}
\end{figure}

\begin{figure}[htb]
  \centering
  \includegraphics{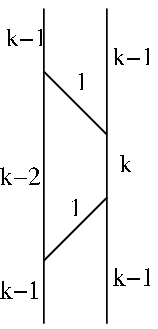}
  \caption{}
  \label{G}
\end{figure}

Let $ G = L\widetilde{Z}_{\mathfrak{r}_i+\mathfrak{q}_{i+k-2}}^{\mathfrak{q}_{i-1}+\mathfrak{q}_{i+k-2}} \widetilde{\epsilon}_{\mathfrak{r}_{i}+\mathfrak{p}_{i+k-1}}^{\mathfrak{r}_i+\mathfrak{q}_{i+k-2}}[-(k-1)]
L\widetilde{Z}_{\mathfrak{r}_i+\mathfrak{q}_{i+k-2}}^{\mathfrak{r}_{i}+\mathfrak{p}_{i+k-1}} \widetilde{\epsilon}_{\mathfrak{q}_{i-1}+\mathfrak{q}_{i+k-2}}^{\mathfrak{r}_i+\mathfrak{q}_{i+k-2}}[-(k-2)]. $   It is the functor for the diagram in figure ~\ref{G}.

\begin{corollary}
$ H_1 \cong H_2 \cong G. $
\end{corollary}

\begin{proof}
The two corollaries, imply that $ H_1 \cong H_3 $ and $ H_2 \cong H_4. $  Then, using the first relation, both of these functors are seen to be isomorphic to $ G. $
\end{proof}

Now we are prepared to prove the main relations of this section.
Let $ A $ and $ B $ denote the functors assigned to the diagrams on the left and right respectively of figure ~\ref{A,B}.
\begin{figure}[htb]
  \centering
  \includegraphics{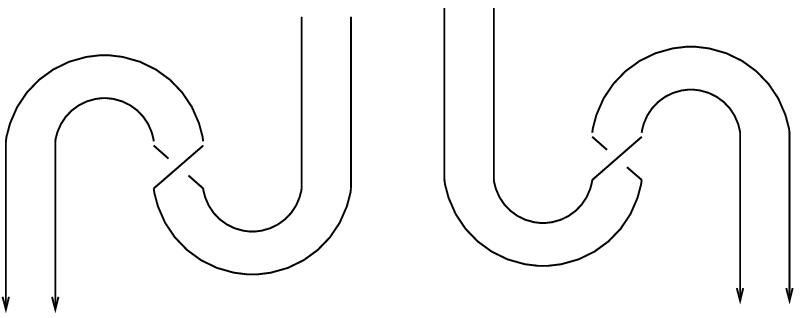}
  \caption{}
  \label{A,B}
\end{figure}

Let $ C $ and $ D $ denote the functors assigned to the diagrams on the left and right respectively of figure ~\ref{C,D}. 
\begin{figure}[htb]
  \centering
  \includegraphics{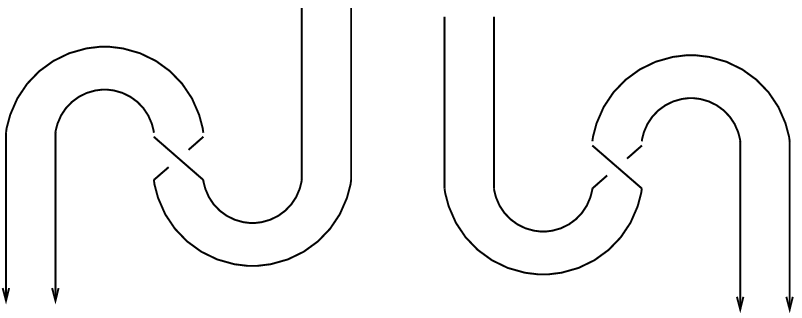}
  \caption{}
  \label{C,D}
\end{figure}

\begin{prop}
\begin{enumerate}
\item The functors $ A $ and $ B $ are isomorphic.
\item The functors $ C $ and $ D $ are isomorphic.
\end{enumerate}
\end{prop}

\begin{proof}
The first pair of functors are cones of morphisms $ G \rightarrow \Id[1]\langle 1 \rangle. $  The second pair are cones of morphisms $ \Id[-1] \langle -1 \rangle \rightarrow G. $ 
Clearly,
$$ \Hom(G, \Id[1] \langle 1 \rangle) \cong \Hom(\Id[-1] \langle -1 \rangle, G) \cong
\End(L\widetilde{Z}_{(k-2,1,k-1)}^{(k-2,k)} \widetilde{\epsilon}_{(k-1,k-1)}^{(k-2,1,k-1)}). $$
By theorem 34 of [MOS], this space is isomorphic to a space of endomorphisms of a projective functor:
$$ \End(\widetilde{\theta}_{(k-2,1,k-1)}^{(k-2,k)} \widetilde{\theta}_{(k-1,k-1)}^{(k-2,1,k-1)}). $$
By proposition 8.7 of [Str2], this is isomorphic to 
$$ \End_{C^{(k-2,k)}-\text{gmod}-C^{(k-1,k-1)}}(C^{(k-2,1,k-1)} \otimes_{C^{(k-1,k-1)}} C^{(k-1,k-1)}). $$
As a graded $ (C^{(k-2,k)},C^{(k-1,k-1)})- $ bimodule, $ C^{(k-2,1,k-1)} \otimes_{C^{(k-1,k-1)}} C^{(k-1,k-1)} $ is generated by $ 1 \otimes 1. $  Thus the degree zero component of this space of endomorphisms is one dimensional.  
Since there is a unique homogenous morphism up to scalar of degree 0 between the functors, the cones must be isomorphic.
\end{proof}

\subsection{Relation 4}
We will show that the functors assigned to the crossings, $ \Pi_i $ and $ \Omega_i $ are auto-equivalences of the derived category.

Let $ \alpha \colon \widetilde{\epsilon}_{\mathfrak{s}_i} L\widetilde{Z}_{\mathfrak{s}_i} \rightarrow \Id[2] \langle 1 \rangle $ and
$ \beta \colon \Id \rightarrow \widetilde{\epsilon}_{\mathfrak{s}_i} L \widetilde{Z}_{\mathfrak{s}_i} \langle 1 \rangle. $

\begin{prop}
$ \Pi_i \circ \Omega_i \cong \Id. $
\end{prop}

\begin{proof}
We would like to show that $ \text{Cone}(\beta) \text{Cone}(\alpha) \cong \Id[r] \langle s \rangle $ for some $ r $ and $ s. $
Consider the following commutative diagram.

\begin{tiny}
\xymatrix
{
&\text{Cone}(\beta) \widetilde{\epsilon}_{\mathfrak{s}_i}L\widetilde{Z}^{\mathfrak{s}_i} \ar[r]^{g}		& \text{Cone}(\beta)\Id[2]\langle 1 \rangle \ar[r]	&\text{Cone}(\beta)\text{Cone}(\alpha)	\\
&\widetilde{\epsilon}_{\mathfrak{s}_i}L\widetilde{Z}^{\mathfrak{s}_i} \langle 1 \rangle \widetilde{\epsilon}_{\mathfrak{s}_i}L\widetilde{Z}^{\mathfrak{s}_i} \ar[r]\ar[u]^{f}		&\widetilde{\epsilon}_{\mathfrak{s}_i}L\widetilde{Z}^{i} \langle 1\rangle \Id[2]\langle 1 \rangle \ar[r]\ar[u]	&\widetilde{\epsilon}_{\mathfrak{s}_i}L\widetilde{Z}^{\mathfrak{s}_i} \langle 1 \rangle \text{Cone}(\alpha) \ar[u] \\
&\Id \circ \widetilde{\epsilon}_{\mathfrak{s}_i}L\widetilde{Z}^{\mathfrak{s}_i} \ar[r]\ar[u]^{h}	&\Id \circ \Id[2] \langle 1 \rangle \ar[u]\ar[r]	&\Id \circ \text{Cone}(\alpha)\ar[u] \\
}
\end{tiny}

By properties of adjoint functors, $ h $ is split inclusion so  $ f $ is surjection onto $ \widetilde{\epsilon}_{\mathfrak{s}_i} L\widetilde{Z}^{\mathfrak{s}_i}[2] \langle 2 \rangle $ 
By commutativity, $ g $ is the map fitting into the distinguished triangle

\xymatrix
{
&\Id[2] \langle 1 \rangle \ar[r]^{\beta[2] \langle 1 \rangle}	&\widetilde{\epsilon}_{\mathfrak{s}_i} L\widetilde{Z}^{\mathfrak{s}_i}[2] \langle 2 \rangle \ar[r]^{g}	
&\text{Cone}(\beta)[2]\langle 1 \rangle.\\
}
Therefore $ \text{Cone}(\beta)\text{Cone}(\alpha) \cong \Id[3] \langle 1 \rangle. $
Then with the shifts built into the definition of $ \Pi_i $ and $ \Omega_i, $ we have $ \Pi_i \circ \Omega_i \cong \Id. $
\end{proof}

\begin{prop}
$ \Omega_i \circ \Pi \cong \Id. $
\end{prop}

\begin{proof}
As in the previous proposition we would like to show that $ \text{Cone}(\alpha) \text{Cone}(\beta) \cong \Id[r] \langle s \rangle $ for some $ r $ and $ s. $
Consider the following commutative diagram.

\begin{tiny}
\xymatrix
{
&\text{Cone}(\alpha)  \ar[r]		& \text{Cone}(\alpha)\widetilde{\epsilon}_{\mathfrak{s}_i}L\widetilde{Z}^{\mathfrak{s}_i} \langle 1 \rangle \ar[r]	&\text{Cone}(\alpha)\text{Cone}(\beta)	\\
&\Id[2]\langle 1 \rangle \ar[r]\ar[u]		&\Id[2]\langle 1 \rangle \widetilde{\epsilon}_{\mathfrak{s}_i}L\widetilde{Z}^{i} \langle 1 \rangle \ar[r]\ar[u]	&\Id[2] \langle 1 \rangle \text{Cone}(\beta) \ar[u] \\
&\widetilde{\epsilon}_{\mathfrak{s}_i}L\widetilde{Z}^{\mathfrak{s}_i} \ar[r]^{f}\ar[u]	&\widetilde{\epsilon}_{\mathfrak{s}_i}L\widetilde{Z}^{\mathfrak{s}_i}\widetilde{\epsilon}_{\mathfrak{s}_i}L\widetilde{Z}^{\mathfrak{s}_i} \langle 1 \rangle \ar[u]\ar[r]	
& \widetilde{\epsilon}_{\mathfrak{s}_i}L\widetilde{Z}^{\mathfrak{s}_i} \text{Cone}(\beta)\ar[u]^{h} \\
}
\end{tiny}

Note that $ f = \widetilde{\epsilon}_{\mathfrak{s}_i} L\widetilde{Z}^{\mathfrak{s}_i} adj $ and $ g \colon \widetilde{\epsilon}_{\mathfrak{s}_i} L\widetilde{Z}^{\mathfrak{s}_i} \widetilde{\epsilon}_{\mathfrak{s}_i} L\widetilde{Z}^{\mathfrak{s}_i} \langle 1 \rangle \rightarrow \widetilde{\epsilon}_{\mathfrak{s}_i} L\widetilde{Z}^{\mathfrak{s}_i} $ is given by $ \Id_{|\widetilde{\epsilon}_{\mathfrak{s}_i}} adj' \Id_{|L\widetilde{Z}^{\mathfrak{s}_i}} \langle 1 \rangle, $
where $ adj \colon \Id \rightarrow \widetilde{\epsilon}_{\mathfrak{s}_i} L\widetilde{Z}^{\mathfrak{s}_i}  $ is the obvious adjunction map and 
$ adj' \colon L\widetilde{Z}^{\mathfrak{s}_i}  \widetilde{\epsilon}_{\mathfrak{s}_i} \rightarrow \Id $ is an adjunction map as well.
By properties of adjunction maps, $ g \circ f = \Id. $ Thus $ f $ is split injective so 
$ \widetilde{\epsilon}_{\mathfrak{s}_i}L\widetilde{Z}^{\mathfrak{s}_i} \text{Cone}(\beta) \cong \widetilde{\epsilon}_{\mathfrak{s}_i}L\widetilde{Z}^{\mathfrak{s}_i}[2] \langle 2 \rangle. $
As in the previous proposition, the map $ h $ fits into the distinguished triangle

\xymatrix
{
&\Id[2] \langle 1 \rangle \ar[r]^{\beta[2] \langle 1 \rangle}	&\widetilde{\epsilon}_{\mathfrak{s}_i} L\widetilde{Z}^{\mathfrak{s}_i}[2] \langle 2 \rangle \ar[r]^{h}	
&\text{Cone}(\beta)[2]\langle 1 \rangle.\\
}

Therefore $ \text{Cone}(\alpha)\text{Cone}(\beta) \cong \Id $ and we have shown that $ \Pi_i $ and $ \Omega_i $ are inverse auto-equivalences.
\end{proof}

\subsection{Relation 5}

\begin{lemma}
\label{sch7.2}
Suppose we have the following morphism of distinguished triangles where $ f $ and $ h $ are split inclusion maps.  

\begin{tiny}
\xymatrix
{
&X\ar[r]^{\alpha}\ar[d]^{f}	&Y\ar[r]^{\beta}\ar[d]^{h}	&Z\ar[d]^{j}\\
&X^{'}\ar[r]^{\alpha^{'}}	&Y^{'}\ar[r]^{\beta^{'}}	&Z^{'}.
}
\end{tiny}
Then $ j $ is also a split inclusion.
\end{lemma}

\begin{proof}
See exercise 7.2 of [Sch].
\end{proof}

\begin{prop}
There is a functorial braid relation $ \Pi_i \circ \Pi_{i+1} \circ \Pi_i \cong \Pi_{i+1} \circ \Pi_i \circ \Pi_{i+1}. $
\end{prop}

\begin{proof}
First we will consider the commutative diagram expressing the functor $ \Pi_i \circ \Pi_{i+1}. $

\begin{tiny}
\xymatrix
{
&\text{Cone}(\alpha_{i+1})\ar[r]	& \widetilde{\epsilon}_{\mathfrak{s}_i}L\widetilde{Z}^{\mathfrak{s}_i}\langle 1 \rangle \text{Cone}(\alpha_{i+1})\ar[r]	&\text{Cone}(\alpha_i)\text{Cone}(\alpha_{i+1})	\\
&\widetilde{\epsilon}_{\mathfrak{s}_{i+1}}L\widetilde{Z}^{\mathfrak{s}_{i+1}} \langle 1 \rangle\ar[r]\ar[u]		&\widetilde{\epsilon}_{\mathfrak{s}_i}L\widetilde{Z}^{\mathfrak{s}_i} \langle 1 \rangle\widetilde{\epsilon}_{\mathfrak{s}_{i+1}}L\widetilde{Z}^{\mathfrak{s}_{i+1}} \langle 1 \rangle \ar[r]\ar[u]	&\text{Cone}(\alpha_{i+1}) \widetilde{\epsilon}_{\mathfrak{s}_{i+1}}L\widetilde{Z}^{\mathfrak{s}_{i+1}} \langle 1 \rangle \ar[u]\\
&\Id\ar[r]^{\alpha_i}\ar[u]^{\alpha_{i+1}}	&\widetilde{\epsilon}_{\mathfrak{s}_i}L\widetilde{Z}^{\mathfrak{s}_i}\langle 1 \rangle\ar[u]\ar[r]	&\text{Cone}(\alpha_i)\ar[u] \\
}
\end{tiny}

Apply $ \widetilde{\epsilon}_{\mathfrak{s}_{i+1}}L\widetilde{Z}^{\mathfrak{s}_{i+1}} $ to $ \Pi_i \circ \Pi_{i+1} $ and to the above diagram to get the following diagram.

\begin{tiny}
\xymatrix
{
&\widetilde{\epsilon}_{\mathfrak{s}_{i+1}}L\widetilde{Z}^{\mathfrak{s}_{i+1}} \langle 1 \rangle \text{Cone}(\alpha_{i+1})\ar[r]	& \widetilde{\epsilon}_{\mathfrak{s}_{i+1}}L\widetilde{Z}^{\mathfrak{s}_{i+1}} \langle 1 \rangle \widetilde{\epsilon}_{\mathfrak{s}_i}L\widetilde{Z}^{\mathfrak{s}_i}\langle 1 \rangle \text{Cone}(\alpha_{i+1})\ar[r]	&\widetilde{\epsilon}_{\mathfrak{s}_{i+1}}L\widetilde{Z}^{\mathfrak{s}_{i+1}} \langle 1 \rangle \text{Cone}(\alpha_i)\text{Cone}(\alpha_{i+1})	\\
&\widetilde{\epsilon}_{\mathfrak{s}_{i+1}}L\widetilde{Z}^{\mathfrak{s}_{i+1}} \langle 1 \rangle\widetilde{\epsilon}_{\mathfrak{s}_{i+1}}L\widetilde{Z}^{\mathfrak{s}_{i+1}} \langle 1 \rangle\ar[r]\ar[u]		&\widetilde{\epsilon}_{\mathfrak{s}_{i+1}}L\widetilde{Z}^{\mathfrak{s}_{i+1}} \langle 1 \rangle\widetilde{\epsilon}_{\mathfrak{s}_i}L\widetilde{Z}^{\mathfrak{s}_i} \langle 1 \rangle\widetilde{\epsilon}_{\mathfrak{s}_{i+1}}L\widetilde{Z}^{\mathfrak{s}_{i+1}} \langle 1 \rangle \ar[r]\ar[u]	&\widetilde{\epsilon}_{i+1}L\widetilde{Z}_{i+1} \langle 1 \rangle \text{Cone}(\alpha_{i+1}) \widetilde{\epsilon}_{\mathfrak{s}_{i+1}}L\widetilde{Z}^{\mathfrak{s}_{i+1}} \langle 1 \rangle \ar[u]\\
&\widetilde{\epsilon}_{\mathfrak{s}_{i+1}}L\widetilde{Z}^{\mathfrak{s}_{i+1}} \langle 1 \rangle\Id\ar[r]^{\widetilde{\epsilon}_{\mathfrak{s}_{i+1}}L\widetilde{Z}^{\mathfrak{s}_{i+1}} \langle 1 \rangle \alpha_i}\ar[u]^{\widetilde{\epsilon}_{\mathfrak{s}_{i+1}}L\widetilde{Z}^{\mathfrak{s}_{i+1}} \langle 1 \rangle \alpha_{i+1}}	&\widetilde{\epsilon}_{\mathfrak{s}_{i+1}}L\widetilde{Z}^{\mathfrak{s}_{i+1}} \langle 1 \rangle\widetilde{\epsilon}_{\mathfrak{s}_i}L\widetilde{Z}^{\mathfrak{s}_i}\langle 1 \rangle\ar[u]\ar[r]	&\widetilde{\epsilon}_{\mathfrak{s}_{i+1}}L\widetilde{Z}^{\mathfrak{s}_{i+1}} \langle 1 \rangle \text{Cone}(\alpha_i)\ar[u] \\
}
\end{tiny}

There is an obvious morphism from the first diagram to the second one.  All of the squares commute.  Now apply $ \text{Cone}(\alpha_{i+1}) $ to the first diagram to get the following commutative diagram.

\begin{tiny}
\xymatrix
{
&\text{Cone}(\alpha_{i+1})\text{Cone}(\alpha_{i+1})\ar[r]	&\text{Cone}(\alpha_{i+1})\widetilde{\epsilon}_{\mathfrak{s}_i}L\widetilde{Z}^{\mathfrak{s}_i}\langle 1 \rangle \text{Cone}(\alpha_{i+1})\ar[r]		&\text{Cone}(\alpha_{i+1})\text{Cone}(\alpha_i)\text{Cone}(\alpha_{i+1})	\\
&\text{Cone}(\alpha_{i+1})\widetilde{\epsilon}_{\mathfrak{s}_{i+1}}L\widetilde{Z}^{\mathfrak{s}_{i+1}} \langle 1 \rangle\ar[r]\ar[u]		&\text{Cone}(\alpha_{i+1})\widetilde{\epsilon}_{\mathfrak{s}_i}L\widetilde{Z}^{\mathfrak{s}_i} \langle 1 \rangle\widetilde{\epsilon}_{\mathfrak{s}_{i+1}}L\widetilde{Z}^{\mathfrak{s}_{i+1}} \langle 1 \rangle \ar[r]\ar[u]	&\text{Cone}(\alpha_{i+1})\text{Cone}(\alpha_{i+1}) \widetilde{\epsilon}_{\mathfrak{s}_{i+1}}L\widetilde{Z}^{\mathfrak{s}_{i+1}} \langle 1 \rangle \ar[u]\\
&\text{Cone}(\alpha_{i+1})\Id\ar[r]^{\text{Cone}(\alpha_{i+1})\alpha_i}\ar[u]^{\text{Cone}(\alpha_{i+1})\alpha_{i+1}}	&\text{Cone}(\alpha_{i+1})\widetilde{\epsilon}_i L\widetilde{Z}_i \langle 1 \rangle\ar[u]\ar[r]&\text{Cone}(\alpha_{i+1})\text{Cone}(\alpha_i)\ar[u] \\
}
\end{tiny}

There is a morphism from the second diagram to this one coming from the distinguished triangle 
$$ \Id \rightarrow \widetilde{\epsilon}_{\mathfrak{s}_{i+1}}L\widetilde{Z}^{\mathfrak{s}_{i+1}} \langle 1 \rangle \rightarrow \text{Cone}(\alpha_{i+1}). $$  
All of the resulting squares commute.

Let $ F = \text{Cone}((\text{Cone}(\Id \rightarrow \widetilde{\epsilon}_{\mathfrak{s}_{i+1}}L \widetilde{Z}^{\mathfrak{s}_{i+1}} \langle 1 \rangle)) \rightarrow \text{Cone}(\widetilde{\epsilon}_{\mathfrak{s}_i} L\widetilde{Z}^{\mathfrak{s}_i} \langle 1 \rangle \rightarrow \widetilde{\epsilon}_{\mathfrak{s}_i} L\widetilde{Z}^{\mathfrak{s}_i} \langle 1 \rangle \widetilde{\epsilon}_{\mathfrak{s}_{i+1}} L\widetilde{Z}^{\mathfrak{s}_{i+1}} \langle 1 \rangle)). $

Then $ \Pi_{i+1} \circ \Pi_i \circ \Pi_{i+1} = 
\text{Cone}(F \rightarrow \widetilde{\epsilon}_{\mathfrak{s}_{i+1}} L\widetilde{Z}^{\mathfrak{s}_{i+1}} \langle 1 \rangle F). $
This is isomorphic to
$$ \text{Cone}(F \rightarrow \text{Cone}(\widetilde{\epsilon}_{\mathfrak{s}_{i+1}} L\widetilde{Z}^{\mathfrak{s}_{i+1}}[2] \langle 3 \rangle \rightarrow \text{Cone}(\widetilde{\epsilon}_{\mathfrak{s}_{i+1}} L\widetilde{Z}^{\mathfrak{s}_{i+1}} \langle 1 \rangle
\widetilde{\epsilon}_{\mathfrak{s}_i} L\widetilde{Z}^{\mathfrak{s}_i} \langle 1 \rangle \rightarrow \widetilde{\epsilon}_{\mathfrak{s}_{i+1}} L\widetilde{Z}^{\mathfrak{s}_{i+1}}[2] \langle 3 \rangle \oplus \widetilde{\epsilon}_{\mathfrak{t}_i} L\widetilde{Z}^{\mathfrak{t}_i}[1] \langle 3 \rangle))). $$

We claim that 
$$ \text{Cone}(\widetilde{\epsilon}_{\mathfrak{s}_{i+1}} L\widetilde{Z}^{\mathfrak{s}_{i+1}}[2] \langle 3 \rangle \rightarrow \text{Cone}(\widetilde{\epsilon}_{\mathfrak{s}_{i+1}} L\widetilde{Z}^{\mathfrak{s}_{i+1}} \langle 1 \rangle
\widetilde{\epsilon}_{\mathfrak{s}_i} L\widetilde{Z}^{\mathfrak{s}_i} \langle 1 \rangle \rightarrow \widetilde{\epsilon}_{\mathfrak{s}_{i+1}} L\widetilde{Z}^{\mathfrak{s}_{i+1}}[2] \langle 3 \rangle \oplus \widetilde{\epsilon}_{\mathfrak{t}_i} L\widetilde{Z}^{\mathfrak{t}_i}[1] \langle 3 \rangle))  $$
is isomorphic to 
$ \text{Cone}(\widetilde{\epsilon}_{\mathfrak{s}_{i+1}} L\widetilde{Z}^{\mathfrak{s}_{i+1}} \langle 1 \rangle
\widetilde{\epsilon}_{\mathfrak{s}_i} L\widetilde{Z}^{\mathfrak{s}_i} \langle 1 \rangle \rightarrow  \widetilde{\epsilon}_{\mathfrak{t}_i} L\widetilde{Z}^{\mathfrak{t}_i}[1] \langle 3 \rangle) . $

We have the following commutative diagram

\begin{tiny}
\xymatrix
{
&\widetilde{\epsilon}_{\mathfrak{s}_{i+1}}L \widetilde{Z}^{\mathfrak{s}_{i+1}} \langle 1 \rangle \ar[r]^{\rho}\ar[d]	&\widetilde{\epsilon}_{\mathfrak{s}_{i+1}}L \widetilde{Z}^{\mathfrak{s}_{i+1}} \langle 1 \rangle\widetilde{\epsilon}_{\mathfrak{s}_{i+1}}L \widetilde{Z}^{\mathfrak{s}_{i+1}} \langle 1 \rangle\ar[r]^{\pi}\ar[d]^{f}	&\text{Cone}(\rho)\ar[d]^{\nu}\\
&\widetilde{\epsilon}_{\mathfrak{s}_{i+1}} L \widetilde{Z}^{\mathfrak{s}_{i+1}} \langle 1 \rangle \widetilde{\epsilon}_{\mathfrak{s}_i} L \widetilde{Z}^{\mathfrak{s}_i} \langle 1 \rangle\ar[r]^{\mu}	&\widetilde{\epsilon}_{\mathfrak{s}_{i+1}} L \widetilde{Z}^{\mathfrak{s}_{i+1}} \langle 1 \rangle \widetilde{\epsilon}_{\mathfrak{s}_i} L \widetilde{Z}^{\mathfrak{s}_i} \langle 1 \rangle \widetilde{\epsilon}_{\mathfrak{s}_{i+1}} L \widetilde{Z}^{\mathfrak{s}_{i+1}} \langle 1 \rangle \ar[r]^{g}	&\text{Cone}(\mu)\\
}
\end{tiny}

By properties of adjunction maps, $ \pi $ maps $ \widetilde{\epsilon}_{\mathfrak{s}_{i+1}}L \widetilde{Z}^{\mathfrak{s}_{i+1}}[2] \langle 3 \rangle $ onto itself and $ \widetilde{\epsilon}_{\mathfrak{s}_{i+1}}L \widetilde{Z}^{\mathfrak{s}_{i+1}} \langle 1 \rangle $ to zero.
Thus,
$$ \nu \pi(\widetilde{\epsilon}_{\mathfrak{s}_{i+1}}L \widetilde{Z}^{\mathfrak{s}_{i+1}}[2] \langle 3 \rangle) \cong \nu(\widetilde{\epsilon}_{\mathfrak{s}_{i+1}}L \widetilde{Z}^{\mathfrak{s}_{i+1}}[2] \langle 3 \rangle) \cong 
g \circ f(\widetilde{\epsilon}_{\mathfrak{s}_{i+1}}L \widetilde{Z}^{\mathfrak{s}_{i+1}}[2] \langle 3 \rangle). $$
Let us examine the map $ f. $   There is an adjunction map $ Id[2] \langle 1 \rangle \rightarrow L\widetilde{Z}^{\mathfrak{s}_i} \widetilde{\epsilon}_{\mathfrak{s}_i}. $ There is another map
$ L\widetilde{Z}^{\mathfrak{s}_i} \widetilde{\epsilon}_{\mathfrak{s}_i} \rightarrow L\widetilde{Z}^{\mathfrak{s}_i} \widetilde{\epsilon}_{\mathfrak{s}_{i+1}} L\widetilde{Z}^{\mathfrak{s}_{i+1}} \langle 1 \rangle \widetilde{\epsilon}_{\mathfrak{s}_i} $ which also comes from adjunction.  From general principles (see [Mac]), the compositon of these two maps is adjunction 
$ \Id[2]\langle 1 \rangle \rightarrow L\widetilde{Z}^{\mathfrak{s}_i} \widetilde{\epsilon}_{\mathfrak{s}_{i+1}} L\widetilde{Z}^{\mathfrak{s}_{i+1}} \langle 1 \rangle \widetilde{\epsilon}_{\mathfrak{s}_i}. $  We know from proposition ~\ref{diagram5} that this map is a split inclusion.  This is how the map $ f $ is constructed. This then gives the following commutative diagram of distinguished triangles.

\begin{tiny}
\xymatrix
{
&0 \ar[r]\ar[d]	&\widetilde{\epsilon}_{\mathfrak{s}_{i+1}}L \widetilde{Z}^{\mathfrak{s}_{i+1}}[2] \langle 3 \rangle\ar[r]\ar[d]	 &\widetilde{\epsilon}_{\mathfrak{s}_{i+1}}L \widetilde{Z}^{\mathfrak{s}_{i+1}}[2] \langle 3 \rangle\ar[d]^{\delta}\\
&\widetilde{\epsilon}_{\mathfrak{s}_{i+1}}L \widetilde{Z}^{\mathfrak{s}_{i+1}} \langle 1 \rangle \widetilde{\epsilon}_{\mathfrak{s}_i}L \widetilde{Z}^{\mathfrak{s}_i} \langle 1 \rangle\ar[r]^{\mu}\ar[d]	&\widetilde{\epsilon}_{\mathfrak{s}_{i+1}}L \widetilde{Z}^{\mathfrak{s}_{i+1}} \langle 1 \rangle \widetilde{\epsilon}_{\mathfrak{s}_i}L \widetilde{Z}^{\mathfrak{s}_i} \langle 1 \rangle \widetilde{\epsilon}_{\mathfrak{s}_{i+1}}L \widetilde{Z}^{\mathfrak{s}_{i+1}} \langle 1 \rangle \ar[r]^{g}\ar[d]	&\text{Cone}(\mu)\ar[d]\\
&\widetilde{\epsilon}_{\mathfrak{s}_{i+1}}L \widetilde{Z}^{\mathfrak{s}_{i+1}} \langle 1 \rangle \widetilde{\epsilon}_{\mathfrak{s}_i}L \widetilde{Z}^{\mathfrak{s}_i} \langle 1 \rangle\ar[r]
&\widetilde{\epsilon}_{\mathfrak{t}_i} L\widetilde{Z}^{\mathfrak{t}_i}[1] \langle 3 \rangle \ar[r]
&\text{Cone}(\delta)
}
\end{tiny}

This implies the claim.

The functor $ \Pi_i \circ \Pi_{i+1} \circ \widetilde{\epsilon}_{\mathfrak{s}_i} L\widetilde{Z}^{\mathfrak{s}_i} \langle 1 \rangle $ is expressed by the following commutative diagram.  There is a morphism from the first diagram to this one such that all the resulting squares commute.

\begin{tiny}
\xymatrix
{
&\text{Cone}(\alpha_{i+1})\ar[r]\widetilde{\epsilon}_{\mathfrak{s}_i} L\widetilde{Z}^{\mathfrak{s}_i} \langle 1 \rangle	& \widetilde{\epsilon}_{\mathfrak{s}_i}L\widetilde{Z}^{\mathfrak{s}_i}\langle 1 \rangle \text{Cone}(\alpha_{i+1})\ar[r]\widetilde{\epsilon}_{\mathfrak{s}_i} L\widetilde{Z}^{\mathfrak{s}_i} \langle 1 \rangle	&\text{Cone}(\alpha_i)\text{Cone}(\alpha_{i+1})\widetilde{\epsilon}_{\mathfrak{s}_i} L\widetilde{Z}^{\mathfrak{s}_i} \langle 1 \rangle	\\
&\widetilde{\epsilon}_{\mathfrak{s}_{i+1}}L\widetilde{Z}^{\mathfrak{s}_{i+1}} \langle 1 \rangle\ar[r]\ar[u]	\widetilde{\epsilon}_{\mathfrak{s}_i} L\widetilde{Z}^{\mathfrak{s}_i} \langle 1 \rangle	&\widetilde{\epsilon}_{\mathfrak{s}_i}L\widetilde{Z}^{\mathfrak{s}_i} \langle 1 \rangle\widetilde{\epsilon}_{\mathfrak{s}_{i+1}}L\widetilde{Z}^{\mathfrak{s}_{i+1}} \langle 1 \rangle \widetilde{\epsilon}_{\mathfrak{s}_i} L\widetilde{Z}^{\mathfrak{s}_i} \langle 1 \rangle\ar[r]\ar[u]	&\text{Cone}(\alpha_i)\widetilde{\epsilon}_{\mathfrak{s}_{i+1}}L\widetilde{Z}^{\mathfrak{s}_{i+1}} \widetilde{\epsilon}_{\mathfrak{s}_i}L\widetilde{Z}^{\mathfrak{s}_i} \langle 1 \rangle \ar[u]\\
&\Id\widetilde{\epsilon}_{\mathfrak{s}_i} L\widetilde{Z}^{\mathfrak{s}_i} \langle 1 \rangle\ar[r]^{\alpha_i \widetilde{\epsilon}_{\mathfrak{s}_i} L\widetilde{Z}^{\mathfrak{s}_i} \langle 1 \rangle}\ar[u]^{\alpha_{i+1}\widetilde{\epsilon}_{\mathfrak{s}_i} L\widetilde{Z}^{\mathfrak{s}_i} \langle 1 \rangle}	&\widetilde{\epsilon}_{\mathfrak{s}_i}L\widetilde{Z}^{\mathfrak{s}_i}\langle 1 \rangle\widetilde{\epsilon}_{\mathfrak{s}_i} L\widetilde{Z}^{\mathfrak{s}_i} \langle 1 \rangle \ar[u]\ar[r]	&\text{Cone}(\alpha_i) \widetilde{\epsilon}_{\mathfrak{s}_i} L\widetilde{Z}^{\mathfrak{s}_i} \langle 1 \rangle \ar[u] \\
}
\end{tiny}

The functor $ \Pi_i \circ \Pi_{i+1} \circ \Pi_i $ is expressed through the following diagram.

\begin{tiny}
\xymatrix
{
&\text{Cone}(\alpha_{i+1})\text{Cone}(\alpha_i) \ar[r]	& \widetilde{\epsilon}_{\mathfrak{s}_i}L\widetilde{Z}^{\mathfrak{s}_i}\langle 1 \rangle \text{Cone}(\alpha_{i+1})\text{Cone}(\alpha_i)\ar[r]	&\text{Cone}(\alpha_i)\text{Cone}(\alpha_{i+1})\text{Cone}(\alpha_i)	\\
&\widetilde{\epsilon}_{\mathfrak{s}_{i+1}}L\widetilde{Z}^{\mathfrak{s}_{i+1}} \langle 1 \rangle \text{Cone}(\alpha_i)\ar[r]\ar[u]		&\widetilde{\epsilon}_{\mathfrak{s}_i}L\widetilde{Z}^{\mathfrak{s}_i} \langle 1 \rangle\widetilde{\epsilon}_{\mathfrak{s}_{i+1}}L\widetilde{Z}^{\mathfrak{s}_{i+1}} \langle 1 \rangle \text{Cone}(\alpha_i) \ar[r]\ar[u]	&\text{Cone}(\alpha_{i}) \widetilde{\epsilon}_{\mathfrak{s}_{i+1}}L\widetilde{Z}^{\mathfrak{s}_{i+1}} \langle 1 \rangle \text{Cone}(\alpha_i) \ar[u]\\
&\Id \text{Cone}(\alpha_i)\ar[r]^{\alpha_i \text{Cone}(\alpha_i)}\ar[u]^{\alpha_{i+1} \text{Cone}(\alpha_i)}	&\widetilde{\epsilon}_{\mathfrak{s}_i}L\widetilde{Z}^{\mathfrak{s}_i}\langle 1 \rangle \text{Cone}(\alpha_i) \ar[u]\ar[r]	&\text{Cone}(\alpha_i)\text{Cone}(\alpha_i)\ar[u] \\
}
\end{tiny}

By definition,
$ \Pi_i \circ \Pi_{i+1} \circ \Pi_i = \text{Cone}(F \rightarrow F \widetilde{\epsilon}_{\mathfrak{s}_i}L \widetilde{Z}^{\mathfrak{s}_i} \langle 1 \rangle). $
We repeat the computation done for $ \Pi_{i+1} \circ \Pi_i \circ \Pi_{i+1}, $ and find this is isomorphic to
$ \text{Cone}(F \rightarrow \text{Cone}(\widetilde{\epsilon}_{\mathfrak{s}_{i+1}} L\widetilde{Z}^{\mathfrak{s}_{i+1}} \langle 1 \rangle \widetilde{\epsilon}_{\mathfrak{s}_i} L\widetilde{Z}^{\mathfrak{s}_i} \langle 1 \rangle \rightarrow 
\widetilde{\epsilon}_{\mathfrak{t}_i} L\widetilde{Z}^{\mathfrak{t}_i}[1] \langle 3 \rangle)). $

Now we just have to verify that the maps in both cones are the same.
Consider the commutative diagram

\xymatrix
{
&F\ar[r]\ar[d]	&F\widetilde{\epsilon}_{i} L\widetilde{Z_i} \langle 1 \rangle \ar[d]\\
&\widetilde{\epsilon}_{i+1} L\widetilde{Z}_{i+1} \langle 1 \rangle F\ar[r]	&\widetilde{\epsilon}_{i+1} L\widetilde{Z}_{i+1} \langle 1 \rangle F 
\widetilde{\epsilon}_{i} L\widetilde{Z}_{i} \langle 1 \rangle\\
}

Since $ \widetilde{\epsilon}_{i+1} L\widetilde{Z}_{i+1} \langle 1 \rangle F \rightarrow \widetilde{\epsilon}_{i+1} L\widetilde{Z}_{i+1} \langle 1 \rangle F \widetilde{\epsilon}_{i} L\widetilde{Z}_i \langle 1 \rangle $ and $ F\widetilde{\epsilon}_{i} L\widetilde{Z}_i \langle 1 \rangle \rightarrow \widetilde{\epsilon}_{i+1} L\widetilde{Z}_{i+1} \langle 1 \rangle F \widetilde{\epsilon}_{i} L\widetilde{Z}_i \langle 1 \rangle $ are split inclusions by lemma ~\ref{sch7.2}, both maps in the cone are the same.

\end{proof}

\begin{corollary}
There is an isomorphism $$ \Omega_{i}^{} \Omega_{i+1}^{} \Omega_{i}^{} \cong \Omega_{i+1}^{} \Omega_{i}^{} \Omega_{i+1}^{}. $$
\end{corollary}

\begin{proof}
This follows immediately from the proposition by taking the adjoint of both sides.
\end{proof}

\subsection{Relation 6}
\begin{prop}
There is an isomorphism $ \cap_{i+1,+,r+2} \circ \Omega_{i} \circ \cup_{i+1,+,r} \cong \Id. $
\end{prop}

\begin{proof}
The left hand sides is
$$ L\widetilde{Z}_{\mathfrak{q}_{i+1}}^{\mathfrak{p}_{i+1}} \text{Cone}(\widetilde{\epsilon}_{\mathfrak{s}_i}[-1]L\widetilde{Z}^{\mathfrak{s}_i} \rightarrow \Id[1] \langle 1 \rangle) \widetilde{\epsilon}_{\mathfrak{p}_{i+1}}^{\mathfrak{q}_{i+1}}[-(k-1)][-k] \langle -k \rangle \cong $$
$$ \text{Cone}(L\widetilde{Z}_{\mathfrak{q}_{i+1}}^{\mathfrak{p}_{i+1}} \widetilde{\epsilon}_{\mathfrak{s}_i}[-1]L\widetilde{Z}^{\mathfrak{s}_i} \widetilde{\epsilon}_{\mathfrak{p}_{i+1}}^{\mathfrak{q}_{i+1}}[-(k-1)] \rightarrow
L\widetilde{Z}_{\mathfrak{q}_{i+1}}^{\mathfrak{p}_{i+1}} \widetilde{\epsilon}_{\mathfrak{p}_{i+1}}^{\mathfrak{q}_{i+1}}[-(k-1)+1] \langle 1 \rangle)[-k] \langle -k \rangle. $$
Now apply proposition 4.19 to the left hand side and corollary 4.5 to the right to get
$$ \text{Cone}( \oplus_{j=0}^{k-2} \Id[k-2-2j] \langle k-2-2j \rangle \rightarrow \oplus_{j=0}^{k-1} \Id[k-2j] \langle k-2j \rangle)[-k]\langle -k \rangle \cong \Id. $$
\end{proof}

\begin{prop}
There is an isomorphism $ \cap_{i+1,+,r+2} \circ \Pi_{i} \circ \cup_{i+1,+,r} \cong \Id. $
\end{prop}

\begin{proof}
The left hand side is
$$ L\widetilde{Z}_{\mathfrak{q}_{i+1}}^{\mathfrak{p}_{i+1}} \text{Cone}(\Id \langle -1 \rangle \rightarrow
\widetilde{\epsilon}_{\mathfrak{s}_i} L\widetilde{Z}^{\mathfrak{s}_i} \langle 1 \rangle)\widetilde{\epsilon}_{\mathfrak{p}_{i+1}}^{\mathfrak{q}_{i+1}}[-(k-1)][k-2] \langle k \rangle \cong $$
$$ \text{Cone}(L\widetilde{Z}_{\mathfrak{q}_{i+1}}^{p_{i+1}} \widetilde{\epsilon}_{\mathfrak{p}_{i+1}}^{\mathfrak{q}_{i+1}}[-(k-1)]\langle -1 \rangle \rightarrow 
L\widetilde{Z}_{\mathfrak{q}_{i+1}}^{\mathfrak{p}_{i+1}} \widetilde{\epsilon}_{\mathfrak{s}_i} L\widetilde{Z}^{\mathfrak{s}_i} \widetilde{\epsilon}_{\mathfrak{p}_{i+1}}^{\mathfrak{q}_{i+1}}[-(k-1)])[k-2] \langle k \rangle. $$
Now apply corollary 10 to the left side and proposition 12 to the right  to get
$$ \text{Cone}(\oplus_{j=0}^{k-1} \Id[k-1-2j] \langle k-2-2j \rangle \rightarrow \oplus_{j=0}^{k-2} \Id[k-1-2j] \langle k-2-2j \rangle)[k-2]\langle k \rangle \cong
\Id. $$
\end{proof}

\subsection{Relation 7}
For purposes of this subsection, we have the following abbreviations for the relevant subalgebras.

\begin{define}
\begin{enumerate}
\item $ \mu = \mathfrak{p}_i + \mathfrak{q}_{i+k}. $
\item $ \nu= \mathfrak{q}_{i-1}+\mathfrak{q}_{i+k}. $
\item $ \gamma = \mathfrak{q}_{i-1}+\mathfrak{s}_{i+k-1}+\mathfrak{q}_{i+k}. $
\item $ \delta = \mathfrak{q}_{i-1}+\mathfrak{p}_{i+k}.$
\end{enumerate}
\end{define}

We let $ \alpha \colon \widetilde{\epsilon}_{\gamma}^{\nu} L\widetilde{Z}_{\nu}^{\gamma} \rightarrow \Id[2] \langle 1 \rangle $ and
$ \beta \colon \Id \rightarrow \widetilde{\epsilon}_{\gamma}^{\nu} L \widetilde{Z}_{\nu}^{\gamma} \langle 1 \rangle. $

\begin{prop}
There are functorial isomorphisms:
\begin{enumerate}
\item $$ \cap_{i,-,r+2} \circ \Omega_{i+k-1} \circ \cup_{i+2,+,r} \circ \cap_{i+2,+,r+2} \circ \Pi_{i+k-1} \circ \cup_{i, -,r} \cong
\Id $$
\item $$ \cap_{i,-,r+2} \circ \Pi_{i+k-1} \circ \cup_{i+2,+,r} \circ \cap_{i+2,+,r+2} \circ \Omega_{i+k-1} \circ \cup_{i, -,r} \cong \Id. $$
\end{enumerate}
\end{prop}

\begin{proof}
\begin{enumerate}
\item The left hand side is the functor given below shifted by $ [-3] \langle -1 \rangle. $
$$ L\widetilde{Z}_{\nu}^{\mu} 
\text{Cone}(\Id \overset\beta\to \widetilde{\epsilon}_{\gamma}^{\nu}L\widetilde{Z}_{\nu}^{\gamma}\langle 1 \rangle)\widetilde{\epsilon}_{\delta}^{\nu}[-(k-1)]
L\widetilde{Z}_{\nu}^{\delta}
\text{Cone}(\widetilde{\epsilon}_{\gamma}^{\nu} L\widetilde{Z}_{\nu}^{\gamma} \overset\alpha\to \Id[2] \langle 1 \rangle) \widetilde{\epsilon}_{\mu}^{\nu}[-(k-1)].  $$

This expression gives rise to the following commutative diagram.

\begin{tiny}
\xymatrix
{
&L\widetilde{Z}_{\nu}^{\mu}\text{Cone}(\beta)\widetilde{\epsilon}_{\delta}^{\nu}L\widetilde{Z}_{\nu}^{\delta}\widetilde{\epsilon}_{\gamma}^{\nu}L\widetilde{Z}_{\nu}^{\gamma}\widetilde{\epsilon}_{\mu}^{\nu}[-2(k-1)] \ar[r]		
&L\widetilde{Z}_{\nu}^{\mu} \text{Cone}(\beta) \widetilde{\epsilon}_{\delta}^{\nu}L\widetilde{Z}_{\nu}^{\delta} \widetilde{\epsilon}_{\mu}^{\nu}[-2k+4]\langle 1 \rangle \ar[r]
&L\widetilde{Z}_{\nu}^{\mu} \text{Cone}(\beta) \widetilde{\epsilon}_{\delta}^{\nu}L\widetilde{Z}_{\nu}^{\delta}\text{Cone}(\alpha)\widetilde{\epsilon}_{\mu}^{\nu}[-2(k-1)]	\\
&L\widetilde{Z}_{\nu}^{\mu} \widetilde{\epsilon}_{\gamma}^{\nu}L\widetilde{Z}_{\nu}^{\gamma} \widetilde{\epsilon}_{\delta}^{\nu}L\widetilde{Z}_{\nu}^{\delta} \widetilde{\epsilon}_{\gamma}^{\nu}L\widetilde{Z}_{\nu}^{\gamma}\widetilde{\epsilon}_{\mu}^{\nu}[-2(k-1)]\langle 1 \rangle \ar[u]\ar[r]
&L\widetilde{Z}_{\nu}^{\mu} \widetilde{\epsilon}_{\gamma}^{\nu} L\widetilde{Z}_{\nu}^{\gamma} \widetilde{\epsilon}_{\delta}^{\nu}L\widetilde{Z}_{\nu}^{\delta}\widetilde{\epsilon}_{\mu}^{\nu}[-2k+4]\langle 2 \rangle \ar[u]\ar[r]
&L\widetilde{Z}_{\nu}^{\mu} \widetilde{\epsilon}_{\gamma}^{\nu}L\widetilde{Z}_{\nu}^{\gamma} \widetilde{\epsilon}_{\delta}^{\nu} L\widetilde{Z}_{\nu}^{\delta} \text{Cone}(\alpha) \widetilde{\epsilon}_{\mu}^{\nu}[-2(k-1)]\langle 1 \rangle \ar[u] \\ 
&L\widetilde{Z}_{\nu}^{\mu} \widetilde{\epsilon}_{\delta}^{\nu}L\widetilde{Z}_{\nu}^{\delta} \widetilde{\epsilon}_{\gamma}^{\nu} L\widetilde{Z}_{\nu}^{\gamma} \widetilde{\epsilon}_{\mu}^{\nu}[-2(k-1)]\ar[u]\ar[r]
&L\widetilde{Z}_{\nu}^{\mu} \widetilde{\epsilon}_{\delta}^{\nu}L\widetilde{Z}_{\nu}^{\delta} \widetilde{\epsilon}_{\mu}^{\nu}[-2k+4] \langle 1 \rangle \ar[u]\ar[r]
&L\widetilde{Z}_{\nu}^{\mu} \widetilde{\epsilon}_{\delta}^{\nu}L\widetilde{Z}_{\nu}^{\delta}\text{Cone}(\alpha)\widetilde{\epsilon}_{\mu}^{\nu}[-2(k-1)]\ar[u].\\
}
\end{tiny}

Using the isomorphisms constructed earlier, the above diagram is the same as

\begin{tiny}
\xymatrix
{
&\widetilde{\epsilon}_{\mathfrak{p}_i}^{\mathfrak{q}_i}L\widetilde{Z}_{\mathfrak{q}_i}^{\mathfrak{p}_i}[-2k+4]\langle -k+2 \rangle \oplus \Id[2]\langle 1 \rangle\ar[r]
&\widetilde{\epsilon}_{\mathfrak{p}_i}^{\mathfrak{q}_i}L\widetilde{Z}_{\mathfrak{q}_i}^{\mathfrak{p}_i}[-2k+5] \langle -k+2 \rangle\ar[r] 
&X\\
&\oplus_{j=0}^{k-3} \widetilde{\epsilon}_{\mathfrak{p}_i}^{\mathfrak{q}_i}L\widetilde{Z}_{\mathfrak{q}_i}^{\mathfrak{p}_i}[-2j]\langle k-2-2j \rangle \oplus \Id[2]\langle 1 \rangle\ar[u]\ar[r]
&\oplus_{j=0}^{k-2} \widetilde{\epsilon}_{p_i}^{q_i}L\widetilde{Z}_{\mathfrak{q}_i}^{\mathfrak{p}_i}[2-2j]\langle k-2j \rangle\ar[r]\ar[u]
&\widetilde{\epsilon}_{\mathfrak{p}_i}^{\mathfrak{q}_i}L\widetilde{Z}_{\mathfrak{q}_i}^{\mathfrak{p}_i}[2]\langle k \rangle \oplus \Id\ar[u]\\
&\oplus_{j=0}^{k-2} \widetilde{\epsilon}_{\mathfrak{p}_i}^{\mathfrak{q}_i}L\widetilde{Z}_{\mathfrak{q}_i}^{\mathfrak{p}_i}[-2j]\langle k-2-2j \rangle\ar[r]\ar[u]
&\oplus_{j=0}^{k-1} \widetilde{\epsilon}_{\mathfrak{p}_i}^{\mathfrak{q}_i}L\widetilde{Z}_{\mathfrak{q}_i}^{\mathfrak{p}_i}[2-2j]\langle k-2j \rangle\ar[r]\ar[u]
&\widetilde{\epsilon}_{\mathfrak{p}_i}^{\mathfrak{q}_i}L\widetilde{Z}_{\mathfrak{q}_i}^{\mathfrak{p}_i}[2]\langle k \rangle\ar[u].\\
}
\end{tiny}

Therefore $ X \cong \Id[3]\langle 1 \rangle $ and so the composition of the functors with the its built in shifts is the identity.
\item The proof is the same.  One just has to reflect the commutative diagrams in the diagonal going from the bottom left to the top right.
\end{enumerate}
\end{proof}

\subsection{Relation 8}

\begin{prop}
There are the following functorial isomorphisms:
\begin{enumerate}
\item $$ \cap_{i,+,r+2} \circ \Omega_{i+k-3} \circ \cup_{i-2,-,r} \circ \cap_{i-2,-,r+2} \circ \Pi_{i+k-3} \circ \cup_{i,+,r} 
\cong \Id $$
\item $$ \cap_{i,+,r+2} \circ \Omega_{i+k-3} \circ \cup_{i-2,-,r} \circ \cap_{i-2,-,r+2} \circ \Pi_{i+k-3} \circ \cup_{i,+,r} \cong \Id. $$
\end{enumerate}
\end{prop}

\begin{proof}
One verifies these isomorphisms as in the proofs of relation seven.  The only differences are the subscripts of the parabolic subalgebras.
\end{proof}

\subsection{Skein Relations}
The images of the functors assigned to the crossings in the Grothendieck group satisfy the skein relations for the HOMFLYPT polynomial.  The following propositions follows directly from the definitions of the functors involved.

\begin{theorem}
\label{skein}
The following relations are satisfied in the Grothendieck group:
\begin{enumerate}
\item $ [\Pi_i]=(-1)^{k}q^k([\widetilde{\epsilon}_{\mathfrak{s}_i}L\widetilde{Z}^{\mathfrak{s}_i}]-q^{-1}[\Id]). $
\item $ [\Omega_i]=(-1)^k q^{-k}(-q[\Id]+[\widetilde{\epsilon}_{\mathfrak{s}_i} L\widetilde{Z}^{\mathfrak{s}_i}]) $
\item $ q^{k}[\Omega_i]-q^{-k}[\Pi_i]=(-1)^{k-1}(q^{1}-q^{-1})[\Id] $
\end{enumerate}
\end{theorem}

\end{document}